\documentclass[letterpaper,leqno,10pt,twoside]{amsart}
\usepackage{amsmath,amsthm,amsfonts,amssymb,euscript,mathrsfs,latexsym,marginnote}
\usepackage{graphicx}
\usepackage{float}
\usepackage{placeins}
\usepackage[margin=1.0in]{geometry}
\usepackage[toc,page]{appendix}
\usepackage{relsize}
\usepackage[shortlabels]{enumitem}
\usepackage{dsfont}

\setlist[enumerate,itemize]{labelindent=\parindent,leftmargin=*,align=left}
\makeatletter
\def\fullwidthdisplay{%
  \ifnum\@listdepth>\z@
    \displayindent\dimexpr\@totalleftmargin-\leftmargin+\labelindent\relax \displaywidth\dimexpr\linewidth+\leftmargin-\labelindent\relax
  \else
    \displayindent\z@ \displaywidth\columnwidth
  \fi
}
\makeatother

\usepackage[all]{xy}
\usepackage{tikz}

\usepackage{booktabs,tabularx,array,longtable}

\usepackage{indentfirst}

\allowdisplaybreaks

\numberwithin{equation}{section}

\usepackage[dvipsnames]{xcolor}

\usepackage{hyperref}
\hypersetup{colorlinks=true, pdfstartview={FitV 0},linkcolor=blue!70!black,citecolor=red!70!black, urlcolor=green!60!black}
\definecolor{labelkey}{rgb}{0.6,0,0}

\newtheorem{theorem}{Theorem}[section]
\newtheorem{corollary}[theorem]{Corollary}
\newtheorem{lemma}[theorem]{Lemma}
\newtheorem{proposition}[theorem]{Proposition}

\theoremstyle{definition}
\newtheorem{definition}{Definition}[section]

\newtheorem{assumption}[theorem]{Assumption}

\theoremstyle{remark}
\newtheorem{remark}{Remark}[section]

\providecommand{\abs}[1]{\left\lvert#1\right\rvert}
\providecommand{\babs}[1]{\big\lvert#1\big\rvert}
\providecommand{\nm}[1]{\left\lVert#1\right\rVert}
\providecommand{\bnm}[1]{\big\lVert#1\big\rVert}
\providecommand{\Bnm}[1]{\Big\lVert#1\Big\rVert}
\newcommand{\jbr}[1]{\left\langle#1\right\rangle}
\newcommand{\dualp}[2]{\left\langle#1,#2\right\rangle}

\def\ud{\mathrm d}
\def\ue{e}%{\mathrm e}
\def\ui{\mathrm i}

\def\e{\varepsilon}

\newcommand{\bulkdist}{r}
\newcommand{\bulkdistat}[1]{r_{#1}}
\newcommand{\sidedist}[1]{r_{#1}}
\newcommand{\sidedistat}[2]{r_{#1,#2}}
\newcommand{\scaledbulk}{R}
\newcommand{\scaledbulkat}[1]{R_{#1}}
\newcommand{\signedsideat}[2]{\xi_{#1,#2}}
\newcommand{\tphys}[2]{t_{#1,#2}}
\newcommand{\nphys}[2]{n^{\mathrm{in}}_{#1,#2}}
\newcommand{\twedge}[2]{t^{Y}_{#1,#2}}
\newcommand{\nwedge}[2]{n^{Y,\mathrm{in}}_{#1,#2}}

\newcommand{\ellvtx}{\ell_{\mathrm{vtx}}}
\newcommand{\ellend}{\ell_{\mathrm{end}}}
\newcommand{\ellnor}{\ell_{\mathrm{nor}}}
\newcommand{\ellsep}{\ell_{\mathrm{sep}}}
\newcommand{\ellOm}{\ell_{\Omega}}
\newcommand{\cfs}{c_{\mathrm{fs}}}
\newcommand{\rhodec}{\rho_{\mathrm{dec}}}

\newcommand{\Uphys}{u}
\newcommand{\Uint}{U}
\newcommand{\Uside}{\mathcal U}
\newcommand{\Uwedge}{\mathscr U}
\newcommand{\Err}{\mathscr R}
\newcommand{\Nvtx}{N_{\rm v}}

\newcommand{\Om}{\Omega}
\newcommand{\Sone}{\mathbb S^1}
\newcommand{\pk}{\mathbf P}
\newcommand{\qk}{\mathbf Q}
\newcommand{\Le}{\mathscr L_\e}
\newcommand{\Kw}{\mathscr K_\omega}
\newcommand{\Pw}{\mathsf T_\omega}
\newcommand{\supp}{\operatorname{supp}}
\newcommand{\Ext}{\operatorname{Ext}}
\newcommand{\gaux}{\widetilde g}
\DeclareMathOperator*{\esssup}{ess\,sup}
\DeclareMathOperator*{\essinf}{ess\,inf}
\DeclareMathOperator{\dist}{dist}

\newcommand{\Uapp}{U_{\mathrm{app}}}
\newcommand{\Ubase}{U_{\mathrm{base}}}
\newcommand{\Lameps}{\Lambda_\e}
\newcommand{\deleps}{\delta_\e}
\newcommand{\Tmatch}{T}
\newcommand{\Tcon}{\widetilde{T}}

\newcommand{\pacc}{p}
\newcommand{\Nexp}{N}
\newcommand{\betaw}{\beta}
\newcommand{\sdata}{s_\ast}
\newcommand{\qrate}{q_\ast}
\newcommand{\aoverlap}{a_\ast}
\newcommand{\Kreg}{K_{\rm reg}}
\newcommand{\Llog}{L_{\ln}}
\newcommand{\lambdaMin}{\lambda_\ast}
\newcommand{\kappaLtwo}{\kappa}

\newcommand{\Fourier}{\mathscr F}
\newcommand{\FourierInv}{\Fourier^{-1}}
\newcommand{\Cvertex}{C}
\newcommand{\EllData}{\mathcal A^{\rm ell}}
\newcommand{\UsideZero}{\mathcal U_{\rm sl,0}}
\newcommand{\UsideZeroAt}[1]{\mathcal U_{{\rm sl},#1,0}}

\title[Kinetic Wedge Layers in Polygonal Domains]{Kinetic Wedge Layers and Diffusive Limit of Neutron Transport in Polygonal Domains}

\author{Zhimeng Ouyang}
\address{School of Mathematical Sciences, Peking University, Beijing 100871, China}
\email{ouyangzm@math.pku.edu.cn}

\author{Lei Wu}
\address{Department of Mathematics, Lehigh University, Bethlehem, PA 18015, USA}
\email{lew218@lehigh.edu}

\date{}

\begin{document}

\begin{abstract}
We study the diffusive limit of the steady neutron transport equation with velocity-dependent in-flow boundary data in polygonal domains. On a bounded convex polygon, we construct an approximation by the interior solution, flat side layers, and wedge layers, and prove an algebraic error estimate in $L^\infty$. We also give a complete formulation and well-posedness theory for the kinetic wedge layer, including its algebraic decay. The proof combines a characteristic stability estimate, a weighted Mellin mapping theorem, a two-depth construction and matching scheme, and a shifted superharmonic barrier for the wedge corrector. We also prove the $L^2$ diffusive limit at rate $O(\e^{1/2})$ on general bounded simple polygons, including those with reentrant vertices. 
\end{abstract}

\maketitle
\enlargethispage{4pt}

\setcounter{tocdepth}{1}
\makeatletter
\def\l@subsection{\@tocline{2}{0pt}{2.5pc}{5pc}{}}
\makeatother
\tableofcontents
\section{Introduction}
\label{sec:intro}

The neutron transport equation describes particles that stream through a fixed domain and interact with a background medium. In applications such equations are usually posed in containers with corners. In reactor physics the computational domain is a lattice cell or an assembly with straight sides. In rarefied gas dynamics the same geometry appears at the sharp edge of a plate and in the corner of a cavity. Numerically, deterministic transport and diffusion codes are discretized on polygonal meshes, so a corner of the container is represented exactly rather than smoothed away.

The flat sides of a polygon admit the usual half-space description, while each vertex couples two incident side layers. A global $L^\infty$ diffusion approximation must account for this interaction together with the singular behavior of the limiting elliptic solution. We establish such an approximation on bounded convex polygons. The main new ingredient is a transport problem on a wedge with in-flow boundary data on both rays. We prove existence, uniqueness in $L^\infty$, and algebraic decay in a weighted $L^\infty$ norm for its wedge layer at every convex opening angle, and match it to the interior solution and the two incident side layers. We also prove an independent $L^2$ diffusive limit on general polygons, including those with reentrant vertices.

\subsection{Problem Setup}

Let $\Om\subset\mathbb R^2$ be a bounded simple polygon. Denote its open sides by $E_j$, its vertices by $V_1,\ldots,V_{\Nvtx}$, and its opening angles by $\omega_1,\ldots,\omega_{\Nvtx}$. We first assume that $\Om$ is convex, so
\begin{align}
 0<\omega_m<\pi\quad\text{for every }1\leq m\leq \Nvtx.
 \label{eq:strict-convexity}
\end{align}

All data and solutions are real-valued. The phase variables are $(x,w)\in\Om\times\Sone$. The normalized velocity measure is $\ud\nu(w)=\frac{\ud\varphi}{2\pi}$ with $w=(\cos\varphi,\sin\varphi)\in\Sone$. For a function $f(x,w)$, we set
\begin{align}
 (\pk f)(x):=\int_{\Sone}f(x,w)\,\ud\nu(w), \qquad \qk f:=f-\pk f,
\end{align}
so $\pk f$ is the velocity average and $\qk f$ is the zero-average part.

Let $n(x)$ be the outward unit normal at a non-vertex boundary point. We define the in-flow boundary $\gamma_-$ and the out-flow boundary $\gamma_+$ by
\begin{align}
 \gamma_\pm :=\big\{(x,w)\in\partial\Om\times\Sone: \ \pm w\cdot n(x)>0\big\}, \qquad \gamma_{j,\pm}:=\gamma_\pm\cap(E_j\times\Sone).
 \label{eq:physical-phase-boundary}
\end{align}

The phase-space $L^\infty$ norm, with respect to $\ud x\ud\nu(w)$, is
\begin{align}
 \nm F_\infty :=\nm F_{L^\infty(\Om\times\Sone)}={\esssup}_{(x,w)\in\Om\times\Sone}\abs{F(x,w)}.
 \label{eq:physical-infinity-norms}
\end{align}
Subscripts identify the variables in partial norms: $L^\infty_w=L^\infty(\Sone;\ud\nu)$, and $L^\infty_v$ denotes the same norm in a rotated velocity frame. We write $L^\infty_{\eta,w}=L^\infty(\mathbb R_+\times\Sone;\ud\eta\ud\nu)$.

The kinetic boundary measure and its norms, for $1\le p\le\infty$, are
\begin{align}
 \ud\gamma_{\rm kin} :=\abs{w\cdot n(x)}\,\ud\sigma(x)\ud\nu(w), \qquad \abs h_{p,\pm} :=\nm h_{L^{p}(\gamma_\pm;\ud\gamma_{\rm kin})}.
 \label{eq:l2-linfty-boundary-norms-comp}
\end{align}
Here $\ud\sigma$ is boundary arclength measure. For $p=\infty$, $\abs h_{\infty,\pm}$ is the kinetic $L^\infty(\gamma_\pm;\ud\gamma_{\rm kin})$ norm.

For a mild solution of the transport equation, $f\big|_{\gamma_\pm}$ and $f\big|_{\gamma_{j,\pm}}$ denote characteristic traces. For $H^1$ functions, the same restriction notation denotes the Sobolev trace.

For $0<\e\le1$, the neutron density $\Uphys^\e(x,w)$ satisfies the steady transport equation with in-flow boundary condition
\begin{align}
 \Le \Uphys^\e:=w\cdot\nabla_{\!x}\Uphys^\e+\e^{-1}\qk\, \Uphys^\e=0 \quad\text{in }\Om\times\Sone, \qquad \Uphys^\e\big|_{\gamma_-}=g^\e.
 \label{eq:physical-model}
\end{align}
We study the behavior of $\Uphys^\e$ as the Knudsen number $\e\to0^+$.

\subsection{Literature and Main Contributions}

\paragraph{\emph{Diffusion approximation}}
The diffusive limit of neutron transport has been studied since the 1970s. Formal interior and boundary layer expansions were constructed in \cite{Larsen.Keller1974,Larsen1974}, and rigorous diffusion limits were proved in \cite{Bensoussan.Lions.Papanicolaou1979,Bardos.Santos.Sentis1984}. The interior solution satisfies a diffusion equation. Its boundary values are determined by the end states of half-space Milne problems with velocity-dependent in-flow data. We refer to \cite{Sone2007} for the half-space theory and to \cite{dEon.Williams2018} for the two-dimensional isotropic kernel used here.

\paragraph{\emph{$L^{\infty}$ theory in smooth domains}}
Convergence in $L^\infty$ requires an approximation that includes boundary layers, which is $O(1)$ in $L^\infty$ norm. The $L^2$--$L^\infty$ framework gives remainder estimates for kinetic equations in bounded domains \cite{Guo2010,Esposito.Guo.Kim.Marra2013,Esposito.Guo.Kim.Marra2015}. For curved boundaries, the singularity near grazing velocities requires a geometric correction to the classical Milne expansion. This correction was developed for disks, annuli, and smooth convex domains in \cite{Wu.Guo2015,Wu.Yang.Guo2016,Guo.Wu2017,Wu2020,Wu2021,Wu2021(=)}. In a polygon, the local problem at a vertex is posed on a wedge. We obtain the required $L^\infty$ remainder estimate directly from the comparison principle in Section~\ref{sec:physical-stability}.

\paragraph{\emph{$L^2$ theory and general domains}}
In $L^2$, localization near the boundary gives an additional small factor. An $O(1)$ function in $L^\infty$ with support in a strip of width $O(\e)$ is $O(\e^{\frac{1}{2}})$ in $L^2$. This permits diffusion limits under weaker geometric assumptions. The smooth, possibly nonconvex case is treated in \cite{Guo.Wu2025,Ouyang2024}. The variational approach of \cite{Egger.Schlottbom2014} applies to bounded Lipschitz domains, heterogeneous coefficients, and general self-adjoint scattering operators. That work uses a subcritical scaling with absorption and small in-flow data. It proves convergence at Lipschitz regularity and a quantitative rate under additional regularity. Our $L^2$ result concerns transport equation with velocity-dependent in-flow data, and gives the rate $O(\e^{\frac{1}{2}})$ on general polygons.

\paragraph{\emph{Wedge layers}}
Boundary layers in a dihedral corner were studied in fluid mechanics in \cite{Carrier1947,Rubin1966}. Ordinary and elliptic corner layers were used to obtain $L^\infty$ expansions on rectangles in \cite{Shih.Kellogg1987}. A semilinear reaction--diffusion problem on a polygon led to a vertex correction on a stretched infinite sector in \cite{Kellogg.Kopteva2010}. In kinetic theory, \cite{Leonard1971} applied a multivariable Wiener--Hopf method to one-speed transport in a right-angle quarter space. The corresponding Milne problem was decomposed into two half-space contributions and a vertex correction in \cite{Lam.Leonard1971}. Numerical studies of BGK flow near edges, temperature jumps, and square-cavity vertices appear in \cite{Taguchi.Aoki2012,Taguchi.Tsuji2020,Hattori2024}. Here we prove existence, uniqueness in $L^\infty$, and decay of the wedge layer at every convex opening angle, and use these estimates to justify its matching with the interior solution and side layers.

\paragraph{\emph{Main contributions}}
Our main results are as follows.
\begin{enumerate}[label=\textup{(\arabic*)},leftmargin=*,itemsep=0.35em]

\item \emph{$L^{\infty}$ diffusion approximation on a convex polygon.}
Theorem~\ref{thm:main-comp} gives an approximation by the interior solution, side layers, and wedge layers, with an algebraic error bound in $L^\infty(\Om\times\Sone)$. The in-flow data may depend on velocity and satisfy the sidewise regularity and endpoint compatibility in Assumption~\ref{ass:primitive-data-comp}.

\item \emph{Wedge layer at every convex angle.}
Proposition~\ref{prop:wedge-wellposedness} proves existence, uniqueness in $L^\infty$, and weighted algebraic decay for the wedge layer. A shifted superharmonic barrier yields the weighted estimate for $1<\betaw<\frac{\pi}{\omega_m}$, and a growing barrier proves uniqueness.

\item \emph{Elliptic estimates near the vertices.}
Proposition~\ref{prop:mellin-mapping-comp} gives weighted vertex expansions for the harmonic functions used in the interior solution, following the elliptic theory of \cite{Kondratiev1967,Grisvard1985,Dauge1988,Kozlov.Mazya.Rossmann1997,Kozlov.Mazya.Rossmann2001}. These expansions contain powers determined by the boundary data and opening angle, together with logarithms at resonance. In a term $\e^k\bulkdist^\lambda$, $\bulkdist$ is the distance from the vertex, $k$ is the order of the transport expansion, and $\lambda$ is the radial exponent. We group the differentiated terms by $\alpha=k+\lambda$. This preserves the transport cancellation and gives the decay required for the wedge source. Appendix~\ref{sec:numerical-vertex-singularity} illustrates the possible singularity of $D^2\rho_0$, where $\rho_0$ is the leading harmonic interior solution.

\item \emph{$L^2$ diffusive limit on general polygons.}
Theorem~\ref{thm:vd-l2-wedge-free} gives $L^2$ convergence at rate $O(\e^{\frac{1}{2}})$ on every bounded simple polygon under Assumption~\ref{ass:vd-l2-data}. Reentrant vertices are allowed. The proof in Appendix~\ref{sec:velocity-dependent-l2} uses the leading side layer cut off near its endpoints. To estimate the velocity average, we combine weak identities with different test functions so that the contribution carrying the factor $\e^{-1}$ cancels, as in \cite{Guo.Wu2025}.
\end{enumerate}

\paragraph{\emph{Formal version of the main theorem}}
Here we first present a formal version of Theorem~\ref{thm:main-comp}. Let $\lambdaMin:=\min_m\frac{\pi}{\omega_m}$, choose $1<\betaw<\lambdaMin$, and take the matching exponent $\Tmatch$ as in \eqref{eq:depth-choices-comp}. Under the assumptions of the theorem, we construct $\Uapp^{\e,[\Nexp]}$ from the interior solution, side layers, and wedge layers through a fixed expansion order $\Nexp$, and prove
\begin{align}
 \bnm{\Uphys^\e-\Uapp^{\e,[\Nexp]}}_\infty\le C_q\e^q, \qquad 0<q<\frac{\Tmatch(\betaw-1)}{\Tmatch+\betaw+1}.
\end{align}

\subsection{Notation and Coordinate Systems}
\label{subsec:intro-coordinate-atlas}

Throughout the paper, $C>0$ and $c>0$ denote constants that may change from one inequality to another. They are independent of $\e$ and the overlap radius $\delta$, but may depend on the polygon, the fixed parameters, and the chosen cutoffs and extensions. We write $A\lesssim B$ if $A\le CB$, and $A\simeq B$ if both inequalities hold. Set $\jbr{z}:=\bigl(1+\abs{z}^2\bigr)^{\frac{1}{2}}$.

We use tangential and normal coordinates near each side to write the Milne equation. Near a vertex, we also need a common coordinate system for the two incident sides. The side layer uses the normal scaling $d=\e\eta$, while the wedge layer requires both spatial variables to be rescaled by $\e$.

The index $j$ denotes a side of the polygon. At $V_m$, the index $i\in\{1,2\}$ denotes an incident ray, and $j(m,i)$ is the side containing that ray. We first relate the side coordinate to distance from $V_m$, then introduce the rescaled wedge coordinates.

\paragraph{\emph{Side coordinates}} Let $L_j$ be the length of $E_j$, fix an orientation, and let $s_j\in[0,L_j]$ be arclength measured from its first endpoint. Denote the point with coordinate $s_j$ by $x_j(s_j)$, so $x_j\colon[0,L_j]\to\overline E_j$ is the arclength parametrization. Its derivative $t_j=x_j'(s_j)$ is the constant unit tangent, and $n_j$ is the outward unit normal. Extend the parametrization to the supporting line by $x_j(s)=x_j(0)+st_j$ for $s\in\mathbb R$. For a point $x$ near the side, define the normal coordinate $d_j(x):=-(x-x_j(0))\cdot n_j$, which is positive inside $\Om$. Extend the tangential coordinate off the side by
\begin{align}
 s_j(x):=(x-x_j(0))\cdot t_j.
 \label{eq:global-affine-side-coordinate}
\end{align}

The position is recovered from $(s_j,d_j)$ by
\begin{align}
 x=x_j(s_j)-d_jn_j.
 \label{eq:side-coordinate-map}
\end{align}
The inward normal and tangential velocity components are
\begin{align}
 \mu_j(w):=-w\cdot n_j, \qquad \tau_j(w):=w\cdot t_j.
\end{align}
On the in-flow boundary $\gamma_{j,-}$ we have $\mu_j>0$. By the chain rule,
\begin{align}
 w\cdot\nabla_x=\tau_j\partial_{s_j}+\mu_j\partial_{d_j}.
 \label{eq:side-sign-check}
\end{align}

\paragraph{\emph{Orientation at the endpoints}} On each side we have fixed an orientation. At a vertex, the local tangential distance increases away from that vertex, so one endpoint requires a reversal of orientation. For each $V_m$, let
\begin{align}
 \bulkdistat{m}(x)&:=\abs{x-V_m}, \qquad \sidedistat{m}{i}:=\text{distance from }V_m\text{ along incident ray }i,
 \label{eq:physical-distance-convention}
\end{align}
and $\scaledbulkat{m}:=\frac{\bulkdistat{m}}{\e}$. The radial distance $\bulkdistat{m}$ agrees with $\sidedistat{m}{i}$ on ray $i$. We omit $m$ and write $\bulkdist$, $\scaledbulk=\frac{\bulkdist}{\e}$, and $\sidedist{i}$ when working on a generic sector. The side coordinate is related to the distance from $V_m$ by
\begin{align}
 \varsigma_{m,i}:= \begin{cases}
 +1,&V_m=x_{j(m,i)}(0),\\
 -1,&V_m=x_{j(m,i)}(L_{j(m,i)}),
 \end{cases}
 \qquad s_{m,i}^{\rm ve}(\sidedistat{m}{i}):= \begin{cases}
 \sidedistat{m}{i},&\varsigma_{m,i}=+1,\\
 L_{j(m,i)}-\sidedistat{m}{i},&\varsigma_{m,i}=-1.
 \end{cases}
 \label{eq:local-global-side-coordinate}
\end{align}
The unit tangent directed away from $V_m$, the inward unit normal, and the corresponding velocity components are
\begin{align}
 \tphys{m}{i}&:=\varsigma_{m,i}t_{j(m,i)}, \qquad \nphys{m}{i}:=-n_{j(m,i)}, \\
 \tau_{m,i}(w)&:=w\cdot\tphys{m}{i} =\varsigma_{m,i}\tau_{j(m,i)}(w), \qquad \mu_{m,i}(w):=w\cdot\nphys{m}{i}=-w\cdot n_{j(m,i)}.
 \label{eq:local-ray-tangent-sign}
\end{align}
For $j=j(m,i)$, define the local normal coordinate by $d_{m,i}(x):=d_j(x)$. Since $V_m$ lies on the supporting line of $E_j$, it also equals $(x-V_m)\cdot\nphys{m}{i}$. The normal velocity component satisfies $\mu_{m,i}=\mu_j$. In the tangential direction, \eqref{eq:local-global-side-coordinate} gives $\partial_{\sidedistat{m}{i}}=\varsigma_{m,i}\partial_{s_j}$. Together with $\tau_{m,i}=\varsigma_{m,i}\tau_j$, this yields
\begin{align}
 \tau_{m,i}\partial_{\sidedistat{m}{i}} &=\tau_j\partial_{s_j}, \qquad j=j(m,i).
 \label{eq:tangential-operator-frame-invariance}
\end{align}

\paragraph{\emph{Coordinates near a vertex}} Fix $V_m$. Number the incident rays counterclockwise and take their polar angles to be $\theta_m=0$ and $\theta_m=\omega_m$. Thus $(\bulkdistat{m},\theta_m)$ are polar coordinates in the sector. Relative to ray $i$, we have
\begin{align}
 \signedsideat{m}{i}(x)&:=(x-V_m)\cdot\tphys{m}{i}, \qquad d_{m,i}(x)=(x-V_m)\cdot\nphys{m}{i},
 \label{eq:signed-endpoint-coordinate}\\
 x-V_m &=\signedsideat{m}{i}\,\tphys{m}{i} +d_{m,i}\,\nphys{m}{i}, \qquad \bulkdistat{m}^2=\signedsideat{m}{i}^{2}+d_{m,i}^2.
 \label{eq:bulk-signed-side-coordinate-relation}
\end{align}
The coordinate $\signedsideat{m}{i}$ is signed. On the forward ray it equals $\sidedistat{m}{i}\ge0$, while behind the vertex it is negative. The extension of the side layer is defined to be zero there.

Choose a fixed $\ellvtx>0$ such that the closed balls $\overline B(V_m,8\ellvtx)$ are disjoint and each $\Om\cap B(V_m,8\ellvtx)$ agrees with its tangent sector. We also require \eqref{eq:fixed-realization-separation}. The fixed vertex neighborhoods used for the elliptic estimates and side extensions have radii given by fixed multiples of $\ellvtx$. Within them, the vertex regularization uses radius $O(\e)$ and the wedge matching uses the shrinking overlap radius $\delta$, chosen below with $\e\ll\delta\ll\ellvtx$.

\paragraph{\emph{Rescaling near a side}} Set $\eta_j=\frac{d_j}{\e}$ and keep the tangential coordinate $s_j$. For $F(x,w)=B(s_j,\frac{d_j}{\e},w)$, \eqref{eq:side-sign-check} gives
\begin{align}
 \e\Le F=\e\tau_j\partial_{s_j}B+\mu_j\partial_{\eta_j}B+\qk B.
\end{align}
The leading equation is therefore the flat Milne equation with operator $\mu_j\partial_{\eta_j}+\qk$. The tangential derivative appears with the additional factor $\e$.

\paragraph{\emph{Rescaling near a vertex}} At distance $O(\e)$ from $V_m$, rescale both spatial coordinates. Choose a rotation $O_m\in SO(2)$ such that $O_mK_{\omega_m}+V_m$ is the tangent sector at $V_m$, where
\begin{align}
 K_{\omega_m}:=\big\{\varrho(\cos\theta,\sin\theta): \varrho>0,\ 0<\theta<\omega_m\big\}.
 \label{eq:wedge-sector-def}
\end{align}
Define the wedge-frame tangent and inward normal of ray $i$ by $\twedge{m}{i}:=O_m^T\tphys{m}{i}$ and $\nwedge{m}{i}:=O_m^T\nphys{m}{i}$.

For $Y\in K_{\omega_m}$, define the coordinates relative to ray $i$ by
\begin{align}
 \sigma_{m,i}(Y):=Y\cdot\twedge{m}{i}, \qquad \eta_{m,i}(Y):=Y\cdot\nwedge{m}{i}.
 \label{eq:wedge-side-coordinates}
\end{align}
After fixing $m$, write these as $\sigma_i,\eta_i$. The vertex rescaling is
\begin{align}
 Y_m:=\frac{O_m^T(x-V_m)}{\e}, \qquad v_m:=O_m^Tw, \qquad \scaledbulkat{m}=\abs{Y_m}=\frac{\bulkdistat{m}}{\e}.
 \label{eq:corner-scaling-full}
\end{align}

These definitions give $\sigma_i=\frac{\signedsideat{m}{i}}{\e}$, $\eta_i=\frac{d_{m,i}}{\e}$, and $\scaledbulkat{m}^2=\sigma_i^2+\eta_i^2$. In particular, $\sigma_i$ is signed, and $\sidedistat{m}{i}=\e\sigma_i$ on the forward ray. The rotated velocity components satisfy
\begin{align}
 \tau_i^Y(v):=v\cdot\twedge{m}{i}=\tau_{m,i}(O_mv), \qquad \mu_i^Y(v):=v\cdot\nwedge{m}{i}=\mu_{m,i}(O_mv).
 \label{eq:wedge-frame-velocity-components}
\end{align}
If $F(x,w)=\widehat F(Y_m,v_m)$, then
\begin{align}
 \e\Le F=v_m\cdot\nabla_{Y_m}\widehat F+\qk\widehat F=:\mathscr K_{\omega_m}\widehat F,
\end{align}
Both spatial derivatives now appear in the leading equation. As $\e\to0$, a fixed neighborhood of $V_m$ covers arbitrarily large subsets of $K_{\omega_m}$ in the $Y_m$ coordinates. We therefore solve a transport problem on the infinite wedge. The solution describes the interaction of the two side layers for bounded $\scaledbulkat{m}$ and must match the interior solution and side layers as $\scaledbulkat{m}\to\infty$.

\subsection{Main Difficulties and Methods}
\label{subsec:intro-difficulties}

We construct the approximation by solving harmonic Dirichlet problems for the interior solution, Milne problems for the side layers, and a transport problem at each vertex for the wedge layer. To estimate $\Uphys^\e-\Uapp^{\e,[\Nexp]}$ in $L^\infty(\Om\times\Sone)$, we must control the error throughout the regions $\bulkdistat{m}=O(\e)$. Their shrinking size gives no small factor in this norm. By Proposition~\ref{prop:physical-stability-comp}, the remainder is controlled by $\e^{-1}\bnm{\Le\Uapp^{\e,[\Nexp]}}_\infty$ and $\babs{g^\e-\Uapp^{\e,[\Nexp]}}_{\infty,-}$. Thus the approximation must satisfy both the equation and the in-flow boundary condition to the required accuracy near every vertex.

Two difficulties govern the construction near each vertex. Singular derivatives of the harmonic interior solution enter the higher-order interior solution and side layers. In addition, the two incident side layers interact within distance $O(\e)$ of the vertex. We regularize the finite vertex expansions and restrict the side extensions near the opposite ray. The wedge layer compensates the resulting errors in the transport equation and the in-flow boundary condition.

\paragraph{\emph{Singular derivatives of the interior solution}}
Write $\Uint_k$ for the interior solution multiplying $\e^k$ and $\rho_k:=\pk\Uint_k$ for its velocity average. The interior equations require $\Delta\rho_k=0$ in $\Om$. The boundary values of $\rho_k$ are determined together with the side layer $\Uside_{j,k}$ by a Milne problem at each fixed tangential coordinate $s\in(0,L_j)$.

We recall what the end-state map does. Consider the equation $\mu_j\partial_\eta f+\qk f=F$ with in-flow boundary condition $f(0,w)=h(w)$ on $\{\mu_j>0\}$. Under the boundedness and exponential-decay hypotheses on $(h,F)$ in Proposition~\ref{prop:flat-milne-comp}, its solution converges to a constant $E$ in $L^\infty_w$ as $\eta\to\infty$. The map $\mathcal E_{\mu_j}(h,F)$ returns this constant, and $\mathfrak M_{\mu_j}(h,F)$ returns the decaying part $f-E$, as defined in \eqref{eq:milne-maps-comp}. At leading order, $h$ is the leading in-flow data and $F=0$. Its end state is the function $e_{j,0}(s)$, which supplies $\rho_0\big|_{E_j}=e_{j,0}$, and the decaying part is $\Uside_{j,0}$. Thus $\Uint_0=\rho_0$ is a harmonic function whose boundary values already depend on the velocity-dependent transport data.

The functions to be expanded near $V_m$ are $\rho_k$, beginning with $\rho_0$. Their local expansions contain harmonic lifts of the prescribed ray data and homogeneous harmonic functions with zero boundary values on both rays. To determine the latter, solve $\Delta Z=0$ in $K_{\omega_m}$ with $Z(\bulkdistat{m},0)=Z(\bulkdistat{m},\omega_m)=0$ for $\bulkdistat{m}>0$. The separated form $Z=\bulkdistat{m}^{\lambda}\Phi(\theta)$ gives $\Phi''+\lambda^2\Phi=0$ and $\Phi(0)=\Phi(\omega_m)=0$. Hence the positive exponents and the two corresponding radial powers are
\begin{align}
 \lambda_{m,n}:=\frac{n\pi}{\omega_m}, \qquad \Phi_{m,n}(\theta):=\sin(\lambda_{m,n}\theta), \qquad Z_{m,n}^{\pm}(\bulkdistat{m},\theta):=\bulkdistat{m}^{\pm\lambda_{m,n}}\Phi_{m,n}(\theta), \quad n\ge1.
 \label{eq:intro-zero-ray-pairs-comp}
\end{align}
Here $\theta$ is the polar angle $\theta_m$ at the fixed vertex, and the superscripts $+$ and $-$ identify the positive and negative radial exponents. The compatible leading end states give an ordinary weak harmonic solution $\rho_0\in H^1(\Om)$. The local energy integrals are
\begin{align}
 \int_{K_{\omega_m}\cap B_1}|\nabla Z_{m,n}^{+}|^2\,\ud x &\simeq\int_0^1r^{2\lambda_{m,n}-1}\,\ud r<\infty,
 \label{eq:intro-zero-ray-energy-test-comp}\\
 \int_{K_{\omega_m}\cap B_1}|\nabla Z_{m,n}^{-}|^2\,\ud x &\simeq\int_0^1r^{-2\lambda_{m,n}-1}\,\ud r=\infty.
 \label{eq:intro-negative-ray-energy-test-comp}
\end{align}
Thus the $H^1$ requirement in the weak Dirichlet problem excludes $Z_{m,n}^{-}$. The terms proportional to $Z_{m,n}^{+}$ have finite energy, and their coefficients are determined by the global boundary data. Their exponents are usually noninteger. For noninteger $\lambda_{m,n}$, a nonzero Cartesian derivative of $Z_{m,n}^{+}$ of total order $a$ is homogeneous of degree $\lambda_{m,n}-a$, so sufficiently high derivatives are still singular. These derivatives enter the transport construction through $\qk\Uint_k=-w\cdot\nabla_x\Uint_{k-1}$ in \eqref{eq:interior-recursion-comp} and the side-layer source $-\tau_j\partial_s\Uside_{j,k-1}$ in \eqref{eq:side-profile-recursion-comp}. Consequently, existence of the weak Laplace solution leaves a further difficulty: the higher-order interior solution and side layers can be singular near the vertices, and we must control the derivatives of their vertex-expansion remainders as well.

At higher orders, the boundary data on the two rays can themselves contain singular powers. Their harmonic lifts may also contain logarithms when a boundary exponent is a root of the Dirichlet pencil. We prescribe these singular terms by the normalization in Subsection~\ref{subsec:harmonic-lifts-comp}, subtract their lifts, and solve a global Dirichlet problem for the remaining $H^1$ function. Definition~\ref{def:finite-part-solution-comp} specifies this construction for $\rho_k$. The energy calculation above excludes an additional free $Z_{m,n}^{-}$ term from the $H^1$ correction. The prescribed lifts retain the singular terms required by the boundary data.

\begin{remark}
\label{rem:zero-ray-mode}
The harmonic function $\bulkdistat{m}^{\frac{\pi}{\omega_m}}\sin(\frac{\pi\theta}{\omega_m})=Z_{m,1}^{+}$ has zero trace on both incident rays. Its coefficient in $\rho_0$ is determined by the global Dirichlet problem and can be nonzero even when the data vanish near that vertex. For an obtuse convex opening, $1<\frac{\pi}{\omega_m}<2$, this function has gradient in $L^\infty$ near the vertex and a Hessian of order $\bulkdistat{m}^{\frac{\pi}{\omega_m}-2}$. Appendix~\ref{sec:numerical-vertex-singularity} illustrates this mechanism.
\end{remark}

More precisely, apply Proposition~\ref{prop:mellin-mapping-comp} to $\rho_k$ with boundary data $e_k$, where $e_k\big|_{E_j}=e_{j,k}$. For a truncation exponent $\sigma>0$ satisfying the hypotheses of that proposition, let $H_{{\rm ray},m}$ be the sum of the normalized harmonic lifts of the retained ray data at this fixed order $k$. Then \eqref{eq:mellin-expansion-comp} reads
\begin{align}
 \rho_k(\bulkdistat{m},\theta) &= H_{{\rm ray},m}(\bulkdistat{m},\theta)+\sum_{0<\lambda_{m,n}<\sigma}c_{m,n}[e_k]Z_{m,n}^{+}(\bulkdistat{m},\theta)+\mathcal R_{\rho_k,m}(\bulkdistat{m},\theta).
 \label{eq:intro-harmonic-expansion}
\end{align}
Thus $\mathcal R_{\rho_k,m}$ is the remainder after subtracting both the harmonic lifts and the retained terms proportional to $Z_{m,n}^{+}$ from $\rho_k$. The weighted estimate controls its Cartesian derivatives $D_x^a\mathcal R_{\rho_k,m}$, for every order $a$ allowed by the data regularity. By \eqref{eq:weighted-sector-pointwise-comp}, their size is bounded by $C\bulkdistat{m}^{\sigma-a}$ times the logarithmic weight in \eqref{eq:mellin-remainder-bound-comp}. This derivative estimate is what allows us to pass from the vertex expansion of $\rho_k$ to that of the next interior solution.

The same estimates also control the data in the successive Milne problems. At order $k\ge1$, the preceding interior solution determines $\widehat{\Uint}_k:=-w\cdot\nabla_x\Uint_{k-1}$, which is the zero-average part of $\Uint_k$. The unshifted Milne problem has in-flow data $\gaux_{j,k}-\widehat{\Uint}_k\big|_{E_j}$ and source $-\tau_j\partial_s\Uside_{j,k-1}$. Here $\gaux_{j,k}$ is the auxiliary in-flow coefficient defined in \eqref{eq:auxiliary-inflow-coefficients}. The end-state map returns $e_{j,k}(s)$ by \eqref{eq:endstate-recursion-comp}. We solve the harmonic Dirichlet problem with $\rho_k\big|_{E_j}=e_{j,k}$, set $\Uint_k=\rho_k+\widehat{\Uint}_k$, and subtract $e_{j,k}(s)$ from the Milne solution to obtain $\Uside_{j,k}$. Its in-flow trace is then $\gaux_{j,k}-\Uint_k\big|_{E_j}$, so the interior solution and side layer satisfy the required boundary condition together. Section~\ref{sec:profiles} derives these equations and proves the Milne estimates. Section~\ref{sec:interior-side-hierarchy-comp} combines them with the elliptic estimates, and Proposition~\ref{prop:data-generated-hierarchy-comp} carries the construction through order $\Nexp$.

\paragraph{\emph{Interaction of the two side layers}}
The two side layers interact near each vertex. Write $Y=Y_m$ and $v=v_m$. Both normal variables $\eta_1,\eta_2$ can be $O(1)$ in a bounded region of the wedge. On ray $i$, at $Y=\sigma\twedge{m}{i}$ with $\sigma>0$, the coordinates relative to the other ray $i'=3-i$ are $\sigma_{i'}=\sigma\cos\omega_m$ and $\eta_{i'}=\sigma\sin\omega_m$. Where $\sigma_{i'}>0$, the extension of the other side layer can therefore have $O(1)$ $L^\infty_v$ norm for $\sigma=O(1)$. Its trace changes the in-flow boundary condition on the incident side within distance $O(\e)$ of the vertex. Exponential decay in the normal variable gives no small factor there, so this contribution must be removed or included in the estimate of $\babs{g^\e-\Uapp^{\e,[\Nexp]}}_{\infty,-}$. We use the narrow extensions in \eqref{eq:forward-sector-extension}. The cutoff for ray $i$ is one on that ray and zero near the other ray. Its transport derivative produces a source in the wedge, with a factor $O(\e^{-1})$ at distance $O(\e)$ from the vertex. The wedge construction must compensate this source while preserving the in-flow boundary condition on both rays.

\paragraph{\emph{Regularization and construction of the wedge layer}}
To obtain an approximation in $L^\infty$, we introduce a cutoff near each vertex to remove the singular terms from the interior solution and side layers. More precisely, we multiply their finite vertex expansions by $\chi_{\rm vtx}\bigl(\frac{\bulkdistat{m}}{\e}\bigr)$, with $\chi_{\rm vtx}$ defined in \eqref{eq:vertex-cutoff-choice}. This cutoff is zero for $\bulkdistat{m}\le\e$ and one for $\bulkdistat{m}\ge2\e$. We leave the construction remainders unchanged, since their weighted estimates give the required regularity. This modification is the regularization in \eqref{eq:Ubase-finite-expanded}. Differentiating the cutoff produces an additional transport source where the cutoff varies. The regularization also changes the in-flow boundary values within distance $2\e$ of the vertex. The wedge layer compensates these errors for the terms retained in the matching expansion, together with the source produced by the narrow side extensions.

The choice of the terms to be matched must respect the transport equations at consecutive orders. In \eqref{eq:interior-recursion-comp}, differentiating a term of size $\e^k\bulkdistat{m}^\lambda$ gives a term of size $\e^{k+1}\bulkdistat{m}^{\lambda-1}$ at the next order. Their ratio is $\frac{\e}{\bulkdistat{m}}$. Thus increasing the expansion order gives a smaller term when $\bulkdistat{m}\gg\e$, but gives the same power of $\e$ when $\bulkdistat{m}=\e\scaledbulkat{m}$. That power is determined by
\begin{align}
 \alpha:=k+\lambda.
 \label{eq:physical-exponent-definition-comp}
\end{align}
For a power-logarithmic term with $b\in\mathbb N_0$, the rescaling is explicit:
\begin{align}
 \e^k\bulkdistat{m}^\lambda\bigl(\ln \bulkdistat{m}\bigr)^b =\e^\alpha \scaledbulkat{m}^\lambda \sum_{c=0}^{b}\binom bc(\ln\e)^{b-c}\bigl(\ln \scaledbulkat{m}\bigr)^c.
 \label{eq:intro-power-log-rescaling}
\end{align}
For fixed $\scaledbulkat{m}>0$, all terms with the same $\alpha$ have the same power of $\e$. We retain these terms at every order $0\le k\le\Nexp$ so that their transport and collision terms cancel by \eqref{eq:interior-coefficient-equations-comp} and \eqref{eq:side-profile-recursion-comp}. Applying the wedge operator to this finite sum leaves the derivative of the order-$\Nexp$ term, with spatial degree $\alpha-\Nexp-1$. Increasing $\Nexp$ therefore improves the decay of the generated source as $\scaledbulkat{m}\to\infty$, with the same prefactor $\e^\alpha$.

We use two truncation exponents. To construct the interior solution and side layers, we retain the terms with $\alpha<\Tcon$ in their vertex expansions. We choose $\Tcon$ large enough that the remainders at order $\Nexp$ have the derivatives needed for the source estimates and endpoint Taylor estimates. For wedge matching, we retain all terms with $\alpha<\Tmatch$ through order $\Nexp$, where $\Tmatch<\Tcon$. The terms with $\Tmatch\le\alpha<\Tcon$ remain in $\Ubase^{\e,[\Nexp]}$, the sum of the interior solution and extended side layers regularized near the vertices, before adding the wedge layers. Its definition is \eqref{eq:Ubase-finite-expanded}. The powers of $\e$ in these unmatched terms give the required smallness after applying $\Le$. Thus $\Tcon$ determines the regularity of the remainders in the construction, and $\Tmatch$ determines the accuracy of matching. The index sets are given in \eqref{eq:construction-ledger-comp}--\eqref{eq:intermediate-ledger-comp}.

After rescaling, we group the vertex expansions by the factors $\e^\alpha(\ln\e)^b$. For each matching index $(m,\alpha,b)$, the corresponding terms from the interior solution and both incident side layers form the unregularized sum $\Cvertex_{m,\alpha,b}$ in \eqref{eq:C-coefficient-expanded}. Set $M_{m,\alpha,b}:=M^{\rm con}_{m,\alpha,b}$ using the vertex regularization in \eqref{eq:Mcon-coefficient-expanded}. The cutoff near $Y=0$ makes $M_{m,\alpha,b}\in L^\infty((K_{\omega_m}\cap B_2)\times\Sone)$, and $M_{m,\alpha,b}=\Cvertex_{m,\alpha,b}$ for $\abs Y\ge2$. Similarly, $G_{m,i,\alpha,b}$ is the coefficient of $\e^\alpha(\ln\e)^b$ in the rescaled endpoint Taylor expansion of the in-flow data, as defined in \eqref{eq:primitive-matching-coefficient-comp}. We solve \eqref{eq:defect-wedge-problem} for the wedge layer $D_{m,\alpha,b}$, with source $-\mathscr K_{\omega_m}M_{m,\alpha,b}$ and boundary value equal to $G_{m,i,\alpha,b}$ minus the in-flow trace of $M_{m,\alpha,b}$. Then $M_{m,\alpha,b}+D_{m,\alpha,b}$ solves the homogeneous transport equation on the wedge with in-flow boundary data $G_{m,i,\alpha,b}$. To match this solution with the retained interior and side expansions as $\abs Y\to\infty$, we require $D_{m,\alpha,b}$ to decay. Thus the same terms used for matching determine both the source and boundary data of the wedge layer.

We must verify that these source and boundary data decay at the rates required by the wedge theorem. Equations \eqref{eq:interior-coefficient-equations-comp} and \eqref{eq:side-profile-recursion-comp} cancel the transport and collision terms through order $\Nexp$. The remaining terms are derivatives of spatial degree $\alpha-\Nexp-1$ and derivatives of the cutoffs. The vertex cutoff produces a source supported near $Y=0$. On the transition region of a ray cutoff, the normal and radial variables are comparable, so the Milne decay estimates apply. On ray $i$, the trace of the interior solution and its side layer reproduce the retained in-flow data, and the extension from ray $3-i$ is zero. Hence \eqref{eq:completed-ray-mismatch-split} confines the boundary error to the vertex regularization region. Lemma~\ref{lem:localized-wedge-data} gives the resulting weighted source and trace estimates.

\paragraph{\emph{Decay of the wedge layer}}
For the wedge at $V_m$, write $\omega=\omega_m$ and $\scaledbulk=\abs Y$. The main difficulty is to estimate the velocity average $\rho_D:=\pk D$ on this unbounded domain. Averaging the characteristic formula gives $\rho_D=\Pw\rho_D+H$, where $H$ is the averaged contribution of the source and in-flow data. The positive operator $\Pw$ is convolution with a whole-space kernel, restricted to the sector. This kernel has unit mass and zero first moment. At a finite point $Y$, $\Pw \mathds 1_{K_\omega}(Y)<1$, but $\Pw \mathds 1_{K_\omega}(Y)\to1$ as the distance to both rays tends to infinity. Thus $\nm{\Pw}_{L^\infty\to L^\infty}=1$, and the equation does not give a uniform contraction in $L^\infty$. The second moment of the kernel gives a Laplacian in the large-distance expansion of $\Pw$. We use this diffusive behavior to construct weighted barriers.

Lemma~\ref{lem:sector-barrier} constructs a positive function $B$ satisfying
\begin{align}
 B(Y)&\simeq\jbr{\scaledbulk}^{-\betaw}, \qquad B(Y)-\Pw B(Y)\ge c\jbr{\scaledbulk}^{-\betaw-2}.
 \label{eq:intro-decaying-barrier-estimate}
\end{align}
The inequality $\betaw<\frac{\pi}{\omega_m}$ allows a positive superharmonic function on the sector. We move its singularity outside the wedge so that $B$ is finite and positive at the vertex. The estimate \eqref{eq:intro-decaying-barrier-estimate} holds both near the vertex and as $\scaledbulk\to\infty$. Positivity of $\Pw$ and this estimate bound the Neumann partial sums and give a solution. A growing barrier gives uniqueness in $L^\infty$. Proposition~\ref{prop:wedge-wellposedness} then yields a unique wedge layer $D_{m,\alpha,b}\in L^\infty(K_{\omega_m}\times\Sone)$ with $\nm{D_{m,\alpha,b}(Y,\cdot)}_{L^\infty_v}\le C\jbr{\abs Y}^{-\betaw}$ for almost every $Y$. This weighted estimate requires the source to decay two powers faster than the wedge layer, which determines $\sdata=\betaw+2$ and the lower bound on $\Nexp$ below.

\paragraph{\emph{Adding the wedge layers}}
With $\mathfrak J_m^{\rm mat}$ the matching set defined in \eqref{eq:finite-index-split-expanded}, form the finite sums
\begin{align}
 \mathcal M_m^\e &:=\sum_{(\alpha,b)\in\mathfrak J_m^{\rm mat}} \e^\alpha(\ln\e)^bM_{m,\alpha,b},
 \label{eq:intro-summed-reference-field}\\
 \Uwedge_m^\e &:=\sum_{(\alpha,b)\in\mathfrak J_m^{\rm mat}} \e^\alpha(\ln\e)^b\bigl(M_{m,\alpha,b}+D_{m,\alpha,b}\bigr).
 \label{eq:intro-summed-wedge-fields}
\end{align}
The difference $\mathcal D_m^\e:=\Uwedge_m^\e-\mathcal M_m^\e$ is the corresponding sum of wedge layers. In the fixed vertex neighborhood of $V_m$, $\Ubase^{\e,[\Nexp]}$ contains $\mathcal M_m^\e$ as part of its local expansion. Adding $\mathcal D_m^\e$ therefore gives
\begin{align}
 \Ubase^{\e,[\Nexp]}+\mathcal D_m^\e =\big(\Ubase^{\e,[\Nexp]}-\mathcal M_m^\e\big)+\Uwedge_m^\e.
 \label{eq:intro-wedge-inclusion-exclusion}
\end{align}
Here the wedge functions are evaluated at $(Y_m,v_m)$, and $\Ubase^{\e,[\Nexp]}$ at $(x,w)$. Applying the transport operator, the source $-\mathscr K_{\omega_m}M_{m,\alpha,b}$ in the equation for $D_{m,\alpha,b}$ cancels the source from the same $M_{m,\alpha,b}$ in $\Ubase^{\e,[\Nexp]}$. This cancellation is exact wherever the wedge cutoff equals one. The matched order-$\Nexp$ derivatives and ray-cutoff terms remain in $\bulkdistat{m}\ge\delta$. We also estimate the terms with $\Tmatch\le\alpha<\Tcon$, the remainders of the vertex expansions, and the commutators from the fixed side cutoffs and the wedge cutoff. The exact source decomposition is given in \eqref{eq:source-pieces-expanded}.

\paragraph{\emph{Cutoff and remainder estimate}}
We multiply $\mathcal D_m^\e$ by a cutoff at radius $\delta$ in the $x$ coordinates, where $\e\ll\delta\ll1$. The cutoff is one near the vertex and zero outside a neighborhood of radius $O(\delta)$. In its transition annulus, $\scaledbulkat{m}\simeq\frac{\delta}{\e}$, so the wedge decay estimate gives an $L^\infty$ norm of order $\bigl(\frac{\e}{\delta}\bigr)^\betaw$ for $\mathcal D_m^\e$, up to a logarithmic factor. Differentiation of the cutoff gives $\delta^{-1}$, and the stability estimate in Proposition~\ref{prop:physical-stability-comp} gives another factor $\e^{-1}$. The endpoint Taylor remainder in the in-flow boundary condition is estimated directly by its kinetic $L^\infty$ norm. Its size depends on the truncation $\alpha<\Tmatch$. Up to the common logarithmic factor, the two contributions that determine the choice of $\delta$ are
\begin{align}
 \delta^{\Tmatch}, \qquad \e^{-1}\delta^{-1}\left(\frac{\e}{\delta}\right)^\betaw =\e^{\betaw-1}\delta^{-\betaw-1}.
 \label{eq:intro-overlap-error-scales}
\end{align}
As $\delta$ decreases, the first contribution decreases and the second increases. We choose $\delta=\e^{\aoverlap}$ to equate the two powers of $\e$, which gives $\qrate=\Tmatch\aoverlap$ in \eqref{eq:overlap-rate-parameters-comp}. Proposition~\ref{prop:full-residual} also estimates the sources from differentiating $\Uint_\Nexp$ and $\Uside_{j,\Nexp}$, and from the terms with $\Tmatch\le\alpha<\Tcon$. After multiplication by $\e^{-1}$, the parameter inequalities make these contributions smaller than the two terms in \eqref{eq:intro-overlap-error-scales}. Applying the stability estimate to $\Uphys^\e-\Uapp^{\e,[\Nexp]}$ proves Theorem~\ref{thm:main-comp}.
\subsection{Main Results}
\label{subsec:main-result-comp}

\subsubsection{Parameter assumptions}
\label{subsec:parameters-data}

The integer $\Nexp$ is the order of the expansion in $\e$: we construct the interior solution $\Uint_k(x,w)$ and the side layer $\Uside_{j,k}(s,\eta,w)$ for $0\le k\le\Nexp$, and use the sums
\begin{align}
 \sum_{k=0}^{\Nexp}\e^k\Uint_k(x,w), \qquad \sum_{k=0}^{\Nexp}\e^k\Uside_{j,k}(s,\eta,w).
 \label{eq:intro-expansion-order}
\end{align}
The side layer is evaluated at $\eta=\frac{d_j(x)}{\e}$ and extended from its side with fixed cutoffs. These extensions and the regularization near the vertices are specified in Subsection~\ref{subsec:main-approximation-notation} below.

At $V_m$, write $\Uside_{m,i,k}(\sidedistat{m}{i},\eta,w):=\Uside_{j(m,i),k}\bigl(s_{m,i}^{\rm ve}(\sidedistat{m}{i}),\eta,w\bigr)$ for the side layer in the endpoint coordinates \eqref{eq:local-global-side-coordinate}. We also expand each of the interior solutions and side layers near the vertex. At order $k$, the truncation exponent $\Tcon$ retains the powers $\bulkdistat{m}^{\alpha-k}$ and $\sidedistat{m}{i}^{\alpha-k}$ with $\alpha<\Tcon$:
\begin{align}
 \Uint_k(\bulkdistat{m},\theta_m,w) &=\sum_{\alpha<\Tcon}\sum_b\bulkdistat{m}^{\alpha-k}\bigl(\ln\bulkdistat{m}\bigr)^b A^U_{m,k,\alpha,b}(\theta_m,w)+\mathcal R^{U,\Tcon}_{m,k},
 \label{eq:intro-interior-vertex-expansion}\\
 \Uside_{m,i,k}(\sidedistat{m}{i},\eta,w) &=\sum_{\alpha<\Tcon}\sum_b\sidedistat{m}{i}^{\alpha-k}\bigl(\ln\sidedistat{m}{i}\bigr)^b B^{\rm sl}_{m,i,k,\alpha,b}(\eta,w)+\mathcal R^{B,\Tcon}_{m,i,k}.
 \label{eq:intro-side-vertex-expansion}
\end{align}
Each sum runs over the finitely many pairs $(\alpha,b)$ generated by the interior and side equations, with $\alpha<\Tcon$ and $b\in\mathbb N_0$. The functions $A^U_{m,k,\alpha,b}$ and $B^{\rm sl}_{m,i,k,\alpha,b}$ are the corresponding coefficients, and each $\mathcal R$ is the difference between the function and its displayed finite sum. Section~\ref{sec:interior-side-hierarchy-comp} constructs these expansions in \eqref{eq:weighted-interior-decomposition-comp}--\eqref{eq:weighted-side-decomposition-comp}. The exponents include Taylor powers and the powers determined by the opening angles. Logarithms occur when the elliptic equations are resonant.

The integer $\pacc\ge1$ specifies the accuracy required of these remainders after the finite expansions have been subtracted. We choose $\Tcon$ large enough that, even at the last order $k=\Nexp$,
\begin{align}
 \bnm{D_x^q\mathcal R^{U,\Tcon}_{m,\Nexp}}_{L^\infty_w} &\le C\bulkdistat{m}^{\pacc-q}, \qquad 0\le q\le\pacc,
 \label{eq:intro-remainder-order-interior}\\
 \bnm{\partial_{\sidedistat{m}{i}}^q\mathcal R^{B,\Tcon}_{m,i,\Nexp}}_{L^\infty_{\eta,w}} &\le C\sidedistat{m}{i}^{\pacc-q}, \qquad 0\le q\le\pacc.
 \label{eq:intro-remainder-order-side}
\end{align}
Thus $\pacc$ controls both the vanishing order of the remainders and the number of derivatives estimated. Increasing $\pacc$ requires retaining more terms in \eqref{eq:intro-interior-vertex-expansion}--\eqref{eq:intro-side-vertex-expansion}. Corollary~\ref{cor:terminal-derivative-bounds} proves these bounds, including exponential decay in $\eta$ for the side remainder.

After rescaling $\bulkdistat{m}=\e\scaledbulkat{m}$, formula \eqref{eq:intro-power-log-rescaling} groups these terms by $\e^\alpha(\ln\e)^b$. Let $\mathfrak J_m^{\rm con}$ be the finite set of pairs obtained by rescaling all the retained interior and side terms with $\alpha<\Tcon$ through order $\Nexp$, as specified in \eqref{eq:finite-index-split-expanded}. The matching exponent $\Tmatch$ specifies which of them receive a wedge layer. The sum $\mathcal D_m^\e$ introduced after \eqref{eq:intro-summed-wedge-fields} is therefore
\begin{align}
 \mathcal D_m^\e(Y,v)&=\sum_{\substack{(\alpha,b)\in\mathfrak J_m^{\rm con}\\\alpha<\Tmatch}}\e^\alpha(\ln\e)^bD_{m,\alpha,b}(Y,v), \qquad \mathfrak J_m^{\rm mat}:=\big\{(\alpha,b)\in\mathfrak J_m^{\rm con}:\alpha<\Tmatch\big\}.
 \label{eq:intro-wedge-expansion}
\end{align}
The functions $D_{m,\alpha,b}$ are the wedge layers described in Subsection~\ref{subsec:main-approximation-notation}. We take $\Tmatch<\Tcon$. The exponent $\Tcon$ ensures the remainder estimates \eqref{eq:intro-remainder-order-interior}--\eqref{eq:intro-remainder-order-side}, while $\Tmatch$ determines the accuracy of the in-flow boundary condition after matching. The remaining parameter $\betaw$ is the decay exponent in $\nm{D_{m,\alpha,b}(Y,\cdot)}_{L^\infty_v}\le C\jbr{\abs Y}^{-\betaw}$.

Recall $\lambdaMin=\min_m\frac{\pi}{\omega_m}>1$. Fix the remainder order $\pacc\ge1$, choose the wedge decay exponent $\betaw$, and then choose the integer expansion order $\Nexp$ so that
\begin{align}
 1<\betaw<\lambdaMin, \qquad \sdata:=\betaw+2, \qquad \Nexp\ge \pacc+\sdata+4.
 \label{eq:N-choice-comp}
\end{align}
Here $\sdata$ is the exponent in the weighted norms of the wedge data.
With $\Nexp$ fixed, let
\begin{align}
 \mathscr E_\Nexp:=\mathbb N_0\cup \bigcup_m\bigcup_{k=0}^\Nexp \Big\{k+\frac{n\pi}{\omega_m}, k-\frac{n\pi}{\omega_m}:n\in\mathbb N\Big\}.
 \label{eq:exceptional-set-comp}
\end{align}
The set is locally finite. Choose the two truncation exponents
\begin{align}
 \Tmatch\in(\pacc+3,\pacc+4)\setminus\mathscr E_\Nexp, \qquad \Tcon\in(\Nexp+\pacc+1,\Nexp+\pacc+2)\setminus\mathscr E_\Nexp.
 \label{eq:depth-choices-comp}
\end{align}
We retain $\alpha<\Tcon$ to construct the interior solution and side layers, and $\alpha<\Tmatch$ to define the wedge data. The choices in \eqref{eq:depth-choices-comp} give $\Tcon>\Tmatch$ and $\Tcon-\Nexp>\pacc+1$. They also keep the Mellin weights away from the roots of the Dirichlet pencil. The positive gap between each truncation exponent and the finitely many retained powers absorbs their logarithmic factors.

The overlap and convergence exponents are
\begin{align}
 \aoverlap:=\frac{\betaw-1}{\betaw+1+\Tmatch}, \qquad \qrate:=\Tmatch \aoverlap =\frac{\Tmatch(\betaw-1)}{\betaw+1+\Tmatch}.
 \label{eq:overlap-rate-parameters-comp}
\end{align}
Thus $\delta=\e^{\aoverlap}$ balances $\delta^{\Tmatch}$ and $\e^{\betaw-1}\delta^{-\betaw-1}$. Finally, the required sidewise regularity is
\begin{align}
 \Kreg:=\Nexp+\pacc+7.
 \label{eq:Kreg-choice-comp}
\end{align}
With \eqref{eq:depth-choices-comp}, $\Kreg$ is two derivatives above the threshold $\max\big\{\lceil\Tcon\rceil+3,\pacc+\Nexp+1,\Nexp+4\big\}=\Nexp+\pacc+5$ in Definition~\ref{def:admissible-hierarchy-depth-comp}.

\begin{remark}
\label{rem:parameter-inequalities}
The expansion order satisfies
\begin{align}
 \Nexp\ge \pacc+\sdata+4>\Tmatch+\sdata.
 \label{eq:N-dominates-matching-decay}
\end{align}
For every matched $\alpha<\Tmatch$, differentiating the order-$\Nexp$ term in the rescaled vertex expansion gives a source of radial degree $\alpha-\Nexp-1<-\sdata-1$. The extra power absorbs its logarithms in Lemma~\ref{lem:localized-wedge-data}. The condition on $\Tcon$ gives the order-$\Nexp$ remainder bounds in Corollary~\ref{cor:terminal-derivative-bounds}.

The final rate comparison uses
\begin{align}
 \frac{\qrate+1+\aoverlap}{1-\aoverlap}=\betaw,
 \label{eq:N-rate-redundancy-comp}
\end{align}
so $(1-\aoverlap)\Nexp-1-\aoverlap>\qrate$ is equivalent to $\Nexp>\betaw$. Also $\lambdaMin\le3$ by the angle sum of a convex polygon. Consequently $\qrate<\betaw-1<2$, while $\Tmatch-2>\pacc+1\ge2$. These inequalities ensure that the order-$\Nexp$ source, the intermediate sources with $\Tmatch\le\alpha<\Tcon$, and the $O(\e^2)$ in-flow remainder contribute errors smaller than every asserted strict rate. All parameters are fixed independently of $\e$.
\end{remark}

\subsubsection{Data assumptions and notation}
\label{subsec:main-approximation-notation}

\begin{assumption}
\label{ass:primitive-data-comp}
\leavevmode
\begin{enumerate}[label=\textup{(\arabic*)},leftmargin=*,itemsep=0.35em]
\item \emph{In-flow data and sidewise regularity.}
In the side coordinates \eqref{eq:side-coordinate-map}, write $g_j^\e(s,w):=g^\e(x_j(s),w)$ for $\mu_j(w)>0$. On each open side, assume
\begin{align}
 g_j^\e(s,w)=g_{j,0}(s,w)+\e g_{j,1}(s,w) +g_{j,\mathrm{rem}}^\e(s,w), \qquad \nm{g_{j,\mathrm{rem}}^\e}_\infty\le C\e^2.
 \label{eq:data-expansion-comp}
\end{align}
The remainder norm is the essential supremum over $[0,L_j]\times\{w\in\Sone:\mu_j(w)>0\}$. For some $0<\gamma<1$, assume also that $g_{j,k}\in C^{\Kreg-k,\gamma}\bigl([0,L_j];L^\infty(\{\mu_j>0\})\bigr)$ for $k=0,1$, with these norms bounded by a fixed constant.

\item \emph{Compatibility of the leading Milne end states.}
At each vertex, use the coordinates oriented away from the vertex in \eqref{eq:local-global-side-coordinate} to write, for $k=0,1$,
\begin{align}
 g_{m,i,k}(\sidedistat{m}{i},w) &:=g_{j(m,i),k}\bigl(s_{m,i}^{\rm ve}(\sidedistat{m}{i}),w\bigr), \qquad i=1,2.
 \label{eq:endpoint-pullbacks}
\end{align}

Require the two leading Milne end states to have the same endpoint value at each vertex. With $\mathcal E_\mu$ defined in \eqref{eq:milne-maps-comp}, this condition is
\begin{align}
 \mathcal E_{\mu_{m,1}}(g_{m,1,0}(0,\cdot),0) =\mathcal E_{\mu_{m,2}}(g_{m,2,0}(0,\cdot),0).
 \label{eq:leading-compatibility-comp}
\end{align}
\end{enumerate}
\end{assumption}

Under \eqref{eq:leading-compatibility-comp}, the sidewise leading end states define a single boundary function
\begin{align}
 e_0(x_j(s)):=e_{j,0}(s) :=\mathcal E_{\mu_j}\bigl(g_{j,0}(s,\cdot),0\bigr), \qquad 0\le s\le L_j.
 \label{eq:leading-boundary-function}
\end{align}

\begin{remark}
\label{rem:leading-end-state-mismatch}
Condition \eqref{eq:leading-compatibility-comp} requires the leading Milne end states on adjacent sides to agree at their common vertex. The condition allows different velocity-dependent in-flow data on the two incident sides. Together with the sidewise regularity, this compatibility gives $e_0\in H^{\frac{1}{2}}(\partial\Om)$ and hence the weak harmonic solution $\rho_0\in H^1(\Om)$. Convexity further gives $\nabla\rho_0\in L^\infty(\Om)$ by Lemma~\ref{lem:leading-field-convex-corner}.
\end{remark}

To construct the interior solution and side layers through order $\Nexp$, we use the auxiliary in-flow coefficients
\begin{align}
 \gaux_{j,k}:= \begin{cases}
 g_{j,k},&0\le k\le1,\\
 0,&2\le k\le \Nexp.
 \end{cases}
 \label{eq:auxiliary-inflow-coefficients}
\end{align}
The equations generate the higher interior solutions and side layers from derivatives of the preceding ones, with $\gaux_{j,k}=0$ for $k\ge2$. Their sum satisfies the in-flow boundary condition with coefficient $\gaux_{j,k}$ at each order.

In the same endpoint coordinates, set
\begin{align}
 \gaux_{m,i,k}(\sidedistat{m}{i},w) &:=\gaux_{j(m,i),k}\bigl(s_{m,i}^{\rm ve}(\sidedistat{m}{i}),w\bigr), \qquad 0\le k\le\Nexp,\quad i=1,2.
 \label{eq:auxiliary-endpoint-pullbacks}
\end{align}

Set
\begin{align}
 \Llog&:=\Nexp+1, \qquad \deleps:=\e^{\aoverlap}, \qquad \Lameps:=\bigl(1+\abs{\ln\e}\bigr)^{\Llog}.
 \label{eq:delta-lambda-comp}
\end{align}

Let $\Uint_k$ and $\Uside_{j,k}$, $0\le k\le\Nexp$, be the interior solution and side layers constructed in Sections~\ref{sec:profiles}--\ref{sec:interior-side-hierarchy-comp}. Fix the harmonic lifts in Subsection~\ref{subsec:harmonic-lifts-comp} and use Definition~\ref{def:finite-part-solution-comp} for the higher $\rho_k$. Let $\Ext_j^\e$ be the side extension in \eqref{eq:physical-side-assembly}: it evaluates the normal variable at $\eta=\frac{d_j(x)}{\e}$ and multiplies the side layer by the fixed normal and endpoint cutoffs. The resulting sum is
\begin{align}
 \mathscr H_\Nexp^\e(x,w)&:=\sum_{k=0}^\Nexp\e^k\Uint_k(x,w)+\sum_j\Ext_j^\e\!\Big[\sum_{k=0}^\Nexp\e^k\Uside_{j,k}\Big](x,w).
 \label{eq:raw-hierarchy-expanded}
\end{align}

For $(\alpha,b)\in\mathfrak J_m^{\rm con}$, let $\Cvertex_{m,\alpha,b}(Y,v)$ be the coefficient of $\e^\alpha(\ln\e)^b$ obtained by rescaling the finite vertex expansions in \eqref{eq:raw-hierarchy-expanded}. Its exact formula is \eqref{eq:C-coefficient-expanded}. We regularize this coefficient by
\begin{align}
 M^{\rm con}_{m,\alpha,b}(Y,v):=\chi_{\rm vtx}(\abs Y)\Cvertex_{m,\alpha,b}(Y,v),
 \label{eq:intro-regularized-vertex-coefficient}
\end{align}
where the fixed smooth cutoff $\chi_{\rm vtx}$ is zero for $\abs Y\le1$ and one for $\abs Y\ge2$, as in \eqref{eq:vertex-cutoff-choice}. Define the regularized sum of the interior solution and side layers by
\begin{align}
 \Ubase^{\e,[\Nexp]}(x,w)&:=\mathscr H_\Nexp^\e(x,w)+\sum_m\sum_{(\alpha,b)\in\mathfrak J_m^{\rm con}}\e^\alpha(\ln\e)^b\bigl(M^{\rm con}_{m,\alpha,b}-\Cvertex_{m,\alpha,b}\bigr)(Y_m,v_m).
 \label{eq:intro-regularized-interior-side-sum}
\end{align}
The difference $M^{\rm con}_{m,\alpha,b}-\Cvertex_{m,\alpha,b}$ vanishes for $\abs Y\ge2$. We evaluate this difference in the corresponding vertex chart and extend it by zero outside that chart. Formula \eqref{eq:intro-regularized-interior-side-sum} therefore modifies the finite vertex expansions only where $\bulkdistat{m}<2\e$. Section~\ref{sec:matching-realization} gives the equivalent local formula \eqref{eq:Ubase-finite-expanded} and proves that the singular terms cancel, so $\Ubase^{\e,[\Nexp]}\in L^\infty(\Om\times\Sone)$ for each sufficiently small $\e>0$.

For each matching pair $(\alpha,b)\in\mathfrak J_m^{\rm mat}$, set $M_{m,\alpha,b}=M^{\rm con}_{m,\alpha,b}$ and let $D_{m,\alpha,b}$ be the unique $L^\infty$ mild solution of \eqref{eq:defect-wedge-problem}, with data given by \eqref{eq:generated-wedge-data}. Its source is $-\mathscr K_{\omega_m}M_{m,\alpha,b}$, and its in-flow boundary data are the retained endpoint Taylor data minus the trace of $M_{m,\alpha,b}$. These functions define the sum $\mathcal D_m^\e$ in \eqref{eq:intro-wedge-expansion}. Fix a smooth cutoff $\chi_{\rm ov}$ equal to one on $[0,1]$ and zero on $[2,\infty)$, and write $\chi_m^\delta(x)=\chi_{\rm ov}(\frac{\bulkdistat{m}}{\delta})$, as in \eqref{eq:corner-cutoff-physical}. The full approximation is
\begin{align}
 \Uapp^{\e,[\Nexp]}(\delta)(x,w)&:=\Ubase^{\e,[\Nexp]}(x,w)+\sum_m\chi_m^\delta(x)\mathcal D_m^\e(Y_m,v_m).
 \label{eq:forward-composite-definition}
\end{align}
Thus \eqref{eq:raw-hierarchy-expanded}, \eqref{eq:intro-regularized-interior-side-sum} and \eqref{eq:intro-wedge-expansion} specify all three parts of \eqref{eq:forward-composite-definition}. All cutoffs, extensions and length scales are fixed as in Section~\ref{sec:matching-realization}. Proposition~\ref{prop:profile-wellposedness-comp} proves existence and uniqueness for the functions in this approximation.

\subsubsection{Main theorem}

\begin{theorem}[$L^{\infty}$ diffusive limit]
\label{thm:main-comp}
Let $\Om$ be a bounded simple convex polygon satisfying \eqref{eq:strict-convexity}, and suppose the in-flow data satisfy Assumption~\ref{ass:primitive-data-comp}. Choose the parameters as in Subsection~\ref{subsec:parameters-data} and the approximation $\Uapp^{\e,[\Nexp]}(\deleps)$ as in Subsection~\ref{subsec:main-approximation-notation}. Let $\Uphys^\e\in L^\infty(\Om\times\Sone)$ be the unique mild solution of \eqref{eq:physical-model}. Then the following hold.
\begin{enumerate}[label=\textup{(\arabic*)},leftmargin=*,itemsep=0.35em]
\item There are $\e_0>0$ and $C>0$, depending only on the polygon, the fixed construction choices, and the data norms, such that for $0<\e\le\e_0$,
\begin{align}
 \bnm{\Uphys^\e-\Uapp^{\e,[\Nexp]}(\deleps)}_\infty \le C\Lameps\left[ \e^{\qrate}+\e^{(1-\aoverlap)\Nexp-1-\aoverlap}+\e^{\Tmatch-2}\right] +C\e^2.
 \label{eq:main-estimate-comp}
\end{align}
Therefore, for every $0<q<\qrate$ there is $C_q$ with
\begin{align}
 \bnm{\Uphys^\e-\Uapp^{\e,[\Nexp]}(\deleps)}_\infty \le C_q\e^q.
 \label{eq:strict-rate-comp}
\end{align}
\item The leading interior solution $\Uint_0=\rho_0$ is the ordinary weak solution of
\begin{align}
 \Delta\rho_0=0\quad\text{in }\Om, \qquad \rho_0\big|_{\partial\Om}=e_0.
 \label{eq:diffusion-problem-comp}
\end{align}
Here $e_0$ is the boundary function in \eqref{eq:leading-boundary-function}.
\end{enumerate}
\end{theorem}

\begin{remark}
\label{rem:main-error-terms}
The three terms in brackets in \eqref{eq:main-estimate-comp} come from the matching and wedge cutoff, the derivatives of $\Uint_\Nexp$ and $\Uside_{j,\Nexp}$, and the terms with $\Tmatch\le\alpha<\Tcon$ in $\Ubase^{\e,[\Nexp]}$, respectively. The final $C\e^2$ comes from the prescribed in-flow remainder. The inequalities in Subsection~\ref{subsec:parameters-data} make every term bounded by $C_q\e^q$ for $q<\qrate$. The strict inequality absorbs the logarithmic factor $\Lameps$.
\end{remark}

\begin{remark}
\label{rem:geometric-rate}
Since $\qrate<\betaw-1$ and $\betaw<\lambdaMin$, the rate supplied by this construction satisfies
\begin{align}
 \qrate<\lambdaMin-1=\min_m\frac{\pi}{\omega_m}-1.
 \label{eq:geometric-rate-cap}
\end{align}
This range follows from the wedge decay and cutoff estimates. It approaches zero as the largest opening angle approaches $\pi$. For fixed $\betaw$, $\qrate$ increases with $\Tmatch$ and tends to $\betaw-1$. Increasing $\pacc$ permits a larger $\Tmatch$, but requires a larger $\Nexp$ and more tangential derivatives of the data. Both contributions in \eqref{eq:intro-overlap-error-scales} remain in the remainder estimate, so an improvement in wedge decay must be combined with sufficient matching accuracy.
\end{remark}

\begin{corollary}[Convergence in the interior]
\label{cor:interior-convergence}
Under the assumptions of Theorem~\ref{thm:main-comp}, for every compact $K_0\Subset\Om$,
\begin{align}
 \nm{\Uphys^\e-\rho_0}_{L^\infty(K_0\times\Sone)}\rightarrow0 \quad \text{as }\e\to0.
 \label{eq:interior-limit-comp}
\end{align}
\end{corollary}

\begin{remark}
\label{rem:convexity-structural}
Convexity enters both the elliptic and kinetic analysis. The inequality $\frac{\pi}{\omega_m}>1$ gives the regularity of $\rho_0$ in Lemma~\ref{lem:leading-field-convex-corner} and permits the choice $1<\betaw<\frac{\pi}{\omega_m}$. Convexity also keeps every chord between points of a sector inside it, yielding the restricted convolution in \eqref{eq:restricted-convolution-operator}, and supports the positive shifted barriers of Subsection~\ref{subsec:shifted-barrier}. The Mellin mapping statement of Proposition~\ref{prop:mellin-mapping-comp} itself applies to every opening $0<\omega_m<2\pi$.
\end{remark}

\begin{remark}
\label{rem:leading-composite-suffices}
The leading-order approximation consists of $\rho_0$, the side layers $\Uside_{j,0}$ and the wedge layers $D_{m,0,0}$, with the same extensions and cutoffs. Its formula is \eqref{eq:leading-composite-def}. The estimate \eqref{eq:leading-composite-gap-bound} bounds the contribution of the remaining terms by $C\Lameps\e$. Combining it with \eqref{eq:strict-rate-comp} yields
\begin{align}
 \bnm{\Uphys^\e-\Uapp^{\e,[0]}(\deleps)}_\infty &\le C_q\e^q, \qquad 0<q<\min\{\qrate,1\}.
 \label{eq:leading-composite-rate}
\end{align}

\end{remark}

\begin{remark}[The evolutionary problem]
\label{rem:evolutionary-problem}
For the evolutionary equation $\e^2\partial_t\Uphys^\e+\e w\cdot\nabla_x\Uphys^\e+\qk\Uphys^\e=0$, the formal limit is $\partial_t\rho_0-\frac{1}{2}\Delta\rho_0=0$. The construction would require a time-dependent interior solution, side layers, and an initial layer on the scale $\tau=\frac{t}{\e^2}$, as in \cite{Wu2017,Wu2019}. The wedge estimates must also be extended away from the initial time. Suppose a time-dependent family of stationary wedge layers satisfies $\nm{\partial_tD(t,Y,\cdot)}_{L^\infty_v}=O(\scaledbulk^{-\betaw})$ for almost every $Y$ as $\scaledbulk=\abs Y\to\infty$. The correction $D_2$ at order $\e^2$ would then satisfy the formal wedge equation $\Kw D_2=-\partial_tD$. This source lacks the two extra decay powers required by Proposition~\ref{prop:wedge-wellposedness}. A parabolic wedge construction, matching of its velocity average with the interior solution, and stability estimates that use spatial localization are possible approaches.

The interaction of the initial layer with the side and wedge layers requires evolution problems on a half-space and a wedge. Their data meet at $\tau=0$ and on one or both in-flow boundaries. Grazing trajectories obstruct regularity, and the collision operator does not damp the velocity average. The resulting spatial diffusion may therefore prevent exponential decay in $\tau$. One would need compatible initial and in-flow traces, spatially weighted $L^\infty$ estimates, and decay and matching in both $(\tau,\eta)$ and $(\tau,Y)$, with constants independent of $\e$ in the common regions of validity. A possible method is to separate the velocity average from the zero-average part, use parabolic estimates for the former, and combine them with weighted transport estimates for the latter.
\end{remark}

\subsubsection{Secondary theorem}

We also prove the $L^2$ diffusive limit on general polygons, including polygons with reentrant vertices. The approximation uses $\rho_0$ and the leading side layer cut off near its endpoints. Exponential Milne decay gives an $O(\e^{\frac{1}{2}})$ bound for the side layer in $L^2$. The key step in the remainder estimate is the cancellation between the weak identities obtained with different test functions, as in \cite{Guo.Wu2025}.

For this subsection and Appendix~\ref{sec:velocity-dependent-l2}, assume
\begin{align}
 \Om &\text{ is a bounded simple polygon with } 0<\omega_m<2\pi,\quad\omega_m\ne\pi \qquad(1\leq m\leq\Nvtx).
 \label{eq:vd-l2-angle-assumption}
\end{align}
The sides $E_j$ are the maximal straight sides of $\Om$. Collinear segments are merged before imposing the sidewise regularity in Assumption~\ref{ass:vd-l2-data}.

Let $g^\e$ be the in-flow boundary data in \eqref{eq:physical-model}. Define $g_0(x_j(s),w):=g_{j,0}(s,w)$ for $0<s<L_j$ and $\mu_j(w)>0$. Under Assumption~\ref{ass:primitive-data-comp}, $g_{j,0}$ is the leading coefficient in the data expansion \eqref{eq:data-expansion-comp}. For the $L^2$ theorem, we impose the following assumptions instead.

\begin{assumption}
\label{ass:vd-l2-data}
\leavevmode
\begin{enumerate}[label=\textup{(\arabic*)},leftmargin=*,itemsep=0.35em]
\item \emph{In-flow data and sidewise regularity.}
On each side, assume
\begin{align}
 g_{j,0} &\in C^1\bigl([0,L_j]; L^\infty(\{w\in\Sone:\mu_j(w)>0\})\bigr).
 \label{eq:vd-l2-leading-regularity}
\end{align}
Assume also that the in-flow data $g^\e$ are uniformly bounded in $L^\infty(\gamma_-)$ and satisfy $\abs{g^\e-g_0}_{2,-}=O(\e^{\frac{1}{2}})$. Denote the combined bound by
\begin{align}
 \mathcal G_{\rm vd} &:={} \sum_j\nm{g_{j,0}}_{C^1_sL^\infty_w} + \sup_{0<\e\leq1} \abs{g^\e}_{\infty,-} + \sup_{0<\e\leq1} \e^{-\frac{1}{2}}\abs{g^\e-g_0}_{2,-} <\infty.
 \label{eq:vd-l2-data-norm}
\end{align}
Here $\nm{g_{j,0}}_{C^1_sL^\infty_w}$ uses the interval $[0,L_j]$ and the velocity set $\{\mu_j>0\}$ from \eqref{eq:vd-l2-leading-regularity}.

\item \emph{Compatibility of the leading Milne end states.}
Let $e_{j,0}$ be the leading Milne end state on $E_j$, defined by \eqref{eq:leading-boundary-function}. We require the two leading end states to have the same endpoint value at each vertex:
\begin{align}
 e_{j(m,1),0}\bigl(s_{m,1}^{\rm ve}(0)\bigr) &= e_{j(m,2),0}\bigl(s_{m,2}^{\rm ve}(0)\bigr).
 \label{eq:vd-l2-end-state-compatibility}
\end{align}
\end{enumerate}
\end{assumption}

Assumption~\ref{ass:primitive-data-comp} implies Assumption~\ref{ass:vd-l2-data}: the expansion $g_j^\e-g_{j,0}=\e g_{j,1}+g_{j,\mathrm{rem}}^\e$ gives $\abs{g^\e-g_0}_{2,-}\le C\e$ on the finite boundary, and it includes the required leading regularity and endpoint compatibility.

By \eqref{eq:vd-l2-leading-regularity} and boundedness of the end-state operator, the map $s\mapsto e_{j,0}(s)$ is $C^1$ on each closed side. Condition \eqref{eq:vd-l2-end-state-compatibility} therefore makes the piecewise function $e_0\big|_{E_j}:=e_{j,0}$ Lipschitz in boundary arclength. Hence $e_0\in H^1(\partial\Om)\subset H^{\frac{1}{2}}(\partial\Om)$.

Let $\rho_0\in H^1(\Om)$ be the unique weak harmonic function with trace $e_0$. A bounded right inverse of the Sobolev trace map, energy minimization, and the maximum principle give
\begin{align}
 \rho_0\big|_{\partial\Om} &=e_0, \qquad \int_\Om \nabla_x\rho_0\cdot\nabla_x\varphi\,\ud x =0 \quad \bigl(\varphi\in H^1_0(\Om)\bigr),
 \label{eq:vd-l2-harmonic-limit}
\end{align}
and
\begin{align}
 \nm{e_0}_{W^{1,\infty}(\partial\Om)} + \nm{\rho_0}_{H^1(\Om)} + \nm{\rho_0}_{L^\infty(\Om)} &\leq C_\Om\mathcal G_{\rm vd}.
 \label{eq:vd-l2-harmonic-estimate}
\end{align}
Thus the Milne end-state Dirichlet problem \eqref{eq:leading-variational-Dirichlet}, \eqref{eq:diffusion-problem-comp} also defines $\rho_0$ on a nonconvex polygon. For the leading side layers with endpoint cutoffs, write
\begin{align}
 \UsideZero^\e&:=\sum_j\UsideZeroAt{j}^\e.
 \label{eq:vd-l2-global-leading-layer}
\end{align}
Each $\UsideZeroAt{j}^\e$ is the leading Milne solution minus its end state, evaluated at the normal coordinate $\frac{d_j}{\e}$ and multiplied by the spatial cutoff in \eqref{eq:vd-l2-side-cutoff}. The precise definition is \eqref{eq:vd-l2-localized-side-layer}.

Write $\Err^\e=\Uphys^\e-\rho_0-\UsideZero^\e$, set $f^\e=\Uphys^\e-\rho_0$, and solve $-\Delta\phi^\e=\pk\Err^\e$ with $\phi^\e\in H^1_0(\Om)$. The test $w\cdot\nabla_x\phi^\e$ produces the desired $L^2$ norm of $\pk\Err^\e$, together with the term $\e^{-1}\int_{\Om\times\Sone}\qk f^\e\,w\cdot\nabla_x\phi^\e\,\ud x\ud\nu$. Adding $\e^{-1}$ times the identity obtained with $\phi^\e$ cancels this term exactly, as shown in \eqref{eq:vd-l2-first-test-cancellation} and \eqref{eq:vd-l2-punctured-second-test}. Appendix~\ref{sec:velocity-dependent-l2} justifies these tests at reentrant vertices and combines the resulting estimate with harmonicity of $\rho_0$, the Milne moment identities, and the energy estimate.

\begin{theorem}[$L^{2}$ diffusive limit]
\label{thm:vd-l2-wedge-free}
Let $\Om$ be a bounded simple polygon satisfying \eqref{eq:vd-l2-angle-assumption}, and suppose the in-flow data satisfy Assumption~\ref{ass:vd-l2-data}. Define $\rho_0$ and $\UsideZero^\e$ by \eqref{eq:vd-l2-harmonic-limit} and \eqref{eq:vd-l2-global-leading-layer}. Let $\Uphys^\e\in L^\infty(\Om\times\Sone)$ be the unique mild solution of \eqref{eq:physical-model}. For $0<\e\leq1$, the following estimates hold with a constant $C_\Om>0$ depending only on the polygon and the fixed parameters and cutoffs.
\begin{enumerate}[label=\textup{(\arabic*)},leftmargin=*,itemsep=0.35em]
\item The diffusion approximation satisfies
\begin{align}
 \nm{\Uphys^\e-\rho_0}_{L^2(\Om\times\Sone)} + \nm{\pk\Uphys^\e-\rho_0}_{L^2(\Om)} &\leq C_\Om\mathcal G_{\rm vd}\e^{\frac{1}{2}},
 \label{eq:vd-l2-main-rate}
\end{align}
\item The side layer and the remainder satisfy
\begin{align}
 \nm{\UsideZero^\e}_{L^2(\Om\times\Sone)}+\nm{\Uphys^\e-\rho_0-\UsideZero^\e} _{L^2(\Om\times\Sone)} &\leq C_\Om\mathcal G_{\rm vd}\e^{\frac{1}{2}},
 \label{eq:vd-l2-full-residual-rate}\\
 \nm{\qk(\Uphys^\e-\rho_0-\UsideZero^\e)} _{L^2(\Om\times\Sone)} &\leq C_\Om\mathcal G_{\rm vd}\e,
 \label{eq:vd-l2-residual-rates}\\
 \abs{\Uphys^\e-\rho_0-\UsideZero^\e}_{2,+} &\leq C_\Om\mathcal G_{\rm vd}\e^{\frac{1}{2}}.
 \label{eq:vd-l2-outgoing-residual-rate}
\end{align}
\end{enumerate}
\end{theorem}

\begin{remark}
\label{rem:vd-l2-nonprepared}
The leading in-flow data may have different velocity-dependent values at the two incident side endpoints. Endpoint compatibility \eqref{eq:vd-l2-end-state-compatibility} gives $e_0\in H^{\frac{1}{2}}(\partial\Om)$ and $\rho_0\in H^1(\Om)$. The remainder after subtracting the side layer is $\Uphys^\e-\rho_0-\UsideZero^\e$. Its zero-average part is $O(\e)$ in $L^2(\Om\times\Sone)$, and both its full $L^2$ norm and its out-flow kinetic $L^2$ norm are $O(\e^{\frac{1}{2}})$. For each fixed $\e>0$, the mild solution and this remainder belong to the kinetic graph space \eqref{eq:vd-l2-graph-space}. This regularity justifies the energy identity and boundary traces used in Appendix~\ref{sec:velocity-dependent-l2}.
\end{remark}

\subsection{Organization of the Paper}
\label{subsec:organization}

Section~\ref{sec:physical-stability} proves the characteristic formulation, comparison principle, and stability estimate for $\Le$. Section~\ref{sec:profiles} derives the equations for the interior solution and side layers and proves the Milne estimates. Section~\ref{sec:polygonal-profiles-comp} gives the normalized vertex expansions and weighted remainder estimates for the harmonic functions. Using these results, Section~\ref{sec:interior-side-hierarchy-comp} constructs $\rho_k$, $\Uint_k$, and $\Uside_{j,k}$ for $0\le k\le\Nexp$, with estimates for their vertex expansions and derivatives.

Section~\ref{sec:wedge-layer} proves existence, uniqueness, and decay for the wedge problem with a prescribed source and in-flow boundary data. Section~\ref{sec:matching-realization} regularizes the vertex expansions of the interior solution and side layers, verifies the required estimates on the wedge data, and defines $\Uapp^{\e,[\Nexp]}$. The same $M_{m,\alpha,b}$ occurs in this approximation and in the wedge equation, giving the cancellation described above. Section~\ref{sec:main-proof} estimates the remaining source and boundary errors and applies the stability estimate to prove Theorem~\ref{thm:main-comp} and Corollary~\ref{cor:interior-convergence}.

Appendix~\ref{sec:velocity-dependent-l2} proves the $L^2$ theorem on general polygons using only the leading Milne theory and harmonic energy estimates. Appendix~\ref{sec:numerical-vertex-singularity} gives a numerical example of the possible singularity of $D^2\rho_0$ at a vertex.

\subsection*{Acknowledgements}
Lei Wu was supported in part by the National Science Foundation under grant DMS-2405161.

\section{Transport Well-Posedness and Stability}
\label{sec:physical-stability}

The remainder problem \eqref{eq:error-equation}--\eqref{eq:error-incoming-data} has the original transport operator with a source term and an error in the in-flow boundary condition. We prove well-posedness and an $L^\infty(\Om\times\Sone)$ estimate with a factor $\e^{-1}$ multiplying the source norm and a factor one multiplying the $L^\infty(\gamma_-)$ norm of the in-flow boundary data. The source must therefore be one power of $\e$ smaller than the required remainder bound. Remark~\ref{rem:l2-linfty-comparison} compares this normalization with the $L^2$--$L^\infty$ framework of \cite{Guo2010,Esposito.Guo.Kim.Marra2013,Esposito.Guo.Kim.Marra2015} and the neutron-transport estimates in \cite{Wu.Guo2015,Guo.Wu2017,Wu2020}.

Positivity of $\pk$ and preservation of constants give a comparison principle from the characteristic formula. We then construct a nonnegative barrier satisfying $\Le\Psi^\e=1$ and $\nm{\Psi^\e}_\infty=O(\e^{-1})$. A multiple of this barrier controls the source, while constant comparison functions control the in-flow boundary data. Remark~\ref{rem:physical-stability-sharp-comp} proves that the factor $\e^{-1}$ cannot be improved. The argument applies to every bounded polygon, including nonconvex polygons.

\subsection{Mild Formulation and Comparison Principle}

We use backward characteristics. For $(x,w)\in\Om\times\Sone$ define the backward exit time and the corresponding boundary point by
\begin{align}
 t_b(x,w)&:=\inf\{t>0:x-tw\notin\Om\},\qquad x_b(x,w):=x-t_b(x,w)w.
 \label{eq:physical-backward-exit}
\end{align}
Because $\Om$ is bounded, $0<t_b(x,w)\le D_\Om:=\operatorname{diam}\Om$ for almost every $(x,w)$.

Consider the transport equation with source $S$ and in-flow boundary data $h$,
\begin{align}
 \Le f&=S\quad\text{in }\Om\times\Sone,
 \label{eq:forced-physical-comp} \\
 f\big|_{\gamma_-}&=h.
 \label{eq:forced-physical-incoming-comp}
\end{align}
Restricting \eqref{eq:forced-physical-comp} to the backward characteristic $x-tw$ gives an ordinary differential equation. Multiplying by the integrating factor $\ue^{-\frac t\e}$ and integrating from $0$ to $t_b(x,w)$, we obtain
\begin{align}
 f(x,w) ={}&\ue^{-\frac{t_b(x,w)}{\e}}h\big(x_b(x,w),w\big)+\int_0^{t_b(x,w)}\ue^{-\frac{t}{\e}} \Big( \e^{-1}\pk f(x-tw)+S(x-tw,w) \Big)\ud t.
 \label{eq:first-event-physical-comp}
\end{align}

We use the phase-space norm in \eqref{eq:physical-infinity-norms} and the kinetic boundary measure and norms in \eqref{eq:l2-linfty-boundary-norms-comp}.

For $S\in L^\infty(\Om\times\Sone)$ and $h\in L^\infty(\gamma_-)$, we call $f\in L^\infty(\Om\times\Sone)$ an \emph{$L^\infty$ mild solution} of \eqref{eq:forced-physical-comp}--\eqref{eq:forced-physical-incoming-comp} if \eqref{eq:first-event-physical-comp} holds for almost every $(x,w)$. Almost every backward ray meets $\gamma_-$ at a non-vertex, non-grazing point, so the boundary term is well defined with respect to the kinetic boundary measure. We may use measurable representatives in this formula. Changes on null sets preserve the identity as an equality of $L^\infty$ classes.

Every backward characteristic has length at most $D_\Om$, so the coefficient of the in-flow boundary data in \eqref{eq:first-event-physical-comp} is at least $\ue^{-\frac{D_\Om}{\e}}>0$. Since the averaging measure is positive and normalized, velocity averaging cannot lower the essential infimum. These facts give the comparison principle below.

\begin{lemma}
\label{lem:physical-comparison-comp}
Let $g\in L^\infty(\Om\times\Sone)$ be a real-valued mild solution of $\Le g=F$ with $g\big|_{\gamma_-}=k$. If $F\ge0$ almost everywhere and $k\ge0$ almost everywhere on $\gamma_-$, then $g\ge0$ almost everywhere in $\Om\times\Sone$. Consequently, if $g_1,g_2$ are real-valued $L^\infty$ mild solutions with sources $F_1\le F_2$ and in-flow boundary data $k_1\le k_2$ almost everywhere, then $g_1\le g_2$ almost everywhere.
\end{lemma}

\begin{proof}
Set $m_g:=\essinf_{\Om\times\Sone}g$. Since $\ud\nu$ is positive and normalized, $\pk g(x)\ge m_g$ for almost every $x$, so \eqref{eq:first-event-physical-comp} and $F,k\ge0$ imply
\begin{align}
 g(x,w) &\ge m_g\int_0^{t_b(x,w)} \e^{-1}\ue^{-\frac{t}{\e}}\ud t=m_g\left(1-\ue^{-\frac{t_b(x,w)}{\e}}\right).
\end{align}

Suppose that $m_g<0$. Since $t_b(x,w)\le D_\Om$, multiplication by the negative number $m_g$ gives $g(x,w)\ge m_g\bigl(1-\ue^{-\frac{D_\Om}{\e}}\bigr)$ almost everywhere. Taking the essential infimum gives $m_g\ge m_g\bigl(1-\ue^{-\frac{D_\Om}{\e}}\bigr)>m_g$, a contradiction. Hence $m_g\ge0$. By linearity of \eqref{eq:first-event-physical-comp}, $g_2-g_1$ is an $L^\infty$ mild solution with source $F_2-F_1$ and in-flow boundary data $k_2-k_1$, both nonnegative. Applying the first assertion to this difference proves the second.
\end{proof}

\begin{lemma}
\label{lem:physical-wellposed-comp}
Let $S\in L^\infty(\Om\times\Sone)$ and $h\in L^\infty(\gamma_-)$. The problem \eqref{eq:forced-physical-comp}--\eqref{eq:forced-physical-incoming-comp} has a unique mild solution in $L^\infty(\Om\times\Sone)$, defined by \eqref{eq:first-event-physical-comp}.
\end{lemma}

\begin{proof}
Let $\mathcal T_{\rm char}f$ denote the right-hand side of \eqref{eq:first-event-physical-comp}, regarded as a map of $L^\infty(\Om\times\Sone)$ into itself. Positivity and normalization of $\pk$ give $\nm{\pk f_1-\pk f_2}_\infty\le\nm{f_1-f_2}_\infty$, so, using $t_b\le D_\Om$,
\begin{align}
 \bnm{\mathcal T_{\rm char} f_1-\mathcal T_{\rm char} f_2}_\infty &\le \esssup_{(x,w)} \int_0^{t_b(x,w)} \e^{-1}\ue^{-\frac{t}{\e}}\ud t\, \nm{f_1-f_2}_\infty \le \left(1-\ue^{-\frac{D_\Om}{\e}}\right) \nm{f_1-f_2}_\infty.
\end{align}
The coefficient is strictly smaller than one for each fixed $\e>0$, so Banach's fixed-point theorem gives a unique fixed point of $\mathcal T_{\rm char}$ in $L^\infty(\Om\times\Sone)$. It is the required mild solution.

For almost every non-grazing, non-vertex incoming pair $(y,w)\in\gamma_-$, let $t_+(y,w)$ be the length of the connected forward characteristic segment starting at $y$. For $0<t<t_+(y,w)$, choose the representative
\begin{align}
 f^\sharp(y+tw,w)&=\ue^{-\frac t\e}h(y,w)+\int_0^t\ue^{-\frac{t-s}{\e}}\bigl[\e^{-1}\pk f(y+sw)+S(y+sw,w)\bigr]\,\ud s.
 \label{eq:physical-characteristic-representative}
\end{align}
The fixed-point identity and Fubini's theorem in the characteristic coordinates, whose Jacobian is $\abs{w\cdot n(y)}$, show that $f^\sharp=f$ almost everywhere. This representative is absolutely continuous in $t$, tends to $h(y,w)$ as $t\downarrow0$, and satisfies $\partial_t f^\sharp=S-\e^{-1}\qk f^\sharp$ almost everywhere. It extends to the outgoing endpoint as well. Continuity along each such segment transfers the bulk essential bound to both endpoint traces, giving $\abs{f}_{\infty,\pm}\le\nm f_\infty$. Integrating the one-dimensional identity against test functions gives $\Le f=S$ in distributions. All trace assertions hold outside a set of zero kinetic boundary measure.
\end{proof}

The contraction constant above approaches one as $\e\to0$, so the fixed-point estimate contains the factor $\ue^{\frac{D_\Om}{\e}}$. To obtain the required factor $\e^{-1}$, we construct a barrier whose $L^\infty$ norm and image under $\Le$ can be computed explicitly.

\subsection{Barrier Argument and \texorpdfstring{$L^\infty$}{L-infinity} Estimate}

The barrier combines a quadratic function of $x$ with a function linear in $w$. The collision term of the latter cancels the streaming derivative of the quadratic function. Choose $x_0\in\mathbb R^2$ and $\ellOm>0$ so that $\overline\Om\subset B(x_0,\ellOm)$, and write $z=x-x_0$. We define
\begin{align}
 \Psi^\e(x,w) :={}& \frac{\ellOm^2-\abs{z}^2}{2\e}+z\cdot w+\ellOm.
 \label{eq:physical-lyapunov-comp}
\end{align}
The constant $\ellOm$ makes $\Psi^\e$ nonnegative. Since $\abs{z}\le \ellOm$ and $\abs{z\cdot w}\le \ellOm$, we have
\begin{align}
 0\le\Psi^\e(x,w) \le{}&\frac{\ellOm^2}{2\e}+2\ellOm \le C_\Om\e^{-1}, \qquad 0<\e\le1.
 \label{eq:physical-lyapunov-size-comp}
\end{align}

The quadratic term and the constant $\ellOm$ are velocity independent, so $\qk$ annihilates them. Also $\pk(z\cdot w)=z\cdot\int_{\Sone}w\,\ud\nu(w)=0$, so that $\qk(z\cdot w)=z\cdot w$.

Using $\abs w=1$, we compute
\begin{align}
 w\cdot\nabla_x \Big(\frac{\ellOm^2-\abs{z}^2}{2\e}\Big)&=-\e^{-1}z\cdot w,\qquad
 w\cdot\nabla_x(z\cdot w)=1,\qquad
 \e^{-1}\qk\Psi^\e=\e^{-1}z\cdot w.
\end{align}
The first and third terms cancel. Hence
\begin{align}
 \Le\Psi^\e&=w\cdot\nabla_x\Psi^\e+\e^{-1}\qk\Psi^\e=1.
 \label{eq:physical-lyapunov-identity-comp}
\end{align}

\begin{proposition}[Stability estimate]
\label{prop:physical-stability-comp}
Let $0<\e\le1$, and let $S\in L^\infty(\Om\times\Sone)$ and $h\in L^\infty(\gamma_-)$ be real-valued. The mild solution $f\in L^\infty(\Om\times\Sone)$ of \eqref{eq:forced-physical-comp}--\eqref{eq:forced-physical-incoming-comp} satisfies
\begin{align}
 \nm f_\infty \le \abs h_{\infty,-} +C_\Om\e^{-1}\nm S_\infty.
 \label{eq:physical-stability-bound-comp}
\end{align}
The constant depends only on $\Om$, and not on $\e$, $S$, or $h$.
\end{proposition}

\begin{proof}
Define the comparison functions $v_\pm:=\abs h_{\infty,-}+\nm S_\infty\Psi^\e\mp f$.

Integrating the classical identity \eqref{eq:physical-lyapunov-identity-comp} along a backward characteristic with the integrating factor $\ue^{-\frac{t}{\e}}$ gives \eqref{eq:first-event-physical-comp} for $\Psi^\e$, with source $1$ and in-flow boundary data $\Psi^\e\big|_{\gamma_-}$. Constants solve the same mild equation with zero source and constant in-flow boundary data. Linearity then shows that $v_\pm$ are $L^\infty$ mild solutions.

The identity \eqref{eq:physical-lyapunov-identity-comp} gives $\Le v_\pm=\nm S_\infty\mp S\ge0$. On the in-flow boundary, nonnegativity of $\Psi^\e$ gives $v_\pm\big|_{\gamma_-}\ge\abs h_{\infty,-}\mp h\ge0$.

Lemma~\ref{lem:physical-comparison-comp} therefore gives $v_\pm\ge0$ almost everywhere, hence $\abs{f(x,w)}\le\abs h_{\infty,-}+\nm S_\infty\Psi^\e(x,w)$ for almost every $(x,w)$. Taking the essential supremum and using \eqref{eq:physical-lyapunov-size-comp} proves \eqref{eq:physical-stability-bound-comp}.
\end{proof}

\begin{remark}
\label{rem:physical-stability-sharp-comp}
The factor $\e^{-1}$ cannot be improved for general sources in $L^\infty(\Om\times\Sone)$. Choose a ball $\overline{B(x_c,\ell)}\subset\Om$ and let $f^\e$ solve $\Le f^\e=1$ with $f^\e\big|_{\gamma_-}=0$, so that $f^\e\ge0$ almost everywhere by Lemma~\ref{lem:physical-comparison-comp}.

The lower bound comes from a barrier for the ball. In $B(x_c,\ell)$ set $\Psi_B^\e(x,w):=\frac{\ell^2-\abs{x-x_c}^2}{2\e}+(x-x_c)\cdot w$. The computation giving \eqref{eq:physical-lyapunov-identity-comp} yields $\Le\Psi_B^\e=1$. Write $n_B$ for the outward unit normal and $\gamma_{B,-}:=\{(y,w)\in\partial B(x_c,\ell)\times\Sone:w\cdot n_B(y)<0\}$ for the in-flow boundary of the ball. There the quadratic term vanishes, so $\Psi_B^\e\big|_{\gamma_{B,-}}(y,w)=\ell\,w\cdot n_B(y)\le0$.

To identify the in-flow trace on $\gamma_{B,-}$, let $t_B(x,w)$ be the backward exit time from the ball and set $y_B=x-t_B(x,w)w$. Splitting \eqref{eq:first-event-physical-comp} at $t_B$ gives the remaining integral as $\ue^{-\frac{t_B}{\e}}f^\e\big|_{\gamma_{B,-}}(y_B,w)$, where
\begin{align}
 f^\e\big|_{\gamma_{B,-}}(y_B,w)&=\int_0^{t_b(y_B,w)}\ue^{-\frac{t}{\e}}\bigl(\e^{-1}\pk f^\e(y_B-tw)+1\bigr)\ud t.
\end{align}
For almost every ball chord, this expression is the characteristic trace of $f^\e$ on $\gamma_{B,-}$. It belongs to the kinetic $L^\infty$ space on $\gamma_{B,-}$ and is nonnegative almost everywhere because $\pk f^\e\ge0$. Thus $f^\e-\Psi_B^\e$ has zero source and nonnegative in-flow boundary data on the ball. The proof of Lemma~\ref{lem:physical-comparison-comp} uses only boundedness of the domain, so it applies to the ball and gives $f^\e\ge\Psi_B^\e$ almost everywhere there.

For $\abs{x-x_c}\le \frac{\ell}{2}$ and $0<\e\le \frac{\ell}{4}$, we have $\Psi_B^\e(x,w)\ge\frac{3\ell^2}{8\e}-\frac{\ell}{2}\ge\frac{\ell^2}{4\e}$. Thus $\nm{f^\e}_\infty\ge\frac{\ell^2}{4\e}$. With Proposition~\ref{prop:physical-stability-comp}, this proves $\nm{\Le^{-1}}_{L^\infty\to L^\infty}\simeq\e^{-1}$ for the inverse with zero in-flow trace as $\e\to0$.
\end{remark}

\begin{remark}
\label{rem:l2-linfty-comparison}
In the normalization $\e w\cdot\nabla_x f+\qk f=F$ with in-flow boundary data $h$, one has $F=\e S$, so \eqref{eq:physical-stability-bound-comp} becomes
\begin{align}
 \nm f_\infty \leq \abs h_{\infty,-}+C_\Om\e^{-2}\nm F_\infty.
\end{align}
Remark~\ref{rem:physical-stability-sharp-comp}, applied with $S=\e^{-1}F$, shows that this power is sharp for a general $L^\infty$ source in that normalization.

For comparison, the energy estimate and double iteration of the characteristic formula in the disk setting of \cite[Section~2]{Wu.Guo2015}, written with $F=\e S$ and the boundary norms \eqref{eq:l2-linfty-boundary-norms-comp}, give
\begin{align}
 \nm f_\infty \leq C_\Om\left( \e^{-2}\nm S_{L^2(\Om\times\Sone)} +\e\nm S_\infty +\e^{-\frac32}\abs h_{2,-} +\abs h_{\infty,-} \right).
\end{align}
Using $\nm S_{L^2(\Om\times\Sone)}\le\abs\Om^{\frac12}\nm S_\infty$ in this bound gives the power $\e^{-2}$. The $L^{2m}$--$L^\infty$ refinement in the smooth convex-domain analyses \cite{Guo.Wu2017,Wu2020} gives the source power $\e^{-1-\frac{1}{m}}$ for fixed $m>2$. Proposition~\ref{prop:physical-stability-comp} gives the power $\e^{-1}$ directly on bounded polygons.

\end{remark}

\subsection{From the Distributional Equation to the Mild Formulation}

The side layers and wedge layers in \eqref{eq:forward-composite-definition} solve their equations in the mild sense. We will first verify the equation for the approximate solution in distributions. The following lemma then allows us to apply Proposition~\ref{prop:physical-stability-comp} to its remainder.

\begin{lemma}
\label{lem:admissible-mild-comp}
Let $U,S\in L^\infty(\Om\times\Sone)$ be real-valued and suppose that
\begin{align}
 w\cdot\nabla_xU&=S-\e^{-1}\qk U \quad\text{in distributions on }\Om\times\Sone.
 \label{eq:admissible-transport-identity}
\end{align}
Then the following statements hold.
\begin{enumerate}[label=\textup{(\roman*)},leftmargin=*]
\item $U$ has a representative that is absolutely continuous along almost every connected characteristic segment. Its in-flow endpoint limits define $U\big|_{\gamma_-}\in L^\infty(\gamma_-)$ with $\abs U_{\infty,-}\le\nm U_\infty$.
\item $U$ is an $L^\infty$ mild solution of \eqref{eq:forced-physical-comp}--\eqref{eq:forced-physical-incoming-comp} with source $S$ and in-flow boundary data $U\big|_{\gamma_-}$.
\end{enumerate}
\end{lemma}

\begin{proof}
Slice \eqref{eq:admissible-transport-identity} first in velocity and then along lines parallel to that fixed velocity. Fubini's theorem and the one-dimensional characterization of weak derivatives give a $W^{1,\infty}$ representative on almost every connected segment, with derivative bounded by $\nm S_\infty+2\e^{-1}\nm U_\infty$. This representative is Lipschitz and extends to both endpoints of the finite segment.

For a non-vertex point $(y,w)\in\gamma_-$, set $t_+(y,w):=\sup\{t>0:y+\tau w\in\Om\text{ for every }0<\tau<t\}$. This is the length of the connected forward chord starting at $y$. The characteristic coordinates are $x=y+\tau w$ with $0<\tau<t_+(y,w)$, and their Jacobian is $\abs{w\cdot n(y)}$. On a nonconvex polygon, each connected segment has its own endpoint on the in-flow boundary. Fubini's theorem in these coordinates gives measurable limits on $\gamma_-$ almost everywhere with respect to the kinetic boundary measure. Their bound follows from the $L^\infty$ bound for $U$ along the segments. Vertex and grazing segments form null sets.

For almost every $(x,w)$, put $\phi(t):=\ue^{-\frac{t}{\e}}U(x-tw,w)$. The sliced identity gives
\begin{align}
 \phi'(t)&=-\ue^{-\frac{t}{\e}}\bigl(\e^{-1}\pk U(x-tw)+S(x-tw,w)\bigr) \quad\text{for almost every }0<t<t_b(x,w).
\end{align}
Integration to $t_b(x,w)$ gives \eqref{eq:first-event-physical-comp} with $f=U$ and $h=U\big|_{\gamma_-}$, proving the assertion.
\end{proof}

The hypotheses of Lemma~\ref{lem:admissible-mild-comp} are preserved by finite sums. If $a\in W^{1,\infty}(\Om)$ is independent of velocity, the product rule and commutation with $\qk$ give
\begin{align}
 \Le(aU)&=aS+(w\cdot\nabla_xa)U.
 \label{eq:admissible-cutoff-product}
\end{align}
Both terms belong to $L^\infty(\Om\times\Sone)$, so the lemma also applies to $aU$. We use these facts for \eqref{eq:forward-composite-definition} at each fixed $\e$ and overlap radius. Proposition~\ref{prop:full-residual} gives the dependence of the residual bounds on these parameters.
\section{Interior Solution and Side Layer}
\label{sec:profiles}

The interior expansion with a Milne boundary layer originates in \cite{Larsen1974,Larsen.Keller1974}, with rigorous analyses in \cite{Bensoussan.Lions.Papanicolaou1979,Bardos.Santos.Sentis1984}. We use the construction developed for neutron transport in \cite{Wu.Guo2015,Guo.Wu2017,Wu2020,Ouyang2024}. To apply it on a polygon, we prove weighted estimates for the Milne solution and its end state, including their dependence on the tangential parameter near an endpoint.

We first derive the equations for the interior solutions $\Uint_k$ and their velocity averages $\rho_k$. The solvability condition for the interior equation gives the Laplace equation for $\rho_k$ by Lemma~\ref{lem:angular-moments-comp}. Its boundary condition is determined by a Milne problem on each side. The prescribed data, together with the interior solutions and side layers from lower orders, determine the in-flow boundary data and source of that problem at order $k$, as specified in \eqref{eq:microscopic-recursion-comp}--\eqref{eq:endstate-recursion-comp}. The end state $e_{j,k}(s)$ is the limit of the Milne solution as $\eta\to\infty$, in the sense of \eqref{eq:milne-maps-comp}. We impose $\rho_k\big|_{E_j}=e_{j,k}$ and subtract $e_{j,k}(s)$ from the Milne solution to obtain the decaying side layer $\Uside_{j,k}$.

\subsection{Interior Solution}

Insert the finite series $\sum_{k=0}^\Nexp\e^k\Uint_k$ into the scaled equation $\e w\cdot\nabla_xu+\qk u=0$. Comparing coefficients through order $\Nexp$ gives
\begin{align}
 \qk\Uint_0=0, \qquad &\qk\Uint_k=-w\cdot\nabla_x\Uint_{k-1} \quad(1\le k\le\Nexp).
 \label{eq:interior-coefficient-equations-comp}
\end{align}

The equation $\qk\Uint_0=0$ gives $\Uint_0=\rho_0(x)$. Since $\qk$ is the identity on functions with zero velocity average and vanishes on functions independent of $w$, we define the interior solution at each subsequent order by
\begin{align}
 &\Uint_k=\rho_k-w\cdot\nabla_x\Uint_{k-1}, \qquad \rho_k:=\pk\Uint_k.
 \label{eq:interior-recursion-comp}
\end{align}
The recursion \eqref{eq:interior-recursion-comp} solves \eqref{eq:interior-coefficient-equations-comp} if and only if $\pk(w\cdot\nabla_x\Uint_{k-1})=0$. Under this condition its velocity average is indeed $\rho_k$, and its microscopic part is $\widehat{\Uint}_k:=\qk\Uint_k=\Uint_k-\rho_k=-w\cdot\nabla_x\Uint_{k-1}$, with $\widehat{\Uint}_0:=0$. We use the following angular identities to derive the Laplace equations from these solvability conditions. Odd velocity moments vanish, and every even moment of a directional derivative is a multiple of an iterated Laplacian.

\begin{lemma}
\label{lem:angular-moments-comp}
Let $\phi$ be a smooth function on $\Om$ and let $q\geq0$ be an integer. Then the angular moments of the directional derivatives of $\phi$ satisfy, pointwise in $\Om$,
\begin{align}
 \pk\bigl(\bigl(w\cdot\nabla\bigr)^{2q+1}\phi\bigr) &=0,
 \label{eq:odd-angular-moment-comp}\\
 \pk\bigl(\bigl(w\cdot\nabla\bigr)^{2q}\phi\bigr) &=c_q\Delta^q\phi, \qquad c_q=\frac{(2q)!}{2^{2q}\bigl(q!\bigr)^2}.
 \label{eq:even-angular-moment-comp}
\end{align}
Let now $k\geq0$, and let $\rho_0,\ldots,\rho_k$ be smooth and harmonic in $\Om$. Set $\Uint_0:=\rho_0$ and define $\Uint_\ell$ by the recursion \eqref{eq:interior-recursion-comp} for $1\leq\ell\leq k$. Then
\begin{align}
 \Uint_\ell &= \sum_{j=0}^{\ell} \bigl(-1\bigr)^j\bigl(w\cdot\nabla\bigr)^j\rho_{\ell-j}, \qquad 0\leq\ell\leq k,
 \label{eq:explicit-interior-profile-comp}\\
 \pk\Uint_\ell &= \rho_\ell, \qquad \pk(w\cdot\nabla\Uint_\ell)=0, \qquad 0\leq\ell\leq k.
 \label{eq:interior-mean-solvability-comp}
\end{align}
Therefore, these interior solutions satisfy \eqref{eq:interior-coefficient-equations-comp}. The identity $\pk(w\cdot\nabla\Uint_k)=0$ is the compatibility condition for constructing $\Uint_{k+1}$. Only $\rho_0,\ldots,\rho_{k-1}$ are required to be harmonic for these identities at order $k$.
\end{lemma}

\begin{proof}
Write $D_w:=w\cdot\nabla_x$. Since $\ud\nu$ is invariant under $w\mapsto-w$, every odd angular moment vanishes, which is \eqref{eq:odd-angular-moment-comp}.

For the even moments, rotational invariance of $\ud\nu$ gives, for $\zeta\in\mathbb R^2$, the polynomial identity $\int_{\Sone}\bigl(w\cdot\zeta\bigr)^{2q}\,\ud\nu(w)=\abs{\zeta}^{2q}\frac1{2\pi}\int_0^{2\pi}\cos^{2q}\varphi\,\ud\varphi=c_q\bigl(\zeta_1^2+\zeta_2^2\bigr)^q$ in $\zeta_1,\zeta_2$, with $c_q=\frac{(2q)!}{2^{2q}\bigl(q!\bigr)^2}$. Comparing coefficients and replacing $\zeta_i$ by the commuting derivatives $\partial_{x_i}$ gives $\pk(D_w^{2q}\phi)=c_q\bigl(\partial_{x_1}^2+\partial_{x_2}^2\bigr)^q\phi=c_q\Delta^q\phi$, which is \eqref{eq:even-angular-moment-comp}.

Iterating \eqref{eq:interior-recursion-comp} gives \eqref{eq:explicit-interior-profile-comp} for each interior solution. Taking the velocity average, the term with $j=0$ gives $\rho_\ell$, and every term with odd $j$ vanishes. For each positive even index $j=2q$, harmonicity gives $\pk(D_w^{2q}\rho_{\ell-2q})=c_q\Delta^q\rho_{\ell-2q}=0$. Hence $\pk\Uint_\ell=\rho_\ell$.

Applying one more $D_w$ gives $\pk(D_w\Uint_\ell)=\sum_{j=0}^{\ell}\bigl(-1\bigr)^j\pk(D_w^{j+1}\rho_{\ell-j})$. Terms with $j+1$ odd vanish by angular symmetry. If $j+1=2q$ is even, then $q\geq1$ and the term equals $-c_q\Delta^q\rho_{\ell-j}=0$. Thus $\pk(D_w\Uint_\ell)=0$, completing \eqref{eq:interior-mean-solvability-comp}.

Finally, $\qk=I-\pk$ gives $\qk\Uint_\ell=\Uint_\ell-\rho_\ell=-w\cdot\nabla\Uint_{\ell-1}$ for $1\le\ell\le k$, while $\qk\Uint_0=0$. Thus the interior solutions satisfy \eqref{eq:interior-coefficient-equations-comp}, and $\pk(w\cdot\nabla\Uint_k)=0$ is the compatibility condition for the next order.
\end{proof}

The Laplace equation for $\rho_k$ first enters at expansion order $k+2$. If $\rho_0,\ldots,\rho_{k-1}$ are harmonic, applying the same angular-moment computation to \eqref{eq:interior-recursion-comp} gives
\begin{align}
 \pk(w\cdot\nabla_x\Uint_{k+1}) &=-\tfrac12\Delta\rho_k.
 \label{eq:interior-laplace-solvability}
\end{align}
Thus $\Delta\rho_k=0$ is the solvability condition for $\Uint_{k+2}$. The finite construction imposes these harmonicity conditions for every $0\le k\le\Nexp$.

\subsection{Side Layer}

Near an open side $E_j$, write $s=s_j$ for the arclength coordinate and introduce the stretched normal variable $\eta=\frac{d_j}{\e}$. In the side coordinates of \eqref{eq:side-coordinate-map}, the streaming operator is $\e w\cdot\nabla_x=\mu_j\partial_\eta+\e\tau_j\partial_s$. Adding the side layer $\sum_{k=0}^\Nexp\e^k\Uside_{j,k}(s,\eta,w)$ to the interior solution and comparing powers of $\e$ gives
\begin{align}
 \big(\mu_j\partial_\eta+\qk\big)\Uside_{j,k} &=-\tau_j\partial_s\Uside_{j,k-1}, \qquad 0\le k\le\Nexp,\quad \Uside_{j,-1}:=0,\\
 \Uside_{j,k}(s,0,w) &=\gaux_{j,k}(s,w)-\Uint_k\big|_{E_j}(s,w), \qquad \mu_j(w)>0.
\end{align}
Here $\gaux_{j,k}$ is the auxiliary in-flow data in \eqref{eq:auxiliary-inflow-coefficients}, $\Uside_{j,k}(s,0,w)$ denotes the characteristic trace for $\mu_j(w)>0$, and $\Uint_k\big|_{E_j}$ is the restriction to the open side $E_j$. At each fixed tangential parameter, we require the side layer to decay in $L^\infty_w$ as $\eta\to\infty$. Its in-flow boundary condition makes the sum of the interior solution and side layer agree with $\gaux_{j,k}$ at that order. The leading side layer $\Uside_{j,0}$ has zero source. For $k\ge1$, the source is $-\tau_j\partial_s\Uside_{j,k-1}$, which is known from the preceding side layer.

To determine the boundary value of $\rho_k$, first solve the Milne problem with in-flow boundary data $\gaux_{j,k}-\widehat{\Uint}_k\big|_{E_j}$ and source $-\tau_j\partial_s\Uside_{j,k-1}$. For each fixed $s$, its limit in $L^\infty_w$ as $\eta\to\infty$ is the Milne end state $e_{j,k}(s)$ in \eqref{eq:endstate-recursion-comp}. Set $\rho_k\big|_{E_j}=e_{j,k}$ and subtract $e_{j,k}(s)$ from the Milne solution. The resulting side layer $\Uside_{j,k}$ satisfies the in-flow boundary condition above. We give this construction at each order in \eqref{eq:microscopic-recursion-comp}--\eqref{eq:side-profile-inflow-decay-comp}.

We next prove weighted estimates for the Milne solution and its end state, including their dependence on the tangential parameter. We first solve the equation for the velocity average in Lemma~\ref{lem:milne-factorization-comp}, using Wiener--Hopf factorization. We then construct the Milne solution in Proposition~\ref{prop:flat-milne-comp}.

The half-space problem is classical. Existence, uniqueness and exponential convergence to an end state were proved for the Boltzmann equation in \cite{Bardos.Caflisch.Nicolaenko1986}. The admissible number of in-flow boundary conditions was settled in \cite{Coron.Golse.Sulem1988}. The stationary linearized theory in a half-space is developed in \cite{Golse.Poupaud.Neunzert1989}. For the one-speed isotropic kernel used here, the radiative-transfer theory is developed in \cite{Chandrasekhar1950}, and the half-space theory is summarized in \cite{Sone2007}. The $L^2$--$L^{\infty}$ approach appears in \cite{Cercignani.Marra.Esposito1998, Wu.Guo2015}.

Fix a unit inward normal $n_{\rm in}$ and put $\mu(w)=w\cdot n_{\rm in}$. We take the in-flow boundary data in $\mathcal H_+:=L^\infty(\{w\in\Sone:\mu(w)>0\},\ud\nu)$. For $0<\kappa<1$, introduce the exponentially weighted spaces
\begin{align}
 X_\kappa &:=\Big\{F\colon\mathbb R_+\times\Sone\to\mathbb R: \nm F_{X_\kappa}:= \esssup_{\eta>0,w\in\Sone} \ue^{\kappa\eta}\abs{F(\eta,w)}<\infty\Big\},
 \label{eq:Milne-exponential-space}\\
 Y_\kappa &:=\Big\{\rhodec\colon\mathbb R_+\to\mathbb R: \nm{\rhodec}_{Y_\kappa}:= \esssup_{\eta>0} \ue^{\kappa\eta}\abs{\rhodec(\eta)}<\infty\Big\}.
\end{align}
These spaces consist of measurable functions modulo equality almost everywhere. The norm on $X_\kappa$ is the weighted $L^\infty(\mathbb R_+\times\Sone)$ norm, and $Y_\kappa$ is the corresponding weighted space for functions of $\eta$ alone. We construct the Milne solution as an end state $E\in\mathbb R$ plus a part in $X_\kappa$. The equation and in-flow boundary data determine both $E$ and the decaying part.

For $h\in\mathcal H_+$ and $F\in L^\infty(\mathbb R_+\times\Sone)$, consider the flat Milne problem
\begin{align}
 \mu\partial_\eta f+\qk f &=F, \qquad f\big|_{\eta=0}(w)=h(w)\quad(\mu>0).
 \label{eq:milne-general-comp}
\end{align}

The restriction at $\eta=0$ denotes the characteristic trace. We prescribe it by $h$ on $\{\mu>0\}$. Put $\rho(\eta)=\pk f(\eta)$. Integrating along characteristics, we obtain
\begin{align}
 f(\eta,w) &=\ue^{-\frac{\eta}{\mu}}h(w) +\frac1\mu\int_0^\eta \ue^{-\frac{\eta-y}{\mu}} \Big(\rho(y)+F(y,w)\Big)\,\ud y, \qquad\mu>0,
 \label{eq:milne-char-plus-comp}\\
 f(\eta,w) &=\frac1{\abs\mu}\int_\eta^\infty \ue^{-\frac{y-\eta}{\abs\mu}} \Big(\rho(y)+F(y,w)\Big)\,\ud y, \qquad\mu<0,
 \label{eq:milne-char-minus-comp}
\end{align}
The formula for $\mu<0$ follows by integrating from infinity back to $\eta$, since membership in $L^\infty(\mathbb R_+\times\Sone)$ excludes the exponentially growing homogeneous term.

A \emph{mild solution in $L^\infty$} of \eqref{eq:milne-general-comp} is an $f\in L^\infty(\mathbb R_+\times\Sone)$ satisfying \eqref{eq:milne-char-plus-comp}--\eqref{eq:milne-char-minus-comp} with $\rho=\pk f$ almost everywhere on $\{\mu(w)\ne0\}$. The grazing set $\{\mu=0\}$ is $\nu$-null and needs no separate prescription.

At fixed $\eta$, we use the representative given by the characteristic formulas. For almost every nongrazing velocity, this representative is locally absolutely continuous in $\eta$. If $f-E\in X_\kappa$, Fubini's theorem gives $\ue^{\kappa\eta}\abs{f(\eta,w)-E}\le\nm{f-E}_{X_\kappa}$ for almost every $\eta$ outside a fixed $\nu$-null set of velocities. Continuity along each such characteristic extends this bound to every $\eta\ge0$. The same argument applies to the unweighted $L^\infty$ bound. Thus the estimates below in $L^\infty_w$ at a fixed $\eta$ refer to this representative.

We will use two elementary identities to prove the weighted estimates. Since $0<\abs\mu\le1$ for $\nu$-almost every $w\in\Sone$, direct integration gives, for every $0<\kappa_*<1$,
\begin{align}
 \frac1\mu\int_0^\eta \ue^{-\frac{\eta-y}{\mu}}\ue^{-\kappa_* y}\,\ud y &= \frac{\ue^{-\kappa_*\eta}-\ue^{-\frac{\eta}{\mu}}}{1-\kappa_*\mu} \quad(\mu>0),
 \label{eq:exponential-convolution-identities-comp}\\
 \frac1{\abs\mu}\int_\eta^\infty \ue^{-\frac{y-\eta}{\abs\mu}}\ue^{-\kappa_* y}\,\ud y &= \frac{\ue^{-\kappa_*\eta}}{1+\kappa_*\abs\mu} \quad(\mu<0),
 \label{eq:exponential-convolution-minus-comp}
\end{align}
where $1-\kappa_*\mu\ge1-\kappa_*>0$.

Averaging \eqref{eq:milne-char-plus-comp}--\eqref{eq:milne-char-minus-comp} in velocity gives $(I-K_+)\rho=H_{h,F}$. The operator on the velocity average is $(K_+\phi)(\eta):=\int_0^\infty k_0(\abs{\eta-y})\phi(y)\,\ud y$, where $k_0(t):=\int_{\{\mu>0\}}\mu^{-1}\ue^{-\frac{t}{\mu}}\,\ud\nu(w)$ for $t>0$. This defines a bounded operator on $L^\infty(\mathbb R_+)$, because $k_0\ge0$ and $2\int_0^\infty k_0(t)\,\ud t=1$. The in-flow data $h$ and source $F$ give
\begin{align}
 \label{eq:milne-H-comp} H_{h,F}(\eta) :={}&\int_{\{\mu>0\}}\ue^{-\frac{\eta}{\mu}}h(w)\,\ud\nu(w)\\
 &+\int_{\{\mu>0\}}\frac1\mu\int_0^\eta \ue^{-\frac{\eta-y}{\mu}}F(y,w)\,\ud y\ud\nu(w)+\int_{\{\mu<0\}}\frac1{\abs\mu}\int_\eta^\infty \ue^{-\frac{y-\eta}{\abs\mu}}F(y,w)\,\ud y\ud\nu(w).\notag
\end{align}

We therefore first solve $(I-K_+)\rho=H_{h,F}$ for the velocity average. The characteristic formulas then give $f$ and the identity $\pk f=\rho$. The following lemma gives the required solution and an estimate for the exponential convergence of $\rho$ to its end state.

\begin{lemma}
\label{lem:milne-factorization-comp}
For every $0<\kappa<\kappa_1<1$ and every $H\in Y_{\kappa_1}$, there is a unique pair $(E,\rhodec)\in\mathbb R\times Y_\kappa$ such that $(I-K_+)(E+\rhodec)=H$ almost everywhere on $\mathbb R_+$. Moreover the map $H\mapsto(E,\rhodec)$ is linear and
\begin{align}
 \abs E+\nm{\rhodec}_{Y_\kappa} \le C_{\kappa,\kappa_1}\nm H_{Y_{\kappa_1}}.
 \label{eq:abstract-halfline-bound-comp}
\end{align}

\end{lemma}

\begin{proof}
We seek $\rho=E+\rhodec$ by Wiener--Hopf factorization. The Fourier symbol has a double zero at the origin. We split it between factors whose inverse Fourier kernels are supported on opposite half-lines and invert each simple zero by one-sided integration. The limit of the positive primitive determines $E$, and weighted estimates give $\rhodec\in Y_\kappa$. The kinetic energy identity proves uniqueness. All Fourier transforms below are taken in the normal variable.

\paragraph{\underline{Extension to the full line}} Set $K(t):=k_0(\abs t)$, and let $\rho_+:=\mathds 1_{\{\eta>0\}}\rho$, $H_+:=\mathds 1_{\{\eta>0\}}H$ be the zero extensions to $\mathbb R$. For $\eta>0$, we have $(K*\rho_+)(\eta)=\int_0^\infty K(\eta-y)\rho(y)\,\ud y=(K_+\rho)(\eta)$. Write $\delta_0$ for the Dirac mass at the origin. Thus $(I-K_+)\rho=H$ on $\mathbb R_+$ is equivalent to requiring
\begin{align}
 Q_-:=(\delta_0-K)*\rho_+-H_+ \quad\text{to satisfy}\quad \supp Q_-\subset(-\infty,0].
 \label{eq:WH-support-condition-comp}
\end{align}

\paragraph{\underline{Computation of the Fourier symbol}} We use the Fourier transform $\widehat g(z):=\int_{\mathbb R}\ue^{\ui z\eta}g(\eta)\,\ud\eta$ and inverse $g(\eta)=\frac1{2\pi}\int_{\mathbb R}\ue^{-\ui z\eta}\widehat g(z)\,\ud z$. We denote this inverse transform by $\FourierInv$. For $\abs{\operatorname{Im}z}<1$, Fubini's theorem and the symmetry $w\mapsto-w$ combine the two contributions $\int_{\{\mu>0\}}\bigl(\bigl(1-\ui z\mu\bigr)^{-1}+\bigl(1+\ui z\mu\bigr)^{-1}\bigr)\ud\nu$ into an integral over $\Sone$. Rotate coordinates so that $\mu(w)=\cos\varphi$ and take the square root in the branch positive at $z=0$. We obtain
\begin{align}
 \widehat K(z) &= \int_{\Sone} \frac{\ud\nu(w)}{1-\ui z\mu(w)} = \frac1{2\pi}\int_0^{2\pi} \frac{\ud\varphi}{1-\ui z\cos\varphi} = \frac1{\sqrt{1+z^2}}.
 \label{eq:milne-kernel-symbol-comp}
\end{align}
Subtracting this from $1$ and extracting the elementary quotient $\frac{z^2}{1+z^2}$ defines the normalized symbol $\mathfrak m$:
\begin{align}
 \mathfrak a(z) &:= 1-\widehat K(z) = \frac{z^2}{1+z^2}\mathfrak m(z), \qquad \mathfrak m(z):= \frac{\sqrt{1+z^2}} {\sqrt{1+z^2}+1}.
 \label{eq:normalized-milne-symbol-comp}
\end{align}
At the origin the symbol $\mathfrak a$ has a double zero,
\begin{align}
 \mathfrak a(z) &= \tfrac12z^2+O(z^4) \qquad\text{as }z\to0.
 \label{eq:milne-symbol-double-zero-comp}
\end{align}
The identity $\widehat K(0)=1$ expresses conservation: $I-K*$ annihilates constants on the full line. Evenness of $K$ gives $\widehat K'(0)=0$, and the nonzero quadratic coefficient in \eqref{eq:milne-symbol-double-zero-comp} makes the zero exactly double.

\paragraph{\underline{Factorization of $\mathfrak m$}} Let $\mathcal S_1:=\{z\in\mathbb C:\abs{\operatorname{Im}z}<1\}$ and let $\mathfrak s(z):=\sqrt{1+z^2}$ be the branch with $\mathfrak s(0)=1$. Write $z=\xi+\ui y$, where $\xi$ denotes a real Fourier variable. For $\abs y<1$, we have $\operatorname{Re}(1+z^2)=1+\xi^2-y^2>0$. Thus $\mathfrak s$ is holomorphic on $\mathcal S_1$ and $\operatorname{Re}\mathfrak s>0$. It follows that $\mathfrak m=\frac{\mathfrak s}{\mathfrak s+1}$ is holomorphic and nonvanishing on this strip. Since $\mathcal S_1$ is simply connected, we may define $\mathfrak q(z):=\ln\mathfrak m(z)$ by $\mathfrak q(0)=-\ln2$. On the real axis, $\mathfrak q(\xi)=\ln\left(\frac{\sqrt{1+\xi^2}}{\sqrt{1+\xi^2}+1}\right)$ is real and tends to $0$ as $\abs\xi\to\infty$.

We factor $\mathfrak m$ by decomposing the inverse Fourier transform of $\mathfrak q$ on the two half-lines and then exponentiating. To justify this construction and the weighted estimates for the factors and their inverses, we first prove exponential $L^1$ integrability of this inverse Fourier transform.

Fix $0<\vartheta<\vartheta'<\kappa_1$. We claim that, uniformly on the closed substrip $\abs{\operatorname{Im}z}\le\vartheta'$,
\begin{align}
 \mathfrak q^{(j)}(z) &= O\!\left(\jbr{\operatorname{Re}z}^{-1-j}\right), \qquad j=0,1,2.
 \label{eq:log-symbol-shifted-bounds-comp}
\end{align}
For $\abs y\le\vartheta'$, we have $\abs{\mathfrak s(z)}^2\ge\operatorname{Re}(1+z^2)=1+\xi^2-y^2\ge\xi^2+1-\bigl(\vartheta'\bigr)^2$ and $\abs{\mathfrak s(z)}^2\le1+\abs z^2\le1+\xi^2+\bigl(\vartheta'\bigr)^2$. Hence
\begin{align}
 c_{\vartheta'}\jbr{\xi} &\leq \abs{\mathfrak s(z)} \leq C_{\vartheta'}\jbr{\xi}.
 \label{eq:s-comparability-comp}
\end{align}
Since $\operatorname{Re}\mathfrak s>0$, we also have $\abs{\mathfrak s(z)+t}\ge\abs{\mathfrak s(z)}$ for $t\in[0,1]$. We now write
\begin{align}
 \mathfrak q(z) &= \ln \mathfrak s(z)-\ln\bigl(\mathfrak s(z)+1\bigr)= -\int_0^1\frac{\ud t}{\mathfrak s(z)+t}.
 \label{eq:q-integral-representation-comp}
\end{align}

Hence \eqref{eq:s-comparability-comp} and \eqref{eq:q-integral-representation-comp} give $\abs{\mathfrak q(z)}\le\abs{\mathfrak s(z)}^{-1}\le C_{\vartheta'}\jbr\xi^{-1}$. Differentiating \eqref{eq:q-integral-representation-comp} with $\mathfrak s'=\frac{z}{\mathfrak s}$, $\mathfrak s''=\mathfrak s^{-3}$ gives $\mathfrak q'(z)=\int_0^1\mathfrak s'\bigl(\mathfrak s+t\bigr)^{-2}\ud t$ and $\mathfrak q''(z)=\int_0^1\{\mathfrak s''\bigl(\mathfrak s+t\bigr)^{-2}-2\bigl(\mathfrak s'\bigr)^2\bigl(\mathfrak s+t\bigr)^{-3}\}\ud t$. Since $\abs{\mathfrak s'}\le C_{\vartheta'}$ and $\abs{\mathfrak s''}\le C_{\vartheta'}\jbr\xi^{-3}$, the same two bounds yield $\abs{\mathfrak q'}\le C_{\vartheta'}\jbr\xi^{-2}$ and $\abs{\mathfrak q''}\le C_{\vartheta'}(\jbr\xi^{-5}+\jbr\xi^{-3})\le C_{\vartheta'}\jbr\xi^{-3}$, which is \eqref{eq:log-symbol-shifted-bounds-comp}.

Let $\ell:=\FourierInv(\mathfrak q\big|_{\mathbb R})$, formally $\ell(\eta)=\frac1{2\pi}\int_{\mathbb R}\ue^{-\ui\xi\eta}\mathfrak q(\xi)\,\ud\xi$. The integral is understood in the $L^2$ sense by \eqref{eq:log-symbol-shifted-bounds-comp}, or, for $\eta\ne0$, as an improper integral convergent by the bound on $\mathfrak q'$.

We first estimate $\ell$ near the origin. Since $\mathfrak q=O(\jbr\xi^{-1})$ on the real axis, we have $\mathfrak q\big|_{\mathbb R}\in L^2(\mathbb R)$. By Plancherel's theorem, $\ell\in L^2(\mathbb R)$. Cauchy--Schwarz then gives $\ell\in L^1(-1,1)$.

To estimate $\ell$ for $\abs\eta\ge1$, we shift the inversion contour. For each fixed $\eta\ne0$, move it to
\begin{align}
 \operatorname{Im}z=-\operatorname{sgn}(\eta)\vartheta',
 \label{eq:ell-shifted-line-comp}
\end{align}
downward when $\eta>0$ and upward when $\eta<0$. On this line, $z=\xi-\ui\operatorname{sgn}(\eta)\vartheta'$ and $\ue^{-\ui z\eta}=\ue^{-\ui\xi\eta}\ue^{-\vartheta'\abs\eta}$, which gives the required exponential decay.

Apply Cauchy's theorem to the holomorphic integrand $\ue^{-\ui z\eta}\mathfrak q(z)$ on the rectangle between the real axis and \eqref{eq:ell-shifted-line-comp}. No singularity is crossed because $\mathfrak q$ is holomorphic on $\mathcal S_1$ and $\vartheta'<1$. On each intermediate horizontal line, $\operatorname{Im}z$ has sign opposite to $\eta$. Therefore $\abs{\ue^{-\ui z\eta}}=\ue^{(\operatorname{Im}z)\eta}\le1$. On the vertical sides $\operatorname{Re}z=\pm A$, estimate \eqref{eq:log-symbol-shifted-bounds-comp} gives $\abs{\mathfrak q(z)}\le C_{\vartheta'}\jbr{A}^{-1}$. The sides have length $\vartheta'$, so their contribution tends to zero as $A\to\infty$. Hence $\ell(\eta)=\frac{\ue^{-\vartheta'\abs\eta}}{2\pi}\int_{\mathbb R}\ue^{-\ui\xi\eta}\mathfrak q\bigl(\xi-\ui\operatorname{sgn}(\eta)\vartheta'\bigr)\,\ud\xi$.

To estimate this oscillatory integral, integrate by parts twice in $\xi$. The bounds on $\mathfrak q$ and $\mathfrak q'$ in \eqref{eq:log-symbol-shifted-bounds-comp} make all boundary terms vanish, giving
\begin{align}
 \ell(\eta) &=\frac{\ue^{-\vartheta'\abs\eta}}{2\pi\bigl(\ui\eta\bigr)^2}\int_{\mathbb R}\ue^{-\ui\xi\eta}\mathfrak q''\bigl(\xi-\ui\operatorname{sgn}(\eta)\vartheta'\bigr)\,\ud\xi, \qquad \eta\ne0.
 \label{eq:log-symbol-kernel-representation}
\end{align}
The bound $\abs{\mathfrak q''}\le C_{\vartheta'}\jbr\xi^{-3}$ is integrable and yields
\begin{align}
 \abs{\ell(\eta)} &\leq C_{\vartheta'} \ue^{-\vartheta'\abs\eta}\abs\eta^{-2}, \qquad \abs\eta\ge1,
 \label{eq:log-symbol-kernel-decay-comp}
\end{align}
and combining it on $\abs\eta\ge1$ with the Plancherel bound $\ell\in L^1(-1,1)$ gives the weighted integrability
\begin{align}
 \ue^{\vartheta\abs{\cdot}}\ell &\in L^1(\mathbb R).
 \label{eq:log-symbol-weighted-L1-comp}
\end{align}

Split $\ell_+:=\mathds 1_{[0,\infty)}\ell$, $\ell_-:=\mathds 1_{(-\infty,0]}\ell$ and set $\mathfrak m_\pm:=\exp(\widehat{\ell_\pm})$, so that $\mathfrak m=\mathfrak m_+\mathfrak m_-$ because $\ell=\ell_++\ell_-$. The kernels of $\mathfrak m_\pm$ and $\mathfrak m_\pm^{-1}$ are given explicitly by the exponential series
\begin{align}
 \FourierInv(\mathfrak m_\pm) &= \delta_0+ \sum_{j=1}^\infty\frac{\ell_\pm^{*j}}{j!},\qquad \FourierInv(\mathfrak m_\pm^{-1}) = \delta_0+ \sum_{j=1}^\infty\frac{\bigl(-\ell_\pm\bigr)^{*j}}{j!},
 \label{eq:factor-convolution-series-comp}
\end{align}
$\ell_\pm^{*j}$ denoting the $j$-fold convolution. Since $\ue^{\vartheta\abs\eta}$ is submultiplicative, \eqref{eq:log-symbol-weighted-L1-comp} makes both series converge in the unitized weighted convolution algebra $\mathbb C\delta_0+L^1(\ue^{\vartheta\abs\eta}\,\ud\eta)$. Convolution preserves each half-line, so the four kernels have one-sided supports:
\begin{align}
 \supp\FourierInv(\mathfrak m_+^{\pm1}) &\subset[0,\infty),\qquad \supp\FourierInv(\mathfrak m_-^{\pm1}) \subset(-\infty,0].
 \label{eq:factor-kernel-support}
\end{align}
This is the weighted one-sided Wiener--Hopf factorization used below. The classical construction is described in \cite[Chapters~1--2]{Noble1988} and \cite{dEon.Williams2018}.

To split the double zero in \eqref{eq:milne-symbol-double-zero-comp}, set $\mathfrak a_+(z):=\frac{z}{z+\ui}\mathfrak m_+(z)$ and $\mathfrak a_-(z):=\frac{z}{z-\ui}\mathfrak m_-(z)$. Then $\mathfrak a=\mathfrak a_+\mathfrak a_-$, and each factor has a simple zero at $z=0$. We invert these zeros below by one-sided integration. Taking the Fourier transform of \eqref{eq:WH-support-condition-comp}, we obtain, for real $z$ and in the space $\mathcal S'(\mathbb R)$ of tempered distributions,
\begin{align}
 \widehat{Q_-}(z)&=\mathfrak a_+(z)\mathfrak a_-(z)\widehat{\rho_+}(z)-\widehat{H_+}(z).
 \label{eq:spectral-fourier}
\end{align}
We use this factorization to solve the equation by convolution in $\eta$.

\paragraph{\underline{Inversion of $\mathfrak a_-$}} We use the exponent $\kappa$ in the conclusion of the lemma. Let $\mathcal W_\kappa$ be the space of measurable functions on $\mathbb R$ with
\begin{align}
 \nm g_{\mathcal W_\kappa}:=\bnm{\ue^{\kappa\abs{\cdot}}g}_{L^1(\mathbb R)}+\bnm{\ue^{\kappa\abs{\cdot}}g}_{L^\infty(\mathbb R)}<\infty.
\end{align}
We need the $L^1$ norm to estimate the one-sided primitives and their limits, and the $L^\infty$ norm to estimate the velocity average. Define $\mathcal B_\pm g:=\FourierInv(\mathfrak m_\pm^{-1})*g$ for $g\in\mathcal W_\kappa$. The inverse kernels were constructed in \eqref{eq:factor-convolution-series-comp}. We claim that $\mathcal B_\pm$ maps $\mathcal W_\kappa$ into itself, with
\begin{align}
 \nm{\mathcal B_\pm g}_{\mathcal W_{\kappa}} &\leq C_{\kappa,\kappa_1} \nm g_{\mathcal W_{\kappa}}.
 \label{eq:factor-convolution-bound}
\end{align}

For $\vartheta>0$, put $\omega_\vartheta(\eta):=\ue^{\vartheta\abs\eta}$. By the triangle inequality, $\omega_\vartheta(\eta)\le\omega_\vartheta(\eta-y)\omega_\vartheta(y)$. Hence $\omega_\vartheta\abs{f*g}\le(\omega_\vartheta\abs f)*(\omega_\vartheta\abs g)$ almost everywhere. Applying Young's inequality, we obtain
\begin{align}
 \nm{\omega_\vartheta(f*g)}_{L^{p_{\rm L}}} &\leq \nm{\omega_\vartheta f}_{L^1} \nm{\omega_\vartheta g}_{L^{p_{\rm L}}}, \qquad p_{\rm L}\in\{1,\infty\}.
 \label{eq:weighted-Young-factor}
\end{align}
Write the kernel of $\mathfrak m_\pm^{-1}$ as $\delta_0+\beta_\pm$, where $\beta_\pm:=\sum_{j\ge1}\frac{\bigl(-\ell_\pm\bigr)^{*j}}{j!}$ by \eqref{eq:factor-convolution-series-comp}. Repeated application of \eqref{eq:weighted-Young-factor} with $p_{\rm L}=1$ gives $\bnm{\omega_\vartheta\ell_\pm^{*j}}_{L^1}\le\nm{\omega_\vartheta\ell_\pm}_{L^1}^{\,j}$. Summing the series yields $1+\nm{\omega_\vartheta\beta_\pm}_{L^1}\le\exp(\nm{\omega_\vartheta\ell_\pm}_{L^1})$, where the initial $1$ is the mass of the Dirac term. Since $\mathcal B_\pm g=g+\beta_\pm*g$, one further application of \eqref{eq:weighted-Young-factor} gives
\begin{align}
 \nm{\omega_\vartheta\mathcal B_\pm g}_{L^{p_{\rm L}}} &\leq \bigl( 1+\nm{\omega_\vartheta\beta_\pm}_{L^1} \bigr) \nm{\omega_\vartheta g}_{L^{p_{\rm L}}}\leq \exp\bigl( \nm{\omega_\vartheta\ell_\pm}_{L^1} \bigr) \nm{\omega_\vartheta g}_{L^{p_{\rm L}}}, \qquad p_{\rm L}\in\{1,\infty\}.
 \label{eq:Bpm-weighted-Lp-bound}
\end{align}
Take $\vartheta=\kappa$ and choose $\kappa<\vartheta'<\kappa_1$. Estimate \eqref{eq:log-symbol-weighted-L1-comp} gives $\nm{\omega_{\kappa}\ell_\pm}_{L^1}\le\nm{\omega_{\kappa}\ell}_{L^1}\le C_{\kappa,\kappa_1}$. Adding \eqref{eq:Bpm-weighted-Lp-bound} for $p_{\rm L}=1$ and $p_{\rm L}=\infty$ proves \eqref{eq:factor-convolution-bound}.

The zero extension $H_+$ belongs to $\mathcal W_\kappa$. Indeed, $H\in Y_{\kappa_1}$ gives $\abs{H(\eta)}\le\ue^{-\kappa_1\eta}\nm H_{Y_{\kappa_1}}$ for almost every $\eta>0$, so $\nm{\ue^{\kappa\abs{\cdot}}H_+}_{L^\infty(\mathbb R)}\le\nm H_{Y_{\kappa_1}}$ and $\nm{\ue^{\kappa\abs{\cdot}}H_+}_{L^1(\mathbb R)}\le\bigl(\kappa_1-\kappa\bigr)^{-1}\nm H_{Y_{\kappa_1}}$. Therefore,
\begin{align}
 \nm{H_+}_{\mathcal W_{\kappa}} &\leq \left( 1+\frac1{\kappa_1-\kappa} \right) \nm H_{Y_{\kappa_1}}.
 \label{eq:Hplus-Wiener-bound}
\end{align}

To invert the simple zeros of $\mathfrak a_\pm$ at $z=0$, define $(J_-g)(\eta):=\int_\eta^\infty g(t)\,\ud t$ and $(J_+g)(\eta):=\mathds 1_{\{\eta>0\}}\int_0^\eta g(t)\,\ud t$. We apply $J_+$ only to functions supported in $[0,\infty)$. These operators are convolutions with one-sided indicator functions in $L^\infty(\mathbb R)$:
\begin{align}
 J_-g=\mathds 1_{(-\infty,0]}*g, \qquad J_+g=\mathds 1_{[0,\infty)}*g \quad\text{when }\supp g\subset[0,\infty),
 \label{eq:one-sided-primitive-convolution}
\end{align}
For $g\in\mathcal W_\kappa$, the Fourier identities for $J_-$ hold in the strip $-\kappa<\operatorname{Im}z<0$, where both transforms converge. For $J_+$ and $g$ supported in $[0,\infty)$, they hold for $\operatorname{Im}z>0$:
\begin{align}
 \widehat{J_-g}(z) &= \frac1{\ui z}\widehat g(z), \qquad \widehat{(I+J_-)g}(z) = \frac{z-\ui}{z}\widehat g(z),
 \label{eq:Jminus-Fourier-multiplier-comp}\\
 \widehat{J_+g}(z) &= -\frac1{\ui z}\widehat g(z), \qquad \widehat{(I+J_+)g}(z) = \frac{z+\ui}{z}\widehat g(z).
 \label{eq:Jplus-Fourier-multiplier-comp}
\end{align}
Inversion by convolution uses two elementary identities. Write
\begin{align}
 c_-(\eta):=\ue^{\eta}\mathds 1_{\{\eta<0\}}, \qquad c_+(\eta):=\ue^{-\eta}\mathds 1_{\{\eta>0\}},
 \label{eq:simple-zero-kernels}
\end{align}
so that $\FourierInv\left(\frac{z}{z-\ui}\right)=\delta_0-c_-$ and $\FourierInv\left(\frac{z}{z+\ui}\right)=\delta_0-c_+$, each supported in the corresponding half-line. For every $g\in\mathcal W_{\kappa}$, and in the second identity for every such $g$ supported in $[0,\infty)$,
\begin{align}
 (\delta_0-c_-)*\bigl[(I+J_-)g\bigr]=g, \qquad (\delta_0-c_+)*\bigl[(I+J_+)g\bigr]=g.
 \label{eq:one-sided-inverse-identities}
\end{align}
Both follow from Fubini's theorem. Indeed,
\begin{align}
 (c_-*J_-g)(\eta)&=\int_\eta^\infty g(s)\bigl(1-\ue^{-(s-\eta)}\bigr)\,\ud s=(J_-g)(\eta)-(c_-*g)(\eta).
\end{align}
For $g$ supported in $[0,\infty)$, the corresponding integral for $c_+*J_+g$ runs from $0$ to $\eta>0$ and gives $J_+g-c_+*g$. Both sides vanish for $\eta<0$. The weighted $L^1$ bound on $g$ justifies these interchanges.

Write $a_\pm:=\FourierInv(\mathfrak a_\pm)$ and $a:=\delta_0-K=\FourierInv(\mathfrak a)$, so that $a=a_-*a_+$ and, by \eqref{eq:simple-zero-kernels} and \eqref{eq:factor-kernel-support},
\begin{align}
 a_-=(\delta_0-c_-)*\FourierInv(\mathfrak m_-), \quad a_+=(\delta_0-c_+)*\FourierInv(\mathfrak m_+), \quad \supp a_\mp\subset\mp[0,\infty).
 \label{eq:factor-kernels-one-sided}
\end{align}
The kernels $a$, $c_\pm$ and $\FourierInv(\mathfrak m_\pm^{\pm1})$ belong to $\mathbb C\delta_0+L^1(\ue^{\kappa\abs\eta}\,\ud\eta)$. For $a=\delta_0-K$, this follows from $\int_{\mathbb R}\ue^{\kappa\abs t}K(t)\,\ud t=2\int_{\{\mu>0\}}\frac{\ud\nu}{1-\kappa\mu}\le\frac1{1-\kappa}$. Each iterated convolution below has at most one factor that is only assumed bounded, either a one-sided indicator or the function on which the kernels act. The remaining factors are finite measures. Their total variations bound the absolute integral, so Fubini's theorem permits regrouping. In particular, $\FourierInv(\mathfrak m_\mp)*\FourierInv(\mathfrak m_\mp^{-1})=\delta_0$.

We first invert $\mathfrak a_-$ on the full line. The resulting function can be nonzero on both half-lines. We restrict it to $(0,\infty)$, extend this restriction by zero, and then invert $\mathfrak a_+$. The omitted part on $(-\infty,0]$, after convolution with $a_-$, will give the additional term $Q_-$ in \eqref{eq:WH-support-condition-comp}.

Put $u_-:=\mathcal B_-H_+$, $G:=(I+J_-)u_-$ and $G_+:=\mathds 1_{(0,\infty)}G$. Regrouping the convolutions and applying the first identity in \eqref{eq:one-sided-inverse-identities} to $H_+$,
\begin{align}
 a_-*G &= (\delta_0-c_-)*\bigl[(I+J_-)H_+\bigr]=H_+.
 \label{eq:minus-factor-inversion}
\end{align}
We estimate $G_+$ using $\ue^{\kappa\eta}\le\ue^{\kappa t}$ for $t\ge\eta$ and Fubini's theorem:
\begin{align}
 \nm{\ue^{\kappa\,\cdot}J_-u_-}_{L^\infty(0,\infty)}+\kappa\nm{\ue^{\kappa\,\cdot}J_-u_-}_{L^1(0,\infty)} &\le2\nm{\ue^{\kappa\,\cdot}u_-}_{L^1(0,\infty)}.
\end{align}
Together with \eqref{eq:factor-convolution-bound} and \eqref{eq:Hplus-Wiener-bound}, this gives
\begin{align}
 \nm{G_+}_{\mathcal W_{\kappa}} &\leq C_{\kappa,\kappa_1}\nm H_{Y_{\kappa_1}}.
 \label{eq:Gplus-Wiener-bound}
\end{align}
Define $G_-:=G-G_+$. Since $\supp G_-\subset(-\infty,0]$, its contribution below vanishes when the equation is restricted to $\eta>0$.

\paragraph{\underline{Inversion of $\mathfrak a_+$}} Put $u_+:=\mathcal B_+G_+$ and $\rho:=(I+J_+)u_+$ on $\eta>0$. Since $G_+$ and the kernel of $\mathcal B_+$ both have support in $[0,\infty)$, so does $u_+$, and the zero extension $\rho_+$ of $\rho$ equals $(I+J_+)u_+$ on all of $\mathbb R$. The same regrouping, now with the second identity in \eqref{eq:one-sided-inverse-identities} applied to $G_+$, gives $a_+*\rho_+=G_+$. With $a=a_-*a_+$, $G=G_++G_-$ and \eqref{eq:minus-factor-inversion} this yields
\begin{align}
 (\delta_0-K)*\rho_+-H_+ = a_-*G_+-a_-*G = -a_-*G_-.
 \label{eq:WH-negative-support-identity}
\end{align}
The convolutions are absolutely convergent, so this is an identity of locally integrable functions.

The right-hand side is supported in $(-\infty,0]$. Indeed $\supp a_-\subset(-\infty,0]$ by \eqref{eq:factor-kernels-one-sided}, and $G_-$ has the same support condition, so their convolution does too. Comparing \eqref{eq:WH-negative-support-identity} with \eqref{eq:WH-support-condition-comp} identifies $Q_-=-a_-*G_-$. Restricting \eqref{eq:WH-negative-support-identity} to $\eta>0$ therefore gives $(I-K_+)\rho=H$ there.

\paragraph{\underline{Estimate of the end state}} We now write $\rho=(I+J_+)u_+$ as a constant plus an exponentially decaying function. The constant is the limit of the primitive of $u_+$, which we estimate using its weighted $L^1$ norm.

By \eqref{eq:factor-convolution-bound} and \eqref{eq:Gplus-Wiener-bound}, $u_+\in\mathcal W_{\kappa}$. Its support is contained in $[0,\infty)$. Define $E:=\int_0^\infty u_+(t)\,\ud t$ and $\rhodec(\eta):=u_+(\eta)-\int_\eta^\infty u_+(t)\,\ud t$. Splitting the integral for $E$ at $\eta$ gives $\rho(\eta)=u_+(\eta)+\int_0^\eta u_+=E+\rhodec(\eta)$. This formula determines the end state $E$ from $H$.

We have $\abs E\le\nm{u_+}_{L^1(0,\infty)}$. For almost every $\eta>0$, we also have $\ue^{\kappa\eta}\abs{\rhodec(\eta)}\le\ue^{\kappa\eta}\abs{u_+(\eta)}+\int_\eta^\infty\ue^{\kappa t}\abs{u_+(t)}\,\ud t\le\nm{u_+}_{\mathcal W_\kappa}$. Combining these bounds with \eqref{eq:factor-convolution-bound} and \eqref{eq:Gplus-Wiener-bound} proves \eqref{eq:abstract-halfline-bound-comp}. Each step is linear. For real $H$, taking real parts gives a real pair $(E,\rhodec)$ with the same bound. This proves existence and linearity.

\paragraph{\underline{Uniqueness}} We reconstruct the kinetic solution from its velocity average and use the energy identity. The zero in-flow boundary data and constant limit at infinity will make the boundary term nonnegative and the flux at infinity vanish.

Suppose $(I-K_+)(E+\rhodec)=0$ with $\rhodec\in Y_\kappa$, and set $\rho:=E+\rhodec$. Reconstruct $f$ by \eqref{eq:milne-char-plus-comp}--\eqref{eq:milne-char-minus-comp} with $h=0$, $F=0$. The symmetry of $\nu$ and Fubini's theorem give $\pk f=K_+\rho=\rho$. For almost every $w$ with $\mu(w)\ne0$, the characteristic representative is locally absolutely continuous in $\eta$ and satisfies $\mu\partial_\eta f+f-\pk f=0$ for almost every $\eta$. Its in-flow trace is $f\big|_{\eta=0}(w)=0$ for almost every $w$ with $\mu(w)>0$.

Then $\nm{f(\eta,\cdot)-E}_{L^\infty_w}\to0$ as $\eta\to\infty$. In \eqref{eq:exponential-convolution-identities-comp}--\eqref{eq:exponential-convolution-minus-comp}, take $\kappa_*=\kappa$. Their denominators stay bounded away from zero because $\kappa<1$. The part $\rhodec$ contributes $O(\ue^{-\kappa\eta})$ in $L^\infty_w$. The constant part of $\rho$ contributes $-E\ue^{-\frac{\eta}{\mu}}$ to $f-E$ for $\mu>0$ and zero for $\mu<0$. Hence $\nm{f(\eta,\cdot)-E}_{L^\infty(\Sone)}\le C_\kappa(\abs E+\nm{\rhodec}_{Y_\kappa})\ue^{-\kappa\eta}$.

We now prove that every mild solution in $L^\infty(\mathbb R_+\times\Sone)$ of the homogeneous equation with zero in-flow boundary data and convergence in $L^\infty_w$ to a constant must vanish. We also use this uniqueness result in Proposition~\ref{prop:flat-milne-comp}. Multiply the equation by $f$ and first integrate over $\{\abs\mu\geq\delta_\mu\}\times(0,L)$. Since $\mu f\partial_\eta f=f(\rho-f)\in L^\infty((0,L)\times\Sone)$, we may let $\delta_\mu\downarrow0$. Using $\pk f=\rho$ and normalization of $\nu$, we have $\int_{\Sone}(f-\rho)f\,\ud\nu=\int_{\Sone}\abs{f-\rho}^2\,\ud\nu$. Hence
\begin{align}
 \frac12\frac{\ud}{\ud\eta} \int_{\Sone}\mu f^2\,\ud\nu + \int_{\Sone}\abs{f-\rho}^2\,\ud\nu &= 0.
 \label{eq:homogeneous-Milne-energy-identity}
\end{align}
Integrating that identity in $\eta$ and letting $L\to\infty$, the flux at infinity vanishes because $\nm{f(L,\cdot)-E}_{L^\infty_w}\to0$ and $\int_{\Sone}\mu\,\ud\nu=0$, which gives
\begin{align}
 -\frac12 \int_{\Sone}\mu\bigl(f\big|_{\eta=0}(w)\bigr)^2\,\ud\nu(w) + \int_0^\infty\!\int_{\Sone} \abs{f-\rho}^2\,\ud\nu\ud\eta &= 0.
 \label{eq:homogeneous-Milne-energy-conclusion}
\end{align}

The in-flow trace vanishes for $\mu>0$, so the first term of \eqref{eq:homogeneous-Milne-energy-conclusion} equals $\frac12\int_{\{\mu<0\}}\abs\mu\bigl(f\big|_{\eta=0}(w)\bigr)^2\,\ud\nu(w)\ge0$. Both terms there are therefore nonnegative, so both vanish and $f=\rho$ almost everywhere. The kinetic equation then gives $\mu\partial_\eta\rho=0$ for almost every nongrazing velocity, so $\rho$ is constant in $\eta$, and that constant is zero since $f\big|_{\eta=0}=0$ on $\{\mu>0\}$. Hence $f=\rho=0$, and letting $\eta\to\infty$ gives $E=0$ and then $\rhodec=0$, proving uniqueness in $\mathbb R\oplus Y_\kappa$.
\end{proof}

For a Banach space $X$, an interval $I$, $q\in\mathbb N_0$ and $0<\gamma<1$, we use
\begin{align}
 \nm g_{C^{q,\gamma}(I;X)}&:=\sum_{a=0}^q\sup_{s\in I}\nm{\partial_s^ag(s)}_X+\sup_{\substack{s,t\in I\\s\ne t}}\frac{\nm{\partial_s^qg(s)-\partial_s^qg(t)}_X}{\abs{s-t}^\gamma}.
\end{align}
The $C^q(I;X)$ norm is the first sum.

\begin{proposition}[Milne problem]
\label{prop:flat-milne-comp}
Fix $0<\kappa_{\rm out}<\kappa_{\rm in}<1$, and let $h\in\mathcal H_+$ and $F\in X_{\kappa_{\rm in}}$. The following statements hold.
\begin{enumerate}[label=\textup{(\roman*)},leftmargin=*]
\item \emph{Existence and decay.} There is a unique pair
\begin{align}
 E&\in\mathbb R, \qquad f^{\rm dec}\in X_{\kappa_{\rm out}}, \qquad f:=E+f^{\rm dec},
 \label{eq:Milne-output-class-rigorous}
\end{align}
for which $f$ is a mild solution in $L^\infty(\mathbb R_+\times\Sone)$ of \eqref{eq:milne-general-comp}. 

In particular, $\nm{f(\eta,\cdot)-E}_{L^\infty_w}\le\nm{f^{\rm dec}}_{X_{\kappa_{\rm out}}}\ue^{-\kappa_{\rm out}\eta}$. The solution is unique among mild solutions in $L^\infty(\mathbb R_+\times\Sone)$ that converge in $L^\infty_w$ to a constant. Hence $(E,f^{\rm dec})$ is independent of $\kappa_{\rm out}$, although the constant in the decay estimate depends on this rate.

\item \emph{Linearity and estimate.} The Milne end state and the decaying part define the operators
\begin{align}
 \mathcal E_\mu(h,F)&:=E, \qquad \mathfrak M_\mu(h,F):=f^{\rm dec}.
 \label{eq:milne-maps-comp}
\end{align}
The function $\mathfrak M_\mu(h,F)$ solves \eqref{eq:milne-general-comp} with source $F$ and in-flow boundary data $h-\mathcal E_\mu(h,F)$, and tends to zero in $L^\infty_w$ as $\eta\to\infty$. Both operators are linear and satisfy
\begin{align}
 \abs{\mathcal E_\mu(h,F)} +\nm{\mathfrak M_\mu(h,F)}_{X_{\kappa_{\rm out}}} &\le C_{\kappa_{\rm out},\kappa_{\rm in}} \left(\nm h_{\mathcal H_+}+\nm F_{X_{\kappa_{\rm in}}}\right).
 \label{eq:milne-mapping-estimate-comp}
\end{align}

The constant is independent of the unit normal used to define $\mu$.

\item \emph{Tangential regularity.} Let the data depend on $s\in I$. For $q\in\mathbb N_0$, assume that
\begin{align}
 h&\in C^q(I;\mathcal H_+), \qquad F\in C^q(I;X_{\kappa_{\rm in}}).
 \label{eq:Milne-parameter-input-rigorous}
\end{align}
For these data, $\mathcal E_\mu(h,F)(s):=\mathcal E_\mu(h(s),F(s))$, and similarly for $\mathfrak M_\mu$. Then $\mathcal E_\mu(h,F)\in C^q(I)$ and $\mathfrak M_\mu(h,F)\in C^q(I;X_{\kappa_{\rm out}})$. For $0\le a\le q$, differentiation in $s$ commutes with these operators:
\begin{align}
 \partial_s^a\mathcal E_\mu(h,F) &=\mathcal E_\mu(\partial_s^ah,\partial_s^aF),\qquad \partial_s^a\mathfrak M_\mu(h,F) =\mathfrak M_\mu(\partial_s^ah,\partial_s^aF).
 \label{eq:milne-parameter-commutation-comp}
\end{align}
If, in addition, $0<\gamma<1$ and
\begin{align}
 h&\in C^{q,\gamma}(I;\mathcal H_+), \qquad F\in C^{q,\gamma}(I;X_{\kappa_{\rm in}})
 \label{eq:Milne-Holder-input-rigorous}
\end{align}
then the end state and decaying part belong respectively to $C^{q,\gamma}(I)$ and $C^{q,\gamma}(I;X_{\kappa_{\rm out}})$, with
\begin{align}
 \nm{\mathcal E_\mu(h,F)}_{C^{q,\gamma}(I)} &+\nm{\mathfrak M_\mu(h,F)}_{C^{q,\gamma}(I;X_{\kappa_{\rm out}})}
 \label{eq:milne-Holder-parameter-estimate}\le C_{\kappa_{\rm out},\kappa_{\rm in}}\left(\nm h_{C^{q,\gamma}(I;\mathcal H_+)}+\nm F_{C^{q,\gamma}(I;X_{\kappa_{\rm in}})}\right),
\end{align}
with a constant depending only on $\kappa_{\rm out}$ and $\kappa_{\rm in}$, independently of $I$, $q$ and $\gamma$. Under \eqref{eq:Milne-parameter-input-rigorous}, the same estimate holds with the $C^q$ norms. Both conclusions hold with one-sided parameter regularity at an endpoint.

Since the operators act at fixed $s$, they commute with dilations of the parameter. For $\varrho>0$, set $h_\varrho(t):=h(\varrho t)$ and $F_\varrho(t):=F(\varrho t)$. Then $\mathcal E_\mu(h,F)(\varrho t)=\mathcal E_\mu(h_\varrho,F_\varrho)(t)$ and $\mathfrak M_\mu(h,F)(\varrho t)=\mathfrak M_\mu(h_\varrho,F_\varrho)(t)$. Consequently, \eqref{eq:milne-Holder-parameter-estimate} on a reference interval holds with a constant independent of $\varrho$.
\end{enumerate}
\end{proposition}

\begin{proof}\leavevmode
\paragraph{\underline{Existence and decay}} We estimate $H_{h,F}$ in \eqref{eq:milne-H-comp}, solve for the velocity average, and reconstruct $f$ along characteristics. Choose $\kappa_{\rm out}<\kappa_1<\kappa_{\rm in}$. We estimate the two source terms in \eqref{eq:milne-H-comp} by \eqref{eq:exponential-convolution-identities-comp}--\eqref{eq:exponential-convolution-minus-comp} with $\kappa_*=\kappa_{\rm in}$. For the boundary term, $0<\mu\leq1$ gives $\ue^{-\frac{\eta}{\mu}}\leq\ue^{-\eta}$. Multiplying by $\ue^{\kappa_1\eta}$ and taking the essential supremum in $\eta>0$, we obtain $\nm{H_{h,F}}_{Y_{\kappa_1}}\le C_{\kappa_1,\kappa_{\rm in}}(\nm h_{\mathcal H_+}+\nm F_{X_{\kappa_{\rm in}}})$ since $\kappa_1<\kappa_{\rm in}<1$.

Lemma~\ref{lem:milne-factorization-comp} gives a unique $\rho=E+\rhodec$ with $\rhodec\in Y_{\kappa_{\rm out}}$ satisfying $(I-K_+)\rho=H_{h,F}$. Define $f$ by substituting this $\rho$ into \eqref{eq:milne-char-plus-comp}--\eqref{eq:milne-char-minus-comp}. Since $\rhodec\in Y_{\kappa_{\rm out}}$ and $F\in X_{\kappa_{\rm in}}\subset X_{\kappa_{\rm out}}$, the identities \eqref{eq:exponential-convolution-identities-comp}--\eqref{eq:exponential-convolution-minus-comp} with $\kappa_*=\kappa_{\rm out}$ bound the contribution of $\rhodec+F$. The constant part of $\rho$ contributes $-E\ue^{-\frac{\eta}{\mu}}$ to $f-E$ for $\mu>0$ and zero for $\mu<0$. Combining these bounds with the boundary term gives $f-E\in X_{\kappa_{\rm out}}$ and \eqref{eq:milne-mapping-estimate-comp}. Taking the velocity average of the characteristic formulas gives $\pk f=K_+\rho+H_{h,F}=\rho$. The first characteristic formula also gives the trace $f\big|_{\eta=0}=h$ on $\{\mu>0\}$. Thus $f$ is a mild solution of the required Milne problem.

\paragraph{\underline{Uniqueness}} Let $z$ be the difference of two $L^\infty$ mild solutions with the same data, each converging in $L^\infty_w$ to a constant. Then $z$ solves the homogeneous equation with zero in-flow trace and converges to the difference of the two constants. The energy argument in the proof of Lemma~\ref{lem:milne-factorization-comp} gives $z=0$. This proves uniqueness in particular in the class \eqref{eq:Milne-output-class-rigorous}, and shows that the solution does not depend on the choice of $\kappa_{\rm out}$.

\paragraph{\underline{Tangential regularity}} The map $(h,F)\mapsto(\mathcal E_\mu(h,F),\mathfrak M_\mu(h,F))$ is bounded and linear from $\mathcal H_+\times X_{\kappa_{\rm in}}$ to $\mathbb R\times X_{\kappa_{\rm out}}$, independent of the parameter. It commutes with Banach-valued differentiation, which gives \eqref{eq:milne-parameter-commutation-comp} under \eqref{eq:Milne-parameter-input-rigorous}. For the H\"older estimate, apply the map to $\partial_s^q(h,F)(s)-\partial_s^q(h,F)(t)$, divide by $\abs{s-t}^\gamma$, and take the supremum. This proves the $C^{q,\gamma}$ assertion under \eqref{eq:Milne-Holder-input-rigorous}. Both arguments apply to one-sided difference quotients at an endpoint of $I$.
\end{proof}

We will use expansions in fractional powers and logarithms of the tangential parameter as $s\downarrow0$. Since the Milne operators act at each fixed $s$, linearity gives expansions of the same form for the end state and decaying part. The following corollary also gives estimates for the remainders.

\begin{corollary}
\label{cor:milne-polyhom-comp}
Fix $0<\kappa_{\rm in}<1$, $s_0>0$, and a finite index set $\mathcal I\subset\mathbb R\times\mathbb N_0$. For each $(\lambda,b)\in\mathcal I$, let $h_{\lambda,b}\in\mathcal H_+$ and $F_{\lambda,b}\in X_{\kappa_{\rm in}}$ be independent of $s$. Suppose that the in-flow data and source have the expansions
\begin{align}
 h(s)&=\sum_{(\lambda,b)\in\mathcal I}s^\lambda(\ln s)^bh_{\lambda,b}+\mathcal R_h(s)\quad\text{in }\mathcal H_+,\\
 F(s)&=\sum_{(\lambda,b)\in\mathcal I}s^\lambda(\ln s)^bF_{\lambda,b}+\mathcal R_F(s)\quad\text{in }X_{\kappa_{\rm in}}.
\end{align}
For some $q,L\in\mathbb N_0$ and $\sigma\in\mathbb R$, the remainders satisfy $\mathcal R_h\in C^q_{\rm loc}\bigl((0,s_0);\mathcal H_+\bigr)$, $\mathcal R_F\in C^q_{\rm loc}\bigl((0,s_0);X_{\kappa_{\rm in}}\bigr)$ and, for every $0\le a\le q$,
\begin{align}
 \nm{\partial_s^a\mathcal R_h(s)}_{\mathcal H_+} &+\nm{\partial_s^a\mathcal R_F(s)}_{X_{\kappa_{\rm in}}} \le C s^{\sigma-a}\bigl(1+\abs{\ln s}\bigr)^L, \qquad 0<s<s_0.
 \label{eq:Milne-endpoint-remainder-hypothesis}
\end{align}
Then the following statements hold.
\begin{enumerate}[label=\textup{(\roman*)},leftmargin=*]
\item \emph{Expansions and derivative estimates.} The Milne maps give
\begin{align}
 \mathcal E_\mu(h(s),F(s)) &=\sum_{(\lambda,b)\in\mathcal I}s^\lambda(\ln s)^b \mathcal E_\mu(h_{\lambda,b},F_{\lambda,b}) +\mathcal E_\mu(\mathcal R_h,\mathcal R_F),
 \label{eq:endstate-polyhom-comp}\\
 \mathfrak M_\mu(h(s),F(s)) &=\sum_{(\lambda,b)\in\mathcal I}s^\lambda(\ln s)^b \mathfrak M_\mu(h_{\lambda,b},F_{\lambda,b}) +\mathfrak M_\mu(\mathcal R_h,\mathcal R_F),
 \label{eq:milne-polyhom-comp}
\end{align}
The remainder terms are evaluated at the same parameter $s$. For every $0<\kappa_{\rm out}<\kappa_{\rm in}$ and every $0\le a\le q$, their tangential derivatives satisfy
\begin{align}
 \abs{\partial_s^a\mathcal E_\mu(\mathcal R_h,\mathcal R_F)} &+\nm{\partial_s^a\mathfrak M_\mu(\mathcal R_h,\mathcal R_F)}_{X_{\kappa_{\rm out}}} \le C' s^{\sigma-a}\bigl(1+\abs{\ln s}\bigr)^L.
 \label{eq:milne-remainder-polyhom-comp}
\end{align}
\item \emph{H\"older estimates.} If, for some $0<\gamma<1$ and every $0<\varrho<\frac{s_0}{2}$,
\begin{align}
 &\nm{\mathcal R_h(\varrho\,\cdot)}_{C^{q,\gamma}((\frac12,2);\mathcal H_+)}
 +\nm{\mathcal R_F(\varrho\,\cdot)}_{C^{q,\gamma}((\frac12,2);X_{\kappa_{\rm in}})}\le C\varrho^{\sigma}\bigl(1+\abs{\ln\varrho}\bigr)^L,
 \label{eq:Milne-dyadic-remainder-hypothesis}
\end{align}
then, for every $0<\kappa_{\rm out}<\kappa_{\rm in}$ and every $0<\varrho<\frac{s_0}{2}$,
\begin{align}
 &\bnm{\mathcal E_\mu(\mathcal R_h,\mathcal R_F)(\varrho\,\cdot)}_{C^{q,\gamma}((\frac12,2))}
 +\bnm{\mathfrak M_\mu(\mathcal R_h,\mathcal R_F)(\varrho\,\cdot)}_{C^{q,\gamma}((\frac12,2);X_{\kappa_{\rm out}})}\le C_{\kappa_{\rm out},\kappa_{\rm in}}C\varrho^{\sigma}\bigl(1+\abs{\ln\varrho}\bigr)^L.
 \label{eq:Milne-dyadic-remainder-conclusion}
\end{align}
\end{enumerate}
\end{corollary}

\begin{proof}
For each fixed $s$, the factors $s^\lambda(\ln s)^b$ are independent of $\eta,w$. Linearity of the Milne maps in \eqref{eq:milne-maps-comp} therefore gives \eqref{eq:endstate-polyhom-comp}--\eqref{eq:milne-polyhom-comp} term by term. The remainders are the images of $(\mathcal R_h,\mathcal R_F)$. To prove \eqref{eq:milne-remainder-polyhom-comp}, apply \eqref{eq:milne-parameter-commutation-comp} on compact subintervals of $(0,s_0)$ to move each $\partial_s^a$ inside the operators. Estimate \eqref{eq:milne-mapping-estimate-comp} then bounds the result by $C_{\kappa_{\rm out},\kappa_{\rm in}}$ times the left-hand side of \eqref{eq:Milne-endpoint-remainder-hypothesis}.

For the dyadic estimate, the fixed-parameter action gives $\mathcal E_\mu(\mathcal R_h,\mathcal R_F)(\varrho t)=\mathcal E_\mu(\mathcal R_h(\varrho t),\mathcal R_F(\varrho t))$, and likewise for $\mathfrak M_\mu$. Apply \eqref{eq:milne-Holder-parameter-estimate} on $I=(\frac12,2)$ to $(\mathcal R_h(\varrho\,\cdot),\mathcal R_F(\varrho\,\cdot))$. Its constant is independent of $\varrho$, giving \eqref{eq:Milne-dyadic-remainder-conclusion}.

In the dyadic norms \eqref{eq:weighted-boundary-norm-comp} introduced in Section~\ref{sec:polygonal-profiles-comp}, hypothesis \eqref{eq:Milne-dyadic-remainder-hypothesis} states that the remainders in the data belong to $\mathcal C^{q,\gamma}_{\sigma,L}$. By \eqref{eq:Milne-dyadic-remainder-conclusion}, their images belong to the same class, with a constant independent of $\varrho$. We use this estimate to construct the interior solutions and side layers in Section~\ref{sec:interior-side-hierarchy-comp}.
\end{proof}
\section{Elliptic Estimates near the Vertices}
\label{sec:polygonal-profiles-comp}

The Milne end states $e_{j,k}(s)$ give the boundary data for the harmonic functions $\rho_k$. The map $\mathcal E_\mu(h,F)$ in \eqref{eq:milne-maps-comp} takes the in-flow boundary value $h$ and source term $F$ of the Milne problem to the constant limit of its solution as $\eta\to\infty$. In \eqref{eq:endstate-recursion-comp}, the inputs depend on the preceding interior solution and side layer, and the output depends on the side parameter $s$. We solve the harmonic Dirichlet problem with this output as boundary data.

Differentiation in the construction can make the higher-order boundary data singular at a vertex. Their finite expansions on the incident sides determine harmonic lifts of the prescribed powers and logarithms. We subtract these lifts before solving for an $H^1$ correction. The global Dirichlet problem also determines the coefficients of harmonic modes with positive radial exponent and zero boundary values on the two rays. The estimates below control the coefficients of both kinds of terms and the derivatives of the remainder. This permits further differentiation when constructing the next interior solution and side layer.

We define the weighted spaces and harmonic lifts in Subsections~\ref{subsec:jets-germs-comp}--\ref{subsec:harmonic-lifts-comp}, and the solution with prescribed singular part in Subsection~\ref{subsec:finite-part-comp}. Proposition~\ref{prop:mellin-mapping-comp} gives the expansion used to construct $\rho_k$. We then apply it to $\rho_0$ in Subsection~\ref{subsec:leading-laplace-corner-structure}.

Weighted spaces, polyhomogeneous expansions and Mellin analysis at conical points are classical \cite{Kondratiev1967,Grisvard1985,Dauge1988,Kozlov.Mazya.Rossmann1997}. We formulate the finite expansions with Banach-valued coefficients and remainders in weighted H\"older spaces, so that the same notation applies to the interior solutions and side layers in Section~\ref{sec:interior-side-hierarchy-comp}. For the harmonic problem, we specify the normalization of the resonant lifts and the class of solutions with prescribed singular part. Transforming the sector to a cylinder then gives the coefficient formulas and the weighted remainder estimate.

\subsection{Weighted Expansions and Angular Basis}
\label{subsec:jets-germs-comp}

Fix $V_m$ and write $\omega:=\omega_m$. In the local weighted norms below, $x$ denotes a Cartesian coordinate centered at $V_m$ and rotated so that the incident rays have angles $0$ and $\omega$. Thus $\abs{x}=\bulkdistat{m}$. We use $\bulkdistat{m}\in(0,\ellvtx)$ for radial distance and $\sidedistat{m}{i}$ for distance on incident ray $i\in\{1,2\}$. The symbol $s\in(0,\ellvtx)$ denotes a boundary parameter. On a ray, $\bulkdistat{m}=\sidedistat{m}{i}$. We retain the two symbols to distinguish interior and side coordinates as in Subsection~\ref{subsec:intro-coordinate-atlas}. The components $\mu_{m,i},\tau_{m,i}$ of \eqref{eq:local-ray-tangent-sign} retain their vertex index. Endpoint derivatives and Taylor coefficients use the orientation away from the vertex in \eqref{eq:local-global-side-coordinate}.

We allow the coefficients of a finite asymptotic expansion to take values in a Banach space $X$. We use $X=\mathbb R$ for $e_{j,k}$, $X=L^\infty(\Sone)$ for the interior solution $\Uint_k$, and the Milne space $X_\kappa$ of \eqref{eq:Milne-exponential-space} for the side layer $\Uside_{j,k}$.

Fix $q\in\mathbb N_0$, $0<\gamma<1$, $\sigma\in\mathbb R$ and $B\in\mathbb N_0$. For an $X$-valued function $\mathcal R$ on a ray or on the sector $K_\omega$ of \eqref{eq:wedge-sector-def}, define the weighted H\"older norms by
\begin{align}
 \nm{\mathcal R}_{\mathcal C^{q,\gamma}_{\sigma,B}(0,\ellvtx;X)} &:=\sup_{0<\varrho<\frac{\ellvtx}{2}} \varrho^{-\sigma}\bigl(1+\abs{\ln\varrho}\bigr)^{-B} \nm{\mathcal R(\varrho\,\cdot)}_{C^{q,\gamma}((\frac12,2);X)},
 \label{eq:weighted-boundary-norm-comp}\\
 \nm{\mathcal R}_{\mathcal C^{q,\gamma}_{\sigma,B}(K_\omega\cap B_{\ellvtx};X)} &:=\sup_{0<\varrho<\frac{\ellvtx}{2}} \varrho^{-\sigma}\bigl(1+\abs{\ln\varrho}\bigr)^{-B} \nm{\mathcal R(\varrho\,\cdot)}_{C^{q,\gamma}(K_\omega\cap\{\frac12<\abs y<2\};X)},
 \label{eq:weighted-sector-norm-comp}
\end{align}
The spaces denoted by $\mathcal C^{q,\gamma}_{\sigma,B}$ consist of the functions with finite corresponding norm. We omit the value space when $X=\mathbb R$. The change of variables $x=\varrho y$ sends each annulus to a fixed annulus. The factors in front of the norms account for the radial weight $\varrho^\sigma$ and the logarithmic weight $\bigl(1+\abs{\ln\varrho}\bigr)^B$. Thus the rescaled $C^{q,\gamma}$ estimates are uniform in $\varrho$, including the H\"older seminorm of the derivatives of order $q$. Here $D_x^a$ denotes any Cartesian derivative of total order $a\in\mathbb N_0$. The chain rule gives $D_y^a[\mathcal R(\varrho y)]=\varrho^a(D_x^a\mathcal R)(\varrho y)$. Consequently a Cartesian derivative lowers the radial weight by one. In polar coordinates the same loss occurs through $\partial_{\bulkdistat{m}}$ and $\bulkdistat{m}^{-1}\partial_\theta$.

By rescaling, the continuous $X$-valued representatives satisfy the following estimates at every ray point $s$ and sector point $x$. Here $C$ is bounded by a fixed multiple of the corresponding norm in \eqref{eq:weighted-boundary-norm-comp}--\eqref{eq:weighted-sector-norm-comp}:
\begin{align}
 \nm{\partial_s^a\mathcal R(s)}_X &\le C s^{\sigma-a}\bigl(1+\abs{\ln s}\bigr)^B, \qquad 0\le a\le q,
 \label{eq:weighted-boundary-pointwise-comp}\\
 \nm{D_x^a\mathcal R(x)}_X &\leq C \abs x^{\sigma-a}\bigl(1+\abs{\ln \abs x}\bigr)^B, \qquad 0\le a\le q.
 \label{eq:weighted-sector-pointwise-comp}
\end{align}

We write each solution and its boundary data near a vertex as a finite sum of powers and logarithms plus a remainder in these spaces. The powers may be noninteger or negative. Differentiating the finite sum lowers each spatial exponent, while \eqref{eq:weighted-boundary-pointwise-comp}--\eqref{eq:weighted-sector-pointwise-comp} give the corresponding loss for the remainder. Both parts therefore remain controlled at the next order of the construction.

\begin{definition}
\label{def:polyhom-germ-comp}
Let $\sigma\in\mathbb R$. An $X$-valued function $F$ has a finite polyhomogeneous expansion below exponent $\sigma$ if $F(s)=\sum_{(\lambda,b)\in I_\sigma}s^\lambda(\ln s)^bF_{\lambda,b}+\mathcal R_\sigma(s)$, where
\begin{enumerate}[label=\textup{(\roman*)},leftmargin=*]
\item $I_\sigma\subset\{(\lambda,b):\lambda<\sigma,\ b\in\mathbb N_0\}$ is finite and $F_{\lambda,b}\in X$.
\item $\mathcal R_\sigma\in\mathcal C^{q,\gamma}_{\sigma,B}(0,\ellvtx;X)$ for the stated $q,\gamma,B$.
\end{enumerate}
Terms with the same $\lambda$ form a polynomial in $\ln s$. We identify two expansions when their coefficients agree, with absent coefficients read as zero. By Lemma~\ref{lem:finite-polyhom-uniqueness}, this is equivalent to the difference of the represented functions belonging to $\mathcal C^{q,\gamma}_{\sigma,B}(0,\ellvtx;X)$.
\end{definition}

A Taylor polynomial corresponds to $\lambda\in\mathbb N_0$ and $b=0$. We apply Taylor's theorem to the prescribed in-flow boundary data $g_{j,0}$ and $g_{j,1}$. We differentiate singular powers and logarithms directly on the open rays.

In these local coordinates, suppress $m$ and write $\bulkdist=\bulkdistat{m}$. The local domain is $K_\omega\cap B_{\ellvtx}=\{(\bulkdist,\theta):0<\bulkdist<\ellvtx,\ 0<\theta<\omega\}$. The main construction has $0<\omega<\pi$. The angular calculations and weighted expansion below apply to $0<\omega<2\pi$, as stated in Proposition~\ref{prop:mellin-mapping-comp}.

For a harmonic function with zero trace on the two rays, set $u=\bulkdist^\lambda\Phi(\theta)$. Direct computation gives
\begin{align}
 \Delta(\bulkdist^\lambda\Phi) &=\bulkdist^{\lambda-2}\bigl(\Phi''+\lambda^2\Phi\bigr).
 \label{eq:elementary-vertex-pencil}
\end{align}
Thus the angular boundary-value problem is $\Phi''+\lambda^2\Phi=0$ on $(0,\omega)$ with $\Phi(0)=\Phi(\omega)=0$. For $\lambda>0$, it has a nonzero solution precisely at the following exponents:
\begin{align}
 \lambda=\lambda_n&=\frac{n\pi}{\omega}, \qquad \Phi_n(\theta)=\sin\left(\frac{n\pi\theta}{\omega}\right), \qquad n\in\mathbb N.
 \label{eq:elementary-dirichlet-exponents}
\end{align}
The functions $\{\Phi_n\}$ form a complete orthogonal basis of $L^2(0,\omega)$, with $\int_0^\omega\Phi_n\Phi_k\,\ud\theta=\frac{\omega}{2}\delta_{nk}$. Define $Z_n^\pm(\bulkdist,\theta):=\bulkdist^{\pm\lambda_n}\Phi_n(\theta)$ for $n\ge1$. These functions solve the homogeneous Dirichlet problem on the sector:
\begin{align}
 \Delta Z_n^\pm&=0\quad\text{in }K_\omega, \qquad Z_n^\pm\big|_{\theta=0}=Z_n^\pm\big|_{\theta=\omega}=0\quad\text{for }\bulkdist>0.
\end{align}
By \eqref{eq:intro-zero-ray-energy-test-comp}--\eqref{eq:intro-negative-ray-energy-test-comp}, $Z_n^+\in H^1$ locally near the vertex and $Z_n^-\notin H^1$.

The exponent $\lambda_n$ is the radial homogeneity, since $Z_n^+(t\bulkdist,\theta)=t^{\lambda_n}Z_n^+(\bulkdist,\theta)$. It need not be an integer. In a truncated sector, the boundary conditions on the two rays leave the values on the circular arc unspecified. They therefore do not determine the coefficients of $Z_n^\pm$. The global Dirichlet problem determines the coefficients of $Z_n^+$ in the $H^1$ correction. The prescribed singular part fixes the coefficients of $Z_n^-$. We define this solution class in Section~\ref{subsec:finite-part-comp} and estimate the coefficients of $Z_n^+$ in Proposition~\ref{prop:mellin-mapping-comp}.

\subsection{Harmonic Lifts of Ray Data}
\label{subsec:harmonic-lifts-comp}

We first construct the harmonic lift of a power-logarithmic term prescribed on each ray. The lift solves $\Delta H=0$ in $K_\omega$ and has the given boundary values on $\theta=0,\omega$ for $\bulkdist>0$. These values are explicit restrictions on the open rays. We will subtract the lifts and an extension of the boundary remainder to obtain zero ray data. Within the prescribed radial exponent and logarithmic form, a nonzero resonant exponent leaves a free harmonic mode with zero ray data. We fix this freedom by an $L^2(0,\omega)$ orthogonality condition below. This distinguishes the prescribed lift from the additional $Z_n^+$ term whose coefficient is determined by the global Dirichlet problem.

Set $\ell_{\log}:=\ln\bulkdist$ and seek a harmonic function $\bulkdist^\lambda P(\ell_{\log},\theta)$, where $P$ is polynomial in $\ell_{\log}$. We allow any real $\lambda$. Direct calculation gives
\begin{align}
 \Delta\bigl(\bulkdist^\lambda P(\ell_{\log},\theta)\bigr) =\bulkdist^{\lambda-2} \left(\bigl(\partial_{\ell_{\log}}+\lambda\bigr)^2+\partial_\theta^2\right) P(\ell_{\log},\theta),
 \label{eq:polar-pencil-calculation-comp}
\end{align}
Comparing powers of $\ell_{\log}$ reduces harmonicity to a finite triangular system of angular two-point problems.

\paragraph{\underline{Exponent zero $\lambda=0$}} For constant ray values $a_1,a_2$,
\begin{align}
 H_0(\bulkdist,\theta)=a_1+(a_2-a_1)\frac{\theta}{\omega}
 \label{eq:degree-zero-lift-comp}
\end{align}
is harmonic with the prescribed ray values. If $a_1\ne a_2$, then $\abs{\nabla H_0}=\frac{\abs{a_2-a_1}}{\omega\bulkdist}$, so $H_0\notin H^1$ near the vertex. This term is included in the prescribed singular part of any harmonic function $\rho_k$ with unequal constant ray values.

For logarithmic ray data at exponent zero, prescribe $h_i^{[0,b]}$ and seek the lift $H_{0,b}$ in the form
\begin{align}
 h_i^{[0,b]}(\sidedist{i}) &:=\sum_{j=0}^b a_{i,j}\bigl(\ln\sidedist{i}\bigr)^j, \qquad i=1,2,
 \label{eq:zero-exponent-logarithmic-ansatz}\\
 H_{0,b}(\bulkdist,\theta) &:=\sum_{j=0}^b\ell_{\log}^jP_j(\theta).
 \label{eq:zero-exponent-logarithmic-lift}
\end{align}
Comparing powers of $\ell_{\log}$ in \eqref{eq:polar-pencil-calculation-comp} and imposing the ray data gives
\begin{align}
 P_j''+(j+2)(j+1)P_{j+2}=0, \qquad P_j(0)=a_{1,j}, \qquad P_j(\omega)=a_{2,j}, \qquad 0\le j\le b,
 \label{eq:zero-exponent-log-lift-system}
\end{align}
with $P_{b+1}\equiv P_{b+2}\equiv0$.

Starting at $j=b$, each equation in \eqref{eq:zero-exponent-log-lift-system} is a two-point problem for $\partial_\theta^2$ on $(0,\omega)$ with known right-hand side and endpoint values. The Dirichlet realization of $\partial_\theta^2$ has trivial kernel, so it determines $P_j$ uniquely. Descending induction therefore gives a unique lift among functions of radial exponent zero that are polynomial in $\ln\bulkdist$. No logarithmic degree above $b$ is created. The case $b=0$ gives \eqref{eq:degree-zero-lift-comp}.

Repeated Cartesian or tangential differentiation lowers the radial exponent of a resonant $\bulkdist^q\ln\bulkdist$ term. If it reaches exponent zero, the resulting logarithmic ray data are lifted by \eqref{eq:zero-exponent-log-lift-system}.

\paragraph{\underline{Nonresonant nonzero powers}} Let $\lambda\ne0$ and $\sin(\lambda\omega)\ne0$. For ray data $a_1\bulkdist^\lambda$ on $\theta=0$ and $a_2\bulkdist^\lambda$ on $\theta=\omega$, set
\begin{align}
 H_\lambda[a_1,a_2](\bulkdist,\theta) =\bulkdist^\lambda\left[ a_1\frac{\sin(\lambda(\omega-\theta))}{\sin(\lambda\omega)} +a_2\frac{\sin(\lambda\theta)}{\sin(\lambda\omega)} \right]
 \label{eq:nonresonant-harmonic-lift-comp}
\end{align}
The function has the prescribed restrictions on both rays. Its angular factor $\Phi$ solves $\Phi''+\lambda^2\Phi=0$, so $\Delta(\bulkdist^\lambda\Phi)=0$ by \eqref{eq:polar-pencil-calculation-comp}. The same formula applies when $\lambda<0$ and is interpreted on $\bulkdist>0$.

\paragraph{\underline{Resonant nonzero powers}} Let $\lambda=\pm\frac{n\pi}{\omega}$, $n\in\mathbb N$, so $\sin(\lambda\omega)=0$. A log-free lift exists precisely when $a_2=(-1)^na_1$. Every such lift is $a_1\bulkdist^\lambda\cos(\lambda\theta)$ plus a multiple of the zero-trace mode $\bulkdist^\lambda\sin(\lambda\theta)$. The normalization below fixes that multiple to be zero. For a general pair, set $d=a_2-(-1)^na_1$. The vanishing denominator in the nonresonant formula reflects the one-dimensional kernel of the angular Dirichlet operator. Differentiating the homogeneous family with respect to its exponent produces a harmonic term whose trace supplies this mismatch. The resulting logarithm is needed to satisfy both boundary conditions.

As $\bulkdist^\zeta\sin(\zeta\theta)$ is harmonic for every real $\zeta$, so is its parameter derivative,
\begin{align}
 \left.\partial_\zeta \bigl(\bulkdist^\zeta\sin(\zeta\theta)\bigr)\right|_{\zeta=\lambda} &=\bulkdist^\lambda\bigl[(\ln \bulkdist)\sin(\lambda\theta) +\theta\cos(\lambda\theta)\bigr],
 \label{eq:resonant-derivative-comp}
\end{align}
which vanishes on $\theta=0$ and equals $\omega(-1)^n\bulkdist^\lambda$ on $\theta=\omega$.

The lift with the required traces is then
\begin{align}
 \widetilde{H_\lambda}[a_1,a_2] &=a_1\bulkdist^\lambda\cos(\lambda\theta) +\frac{d}{\omega(-1)^n}\bulkdist^\lambda \bigl[(\ln \bulkdist)\sin(\lambda\theta) +\theta\cos(\lambda\theta)\bigr],
 \label{eq:resonant-harmonic-lift-comp}
\end{align}
where the logarithmic coefficient follows from \eqref{eq:resonant-derivative-comp}.

\paragraph{\underline{Logarithmic ray data with $\lambda\ne0$}} Fix $\lambda\ne0$ and prescribe the ray traces $\sidedist{i}^\lambda\sum_{c=0}^b a_{i,c}\bigl(\ln \sidedist{i}\bigr)^c$, $i\in\{1,2\}$, with $a_{i,c}:=0$ for $c>b$. Seek $H_{\lambda,b}$ in the form
\begin{align}
 H_{\lambda,b}(\bulkdist,\theta) &:= \bulkdist^\lambda \sum_{c=0}^{b+\iota_{\rm res}}\bigl(\ln \bulkdist\bigr)^cP_c(\theta), \qquad \iota_{\rm res} = \begin{cases}
 0,&\sin(\lambda\omega)\ne0,\\
 1,&\lambda=\pm\frac{n\pi}{\omega},
 \end{cases}
 \label{eq:general-log-lift-comp}
\end{align}

Set $P_c=0$ for $c>b+\iota_{\rm res}$. By \eqref{eq:polar-pencil-calculation-comp}, harmonicity and the two ray boundary conditions are equivalent to the following equations for $0\le c\le b+\iota_{\rm res}$:
\begin{align}
 \left(\partial_\theta^2+\lambda^2\right)P_c &= -2\lambda(c+1)P_{c+1} -(c+1)(c+2)P_{c+2}, \qquad P_c(0)=a_{1,c}, \quad P_c(\omega)=a_{2,c}.
 \label{eq:nonzero-log-lift-recursion-comp}
\end{align}
If $\sin(\lambda\omega)\ne0$, the Dirichlet realization of $\partial_\theta^2+\lambda^2$ is invertible. Solving \eqref{eq:nonzero-log-lift-recursion-comp} downward from $c=b$ therefore gives $P_b,\ldots,P_0$ uniquely. The lift has logarithmic degree at most $b$ and depends linearly on the ray coefficients.

It remains to solve \eqref{eq:nonzero-log-lift-recursion-comp} at a resonant exponent. Let $\lambda=\pm\frac{n\pi}{\omega}$ and $\Phi_\lambda(\theta):=\sin(\lambda\theta)$. The Dirichlet realization of $\partial_\theta^2+\lambda^2$ has kernel $\operatorname{span}\{\Phi_\lambda\}$. We solve each equation on the orthogonal complement of this kernel after imposing the Fredholm compatibility condition. The coefficient of the next logarithmic degree makes that condition possible. For endpoint values $c_1,c_2$, define $\mathscr C_\lambda[c_1,c_2](\theta):=c_1\cos(\lambda\theta)+\bigl(c_2-(-1)^nc_1\bigr)\bigl(\omega(-1)^n\bigr)^{-1}\theta\cos(\lambda\theta)$. Write $P_c:=\mathscr C_\lambda[a_{1,c},a_{2,c}]+Q_c+\beta_c\Phi_\lambda$. Here $Q_c\in H^2(0,\omega)\cap H_0^1(0,\omega)$ is $L^2$-orthogonal to $\Phi_\lambda$, and $P_{b+1}=\beta_{b+1}\Phi_\lambda$.

Starting at level $b$, the Fredholm compatibility condition for $Q_c$ contains the undetermined coefficient $\beta_{c+1}$. Its coefficient is $-2\lambda(c+1)\int_0^\omega\Phi_\lambda\bigl(\theta\bigr)^2\,\ud\theta=-\lambda(c+1)\omega\ne0$, so compatibility fixes $\beta_{c+1}$ uniquely. The orthogonal equation then determines $Q_c$, while $\beta_c$ remains available for the next lower equation. Finite downward induction gives $P_{b+1},\ldots,P_0$, leaving only $\beta_0$ free.

Finally impose $\int_0^\omega\bigl(P_0-\mathscr C_\lambda[a_{1,0},a_{2,0}]\bigr)(\theta)\Phi_\lambda(\theta)\,\ud\theta=0$, that is $\beta_0=0$, since $Q_0\perp\Phi_\lambda$. This gives existence, uniqueness under that normalization and linear dependence on the ray coefficients.

For $b=0$, the first compatibility condition gives $\beta_1=\frac{a_{2,0}-(-1)^na_{1,0}}{\omega(-1)^n}$. Then $Q_0=0$, and the normalization gives $\beta_0=0$. The resulting lift is exactly \eqref{eq:resonant-harmonic-lift-comp}.

A resonant ray term of logarithmic degree $b$ therefore requires degree at most $b+1$. Compatibility at level $b$ fixes
\begin{align}
 \beta_{b+1} &=\frac{a_{2,b}-(-1)^na_{1,b}}{(b+1)\,\omega(-1)^n}.
 \label{eq:resonant-top-log-coefficient}
\end{align}
This coefficient vanishes precisely when $a_{2,b}=(-1)^na_{1,b}$. The term of logarithmic degree $b+1$ has zero trace on both rays. Its normal derivative satisfies, up to orientation,
\begin{align}
 \partial_n\left[\beta_{b+1}\bulkdist^\lambda\bigl(\ln\bulkdist\bigr)^{b+1}\Phi_\lambda\right]\big|_{\rm ray} &=\pm\beta_{b+1}\lambda\bulkdist^{\lambda-1}\bigl(\ln\bulkdist\bigr)^{b+1}.
\end{align}
Thus a resonant lift raises the logarithmic degree by at most one, and its normal derivative attains that degree when $\beta_{b+1}\ne0$. This gives the bound $L_k^{\rm int}=L_k^{\rm ray}+1$ in \eqref{eq:order-dependent-log-budgets}. At negative $\lambda$, the prescribed normalization fixes the otherwise free zero-trace mode with radial exponent $\lambda$. Lemma~\ref{lem:finite-part-uniqueness-comp} proves that the solution is unchanged by cutoff changes or by transferring an $H^1$ lift between the prescribed part and the regular correction, with this normalization fixed.

\subsection{Dirichlet Problem with a Prescribed Singular Part}
\label{subsec:finite-part-comp}

We first subtract enough terms of the ray expansions to obtain a boundary remainder in $H^{\frac12}(\partial\Om)$. We multiply their harmonic lifts by radial cutoffs and denote the resulting sum by $H^{\rm sing}$. The correction to $H^{\rm sing}$ then satisfies a variational Dirichlet problem in $H^1(\Om)$. Its source consists of cutoff commutators supported on fixed annuli. We use the normalization of Section~\ref{subsec:harmonic-lifts-comp} for all zero-trace modes with negative radial exponent.

For boundary data $h$ given on each open side, let $H_{{\rm ray},m}$ be the sum of the harmonic lifts of the displayed terms on the two rays at $V_m$. The sum includes the negative powers and logarithmic terms given by the ray expansions and the resonant lifting formulas.

Choose smooth radial cutoffs $\chi_m^{\rm sing}=1$ near $V_m$ with disjoint supports. Require each support to lie in a vertex neighborhood where $\Om$ agrees with its sector chart. Define the singular part on $\Om\setminus\{V_1,\ldots,V_{\Nvtx}\}$ by
\begin{align}
 H^{\rm sing}&:=\sum_m\chi_m^{\rm sing}H_{{\rm ray},m},
 \label{eq:finite-part-cutoff-realization}
\end{align}
and, as $\Delta H_{{\rm ray},m}=0$ for $\bulkdistat{m}>0$, the commutators
\begin{align}
 f^{\rm cut}_m&:=[\Delta,\chi_m^{\rm sing}]H_{{\rm ray},m} =2\nabla\chi_m^{\rm sing}\cdot\nabla H_{{\rm ray},m} +(\Delta\chi_m^{\rm sing})H_{{\rm ray},m},
 \label{eq:cutoff-commutator-comp}\\
 f^{\rm cut}&:=\sum_m f^{\rm cut}_m =\Delta H^{\rm sing} \qquad\text{classically away from the vertices.}
 \label{eq:finite-part-commutator-source}
\end{align}

Each commutator is smooth and supported in a closed annulus separated from $V_m$. Extending it by zero near the vertex gives $f^{\rm cut}\in C^\infty(\overline\Om)$. We shall solve $-\Delta v^{\rm reg}=f^{\rm cut}$ with boundary data $h-H^{\rm sing}\big|_{\partial\Om}$. Then $H^{\rm sing}+v^{\rm reg}$ is harmonic and has the prescribed boundary data on the open sides.

The following lemma shows that a boundary remainder of positive radial order belongs to $H^{\frac12}$ across a vertex. We estimate the differences both on one ray and between the two incident rays.
\begin{lemma}
\label{lem:finite-part-remainder-trace}
Fix $V_m$ and the atlas radius $\ellvtx$. The set $\partial\Om\cap B_{\ellvtx}(V_m)$ consists of the two incident ray segments, parametrized by $x=V_m+\sidedistat{m}{i}\tphys{m}{i}$ with $\sidedistat{m}{i}=\abs{x-V_m}\in(0,\ellvtx)$, $i=1,2$. Let $\mathcal R_{m,i}\in C^1((0,\ellvtx])$ be the restriction of a boundary remainder to ray $i$. Suppose that, for some $\sigma>0$ and $B\in\mathbb N_0$,
\begin{align}
 M_{\rm tr} &:= \max_{i=1,2}\nm{ \Big( \babs{\mathcal R_{m,i}(\sidedistat{m}{i})} +\sidedistat{m}{i} \babs{\partial_{\sidedistat{m}{i}} \mathcal R_{m,i}(\sidedistat{m}{i})} \Big) \sidedistat{m}{i}^{-\sigma} \bigl(1+\abs{\ln \sidedistat{m}{i}}\bigr)^{-B} }_{L^\infty(0,\ellvtx)} <\infty.
 \label{eq:remainder-trace-Hhalf-hypothesis}
\end{align}
Define $\mathcal R(V_m+\sidedistat{m}{i}\tphys{m}{i}):=\mathcal R_{m,i}(\sidedistat{m}{i})$ and $\mathcal R(V_m):=0$. Then:
\begin{enumerate}[label=\textup{(\roman*)},leftmargin=*]
\item The two ray values tend to zero at $V_m$, and $\nm{\mathcal R}_{H^{\frac{1}{2}}(\partial\Om\cap B_{\ellvtx}(V_m))}\le C_{\sigma,B,\omega_m,\ellvtx}M_{\rm tr}$.
\item If \eqref{eq:remainder-trace-Hhalf-hypothesis} holds at every vertex and $\mathcal R$ is $H^1$ on the portions of the open sides away from the vertices, then $\mathcal R\in H^{\frac{1}{2}}(\partial\Om)$.
\end{enumerate}
\end{lemma}

\begin{proof}
Write $s:=\sidedistat{m}{i}$ and $\Lambda_B(s):=\bigl(1+\abs{\ln s}\bigr)^B$. By continuity and \eqref{eq:remainder-trace-Hhalf-hypothesis}, $\abs{\mathcal R_{m,i}(s)}\le M_{\rm tr}s^\sigma\Lambda_B(s)$ for every $s\in(0,\ellvtx)$. In particular, $\mathcal R_{m,i}(s)\to0$. We first estimate the $H^{\frac12}$ seminorm on one ray. By symmetry, it suffices to consider $0<\varrho<s<\ellvtx$.

Near the diagonal, $\frac{s}{2}<\varrho<s$, all intermediate distances are comparable to $s$. The derivative bound gives $\babs{\mathcal R_{m,i}(s)-\mathcal R_{m,i}(\varrho)}\le C M_{\rm tr}s^{\sigma-1}\Lambda_B(s)(s-\varrho)$. This contributes at most $C M_{\rm tr}^2\int_0^{\ellvtx}s^{2\sigma-1}\Lambda_B\bigl(s\bigr)^2\,\ud s<\infty$ to the squared seminorm. The integral converges because $\sigma>0$.

On the remaining region $0<\varrho\le\frac{s}{2}$, we have $s-\varrho\simeq s$. The value bound in \eqref{eq:remainder-trace-Hhalf-hypothesis} and the integral estimate $\int_0^{\frac{s}{2}}\varrho^{2\sigma}\Lambda_B\bigl(\varrho\bigr)^2\,\ud\varrho\le C s^{2\sigma+1}\Lambda_B\bigl(s\bigr)^2$ give $\int_0^{\frac{s}{2}}\babs{\mathcal R_{m,i}(s)-\mathcal R_{m,i}(\varrho)}^2\bigl(s-\varrho\bigr)^{-2}\,\ud\varrho\le C M_{\rm tr}^2s^{2\sigma-1}\Lambda_B\bigl(s\bigr)^2$. Its integral in $s$ converges for the same reason.

We next estimate the differences between the two rays. Set $s=\sidedistat{m}{1}$ and $\varrho=\sidedistat{m}{2}$. Since the opening angle is fixed, we have
\begin{align}
 \abs{s\,\tphys{m}{1}-\varrho\,\tphys{m}{2}}^2 &= s^2+\varrho^2-2s\varrho\cos\omega_m \ge c_{\omega_m}(s^2+\varrho^2), \qquad c_{\omega_m}:=1-\max\{\cos\omega_m,0\}>0.
\end{align}
Both traces vanish at the vertex. Hence
\begin{align}
 &\int_0^{\ellvtx}\!\int_0^{\ellvtx} \frac{\abs{\mathcal R_{m,1}(s)-\mathcal R_{m,2}(\varrho)}^2} {\abs{s\,\tphys{m}{1}-\varrho\,\tphys{m}{2}}^2} \,\ud\varrho\,\ud s \le C M_{\rm tr}^2\int_0^{\ellvtx}\!\int_0^{\ellvtx} \frac{s^{2\sigma}\Lambda_B\bigl(s\bigr)^2+\varrho^{2\sigma}\Lambda_B\bigl(\varrho\bigr)^2} {s^2+\varrho^2} \,\ud\varrho\,\ud s<\infty
\end{align}
by the same integral estimate. Unequal limiting constants on the two rays would instead give a logarithmically divergent cross-ray seminorm. A common constant cancels in the trace differences.

The $L^2$ estimate follows from \eqref{eq:remainder-trace-Hhalf-hypothesis}. Combining the local estimates with a smooth partition of unity and the assumed $H^1$ regularity away from the vertices proves the global statement.
\end{proof}

Define the boundary remainder by
\begin{align}
 h_{\rm rem}&:={} h-H^{\rm sing}\big|_{\partial\Om}.
 \label{eq:finite-part-residual-trace}
\end{align}
Here $H^{\rm sing}\big|_{\partial\Om}$ means the explicit restriction to the open sides, and boundary equalities are understood almost everywhere. For a function in $H^1(\Om)$, the same restriction bar denotes its Sobolev trace.

Lemma~\ref{lem:finite-part-remainder-trace} places a boundary remainder of positive order in $H^{\frac12}(\partial\Om)$. Negative powers, exponent-zero logarithms and unequal constant limits are included in $H^{\rm sing}$ before solving for the $H^1$ correction.

We now define the harmonic solution with this prescribed singular part.
\begin{definition}
\label{def:finite-part-solution-comp}
A \emph{harmonic solution with prescribed singular part} and boundary data $h$ is a function $u=H^{\rm sing}+v^{\rm reg}$ on $\Om\setminus\{V_1,\ldots,V_{\Nvtx}\}$. Here $H^{\rm sing}$ is fixed by \eqref{eq:finite-part-cutoff-realization}, and $v^{\rm reg}\in H^1(\Om)$ is required to satisfy
\begin{align}
 v^{\rm reg}\big|_{\partial\Om} &= h_{\rm rem},
 \label{eq:finite-part-variational-trace}\\
 \int_\Om\big(\nabla v^{\rm reg}\cdot\nabla\phi\big)\,\ud x &= \dualp{f^{\rm cut}}{\phi}, \qquad \phi\in H_0^1(\Om).
 \label{eq:finite-part-variational-comp}
\end{align}
Here $h_{\rm rem}$ and $f^{\rm cut}$ are defined by \eqref{eq:finite-part-residual-trace} and \eqref{eq:finite-part-commutator-source}, with the lifting normalization fixed as in Section~\ref{subsec:harmonic-lifts-comp}.
\end{definition}

We require $h_{\rm rem}\in H^{\frac12}(\partial\Om)$. A displayed term may be nonintegrable at $V_m$, but its harmonic lift is defined for $\bulkdistat{m}>0$ and has the prescribed restriction to each open side.

The energy estimates for $Z_{m,n}^{\pm}$ in \eqref{eq:intro-zero-ray-energy-test-comp}--\eqref{eq:intro-negative-ray-energy-test-comp} explain the role of this normalization. An additional $Z_{m,n}^-$ term has infinite Dirichlet energy and cannot occur in $v^{\rm reg}$. The functions $Z_{m,n}^+$ have finite energy. Their coefficients are determined by the global Dirichlet problem for $v^{\rm reg}$. Negative powers and zero-trace modes required by a resonant logarithmic lift are already included in $H^{\rm sing}$.

For each $Z_{m,n}^+$, the derivative estimate is $D_x^aZ_{m,n}^+=O\big(\bulkdistat{m}^{\lambda_{m,n}-a}\big)$.

\begin{lemma}
\label{lem:finite-part-uniqueness-comp}
Assume $h_{\rm rem}\in H^{\frac12}(\partial\Om)$.
\begin{enumerate}[label=\textup{(\roman*)},leftmargin=*]
\item The problem \eqref{eq:finite-part-variational-trace}--\eqref{eq:finite-part-variational-comp} has a unique solution $v^{\rm reg}\in H^1(\Om)$, with
\begin{align}
 \nm{v^{\rm reg}}_{H^1(\Om)} &\le C\left( \nm{h_{\rm rem}}_{H^{\frac12}(\partial\Om)} + \nm{f^{\rm cut}}_{H^{-1}(\Om)} \right).
 \label{eq:finite-part-energy-estimate}
\end{align}
\item The harmonic solution $u=H^{\rm sing}+v^{\rm reg}$ is independent of the radial cutoffs. If $\widetilde u=H^{\rm sing}+\widetilde v$ is harmonic on $\Om\setminus\{V_1,\ldots,V_{\Nvtx}\}$, with $\widetilde v\in H^1(\Om)$ and $\widetilde v\big|_{\partial\Om}=h_{\rm rem}$, then $\widetilde u=u$.
\item Adding finitely many normalized lifts in $H^1$ to the prescribed singular part leaves $u$ unchanged. Such lifts include positive-exponent terms with logarithmic factors and a constant common to both rays. If $L$ is the sum of their cutoff lifts, the new singular part and regular correction are $H^{\rm sing}+L$ and $v^{\rm reg}-L$.
\end{enumerate}
\end{lemma}

\begin{proof}\leavevmode
\paragraph{\underline{Existence and energy estimate}} The trace theorem on the bounded Lipschitz domain $\Om$ gives $v_{\rm ext}\in H^1(\Om)$ with
\begin{align}
 v_{\rm ext}\big|_{\partial\Om} &= h_{\rm rem}, \qquad \nm{v_{\rm ext}}_{H^1(\Om)} \le C\nm{h_{\rm rem}}_{H^{\frac12}(\partial\Om)}.
 \label{eq:finite-part-trace-extension}
\end{align}
Write $v^{\rm reg}=v_{\rm ext}+z$ with $z\in H_0^1(\Om)$. The equation becomes $\int_\Om\nabla z\cdot\nabla\phi\,\ud x=\dualp{f^{\rm cut}}{\phi}-\int_\Om\nabla v_{\rm ext}\cdot\nabla\phi\,\ud x$ for $\phi\in H_0^1(\Om)$. Its right-hand side is bounded by $C\big(\nm{f^{\rm cut}}_{H^{-1}(\Om)}+\nm{v_{\rm ext}}_{H^1(\Om)}\big)\nm{\phi}_{H_0^1(\Om)}$, and the Dirichlet form is coercive by Poincar\'e's inequality. Lax--Milgram gives a unique $z$ with the corresponding bound. Together with \eqref{eq:finite-part-trace-extension}, this proves \eqref{eq:finite-part-energy-estimate}. By \eqref{eq:cutoff-commutator-comp}, $f^{\rm cut}\in L^2(\Om)\subset H^{-1}(\Om)$. Two corrections differ by an element of $H_0^1(\Om)$ orthogonal to $H_0^1(\Om)$ in the Dirichlet form, so they coincide.

Also $\Delta H^{\rm sing}=f^{\rm cut}$ and $-\Delta v^{\rm reg}=f^{\rm cut}$ away from the vertices, so $u$ is harmonic there. Its boundary restriction is defined locally on the open sides by $u\big|_{\partial\Om}:=H^{\rm sing}\big|_{\partial\Om}+v^{\rm reg}\big|_{\partial\Om}=h$ almost everywhere.

\paragraph{\underline{Uniqueness and independence of the cutoffs}} For two admissible cutoff families, the singular parts differ by a smooth function vanishing near every vertex. Their difference lies in $H^1(\Om)$, so the $H^{\frac12}$ condition on the boundary remainder is preserved. The resulting solutions satisfy $\delta u:=u_1-u_2\in H_0^1(\Om)$ and $\Delta\delta u=0$ in $\mathcal D'(\Om)$. Density permits testing with $\delta u$, which gives $\nm{\nabla\delta u}_{L^2(\Om)}^2=0$ and hence $u_1=u_2$. The same argument applies to $\widetilde u-u$.

\paragraph{\underline{Enlargement of the prescribed part}} Let $L_m^{\rm loc}$ be the sum of the added normalized harmonic lifts and $L:=\sum_m\chi_m^{\rm sing}L_m^{\rm loc}$. For each positive-exponent term, its squared gradient is integrable because $\int_0^{\ellvtx}r^{2\lambda-1}\bigl(1+\abs{\ln r}\bigr)^{2c}\,\ud r<\infty$ for $\lambda>0$. A common constant has a smooth cutoff lift. Thus $L\in H^1(\Om)$, and its Sobolev trace $L\big|_{\partial\Om}$ agrees almost everywhere with its explicit restriction to the open sides.

Each $L_m^{\rm loc}$ is harmonic in the open sector. Testing against $C_c^\infty(\Om)$ therefore gives $\Delta L=\sum_m[\Delta,\chi_m^{\rm sing}]L_m^{\rm loc}$ as an interior distribution. The commutator source is supported on fixed annuli, so density extends the identity to
\begin{align}
 \int_\Om\nabla L\cdot\nabla\phi\,\ud x &= -\dualp{\textstyle\sum_m[\Delta,\chi_m^{\rm sing}]L_m^{\rm loc}}{\phi}, \qquad \phi\in H_0^1(\Om).
 \label{eq:enlargement-transfer-identity}
\end{align}
The enlarged problem has singular part $H^{\rm sing}+L$ and boundary remainder $h_{\rm rem}-L\big|_{\partial\Om}$. Its source is $f^{\rm cut}+\Delta L$. By \eqref{eq:enlargement-transfer-identity}, $v^{\rm reg}-L$ has this trace and satisfies its variational equation. Uniqueness identifies it as the enlarged regular correction, so the total solution remains $(H^{\rm sing}+L)+(v^{\rm reg}-L)=u$.
\end{proof}

\subsection{Weighted Vertex Expansion}
\label{subsec:mellin-mapping-comp}

We shall prove that the solution is the sum of the prescribed harmonic lifts, the modes $Z_{m,n}^+$ with $0<\lambda_{m,n}<\sigma$, and a remainder of weight $\sigma$. Proposition~\ref{prop:mellin-mapping-comp} estimates their coefficients and the remainder in terms of the boundary data. Here $\sigma$ is positive and does not belong to the Dirichlet pencil spectrum.

We first give an elliptic estimate on rescaled annuli. It controls the H\"older norm by the $L^2$ norm of the solution and a H\"older norm of the source. We shall also use this estimate in Section~\ref{sec:interior-side-hierarchy-comp}. Write $\mathscr A_\omega:=K_\omega\cap\{\frac12<\abs y<2\}$ and $\mathscr A_\omega^*:=K_\omega\cap\{\frac14<\abs y<4\}$ for the reference annuli of the sector $K_\omega$.

\begin{lemma}
\label{lem:rescaled-annulus}
Let $q\ge2$ be an integer and $0<\gamma<1$. Suppose $f\in C^{q-2,\gamma}(\overline{\mathscr A_\omega^*})$. There is $C=C(\omega,q,\gamma)$ such that every $v\in H^1(\mathscr A_\omega^*)$ with $-\Delta v=f$ in $\mathscr A_\omega^*$ and $v\big|_{\partial K_\omega\cap\overline{\mathscr A_\omega^*}}=0$ in the Sobolev trace sense satisfies
\begin{align}
 \nm v_{C^{q,\gamma}(\overline{\mathscr A_\omega})} &\le C\left( \nm v_{L^2(\mathscr A_\omega^*)} +\nm f_{C^{q-2,\gamma}(\overline{\mathscr A_\omega^*})} \right).
 \label{eq:rescaled-annulus-estimate}
\end{align}
\end{lemma}

\begin{proof}
Cover $\overline{\mathscr A_\omega}$ by finitely many interior balls and boundary half-balls. An interior ball $B(y_0,r)$ is chosen so that $B(y_0,2r)\subset\mathscr A_\omega^*$. A boundary half-ball is centered on $\partial K_\omega\cap\overline{\mathscr A_\omega}$. Its double meets $\partial K_\omega$ in a single flat piece and is otherwise contained in $\mathscr A_\omega^*$. Such a cover exists because $\overline{\mathscr A_\omega}$ stays at positive distance from the vertex and from the artificial circles $\abs y=\frac14,4$.

Work on each chart with three concentric radii $r<\frac32r<2r$. The interior and flat-boundary $W^{2,2}$ estimates for the Dirichlet Laplacian control $\nm v_{W^{2,2}}$ on the middle chart by $\nm v_{L^2}+\nm f_{L^2}$ on the largest. In two dimensions, $W^{2,2}\hookrightarrow C^{0,\gamma}$ for every $0<\gamma<1$. This gives the $C^{0,\gamma}$ bound needed for the next estimate.

The interior and flat-boundary Schauder estimates on the smallest chart then give the $C^{q,\gamma}$ bound in terms of $\nm v_{L^\infty}$ and $\nm f_{C^{q-2,\gamma}}$ on the middle one. Since $C^{q-2,\gamma}(\overline{\mathscr A_\omega^*})\hookrightarrow L^2$ on the bounded set $\mathscr A_\omega^*$, both source norms are controlled by the right-hand side of \eqref{eq:rescaled-annulus-estimate}. Summing over the finite cover gives \eqref{eq:rescaled-annulus-estimate}.
\end{proof}

After rescaling by $y\mapsto\varrho y$, the constant in \eqref{eq:rescaled-annulus-estimate} is independent of $\varrho$.

\begin{proposition}[Weighted expansion of a harmonic solution]
\label{prop:mellin-mapping-comp}
Fix $\Om$, its sector charts, the ray-lifting convention of Section~\ref{subsec:harmonic-lifts-comp}, and finite power-logarithmic index sets $I_{\ell,i}$ at every $V_\ell$, $1\le \ell\le\Nvtx$, $i=1,2$. Let $h$ be boundary data on the open sides. For the prescribed ray terms used to form $H^{\rm sing}$, denote the coefficient of $s^\lambda(\ln s)^b$ at $V_\ell$ on ray $i$ by $a_{\ell,i,\lambda,b}$, $(\lambda,b)\in I_{\ell,i}$. These terms determine the lifts $H_{{\rm ray},\ell}$ at all vertices. At the fixed vertex $V_m$, impose the local remainder condition below, together with the stated global data bounds:

\begin{enumerate}[label=\textup{(A\arabic*)},leftmargin=*,itemsep=0.8em]

\item \emph{Geometry and weight.} Fix the vertex $V_m$. After translation and rotation, $V_m=0$ and $\Om\cap B_{\ellvtx}=K_{\omega_m}\cap B_{\ellvtx}$ with $0<\omega_m<2\pi$, $K_{\omega_m}$ the sector of \eqref{eq:wedge-sector-def}.
Let $q\ge2$ be an integer, $0<\gamma<1$ and $B\in\mathbb N_0$. Let $\sigma>0$ be noncritical, so $\sigma\ne\lambda_{m,n}$ for every $n\in\mathbb N$. Choose $\chi_m^{\rm sing}=1$ on $B_{\ellvtx}(V_m)$, with the other singular-part cutoffs zero there.

\item \emph{Ray expansions at $V_m$.} Assume $I_{m,i}\subset\{(\lambda,b):\lambda<\sigma,\ b\in\mathbb N_0\}$ for $i=1,2$, and write the boundary data on the two incident rays as
\begin{align}
 h_{m,i}(\sidedistat{m}{i}) &= \sum_{(\lambda,b)\in I_{m,i}} a_{m,i,\lambda,b}\, \sidedistat{m}{i}^\lambda \bigl(\ln \sidedistat{m}{i}\bigr)^b + \mathcal R_{m,i}(\sidedistat{m}{i}), \qquad \mathcal R_{m,i} \in \mathcal C_{\sigma,B}^{q,\gamma}(0,\ellvtx),
 \label{eq:mellin-ray-data-comp}
\end{align}
for $i=1,2$, where an exponent occurring on one ray only has coefficient zero on the other.

\item \emph{Global data norms.} Let $H^{\rm sing}$ and $f^{\rm cut}$ be as in \eqref{eq:finite-part-cutoff-realization}--\eqref{eq:finite-part-commutator-source} and $h_{\rm rem}$ as in \eqref{eq:finite-part-residual-trace}, where $H^{\rm sing}\big|_{\partial\Om}$ denotes its explicit restriction to the open sides. Assume that
\begin{align}
 \mathcal N_{\rm glob} :={}& \nm{h}_{C^{q,\gamma} \big( \partial\Om \setminus \bigcup_{\ell=1}^{\Nvtx}B_{\frac{\ellvtx}{4}}(V_\ell) \big)}+ \nm{h_{\rm rem}}_{H^{\frac{1}{2}}(\partial\Om)} + \nm{f^{\rm cut}}_{C^{q-2,\gamma}(\overline\Om)}
 \label{eq:explicit-Mellin-far-norm}
\end{align}
is finite. Define the sum of the absolute values of the expansion coefficients by
\begin{align}
 \mathcal A_{\rm coeff}(h)&:=\sum_{\ell=1}^{\Nvtx}\sum_{i=1}^2\sum_{(\lambda,b)\in I_{\ell,i}}\abs{a_{\ell,i,\lambda,b}}.
 \label{eq:elliptic-germ-data-size}
\end{align}
The sum is taken over all vertices, since the Dirichlet problem uses the data on the whole boundary. Define
\begin{align}
 \EllData_{m,\sigma}(h)&:=\sum_{i=1}^2\nm{\mathcal R_{m,i}}_{\mathcal C_{\sigma,B}^{q,\gamma}(0,\ellvtx)}+\mathcal N_{\rm glob}+\mathcal A_{\rm coeff}(h).
 \label{eq:elliptic-total-data-size}
\end{align}
\end{enumerate}

As in Subsection~\ref{subsec:finite-part-comp}, $H_{{\rm ray},m}$ is the sum of the normalized harmonic lifts of the displayed terms in \eqref{eq:mellin-ray-data-comp}. The formulas in Section~\ref{subsec:harmonic-lifts-comp} make its coefficients linear in $a_{m,i,\lambda,b}$. The conclusions are the following:
\begin{enumerate}[label=\textup{(R\arabic*)},leftmargin=*,itemsep=0.8em]

\item \emph{Vertex expansion.} Let $u$ be the solution with data $h$ and prescribed singular part $H^{\rm sing}$. In $K_{\omega_m}\cap B_{\frac{\ellvtx}{2}}$, with the angular modes of \eqref{eq:intro-zero-ray-pairs-comp}, $u$ has the expansion
\begin{align}
 u(\bulkdistat{m},\theta) &= H_{{\rm ray},m}(\bulkdistat{m},\theta) + \sum_{0<\lambda_{m,n}<\sigma} c_{m,n}[h]\, \bulkdistat{m}^{\lambda_{m,n}}\Phi_{m,n}(\theta) + \mathcal R_{u,m}(\bulkdistat{m},\theta).
 \label{eq:mellin-expansion-comp}
\end{align}

\item \emph{Remainder and coefficient estimate.} The remainder belongs to $\mathcal C_{\sigma,B+1}^{q,\gamma}\bigl(K_{\omega_m}\cap B_{\frac{\ellvtx}{2}}\bigr)$ and satisfies
\begin{align}
 \nm{\mathcal R_{u,m}}_{\mathcal C_{\sigma,B+1}^{q,\gamma}} + \sum_{0<\lambda_{m,n}<\sigma} \abs{c_{m,n}[h]} &\le C\,\EllData_{m,\sigma}(h).
 \label{eq:mellin-remainder-bound-comp}
\end{align}
The constant depends on the fixed polygon, vertex radius $\ellvtx$, finite expansion index sets and exponents, and cutoff and ray-lifting conventions. It also depends on $q,\gamma,\sigma,B$ and the distance of $\sigma$ from the Dirichlet pencil spectrum. It is independent of $h$. The proof gives the remainder in $\mathcal C_{\sigma,B}^{q,\gamma}$. The stated estimate with $B+1$ is the form used by the order-dependent bounds \eqref{eq:order-dependent-log-budgets}.

Here $\mathcal R_{u,m}$ is the difference between $u$ and the two finite sums in \eqref{eq:mellin-expansion-comp}. In particular, \eqref{eq:weighted-sector-pointwise-comp} gives, for every derivative order $0\le a\le q$,
\begin{align}
 \abs{D_x^a\mathcal R_{u,m}(x)}&\le C\EllData_{m,\sigma}(h)\bulkdistat{m}^{\sigma-a}\bigl(1+\abs{\ln\bulkdistat{m}}\bigr)^{B+1},\qquad 0<\bulkdistat{m}<\frac{\ellvtx}{2}.
 \label{eq:mellin-remainder-derivative-bound}
\end{align}
Thus differentiating the solution requires differentiating the displayed harmonic terms and this specific remainder. The exponent of each harmonic term decreases by the derivative order, and the remainder has the controlled weight $\sigma-a$.

\item \emph{The vertex coefficients.} For a fixed polygon and lifting normalization, each map $h\mapsto c_{m,n}[h]$ is a continuous linear functional of the boundary data $h$. Continuity follows from \textup{(R2)}. The coefficients satisfy
\begin{align}
 c_{m,n}[h] &= \frac{2}{\omega_m} \lim_{\bulkdistat{m}\downarrow0} \bulkdistat{m}^{-\lambda_{m,n}} \int_0^{\omega_m} \bigl(u-H_{{\rm ray},m}\bigr)(\bulkdistat{m},\theta) \Phi_{m,n}(\theta)\,\ud\theta, \qquad 0<\lambda_{m,n}<\sigma.
 \label{eq:global-zero-ray-coefficient}
\end{align}
Thus $c_{m,n}[h]$ is independent of the auxiliary cutoffs and of the extension of the ray remainders. It is also the same for two weights strictly greater than $\lambda_{m,n}$, provided that the expansions represent the same boundary data, use the same lifting normalization, and have the same prescribed ray coefficients at every exponent below the smaller weight.
\end{enumerate}
\end{proposition}

\begin{remark}
\label{rem:mellin-hypothesis-usage}
The radius $\ellvtx$ may be replaced throughout by any fixed $\ell'\in[\ellvtx,4\ellvtx]$ for which the sector identity and cutoff conditions in \textup{(A1)} hold on $B_{\ell'}(V_m)$. In that case, the ray remainders are measured on $(0,\ell')$. The sidewise H\"older term in \eqref{eq:explicit-Mellin-far-norm} excludes the balls $B_{\frac{\ell'}{4}}(V_\ell)$, and the conclusion holds on $K_{\omega_m}\cap B_{\frac{\ell'}{2}}$. The constants may depend on this fixed radius. Section~\ref{sec:interior-side-hierarchy-comp} uses $\ell'=4\ellvtx$.

Hypothesis \textup{(A2)} supplies the weighted local trace extension and the decay of its Laplacian. Hypothesis \textup{(A3)} controls the global variational correction and the annular localization source through $\nm{h_{\rm rem}}_{H^{\frac12}(\partial\Om)}$ and $\nm{f^{\rm cut}}_{L^2(\Om)}$ in \eqref{eq:finite-part-global-energy-bound}. These are controlled by the second and third terms of \eqref{eq:explicit-Mellin-far-norm}, and the proof below needs no other global data bound. The H\"older norm of $h$ away from the vertices and the coefficient sum $\mathcal A_{\rm coeff}(h)$ are included in $\EllData_{m,\sigma}(h)$ because the induction in Section~\ref{sec:interior-side-hierarchy-comp} estimates them together with the local remainders. Lemma~\ref{lem:finite-part-remainder-trace} verifies the local contribution to the $H^{\frac12}$ trace condition.
\end{remark}

\begin{remark}
\label{rem:local-versus-global-coefficients}
The two ray expansions at $V_m$ determine $H_{{\rm ray},m}$. The coefficients $c_{m,n}[h]$ also depend on boundary data away from $V_m$, by \textup{(R3)}. Changing data on a side separated from $V_m$ can therefore change $c_{m,n}[h]$ while preserving $H_{{\rm ray},m}$. In particular, boundary data vanishing near the vertex can produce a nonzero coefficient of $Z_{m,n}^+$.
\end{remark}

\begin{proof}
We first use the variational problem to obtain the harmonic solution and its global energy bound. Near $V_m$, subtracting $H_{{\rm ray},m}$ and an angular extension of the ray remainders gives an $H^1$ function with zero trace on both rays. A radial cutoff extends this function to the full sector, and the change of variables $t=-\ln\bulkdistat{m}$ gives an equation on the whole cylinder. Its source has a decaying part from the ray remainder and a part supported in a bounded interval from the cutoff. The latter depends on the global solution. The Green formula identifies the modes with exponents less than $\sigma$. We subtract those modes, estimate the remaining angular sum in $L^2$, and then use Lemma~\ref{lem:rescaled-annulus} to obtain the weighted H\"older estimate.

\paragraph{\underline{Step 1: reduction to homogeneous ray data}} By the $H^{\frac12}$ term in \textup{(A3)}, Lemma~\ref{lem:finite-part-uniqueness-comp} gives a unique $v^{\rm reg}$ and $u=H^{\rm sing}+v^{\rm reg}$, with
\begin{align}
 \nm{v^{\rm reg}}_{H^1(\Om)} &\le C\left( \nm{h_{\rm rem}}_{H^{\frac12}(\partial\Om)} + \nm{f^{\rm cut}}_{H^{-1}(\Om)} \right) \le C\EllData_{m,\sigma}(h).
 \label{eq:finite-part-global-energy-bound}
\end{align}

Set $K_{m,\ellvtx}:=\Om\cap B_{\ellvtx}(V_m)=K_{\omega_m}\cap B_{\ellvtx}$. By the definition of $H_{{\rm ray},m}$, the restrictions of $u-H_{{\rm ray},m}$ to the two rays are $\mathcal R_{m,1}$ and $\mathcal R_{m,2}$. Define $E_m^{\rm ang}(\bulkdistat{m},\theta):=\bigl(1-\frac{\theta}{\omega_m}\bigr)\mathcal R_{m,1}(\bulkdistat{m})+\frac{\theta}{\omega_m}\mathcal R_{m,2}(\bulkdistat{m})$. Since $\bulkdistat{m}=\sidedistat{m}{i}$ on ray $i$, this function has the same ray restrictions. We subtract it to obtain zero boundary data on the rays.

We first estimate this extension. Write $\bulkdistat{m}=\varrho s$ at a dyadic scale $\varrho$. The identity $\bulkdistat{m}^2\Delta=\bigl(\bulkdistat{m}\partial_{\bulkdistat{m}}\bigr)^2+\partial_\theta^2$ shows that both $E_m^{\rm ang}$ and $\bulkdistat{m}^2\Delta E_m^{\rm ang}$ are affine angular combinations of $\mathcal R_{m,i}(\varrho s)$ and $\bigl(s\partial_s\bigr)^2\mathcal R_{m,i}(\varrho s)$. On $\frac12<s<2$, the polar and Cartesian $C^{q,\gamma}$ norms are uniformly equivalent. Multiplying by $\varrho^{-\sigma}\bigl(1+\abs{\ln\varrho}\bigr)^{-B}$ and taking the supremum in $\varrho$, we obtain
\begin{align}
 \nm{E_m^{\rm ang}}_{\mathcal C_{\sigma,B}^{q,\gamma}} + \nm{\bulkdistat{m}^2\Delta E_m^{\rm ang}}_{\mathcal C_{\sigma,B}^{q-2,\gamma}} &\le C \sum_{i=1}^2 \nm{\mathcal R_{m,i}}_{\mathcal C_{\sigma,B}^{q,\gamma}(0,\ellvtx)} \le C\EllData_{m,\sigma}(h).
 \label{eq:weighted-ray-remainder-lift-bound}
\end{align}

The estimate $\abs{E_m^{\rm ang}}+\bulkdistat{m}\abs{\nabla E_m^{\rm ang}}\le C\EllData_{m,\sigma}(h)\, \bulkdistat{m}^\sigma\bigl(1+\abs{\ln \bulkdistat{m}}\bigr)^B$ and $\sigma>0$ then give
\begin{align}
 \nm{E_m^{\rm ang}}_{H^1(K_{m,\ellvtx})}^2 &\le C\EllData_{m,\sigma}\bigl(h\bigr)^2 \int_0^{\ellvtx} \left( \bulkdistat{m}^{2\sigma+1}+\bulkdistat{m}^{2\sigma-1} \right) \bigl(1+\abs{\ln \bulkdistat{m}}\bigr)^{2B}\,\ud \bulkdistat{m} \le C\EllData_{m,\sigma}\bigl(h\bigr)^2
 \label{eq:weighted-ray-lift-H1}
\end{align}
so $E_m^{\rm ang}\in H^1(K_{m,\ellvtx})$ with $\nm{E_m^{\rm ang}}_{H^1(K_{m,\ellvtx})}\le C\EllData_{m,\sigma}(h)$.

By \textup{(A1)}, $H^{\rm sing}=H_{{\rm ray},m}$ in $K_{m,\ellvtx}$, and hence $u-H_{{\rm ray},m}=v^{\rm reg}$ a.e.\ there. Set $W_m:=u-H_{{\rm ray},m}-E_m^{\rm ang}=v^{\rm reg}-E_m^{\rm ang}$. Estimates \eqref{eq:finite-part-global-energy-bound} and \eqref{eq:weighted-ray-lift-H1} give $W_m\in H^1(K_{m,\ellvtx})$. The Sobolev traces of $v^{\rm reg}$ and $E_m^{\rm ang}$ both equal $\mathcal R_{m,i}$ on ray $i$, since the normalized lift reproduces the displayed terms of $h_{m,i}$. Hence $W_m\big|_{\theta=0}=W_m\big|_{\theta=\omega_m}=0$.

Choose $\chi_m\in C_c^\infty([0,\infty))$ with $\chi_m=1$ for $\bulkdistat{m}<\frac{\ellvtx}{2}$ and $\supp\chi_m\subset\{\bulkdistat{m}<\frac{3\ellvtx}{4}\}$. Set $z_m:=\chi_m W_m$ and $\mathcal G_m:=-\bulkdistat{m}^2\Delta z_m$. The function $z_m$ has finite energy and zero trace on $\theta=0,\omega_m$. Extend it by zero radially to the full sector. Harmonicity of $u$ and $H_{{\rm ray},m}$ on the punctured neighborhood gives $\mathcal G_m=\mathcal G_m^{\rm near}-\mathcal G_m^{\rm ann}$. Here $\mathcal G_m^{\rm near}:=\chi_m\bulkdistat{m}^2\Delta E_m^{\rm ang}$ has weight $\bulkdistat{m}^\sigma\bigl(1+\abs{\ln\bulkdistat{m}}\bigr)^B$, and $\mathcal G_m^{\rm ann}:=\bulkdistat{m}^2[\Delta,\chi_m]W_m$ is supported in a fixed annulus separated from $V_m$.

Choose the sector annuli $A_*:=K_{\omega_m}\cap\{\frac{\ellvtx}{3}<\bulkdistat{m}<\frac{7\ellvtx}{8}\}$ and $A^*:=K_{\omega_m}\cap\{\frac{\ellvtx}{4}<\bulkdistat{m}<\frac{15\ellvtx}{16}\}$. Their closures include the straight ray boundaries, and $\supp\nabla\chi_m\cap\overline\Om\subset\overline{A_*}$. The smaller annulus is separated from the vertex and from the two circular boundaries of $A^*$. On $A^*$, $W_m$ has zero trace on the straight sides and satisfies $-\Delta W_m=\Delta E_m^{\rm ang}$. The interior and flat-boundary covering argument in Lemma~\ref{lem:rescaled-annulus} gives $\nm{W_m}_{C^{q,\gamma}(\overline{A_*})}\le C(\nm{W_m}_{L^2(A^*)}+\nm{\Delta E_m^{\rm ang}}_{C^{q-2,\gamma}(\overline{A^*})})\le C\EllData_{m,\sigma}(h)$ by \eqref{eq:finite-part-global-energy-bound}, \eqref{eq:weighted-ray-lift-H1} and \eqref{eq:weighted-ray-remainder-lift-bound}. Thus both parts of the localized source satisfy
\begin{align}
 \nm{\mathcal G_m^{\rm near}} _{\mathcal C_{\sigma,B}^{q-2,\gamma}} + \nm{\mathcal G_m^{\rm ann}} _{C^{q-2,\gamma}} \le C\EllData_{m,\sigma}(h).
 \label{eq:localized-mellin-source-bound}
\end{align}

\paragraph{\underline{Step 2: the equation on the cylinder}} We use the change of variables
\begin{align}
 t:=-\ln\bulkdistat{m}, \qquad z_m^{\rm cyl}(t,\theta) &:=z_m(\ue^{-t},\theta).
 \label{eq:logarithmic-cylinder-change-comp}
\end{align}
Then $\bulkdistat{m}\to0$ corresponds to $t\to\infty$, and a radial power $\bulkdistat{m}^\lambda$ becomes $\ue^{-\lambda t}$. By the choice of cutoff, $z_m^{\rm cyl}=0$ for $t<T_-:=\ln\frac{4}{3\ellvtx}$. Using polar coordinates and the zero trace on the rays, we obtain
\begin{align}
 \int_{\mathbb R}\int_0^{\omega_m} \left(\abs{\partial_tz_m^{\rm cyl}}^2+ \abs{\partial_\theta z_m^{\rm cyl}}^2\right) \,\ud\theta\,\ud t &=\nm{\nabla z_m}_{L^2(K_{\omega_m})}^2,
 \label{eq:cylindrical-gradient-identity-comp}\\
 \nm{z_m^{\rm cyl}}_{L^2(\mathbb R\times(0,\omega_m))}^2 &\le\lambda_{m,1}^{-2} \nm{\partial_\theta z_m^{\rm cyl}}_{L^2(\mathbb R\times(0,\omega_m))}^2.
 \label{eq:cylindrical-poincare-control-comp}
\end{align}
Thus $z_m^{\rm cyl}\in H^1(\mathbb R\times(0,\omega_m))$, with norm bounded by $C\EllData_{m,\sigma}(h)$. Write $\mathcal G_m^{\rm cyl}(t,\theta):=\mathcal G_m(\ue^{-t},\theta)$. Since $\bulkdistat{m}^2\Delta=\partial_t^2+\partial_\theta^2$, the localized equation becomes
\begin{align}
 \left(-\partial_t^2-\partial_\theta^2\right)z_m^{\rm cyl} &=\mathcal G_m^{\rm cyl} \quad\text{on }\mathbb R\times(0,\omega_m),
 \label{eq:cylindrical-source-equation-comp}
\end{align}
with zero trace on $\theta=0,\omega_m$.

The Dirichlet Mellin pencil associated with \eqref{eq:elementary-vertex-pencil} is
\begin{align}
 \mathfrak A_{\omega_m}(\zeta) &:=\zeta^2+\partial_\theta^2: H^2(0,\omega_m)\cap H_0^1(0,\omega_m)\rightarrow L^2(0,\omega_m).
 \label{eq:dirichlet-mellin-pencil-comp}
\end{align}
Its roots are $\pm\lambda_{m,n}$. The corresponding solutions on the cylinder are $\ue^{\pm\lambda_{m,n}t}\Phi_{m,n}(\theta)$. We shall use the $H^1$ condition to determine the Green solution on the whole line and then extract its decaying terms as $t\to\infty$.

Write $z_{m,n}^{\rm cyl}(t):=\frac{2}{\omega_m}\int_0^{\omega_m}z_m^{\rm cyl}(t,\theta)\Phi_{m,n}(\theta)\,\ud\theta$ and $\mathcal G_{m,n}^{\rm cyl}$ likewise. Orthogonality and completeness give $\nm{z_m^{\rm cyl}(t)}_{L^2(0,\omega_m)}^2=\frac{\omega_m}{2}\sum_{n\ge1}\abs{z_{m,n}^{\rm cyl}(t)}^2$, and the same identity holds for $\mathcal G_m^{\rm cyl}$.

By \eqref{eq:localized-mellin-source-bound} and the support of the cutoff, there is $T_0>T_-$ with $T_0\geq1$ such that
\begin{align}
 \nm{\mathcal G_m^{\rm cyl}(t)}_{L^2(0,\omega_m)} \leq C\EllData_{m,\sigma}(h)& \begin{cases}
 0, & t<T_-,\\
 1, & T_-\leq t\leq T_0,\\
 \ue^{-\sigma t}\bigl(1+t\bigr)^B, & t\geq T_0.
 \end{cases}
 \label{eq:cylindrical-source-piecewise-bound-comp}
\end{align}
The interval $[T_-,T_0]$ contains the support of $\mathcal G_m^{\rm ann}$ in the $t$ variable. For $t\geq T_0$, the weighted estimate for $\mathcal G_m^{\rm near}$ gives the last bound.

Since $z_m^{\rm cyl}\in H^1(\mathbb R\times(0,\omega_m))$, each coefficient $z_{m,n}^{\rm cyl}$ belongs to $H^1(\mathbb R)$. Test \eqref{eq:cylindrical-source-equation-comp} against $\phi(t)\Phi_{m,n}(\theta)$, where $\phi\in C_c^\infty(\mathbb R)$. This test is justified by $H_0^1$ density on a bounded slab $(T_1,T_2)\times(0,\omega_m)$ with $\supp\phi\subset(T_1,T_2)$. This gives $-\partial_t^2z_{m,n}^{\rm cyl}+\lambda_{m,n}^2z_{m,n}^{\rm cyl}=\mathcal G_{m,n}^{\rm cyl}$ in $\mathcal D'(\mathbb R)$. The source belongs to $L^1(\mathbb R)\cap L^2(\mathbb R)$ by \eqref{eq:cylindrical-source-piecewise-bound-comp}. Its convolution with the Green kernel $\frac{\ue^{-\lambda_{m,n}\abs{t}}}{2\lambda_{m,n}}$ is another $H^1(\mathbb R)$ solution. The difference solves $-y''+\lambda_{m,n}^2y=0$ and belongs to $H^1(\mathbb R)$. Every solution of this homogeneous equation is a combination of $\ue^{\pm\lambda_{m,n}t}$, so the difference vanishes. Therefore
\begin{align}
 z_{m,n}^{\rm cyl}(t) &=\frac{\ue^{-\lambda_{m,n}t}}{2\lambda_{m,n}} \int_{-\infty}^t\ue^{\lambda_{m,n}s}\mathcal G_{m,n}^{\rm cyl}(s)\,\ud s +\frac{\ue^{\lambda_{m,n}t}}{2\lambda_{m,n}} \int_t^\infty\ue^{-\lambda_{m,n}s}\mathcal G_{m,n}^{\rm cyl}(s)\,\ud s.
 \label{eq:cylindrical-mode-Green-formula}
\end{align}

\paragraph{\underline{Step 3: coefficients and the $L^2$ remainder estimate}} We first consider $0<\lambda_{m,n}<\sigma$. These modes decay more slowly than $\ue^{-\sigma t}$ and appear explicitly in the expansion. Extending the first integral in \eqref{eq:cylindrical-mode-Green-formula} to $\mathbb R$, we define
\begin{align}
 c_{m,n}[h] :=\frac1{2\lambda_{m,n}} \int_{\mathbb R}\ue^{\lambda_{m,n}s} \mathcal G_{m,n}^{\rm cyl}(s)\,\ud s.
 \label{eq:positive-zero-ray-residue}
\end{align}
The integral converges by \eqref{eq:cylindrical-source-piecewise-bound-comp}, since its tail is bounded by $C\EllData_{m,\sigma}(h)\ue^{-(\sigma-\lambda_{m,n})s}(1+s)^B$. There are finitely many such modes, so
\begin{align}
 \sum_{0<\lambda_{m,n}<\sigma}\babs{c_{m,n}[h]} &\le C\EllData_{m,\sigma}(h).
 \label{eq:positive-zero-ray-coefficient-bound}
\end{align}
Subtracting $c_{m,n}[h]\ue^{-\lambda_{m,n}t}$ from \eqref{eq:cylindrical-mode-Green-formula}, we obtain
\begin{align}
 z_{m,n}^{\rm cyl}(t)-c_{m,n}[h]\ue^{-\lambda_{m,n}t}&=-\frac{\ue^{-\lambda_{m,n}t}}{2\lambda_{m,n}}\int_t^\infty\ue^{\lambda_{m,n}s}\mathcal G_{m,n}^{\rm cyl}(s)\,\ud s
 +\frac{\ue^{\lambda_{m,n}t}}{2\lambda_{m,n}}\int_t^\infty\ue^{-\lambda_{m,n}s}\mathcal G_{m,n}^{\rm cyl}(s)\,\ud s,
 \label{eq:low-mode-terminal-tail-formula-comp}
\end{align}
where the two integrals are controlled by $\sigma-\lambda_{m,n}>0$ and $\sigma+\lambda_{m,n}>0$, respectively. Since $B\in\mathbb N_0$, substituting $s=t+u$ and expanding the polynomial gives, for $a>0$ and $t\ge1$,
\begin{align}
 \int_t^\infty\ue^{-as}\bigl(1+s\bigr)^B\,\ud s
 &=\ue^{-at}\sum_{j=0}^B\binom{B}{j}\bigl(1+t\bigr)^{B-j}\frac{j!}{a^{j+1}}\le\ue^{-at}\bigl(1+t\bigr)^B\sum_{j=0}^B\binom{B}{j}\frac{j!}{a^{j+1}}.
\end{align}
Apply this bound with $a=\sigma-\lambda_{m,n}$ and $a=\sigma+\lambda_{m,n}$. After multiplying by the factors outside the integrals, both contributions have size $O\bigl(\ue^{-\sigma t}(1+t)^B\bigr)$ for $t\ge T_0$. This is the reason for subtracting the coefficient in \eqref{eq:positive-zero-ray-residue} before estimating the remainder.

We estimate the sum of the modes with $\lambda_{m,n}>\sigma$ in $L^2(0,\omega_m)$. The distance of $\sigma$ from the spectrum gives a bound uniform in $n$.

Let $\Pi_{m,<\sigma}$ be the orthogonal projection of $L^2(0,\omega_m)$ onto $\operatorname{span}\{\Phi_{m,n}:\lambda_{m,n}<\sigma\}$ and $\Pi_{m,>\sigma}:=I-\Pi_{m,<\sigma}$.

Noncriticality gives the gap $\delta_{m,\sigma}^{\rm spec}:=\dist(\sigma,\{\lambda_{m,n}:n\in\mathbb N\})>0$, hence $\sigma-\lambda_{m,n}\geq\delta_{m,\sigma}^{\rm spec}$ when $\lambda_{m,n}<\sigma$ and $\lambda_{m,n}\ge\sigma+\delta_{m,\sigma}^{\rm spec}$ when $\lambda_{m,n}>\sigma$.

Since there are finitely many modes with $\lambda_{m,n}<\sigma$, \eqref{eq:low-mode-terminal-tail-formula-comp} gives, for $t\geq T_0$,
\begin{align}
 \nm{ \Pi_{m,<\sigma}z_m^{\rm cyl}(t) - \sum_{\lambda_{m,n}<\sigma} c_{m,n}[h] \ue^{-\lambda_{m,n}t}\Phi_{m,n} }_{L^2(0,\omega_m)}&\leq C\EllData_{m,\sigma}(h) \ue^{-\sigma t}\bigl(1+t\bigr)^B.
 \label{eq:low-mode-summed-terminal-bound-comp}
\end{align}

For $\lambda_{m,n}>\sigma$, use \eqref{eq:cylindrical-mode-Green-formula} in its convolution form and Parseval's identity $\nm{\Pi_{m,>\sigma}z_m^{\rm cyl}(t)}_{L^2(0,\omega_m)}^2=\frac{\omega_m}{2}\sum_{\lambda_{m,n}>\sigma}\abs{z_{m,n}^{\rm cyl}(t)}^2$. The convolution kernel is at most $\frac{1}{2\sigma}\ue^{-(\sigma+\delta_{m,\sigma}^{\rm spec})\abs{t-s}}$. Minkowski's integral inequality therefore yields
\begin{align}
 \nm{\Pi_{m,>\sigma}z_m^{\rm cyl}(t)}_{L^2(0,\omega_m)} &\le \frac1{2\sigma}\int_{\mathbb R}\ue^{-(\sigma+\delta_{m,\sigma}^{\rm spec})\abs{t-s}} \nm{\mathcal G_m^{\rm cyl}(s)}_{L^2(0,\omega_m)}\,\ud s.
 \label{eq:high-mode-spectral-gap-comp}
\end{align}
For $t\geq T_0$, split the source estimate \eqref{eq:cylindrical-source-piecewise-bound-comp} at $T_0$. Its two contributions are controlled by
\begin{align}
 \int_{T_-}^{T_0}\ue^{-(\sigma+\delta_{m,\sigma}^{\rm spec})(t-s)}\nm{\mathcal G_m^{\rm cyl}(s)}_{L^2(0,\omega_m)}\,\ud s&\leq C\EllData_{m,\sigma}(h)\ue^{-(\sigma+\delta_{m,\sigma}^{\rm spec})(t-T_0)},\\
 \int_{T_0}^{\infty}\ue^{-(\sigma+\delta_{m,\sigma}^{\rm spec})\abs{t-s}}\ue^{-\sigma s}\bigl(1+s\bigr)^B\,\ud s&\leq C\ue^{-\sigma t}\bigl(1+t\bigr)^B.
\end{align}
For the second inequality, the part with $T_0\le s\le t$ is bounded by $\bigl(\delta_{m,\sigma}^{\rm spec}\bigr)^{-1}\ue^{-\sigma t}\bigl(1+t\bigr)^B$. For the part with $s\ge t$, apply the polynomial tail bound above with $a=2\sigma+\delta_{m,\sigma}^{\rm spec}$, including the factor $\ue^{(\sigma+\delta_{m,\sigma}^{\rm spec})t}$ outside the integral. Inserting the two bounds in \eqref{eq:high-mode-spectral-gap-comp} gives $C\EllData_{m,\sigma}(h)\ue^{-\sigma t}\bigl(1+t\bigr)^B$.

Define $z_{m,\rm rem}^{\rm cyl}(t,\theta):=z_m^{\rm cyl}(t,\theta)-\sum_{\lambda_{m,n}<\sigma}c_{m,n}[h]\ue^{-\lambda_{m,n}t}\Phi_{m,n}(\theta)$. Combining the preceding estimate with \eqref{eq:low-mode-summed-terminal-bound-comp}, we obtain
\begin{align}
 \nm{z_{m,\rm rem}^{\rm cyl}(t)}_{L^2(0,\omega_m)} &\leq C\EllData_{m,\sigma}(h) \ue^{-\sigma t}\bigl(1+t\bigr)^B, \qquad t\geq T_0.
 \label{eq:exact-terminal-L2-bound-comp}
\end{align}

\paragraph{\underline{Step 4: the weighted H\"older estimate}} Write $z_m^{\rm rem}(\bulkdistat{m},\theta):=z_{m,\rm rem}^{\rm cyl}(-\ln\bulkdistat{m},\theta)$. This function belongs to $H^1$ on every fixed annulus, since it is the localized $H^1$ function $z_m$ minus a finite sum of smooth harmonic functions there. We apply Lemma~\ref{lem:rescaled-annulus} using this annular regularity, \eqref{eq:exact-terminal-L2-bound-comp} and the source estimate \eqref{eq:localized-mellin-source-bound}. Choose $\varrho_0\leq\frac14\ue^{-T_0}$ so small that, for $0<\varrho<\varrho_0$, the annulus $\{\frac{\varrho}{4}<\bulkdistat{m}<4\varrho\}$ lies in the exact sector chart and in the region where $\chi_m=1$.

On the reference annuli $\mathscr A_{\omega_m}\subset\mathscr A_{\omega_m}^*$ of Lemma~\ref{lem:rescaled-annulus}, set $z_{m,\varrho}^{\rm rem}(y):=z_m^{\rm rem}(\varrho y)$. The subtracted modes are harmonic, and $\mathcal G_m=-\bulkdistat{m}^2\Delta_xz_m$. Consequently $-\Delta_yz_{m,\varrho}^{\rm rem}(y)=\abs{y}^{-2}\mathcal G_m(\varrho\abs{y},\arg y)=:f_\varrho(y)$ in $\mathscr A_{\omega_m}^*$. The Sobolev trace satisfies $z_{m,\varrho}^{\rm rem}\big|_{\partial K_{\omega_m}\cap\overline{\mathscr A_{\omega_m}^*}}=0$.

The annulus $\big\{\frac{\varrho}{4}<\bulkdistat{m}<4\varrho\big\}$ is covered by the dyadic annuli at scales $\frac{\varrho}{2},\varrho,2\varrho$. Their logarithmic weights are uniformly comparable because the scale ratios are fixed. Estimate \eqref{eq:localized-mellin-source-bound} gives $\nm{f_\varrho}_{C^{q-2,\gamma}(\overline{\mathscr A_{\omega_m}^*})}\leq C\EllData_{m,\sigma}(h)\varrho^\sigma\bigl(1+\abs{\ln\varrho}\bigr)^B$. Polar integration and \eqref{eq:exact-terminal-L2-bound-comp} also give
\begin{align}
 \nm{z_{m,\varrho}^{\rm rem}}_{L^2(\mathscr A_{\omega_m}^*)}^2&=\int_{\frac14}^{4}\nm{z_{m,\rm rem}^{\rm cyl}(-\ln(\varrho s))}_{L^2(0,\omega_m)}^2s\,\ud s\leq C\EllData_{m,\sigma}(h)^2\varrho^{2\sigma}\bigl(1+\abs{\ln\varrho}\bigr)^{2B}.
\end{align}

Apply Lemma~\ref{lem:rescaled-annulus} with $v=z_{m,\varrho}^{\rm rem}$ and $f=f_\varrho$. It bounds $\nm{z_{m,\varrho}^{\rm rem}}_{C^{q,\gamma}(\overline{\mathscr A_{\omega_m}})}$ by $C(\nm{z_{m,\varrho}^{\rm rem}}_{L^2(\mathscr A_{\omega_m}^*)}+\nm{f_\varrho}_{C^{q-2,\gamma}(\overline{\mathscr A_{\omega_m}^*})})$. The preceding estimates bound this sum by $C\EllData_{m,\sigma}(h)\varrho^\sigma\bigl(1+\abs{\ln\varrho}\bigr)^B$ for $0<\varrho<\varrho_0$.

On the compact range $\varrho_0\leq\varrho<\frac{\ellvtx}{4}$, the coefficient bound \eqref{eq:positive-zero-ray-coefficient-bound} controls the subtracted finite harmonic sum in $C^{q,\gamma}$. A finite covering by the charts of Lemma~\ref{lem:rescaled-annulus}, together with the $H^1$ bound from Steps~1--2 and \eqref{eq:localized-mellin-source-bound}, gives the same estimate because the radius is bounded away from zero. Since $\inf_{\varrho_0\leq\varrho<\frac{\ellvtx}{4}}\varrho^\sigma\bigl(1+\abs{\ln\varrho}\bigr)^B>0$, the two ranges combine into $\nm{z_m^{\rm rem}}_{\mathcal C_{\sigma,B}^{q,\gamma}(K_{\omega_m}\cap B_{\frac{\ellvtx}{2}})}\leq C\EllData_{m,\sigma}(h)$.

For $\bulkdistat{m}<\frac{\ellvtx}{2}$, Step~1 gives $u-H_{{\rm ray},m}=z_m+E_m^{\rm ang}$. Thus the remainder in \eqref{eq:mellin-expansion-comp} is $\mathcal R_{u,m}=u-H_{{\rm ray},m}-\sum_{\lambda_{m,n}<\sigma}c_{m,n}[h]\bulkdistat{m}^{\lambda_{m,n}}\Phi_{m,n}=z_m^{\rm rem}+E_m^{\rm ang}$. Combining the bound for $z_m^{\rm rem}$ with \eqref{eq:weighted-ray-remainder-lift-bound} proves the stronger dyadic H\"older estimate
\begin{align}
 \nm{\mathcal R_{u,m}} _{\mathcal C_{\sigma,B}^{q,\gamma} (K_{\omega_m}\cap B_{\frac{\ellvtx}{2}})} &\leq C\EllData_{m,\sigma}(h).
 \label{eq:strong-terminal-remainder-bound-comp}
\end{align}

As $\bigl(1+\abs{\ln\varrho}\bigr)^{-(B+1)}\leq\bigl(1+\abs{\ln\varrho}\bigr)^{-B}$, $\mathcal C_{\sigma,B}^{q,\gamma}\hookrightarrow\mathcal C_{\sigma,B+1}^{q,\gamma}$, so \eqref{eq:strong-terminal-remainder-bound-comp} gives the remainder estimate in \eqref{eq:mellin-remainder-bound-comp}.

The construction is linear in $h$, and the estimates prove continuity in \textup{(R3)}. Boundary data supported away from $V_m$ contribute to the localized source on a bounded $t$-interval. They can change $c_{m,n}[h]$ for $\lambda_{m,n}<\sigma$, while their contributions to modes with $\lambda_{m,n}>\sigma$ decay faster than $\ue^{-\sigma t}$. The exponents $\lambda_{m,n}$ of these zero-trace modes depend only on the angle.

\paragraph{\underline{Step 5: independence of the cutoffs and weight}} The definition \eqref{eq:positive-zero-ray-residue} uses the cutoff and the extension of the ray remainders from Step~1. To prove independence of these choices, we recover the coefficients directly from $u-H_{{\rm ray},m}$. By orthogonality in \eqref{eq:mellin-expansion-comp} and the remainder bound \eqref{eq:strong-terminal-remainder-bound-comp},
\begin{align}
 \frac2{\omega_m}\bulkdistat{m}^{-\lambda_{m,n}} \int_0^{\omega_m} (u-H_{{\rm ray},m})(\bulkdistat{m},\theta) \Phi_{m,n}(\theta)\,\ud\theta &=c_{m,n}[h]+ O\!\left(\EllData_{m,\sigma}(h)\bulkdistat{m}^{\sigma-\lambda_{m,n}} \bigl(1+\abs{\ln\bulkdistat{m}}\bigr)^B\right).
 \label{eq:global-zero-ray-coefficient-from-contour-comp}
\end{align}
The error tends to zero since $\lambda_{m,n}<\sigma$. This proves \eqref{eq:global-zero-ray-coefficient}, whose right-hand side contains only $u$ and the fixed normalized ray lift. By Lemma~\ref{lem:finite-part-uniqueness-comp}, $u$ is independent of the singular-part cutoffs. The extension of the ray remainders does not enter this formula. Hence $c_{m,n}[h]$ is independent of both choices.

Finally, compare expansions of the same boundary data at noncritical weights $\lambda_{m,n}<\sigma<\sigma'$, with identical prescribed ray coefficients below $\sigma$. Terms displayed at weight $\sigma'$ but not at weight $\sigma$ have positive exponents $\lambda_{\rm sp}\in[\sigma,\sigma')$. Their normalized lifts, including the logarithmic terms, lie in $H^1$ near the vertex. The invariance under enlargement in Lemma~\ref{lem:finite-part-uniqueness-comp} therefore gives the same $u$ for both truncations. For some finite $B_{\rm log}$,
\begin{align}
 &\bulkdistat{m}^{-\lambda_{m,n}} \abs{\int_0^{\omega_m} \bigl(H_{{\rm ray},m}^{(\sigma')}-H_{{\rm ray},m}^{(\sigma)}\bigr) (\bulkdistat{m},\theta)\Phi_{m,n}(\theta)\,\ud\theta} \leq C \bulkdistat{m}^{\sigma-\lambda_{m,n}} \bigl(1+\abs{\ln \bulkdistat{m}}\bigr)^{B_{\rm log}}\to0.
\end{align}
Thus \eqref{eq:global-zero-ray-coefficient} gives the same $c_{m,n}[h]$ for both weights.
\end{proof}

\begin{remark}
\label{rem:cutoff-source-no-new-exponent}
The annular source $\mathcal G_m^{\rm ann}$ is supported away from the vertex and enters the cylinder equation on a bounded interval of $t$. Through $W_m$, it depends on the global Dirichlet solution. Changing boundary data away from $V_m$ can therefore change the coefficient integral \eqref{eq:positive-zero-ray-residue}, even when the local ray expansion is fixed. For fixed boundary data, changing the auxiliary cutoff leaves $c_{m,n}[h]$ unchanged by \eqref{eq:global-zero-ray-coefficient}. The possible exponents of these zero-trace modes are always $\lambda_{m,n}=\frac{n\pi}{\omega_m}$.
\end{remark}

\subsection{Regularity of the Leading Harmonic Function}
\label{subsec:leading-laplace-corner-structure}

We now apply the cylinder estimates in Proposition~\ref{prop:mellin-mapping-comp} to the leading harmonic function $\rho_0$. The following regularity result is used in Remarks~\ref{rem:leading-end-state-mismatch} and \ref{rem:convexity-structural}.

\begin{lemma}
\label{lem:leading-field-convex-corner}
Let $V_m$ be a convex vertex, $0<\omega_m<\pi$, and let $\rho_0\in H^1(\Om)$ be the weak harmonic solution with boundary data $e_0$. Suppose its two incident traces are sidewise $C^{2,\gamma}$ and have the same value at $V_m$, denoted by $\rho_0(V_m):=e_0(V_m)$. Let $0<\alpha_c<\min\bigl\{1,\frac{\pi}{\omega_m}-1\bigr\}$. Then there is $p_m\in\mathbb R^2$ such that
\begin{align}
 &\rho_0(x) = \rho_0(V_m)+p_m\cdot(x-V_m) +\mathcal R_m^{\rho_0}(x),
 \label{eq:leading-field-C1-corner-expansion}\\
 &\abs{\mathcal R_m^{\rho_0}(x)} +\bulkdistat{m}\abs{\nabla\mathcal R_m^{\rho_0}(x)} \leq C\bulkdistat{m}^{1+\alpha_c}.
 \label{eq:leading-field-C1-corner-remainder}
\end{align}
In particular $\rho_0\in W^{1,\infty}(\Om\cap B_\delta(V_m))$ for sufficiently small $\delta>0$.
\end{lemma}

\begin{proof}
Since the vertex is convex, $\lambda_{m,1}>1$. Choose a noncritical $\widetilde\sigma$ with $1+\alpha_c<\widetilde\sigma<\min\{2,\frac{\pi}{\omega_m}\}$. The constant ray terms have the common value $\rho_0(V_m)$. Their linear terms determine a unique vector $p_m$ because the two incident tangents are linearly independent. Thus the affine harmonic function $L_m(x):=\rho_0(V_m)+p_m\cdot(x-V_m)$ has the same constant and linear ray terms as $\rho_0$. In this proof, let $\mathcal R_{m,i}$ denote the restriction of $\rho_0-L_m$ to ray $i$. Taylor's theorem gives $\mathcal R_{m,i}\in\mathcal C_{\widetilde\sigma,0}^{2,\gamma}$, since the constant and linear terms have been subtracted and $\widetilde\sigma<2$.

Use a sufficiently small sector neighborhood contained in the two neighborhoods where the assumed trace regularity holds. Define its angular trace extension $E_m^{\rm ang}$ as in Step~1 of the proof of Proposition~\ref{prop:mellin-mapping-comp}. Then $E_m^{\rm ang}\in H^1$ and its weighted norm is bounded by the local trace norms, by \eqref{eq:weighted-ray-remainder-lift-bound}--\eqref{eq:weighted-ray-lift-H1}. Choose a radial cutoff $\chi$ supported in this neighborhood and equal to one on a smaller one. Set $z:=\chi(\rho_0-L_m-E_m^{\rm ang})$. Then
\begin{align}
 -\bulkdistat{m}^2\Delta z &=\chi\bulkdistat{m}^2\Delta E_m^{\rm ang} -\bulkdistat{m}^2[\Delta,\chi](\rho_0-L_m-E_m^{\rm ang}).
 \label{eq:leading-localized-cylinder-source}
\end{align}
The function $z$ has finite energy and zero ray trace. The first term on the right has weight $\bulkdistat{m}^{\widetilde\sigma}$ in dyadic $C^{0,\gamma}$. The second is supported in a fixed annulus. Interior and flat-boundary elliptic estimates there bound it in $C^{0,\gamma}$ by $\nm{\rho_0}_{H^1(\Om)}$ and the local $C^{2,\gamma}$ trace norms.

Consequently the cylinder Green argument \eqref{eq:cylindrical-source-piecewise-bound-comp}--\eqref{eq:strong-terminal-remainder-bound-comp} applies to this localized equation with $q=2$, $\sigma=\widetilde\sigma$ and $B=0$. Since $\widetilde\sigma<\lambda_{m,1}$, its low-mode sum is empty. On the neighborhood where $\chi=1$, it gives $\rho_0-L_m=z+E_m^{\rm ang}\in\mathcal C_{\widetilde\sigma,0}^{2,\gamma}$. The fixed-radius derivative estimates \eqref{eq:weighted-sector-pointwise-comp} and $\widetilde\sigma>1+\alpha_c$ imply \eqref{eq:leading-field-C1-corner-expansion}--\eqref{eq:leading-field-C1-corner-remainder}.
\end{proof}

\begin{remark}
\label{rem:leading-field-scope}
The compatibility condition makes the constant endpoint values equal. The two tangential slopes determine the affine term in \eqref{eq:leading-field-C1-corner-expansion}. This gives the stated $C^1$ regularity at the vertex. The modes $Z_{m,n}^+$ and resonant lifts describe the power-logarithmic terms that may occur in higher derivatives.
\end{remark}
\section{Construction of the Interior Solution and Side Layers}
\label{sec:interior-side-hierarchy-comp}

We construct the interior solution $\Uint_k$ and the side layer $\Uside_{j,k}$ for $0\le k\le\Nexp$. The construction uses the Milne maps \eqref{eq:milne-maps-comp} and the harmonic Dirichlet problem of Section~\ref{sec:polygonal-profiles-comp}.

At order zero, the Milne problem has in-flow boundary data $g_{j,0}$ and zero source. At order $k\ge1$, suppose $\Uint_{k-1}$ and $\Uside_{j,k-1}$ are known. Their derivatives and $\gaux_{j,k}$ determine the Milne boundary data and source. The limit of the Milne solution as $\eta\to\infty$ is $e_{j,k}(s)$ in \eqref{eq:endstate-recursion-comp}. We use this end state as the boundary value of $\rho_k$, solving the weak harmonic Dirichlet problem at $k=0$ and the problem with prescribed singular part in Definition~\ref{def:finite-part-solution-comp} at $k\ge1$. Then \eqref{eq:total-interior-recursion-comp} determines $\Uint_k$. Subtracting $e_{j,k}(s)$ from the Milne solution gives the side layer $\Uside_{j,k}$, which tends to zero as $\eta\to\infty$. Thus the Milne data are known before the Dirichlet problem is solved.

Near each vertex, these solutions have finite expansions in powers and logarithms of the distance to the vertex. We first determine their index sets. Proposition~\ref{prop:data-generated-hierarchy-comp} then gives existence, uniqueness and the weighted estimates needed for matching. These estimates include spatial derivatives, logarithmic powers and exponential decay of the side layers.

At order $k$, let $\lambda$ denote the radial power in a vertex expansion. We retain the terms with combined exponent $\alpha=k+\lambda<\mathsf d$. In Section~\ref{sec:matching-realization}, we take $\mathsf d=\Tcon$ and use the terms with $\alpha<\Tmatch$ to define the wedge data.

We use the coordinates of Subsection~\ref{subsec:intro-coordinate-atlas} and \eqref{eq:physical-distance-convention}. Thus $\bulkdistat{m}$ is the radial distance from $V_m$. On incident ray $i$, we write $x=V_m+\sidedistat{m}{i}\tphys{m}{i}$, so $\sidedistat{m}{i}=\bulkdistat{m}(x)$. On the whole side $E_j$, we use the oriented arclength $s=s_j$. In the outward coordinate \eqref{eq:local-global-side-coordinate}, define $e_{m,i,k}:=e_{j(m,i),k}\circ s_{m,i}^{\rm ve}$. The notation $\Uside_{m,i,k}$ denotes the same pullback of $\Uside_{j(m,i),k}$ in its tangential variable. We use these pullbacks in the ray expansions and omit $m$ from $\bulkdistat{m}$ when the vertex is unspecified.

\subsection{Indices of the Vertex Expansions}

For wedge matching, we collect terms with the same combined exponent of $\e$. We now determine how these exponents change under \eqref{eq:microscopic-recursion-comp}--\eqref{eq:side-profile-inflow-decay-comp}.

For a coefficient in the vertex expansion at $V_m$, write
\begin{align}
 \mathfrak l:=(m,k,\lambda,b;\mathsf t), \qquad \mathsf t\in\mathscr T:=\{\rho,U,e_1,e_2,\mathrm{SL}_1,\mathrm{SL}_2\}.
\end{align}
Here $k$ is the power of $\e$, and $\lambda,b$ are the radial exponent and logarithmic degree. The index $\mathsf t=\mathsf t(\mathfrak l)$ takes the value $\rho,U,e_i$ or $\mathrm{SL}_i$ according as the coefficient belongs to the expansion of $\rho_k$, $\Uint_k$, $e_{m,i,k}$ or $\Uside_{m,i,k}$. As in \eqref{eq:physical-exponent-definition-comp}, put $\alpha(\mathfrak l)=k+\lambda$.

The differentiation in \eqref{eq:microscopic-recursion-comp} and \eqref{eq:side-profile-recursion-comp} changes $(k,\lambda,b)$ to
\begin{align}
 (k+1,\lambda-1,b) &\quad\hbox{or}\quad (k+1,\lambda-1,b-1),
 \label{eq:generic-recurrence-descendants-comp}
\end{align}
since $\partial\bigl(\bulkdist^\lambda(\ln \bulkdist\bigr)^b)=\bulkdist^{\lambda-1}\big\{\lambda(\ln \bulkdist)^b+b(\ln \bulkdist)^{b-1}\big\}$. Both possibilities preserve the combined exponent:
\begin{align}
 (k+1)+(\lambda-1)&=k+\lambda=\alpha.
 \label{eq:recurrence-degree-invariance-comp}
\end{align}

Under $\bulkdist=\e\scaledbulk$, we have $\e^k\bulkdist^\lambda=\e^\alpha\scaledbulk^\lambda$. By \eqref{eq:recurrence-degree-invariance-comp}, differentiation to the next order preserves this power of $\e$. Hence wedge matching must retain all generated orders corresponding to each retained $\alpha$.

\begin{definition}
\label{def:admissible-hierarchy-depth-comp}
We call $\mathsf d$ an \emph{admissible truncation exponent} if $\mathsf d>\Nexp$, $\mathsf d\notin\mathscr E_\Nexp$ and $\Kreg\ge\max\{\lceil\mathsf d\rceil+3,\pacc+\Nexp+1,\Nexp+4\}$. For the order-$\Nexp$ remainder estimate with up to $\pacc$ derivatives, we also require $\mathsf d-\Nexp>\pacc+1$.
\end{definition}

Since $\mathsf d>\Nexp$ and $\mathsf d\notin\mathscr E_\Nexp$, the weights $\mathsf d-k$ are positive and noncritical for the Dirichlet pencil for $0\le k\le\Nexp$. No generated exponent equals $\mathsf d$. In particular, $\mathsf d\notin\mathbb N_0$.

We define the index sets according to \eqref{eq:microscopic-recursion-comp}--\eqref{eq:side-profile-inflow-decay-comp}. The set $\mathscr D_{m,k}$ contains the indices obtained by differentiating the preceding order. These determine the Milne data before we solve for $\rho_k$. In the harmonic lift, $\iota_{\rm res}(\lambda)$ is given by \eqref{eq:general-log-lift-comp} with opening $\omega_m$, and we set $\iota_{\rm res}(0):=0$.

\begin{definition}
\label{def:singular-exponent-ledger}
Fix $V_m$ and let $\mathsf d$ satisfy Definition~\ref{def:admissible-hierarchy-depth-comp}. We define $\mathscr I_m(\mathsf d)$ and $\mathscr D_{m,k}$ as the smallest pair of families satisfying \textup{(L1)}--\textup{(L10)} below. Their entries have the form $\mathfrak l=(m,k,\lambda,b;\mathsf t)$ with $0\le k\le\Nexp$, $b\in\mathbb N_0$ and $\mathsf t\in\mathscr T$. The intersection of any collection of pairs satisfying these rules also satisfies them, so the smallest pair is well-defined.

\begin{enumerate}[label=\textup{(L\arabic*)},leftmargin=*,itemsep=0.3em]

\item \emph{Taylor coefficients of the in-flow data.} Include $(m,k,n,0;e_i)$ for $0\le k\le1$, $i=1,2$ and every integer $n\ge0$.

\item \emph{Harmonic modes $Z^+_{m,n}$.} Include $(m,k,\lambda_{m,n},0;\rho)$ for $0\le k\le\Nexp$ and every $n\in\mathbb N$ such that $\lambda_{m,n}=\frac{n\pi}{\omega_m}<\mathsf d-k$, with $Z^+_{m,n}$ defined in \eqref{eq:intro-zero-ray-pairs-comp}.

\item \emph{Differentiation.} If $(m,k-1,\lambda,b;\mathsf t)$ is an entry with $\mathsf t\in\{U,\mathrm{SL}_1,\mathrm{SL}_2\}$, include $(m,k,\lambda-1,b;\mathsf t)$ and, if $b\ge1$, $(m,k,\lambda-1,b-1;\mathsf t)$. Also place these entries in $\mathscr D_{m,k}$.

\item \emph{Milne end state.} If $(m,k,\lambda,b;U)\in\mathscr D_{m,k}$, include $(m,k,\lambda,b;e_i)$ for $i=1,2$. If $(m,k,\lambda,b;\mathrm{SL}_i)\in\mathscr D_{m,k}$, include $(m,k,\lambda,b;e_i)$ for that $i$.

\item \emph{Harmonic lift.} If $(m,k,\lambda,b;e_i)$ is an entry, include $(m,k,\lambda,c;\rho)$ for $0\le c\le b+\iota_{\rm res}(\lambda)$.

\item \emph{Interior solution.} If $(m,k,\lambda,b;\rho)$ is an entry, include $(m,k,\lambda,b;U)$.

\item \emph{Side layer.} If $(m,k,\lambda,b;U)$ is an entry, include $(m,k,\lambda,b;\mathrm{SL}_i)$ for $i=1,2$. If $(m,k,\lambda,b;e_i)$ is prescribed by \textup{(L1)}, include $(m,k,\lambda,b;\mathrm{SL}_i)$ for that $i$.

\item \emph{Truncation.} In every rule, retain an entry only if $\alpha(\mathfrak l)=k+\lambda<\mathsf d$.

\item \emph{Repeated indices.} If several rules give the same tuple, use one entry and add the corresponding coefficients.

\item \emph{Zero coefficients.} When comparing the expansions of different functions, orders or rays, assign coefficient zero to missing indices, as in \eqref{eq:mellin-ray-data-comp}. All quantities with index $-1$ are zero.
\end{enumerate}
Only \textup{(L3)} adds entries to $\mathscr D_{m,k}$, and it also places them in $\mathscr I_m(\mathsf d)$. Thus
\begin{align}
 \mathscr D_{m,k} \subseteq \bigl\{\mathfrak l\in\mathscr I_m(\mathsf d): \text{$\mathfrak l$ has order }k,\ \mathsf t(\mathfrak l)\in\{U,\mathrm{SL}_1,\mathrm{SL}_2\}\bigr\}.
 \label{eq:differentiated-slots-inclusion}
\end{align}
The inclusion can be strict, since \textup{(L6)} and \textup{(L7)} add further order-$k$ entries to $\mathscr I_m(\mathsf d)$. The set $\mathscr D_{m,k}$ is determined entirely by differentiation of order $k-1$ and is known before the order-$k$ construction. In particular, $\mathscr D_{m,0}=\varnothing$.
\end{definition}

\begin{remark}
\label{rem:ledger-clauses-explained}
The rules above follow directly from \eqref{eq:microscopic-recursion-comp}--\eqref{eq:side-profile-inflow-decay-comp}.

By \eqref{eq:auxiliary-inflow-coefficients}, $\gaux_{j,k}=0$ for $k\ge2$. Hence \textup{(L1)} contains the endpoint Taylor coefficients from \eqref{eq:auxiliary-endpoint-pullbacks} only at orders zero and one. Rule \textup{(L2)} gives the harmonic modes $Z^+_{m,n}$ of \eqref{eq:intro-zero-ray-pairs-comp} in the expansion \eqref{eq:mellin-expansion-comp} of $\rho_k$. Their coefficients are determined by the global Dirichlet problem. Negative powers arise from differentiation in \textup{(L3)} and the harmonic lifts of the resulting ray data, with the solution class fixed by Definition~\ref{def:finite-part-solution-comp}. The two alternatives in \textup{(L3)} are exactly \eqref{eq:generic-recurrence-descendants-comp}, obtained from $-w\cdot\nabla$ in \eqref{eq:microscopic-recursion-comp} and $-\tau_j\partial_s$ in \eqref{eq:endstate-recursion-comp} and \eqref{eq:side-profile-recursion-comp}.

By Corollary~\ref{cor:milne-polyhom-comp}, the end-state map acts separately on each coefficient. This gives \textup{(L4)} and determines the boundary data for $\rho_k$. The harmonic lift in \textup{(L5)} preserves $\lambda$. At a nonzero root of the Dirichlet pencil, it adds at most one logarithm, with coefficient \eqref{eq:resonant-top-log-coefficient}. At $\lambda=0$, the system \eqref{eq:zero-exponent-log-lift-system} preserves the bound $b$. Rule \textup{(L6)} follows from \eqref{eq:total-interior-recursion-comp}. Rule \textup{(L7)} gives the in-flow boundary trace in \eqref{eq:side-profile-inflow-decay-comp}, while the source term has the differentiated side indices in $\mathscr D_{m,k}$.

Some entries of the index set have zero coefficients. In particular, the resonant term $\beta_{b+1}\bulkdist^\lambda\bigl(\ln\bulkdist\bigr)^{b+1}\Phi_\lambda(\theta)$ vanishes on both rays. Therefore the coefficients $B^{\rm sl}_{m,i,k,\alpha,k+1}$ included by \textup{(L7)} are zero. The logarithmic degrees are at most $k$ for $e_{m,i,k}$ and $\Uside_{m,i,k}$, and at most $k+1$ for $\rho_k$ and $\Uint_k$, as stated in \eqref{eq:order-dependent-log-budgets}.

All solution maps are linear. We multiply the solutions only by prescribed angular functions, powers of $\e$ and spatial cutoffs. Thus products of unknown solutions do not introduce additional exponents.
\end{remark}

\begin{lemma}
\label{lem:ledger-finite-closed}
Let $\mathsf d$ satisfy Definition~\ref{def:admissible-hierarchy-depth-comp}. For each $V_m$, the set $\mathscr I_m(\mathsf d)$ is finite and is closed under \textup{(L3)}--\textup{(L7)} whenever the resulting order is at most $\Nexp$ and \textup{(L8)} retains the entry. Moreover, every combined exponent belongs to the exceptional set \eqref{eq:exceptional-set-comp}:
\begin{align}
 \alpha=k+\lambda\in\mathscr E_\Nexp \qquad\text{for every entry }(m,k,\lambda,b;\mathsf t)\in\mathscr I_m(\mathsf d).
 \label{eq:generated-exponents-in-exceptional-set}
\end{align}
In particular, the generated exponents have positive distance from $\Tmatch$ and $\Tcon$ in \eqref{eq:depth-choices-comp}.
\end{lemma}

\begin{proof}
By Definition~\ref{def:singular-exponent-ledger}, the pair $(\mathscr I_m(\mathsf d),\{\mathscr D_{m,k}\})$ contains every entry required by \textup{(L1)}--\textup{(L7)} and retained by \textup{(L8)}. Rule \textup{(L3)} preserves $\alpha$ by \eqref{eq:recurrence-degree-invariance-comp}, and \textup{(L4)}--\textup{(L7)} change neither $k$ nor $\lambda$. Hence all these rules preserve $\alpha<\mathsf d$. Every differentiated entry of order at most $\Nexp$ is therefore retained.

We prove finiteness order by order. New exponents occur only in \textup{(L1)} and \textup{(L2)}. Rule \textup{(L3)} reduces them by one, and it is applied at most $k$ times before order $k$. Consequently every spatial exponent at that order is $n-j$ for $n\in\mathbb N_0$ or $\frac{n\pi}{\omega_m}-j$ for $n\in\mathbb N$, with $0\le j\le k$. Since $-k\le\lambda<\mathsf d-k$ by \textup{(L8)}, there are only finitely many such exponents.

Differentiation does not increase the logarithmic degree, and the harmonic lift increases it by at most one at each of the $\Nexp+1$ orders. Thus $0\le b\le\Llog$ as in \eqref{eq:delta-lambda-comp}. With $k'=k-j$, every combined exponent lies in $\{k'+n:n\in\mathbb N_0\}\cup\{k'+\frac{n\pi}{\omega_m}:n\in\mathbb N\}$ for some $0\le k'\le k\le\Nexp$. This proves \eqref{eq:generated-exponents-in-exceptional-set}. Finally, $\mathscr T$ has six elements and repeated indices are identified by \textup{(L9)}. It follows that $\mathscr I_m(\mathsf d)$ is finite.
\end{proof}

\subsection{Construction and Estimates}
\label{subsec:data-generated-hierarchy-comp}

We now apply Proposition~\ref{prop:flat-milne-comp} and Proposition~\ref{prop:mellin-mapping-comp} to construct the solutions. The variational estimate and local Schauder estimates give the bounds away from the vertices.

\subsubsection{Recursion for the Interior Solution and Side Layers}

Set $\Uint_{-1}=0$ and $\Uside_{j,-1}=0$. Recall from \eqref{eq:auxiliary-inflow-coefficients} that $\gaux_{j,k}=0$ for $2\le k\le\Nexp$. By Proposition~\ref{prop:flat-milne-comp}, the Milne maps are bounded and linear and commute with tangential derivatives. Corollary~\ref{cor:milne-polyhom-comp} gives their action on the finite vertex expansions, including the weighted H\"older estimate for the remainders.

We define $\widehat{\Uint}_k$, $e_{j,k}$, $\rho_k$, $\Uint_k$ and $\Uside_{j,k}$ successively by the following equations for $0\le k\le\Nexp$. The restrictions to $E_j$ are evaluated at $x_j(s)$. Proposition~\ref{prop:data-generated-hierarchy-comp} below proves that the equations are solvable in the stated classes and that their data have the regularity required at the next order.
\begin{align}
 \widehat{\Uint}_0&=0,\quad \widehat{\Uint}_k=-w\cdot\nabla\Uint_{k-1} \quad (k\ge1),
 \label{eq:microscopic-recursion-comp}\\
 e_{j,k}(s) &=\mathcal E_{\mu_j}\left( \gaux_{j,k}(s,\cdot)-\widehat{\Uint}_k\big|_{E_j}, -\tau_j\partial_s\Uside_{j,k-1} \right),
 \label{eq:endstate-recursion-comp}\\
 \Delta\rho_k&=0\quad\text{in }\Om, \qquad\rho_k\big|_{E_j}=e_{j,k},
 \label{eq:scalar-recursion-comp}\\
 \Uint_k&=\rho_k+\widehat{\Uint}_k,
 \label{eq:total-interior-recursion-comp}\\
 \big(\mu_j\partial_\eta+\qk\big)\Uside_{j,k} &=-\tau_j\partial_s\Uside_{j,k-1},
 \label{eq:side-profile-recursion-comp}\\
 \Uside_{j,k}(s,0,w) &=\gaux_{j,k}(s,w)-\Uint_k\big|_{E_j} \ \ (\mu_j>0), \qquad \nm{\Uside_{j,k}(s,\eta,\cdot)}_{L^\infty(\Sone)}\rightarrow0 \ \ (\eta\to\infty).
 \label{eq:side-profile-inflow-decay-comp}
\end{align}
At successive orders the Milne map uses the decreasing exponential rates
\begin{align}
 \kappa_k:=\frac{\Nexp+1-k}{\Nexp+3}, \qquad -1\le k\le \Nexp, \qquad 1>\kappa_{-1}>\kappa_0>\cdots>\kappa_\Nexp=:\kappa_\ast>0.
 \label{eq:milne-rate-ladder-comp}
\end{align}
The space $X_\kappa$, defined in \eqref{eq:Milne-exponential-space}, measures exponential decay in $\eta$ by the weighted $L^\infty_{\eta,w}$ norm $\nm{\ue^{\kappa\eta}F}_{L^\infty_{\eta,w}}$. At order $k$ and each fixed tangential coordinate on the open side, the Milne estimate takes a source in $X_{\kappa_{k-1}}$ to a side layer in $X_{\kappa_k}$. The gap between these rates absorbs the polynomial factors in $\eta$ produced by the source. All orders can then be compared in the common space $X_{\kappa_\ast}$. Their dependence on the tangential coordinate is controlled by the separate weighted vertex estimates.

For fixed $s\in(0,L_j)$, the two inputs in \eqref{eq:endstate-recursion-comp} are the in-flow boundary function $\gaux_{j,k}(s,\cdot)-\widehat{\Uint}_k(x_j(s),\cdot)$ in $L^\infty(\{\mu_j>0\})$ and the source $-\tau_j\partial_s\Uside_{j,k-1}(s,\cdot,\cdot)$ in $X_{\kappa_{k-1}}$. The maps in \eqref{eq:milne-maps-comp} applied to these same inputs give, respectively, the end state $e_{j,k}(s)$ and the decaying part $\Uside_{j,k}(s,\cdot,\cdot)\in X_{\kappa_k}$. Here the rates are those of \eqref{eq:milne-rate-ladder-comp}. At $k=0$, the input pair is $(g_{j,0}(s,\cdot),0)$.

The Milne equation \eqref{eq:side-profile-recursion-comp} holds in $L^\infty(\Sone)$ for a.e.\ $\eta>0$. The in-flow boundary condition holds in $L^\infty(\{\mu_j>0\})$, and the limit as $\eta\to\infty$ is taken in $L^\infty(\Sone)$. By the energy argument in Proposition~\ref{prop:flat-milne-comp}, these conditions determine a unique mild solution in $L^\infty((0,\infty)\times\Sone)$.

At order zero, \eqref{eq:endstate-recursion-comp} gives $e_0\big|_{E_j}=e_{j,0}$, with $e_0$ defined in \eqref{eq:leading-boundary-function}. By \eqref{eq:leading-compatibility-comp}, these boundary values agree at every vertex. We define $\rho_0\in H^1(\Om)$ as the weak harmonic solution satisfying
\begin{align}
 \rho_0\big|_{\partial\Om}=e_0,\qquad \int_\Om\nabla\rho_0\cdot\nabla\phi\,\ud x=0 \quad \bigl(\phi\in H_0^1(\Om)\bigr).
 \label{eq:leading-variational-Dirichlet}
\end{align}

The boundary value of $\rho_0$ is its Sobolev trace. For $k\ge1$, we use Definition~\ref{def:finite-part-solution-comp}: on each open side, the boundary value of $\rho_k$ is the restriction of its prescribed singular part plus the Sobolev trace of its $H^1$ correction. Equations \eqref{eq:microscopic-recursion-comp} and \eqref{eq:total-interior-recursion-comp} then give the restrictions of $\widehat{\Uint}_k$ and $\Uint_k$, using the regularity away from vertices proved below. The value $\Uside_{j,k}(s,0,\cdot)$ in \eqref{eq:side-profile-inflow-decay-comp} is the characteristic trace of the mild Milne solution.

Subtracting the end state $e_{j,k}$ from the Milne solution with inputs \eqref{eq:endstate-recursion-comp} gives the in-flow boundary trace $\gaux_{j,k}-\widehat{\Uint}_k\big|_{E_j}-e_{j,k}=\gaux_{j,k}-\Uint_k\big|_{E_j}$. Thus \eqref{eq:side-profile-inflow-decay-comp} holds in $L^\infty(\{\mu_j>0\})$. Multiplying by $\e^k$ and summing, we obtain
\begin{align}
 \left(\sum_{k=0}^\Nexp\e^k\Uint_k\right)\bigg|_{\gamma_{j,-}}(x_j(s),w) + \sum_{k=0}^\Nexp\e^k\Uside_{j,k}(s,0,w) &= \sum_{k=0}^\Nexp\e^k\gaux_{j,k}(s,w)= g_{j,0}(s,w) +\e g_{j,1}(s,w),
 \label{eq:profile-boundary-telescope-comp}
\end{align}
where the last equality follows from $\gaux_{j,k}=0$ for $k\ge2$. We use this identity to compute the in-flow boundary error of the approximation \eqref{eq:forward-composite-definition}.

In the estimates below, let $0<\kappa\le\kappa_k$, where $\{\kappa_k\}_{k=-1}^{\Nexp}$ is the decreasing sequence in \eqref{eq:milne-rate-ladder-comp}. We use $\kappa=\kappa_k$ at order $k$. Since $\kappa_\ast\le\kappa_k$, we have $\nm F_{X_{\kappa_\ast}}\le\nm F_{X_{\kappa_k}}$. The same inequality holds after tangential differentiation and in the Banach-valued H\"older norms. Hence these estimates also hold with the common weight $\kappa_\ast$ used in Sections~\ref{sec:matching-realization} and~\ref{sec:main-proof}.

For the spatial derivatives, put
\begin{align}
 q_k&:=\max\{\pacc,3\}+\Nexp-k,\qquad 0\le k\le\Nexp.
 \label{eq:hierarchy-derivative-orders}
\end{align}
Then $q_{k-1}=q_k+1$ for $1\le k\le\Nexp$, and $q_\Nexp=\max\{\pacc,3\}$. The loss of one derivative in \eqref{eq:microscopic-recursion-comp} and \eqref{eq:side-profile-recursion-comp} is therefore included in this choice. The assumption \eqref{eq:Kreg-choice-comp} also supplies the derivatives required for the cutoff estimates.

At the chosen construction depth $\mathsf d=\Tcon$ in \eqref{eq:depth-choices-comp}, the Taylor remainder and highest tangential H\"older estimate for the primitive data at orders $k=0,1$ require at most $\max\{\lceil\Tcon\rceil-k,q_k+1\}$ derivatives of $g_{j,k}$. Since $\lceil\Tcon\rceil=\Nexp+\pacc+2$ and $\Kreg=\Nexp+\pacc+7$, the available derivative margins are
\begin{align}
 (\Kreg-k)-(\lceil\Tcon\rceil-k)&=5,\\
 (\Kreg-k)-(q_k+1)&=\pacc+6-\max\{\pacc,3\}\ge4.
\end{align}
Away from the vertices, $q_{k-1}+1=q_k+2$, so one derivative of the preceding profile leaves the required $q_k+1$ derivatives. The endpoint Taylor remainder at matching depth $\Tmatch\in(\pacc+3,\pacc+4)$ requires only $\pacc+4-k$ primitive derivatives, for $k=0,1$.

\subsubsection{Vertex Expansions and Norms}

Let $Q_k$ denote $\rho_k$ or $\Uint_k$, with value space $\mathcal V_\rho:=\mathbb R$ or $\mathcal V_U:=L^\infty(\Sone)$, respectively. For $Q\in\{\rho,U\}$ and $i\in\{1,2\}$, define $\mathcal J^Q_{m,k}(\mathsf d)$, $\mathcal J^e_{m,i,k}(\mathsf d)$ and $\mathcal J^\Uside_{m,i,k}(\mathsf d)$ by taking the pairs $(\alpha,b)$ for which $(m,k,\alpha-k,b;\mathsf t)\in\mathscr I_m(\mathsf d)$ with $\mathsf t=Q$, $e_i$ and $\mathrm{SL}_i$, respectively.

By Subsection~\ref{subsec:intro-coordinate-atlas}, $\Om\cap B(V_m,8\ellvtx)$ agrees with the tangent sector. We apply the elliptic mapping theorem at radius $4\ellvtx$ and obtain \eqref{eq:mellin-expansion-comp}--\eqref{eq:mellin-remainder-bound-comp} on $K_{\omega_m}\cap B_{2\ellvtx}$. We use this latter radius in the norms below. The same expansions enter $\Ubase^{\e,[\Nexp]}$ in \eqref{eq:Ubase-finite-expanded} and the partition of unity in Section~\ref{sec:main-proof} on $\{\bulkdistat{m}<2\ellvtx\}$.

For the ray expansions, \eqref{eq:local-global-side-coordinate} gives $\partial_{\sidedistat{m}{i}}=\varsigma_{m,i}\partial_s$. Thus differentiating a pullback includes the change of sign at $s=L_j$, as in \eqref{eq:endpoint-pullbacks}--\eqref{eq:auxiliary-endpoint-pullbacks}.

For the fixed truncation exponent $\mathsf d$, we write the vertex expansions as
\begin{align}
 Q_k(\bulkdistat{m},\theta_m,w) &= \sum_{(\alpha,b)\in\mathcal J^Q_{m,k}(\mathsf d)} \bulkdistat{m}^{\alpha-k}\bigl(\ln \bulkdistat{m}\bigr)^b A^Q_{m,k,\alpha,b}(\theta_m,w) +\mathcal R^{Q,\mathsf d}_{m,k}(\bulkdistat{m},\theta_m,w),
 \label{eq:weighted-interior-decomposition-comp}\\
 \Uside_{m,i,k}(\sidedistat{m}{i},\eta,w) &= \sum_{(\alpha,b)\in\mathcal J^\Uside_{m,i,k}(\mathsf d)} \sidedistat{m}{i}^{\alpha-k}\bigl(\ln \sidedistat{m}{i}\bigr)^b B^{\rm sl}_{m,i,k,\alpha,b}(\eta,w) +\mathcal R^{B,\mathsf d}_{m,i,k}(\sidedistat{m}{i},\eta,w),
 \label{eq:weighted-side-decomposition-comp}\\
 e_{m,i,k}(\sidedistat{m}{i}) &= \sum_{(\alpha,b)\in\mathcal J^e_{m,i,k}(\mathsf d)} \sidedistat{m}{i}^{\alpha-k}\bigl(\ln \sidedistat{m}{i}\bigr)^b e_{m,i,k,\alpha,b} +\mathcal R^{e,\mathsf d}_{m,i,k}(\sidedistat{m}{i}).
 \label{eq:weighted-endstate-decomposition-comp}
\end{align}
Here $A^\rho_{m,k,\alpha,b}$ is an angular function, $A^U_{m,k,\alpha,b}$ is an angular function with values in $L^\infty(\Sone)$, $B^{\rm sl}_{m,i,k,\alpha,b}$ is a function of $(\eta,w)$, and $e_{m,i,k,\alpha,b}$ is a real number. The velocity argument is omitted when $Q=\rho$. Each remainder is the difference between the function on the left and its finite sum. Proposition~\ref{prop:data-generated-hierarchy-comp} proves that these remainders have weight $\mathsf d-k$ in the norms below, with $B^{\rm sl}_{m,i,k,\alpha,b}\in X_{\kappa_k}$. By Lemma~\ref{lem:ledger-finite-closed}, the sums are finite. Their indices satisfy $\alpha<\mathsf d$ and $0\le b\le\Llog=\Nexp+1$, as in \eqref{eq:delta-lambda-comp}.

At order $k$, the bounds for the logarithmic degrees are
\begin{align}
 L_k^{\rm ray}:=k,\qquad L_k^{\rm int}:=k+1,\qquad 0\le k\le\Nexp.
 \label{eq:order-dependent-log-budgets}
\end{align}
For an integer $q\ge0$, a logarithmic degree bound $L$ and a weight $0<\kappa\le\kappa_k$, define
\begin{align}
 \mathfrak N^{\mathrm{int},\mathsf d}_{m,k} (Q_k;q,L) :={}& \sum_{(\alpha,b)\in\mathcal J^Q_{m,k}(\mathsf d)} \bnm{A^Q_{m,k,\alpha,b}}_{C^q([0,\omega_m];\mathcal V_Q)}+ \bnm{\mathcal R^{Q,\mathsf d}_{m,k}}_{ \mathcal C^{q,\gamma}_{\mathsf d-k,L} (K_{\omega_m}\cap B_{2\ellvtx};\mathcal V_Q)},
 \label{eq:weighted-interior-germ-norm-comp}\\
 \mathfrak N^{B,\mathsf d}_{m,i,k,\kappa} (\Uside_{m,i,k};q,L) :={}& \sum_{(\alpha,b)\in\mathcal J^\Uside_{m,i,k}(\mathsf d)} \nm{B^{\rm sl}_{m,i,k,\alpha,b}}_{X_\kappa}+ \bnm{\mathcal R^{B,\mathsf d}_{m,i,k}}_{ \mathcal C^{q,\gamma}_{\mathsf d-k,L} (0,2\ellvtx;X_\kappa)},
 \label{eq:weighted-side-germ-norm-comp}\\
 \mathfrak N^{e,\mathsf d}_{m,i,k} (e_{m,i,k};q,L) :={}& \sum_{(\alpha,b)\in\mathcal J^e_{m,i,k}(\mathsf d)} \abs{e_{m,i,k,\alpha,b}}+ \bnm{\mathcal R^{e,\mathsf d}_{m,i,k}}_{ \mathcal C^{q,\gamma}_{\mathsf d-k,L}(0,4\ellvtx)}.
 \label{eq:weighted-endstate-germ-norm-comp}
\end{align}
These norms estimate both the coefficients in the finite sums and the weighted remainders. The remainder norms include derivatives and H\"older seminorms on each dyadic annulus. For the side layer, the $X_\kappa$ norm also gives exponential decay in $\eta$ in $L^\infty(\Sone)$. We retain the H\"older seminorms of the highest derivatives in order to apply the parameterized Milne estimate and the Schauder estimate.

Away from the vertices, put $\Om_{\rm away}:=\overline\Om\setminus\bigcup_m B_{\frac{\ellvtx}{4}}(V_m)$ and $E_{j,\rm away}:=x_j([\frac{\ellvtx}{4},L_j-\frac{\ellvtx}{4}])$. To apply local elliptic estimates at successive orders, we need the preceding solutions and boundary data on a larger set. For this purpose, define
\begin{align}
 \ell_k^{\rm away} &:=\frac{\ellvtx}{8}\left(1+\frac{k+1}{\Nexp+2}\right), \qquad -1\leq k\leq \Nexp,
 \label{eq:hierarchy-buffered-away-sets}\\
 \Om_{\rm away}^{(k)} &:=\overline\Om\setminus\bigcup_m B_{\ell_k^{\rm away}}(V_m),\qquad E_{j,\rm away}^{(k)}:=x_j([\ell_k^{\rm away},L_j-\ell_k^{\rm away}]).
 \label{eq:hierarchy-buffered-domain-set}
\end{align}
We have $\Om_{\rm away}\subset\Om_{\rm away}^{(k)}$ and $E_{j,\rm away}\subset E_{j,\rm away}^{(k)}$. Consecutive radii differ by $\frac{\ellvtx}{8(\Nexp+2)}>0$. Thus each order-$k$ set has a relative neighborhood in $\overline\Om$ or $E_j$ contained in the corresponding order-$(k-1)$ set. These neighborhood radii depend only on $\Nexp$ and $\ellvtx$. Norms of $e_{j,k}$ and $\Uside_{j,k}$ on a subset of $E_j$ are taken in the corresponding arclength interval under $x_j$.

Define $e_k$ on the open sides by $e_k\big|_{E_j}=e_{j,k}$. For $k\ge1$, write $\rho_k=H_k^{\rm sing}+v_k^{\rm reg}$ as in Definition~\ref{def:finite-part-solution-comp}. At order zero, set $H_0^{\rm sing}=0$ and $v_0^{\rm reg}=\rho_0$. The source $f_k^{\rm cut}$ is given by \eqref{eq:finite-part-commutator-source} for the cutoff lifts in $H_k^{\rm sing}$, and $f_0^{\rm cut}=0$. The function $v_k^{\rm reg}\in H^1(\Om)$ is the unique variational solution of $-\Delta v_k^{\rm reg}=f_k^{\rm cut}$ with Sobolev trace $e_k-H_k^{\rm sing}\big|_{\partial\Om}$ a.e.\ on $\partial\Om$, where the singular part is restricted to the open sides.

By \eqref{eq:finite-part-global-energy-bound} and \eqref{eq:mellin-remainder-bound-comp}, the $H^1$ norm of $v_k^{\rm reg}$ and the ray remainders control the coefficients of $Z^+_{m,n}$ in \eqref{eq:mellin-expansion-comp}. We also estimate the boundary data away from the endpoints and the cutoff source terms on their supporting annuli.

For $0<\kappa\le\kappa_k$, define
\begin{align}
 \label{eq:hierarchy-global-profile-norm} \mathfrak F_{k,\kappa}:={}&\nm{v_k^{\rm reg}}_{H^1(\Om)} +\nm{\rho_k}_{C^{q_k+1,\gamma}(\Om_{\rm away}^{(k)})} +\nm{\Uint_k}_{C^{q_k+1,\gamma}(\Om_{\rm away}^{(k)};L^\infty_w)} \\
 &+\sum_j\left\{ \nm{e_{j,k}}_{C^{q_k+1,\gamma}(E_{j,\rm away}^{(k)})} +\nm{\Uside_{j,k}}_{C^{q_k+1,\gamma}(E_{j,\rm away}^{(k)};X_\kappa)} \right\}.\notag
\end{align}

The norms in \eqref{eq:hierarchy-global-profile-norm} contain $q_k+1$ derivatives. This permits differentiation once at the next order in \eqref{eq:microscopic-recursion-comp} and \eqref{eq:side-profile-recursion-comp}. We use the nested sets \eqref{eq:hierarchy-buffered-domain-set} for the local boundary estimates, and estimate the vertex remainders through order $q_k$.

Together with the local expansion norms, $\mathfrak F_{k,\kappa}$ controls $\mathcal N_{\rm glob}$ in \eqref{eq:explicit-Mellin-far-norm}. Indeed, the norm of $e_{j,k}$ gives the boundary data away from the endpoints, the trace theorem bounds the boundary remainder by $\nm{v_k^{\rm reg}}_{H^1(\Om)}$, and the explicit harmonic lifts give the cutoff source term. In the proof below, we first estimate this boundary remainder and the cutoff source from the Milne data. The variational estimate then gives the $H^1$ norm, and the elliptic expansion theorem bounds the new angular coefficients and vertex remainder.

We combine the local coefficient and remainder estimates with \eqref{eq:hierarchy-global-profile-norm} by setting
\begin{align}
 \mathfrak N_k^\sharp(\mathsf d) :={}& \mathfrak F_{k,\kappa_k}+\sum_m\Bigg[\sum_{Q\in\{\rho,U\}}\mathfrak N^{\mathrm{int},\mathsf d}_{m,k}\bigl(Q_k;q_k,L_k^{\rm int}\bigr)
 \label{eq:parameterized-complete-profile-norm}\\
 &\qquad+\sum_{i=1}^2\Bigl\{\mathfrak N^{e,\mathsf d}_{m,i,k}\bigl(e_{m,i,k};q_k,L_k^{\rm ray}\bigr)+\mathfrak N^{B,\mathsf d}_{m,i,k,\kappa_k}\bigl(\Uside_{m,i,k};q_k,L_k^{\rm ray}\bigr)\Bigr\}\Bigg],\notag
\end{align}
where $m$ ranges over the vertices and $i=1,2$ over the incident rays. The local terms control the differentiated data near each endpoint, and $\mathfrak F_{k,\kappa_k}$ controls the global energy and derivatives away from the vertices.

For the prescribed in-flow boundary data, put
\begin{align}
 \mathcal G&:=\sum_j\sum_{\ell=0}^1\nm{g_{j,\ell}}_{C^{\Kreg-\ell,\gamma}([0,L_j];L^\infty(\{\mu_j>0\}))},
 \label{eq:hierarchy-primitive-data-norm}
\end{align}
which will control the solutions and their vertex expansions.

\subsubsection{Existence and Weighted Estimates}

We first prove uniqueness of the coefficients in a finite vertex expansion. We will use it in Step~6 and in Section~\ref{sec:matching-realization}.

\begin{lemma}
\label{lem:finite-polyhom-uniqueness}
Let $X$ be a Banach space and $\lambda_1<\cdots<\lambda_J<d_{\rm rem}$. Suppose
\begin{align}
 &\sum_{j=1}^J \bulkdist^{\lambda_j}P_j(\ln \bulkdist)+\mathcal R_{d_{\rm rem}}(\bulkdist)=0,
 \label{eq:finite-polyhom-zero-expansion}\\
 &\nm{\mathcal R_{d_{\rm rem}}(\bulkdist)}_X \le C\bulkdist^{d_{\rm rem}}\bigl(1+\abs{\ln \bulkdist}\bigr)^B,
\end{align}
where each $P_j$ is an $X$-valued polynomial. Then $P_j=0$ for every $j$. This applies in particular to the sector expansions and the $X_\kappa$-valued side expansions above.
\end{lemma}

\begin{proof}
Suppose $P_1\ne0$, and let $b$ be its degree. Dividing \eqref{eq:finite-polyhom-zero-expansion} by $\bulkdist^{\lambda_1}\bigl(\ln\bulkdist\bigr)^b$ and letting $\bulkdist\downarrow0$, we obtain the leading coefficient of $P_1$ in the norm of $X$. All terms with $j>1$ tend to zero, since a positive power dominates every logarithm. The remainder also tends to zero because $d_{\rm rem}>\lambda_1$. This contradicts \eqref{eq:finite-polyhom-zero-expansion}. Descending in $b$ and then repeating for $\lambda_2,\ldots,\lambda_J$ proves the result. For the sector expansions, take $X=C([0,\omega_m];\mathcal V_Q)$. For the side expansions, take $X=X_\kappa$.
\end{proof}

For the fixed truncation exponent $\mathsf d$, let $C^U_{m,k}$ and $C^\Uside_{m,i,k}$ denote the finite sums in \eqref{eq:weighted-interior-decomposition-comp} and \eqref{eq:weighted-side-decomposition-comp}. On the punctured vertex chart, write
\begin{align}
 \mathcal R^U_{m,k}:=\mathcal R^{U,\mathsf d}_{m,k},\qquad &C^U_{m,k}:=\Uint_k-\mathcal R^U_{m,k},
 \label{eq:exact-construction-projection-comp}\\
 \mathcal R^\Uside_{m,i,k}:=\mathcal R^{B,\mathsf d}_{m,i,k},\qquad &C^\Uside_{m,i,k}:=\Uside_{m,i,k}-\mathcal R^\Uside_{m,i,k}.
 \label{eq:exact-side-construction-projection}
\end{align}
All four quantities are zero at index $-1$. This notation suppresses $\mathsf d$.

\begin{proposition}[Construction of the interior solution and side layers]
\label{prop:data-generated-hierarchy-comp}
Let $\Om$ be a bounded convex polygon satisfying \eqref{eq:strict-convexity}. Suppose Assumption~\ref{ass:primitive-data-comp} holds, with parameters as in \eqref{eq:N-choice-comp}--\eqref{eq:Kreg-choice-comp}, and let $\mathsf d$ satisfy Definition~\ref{def:admissible-hierarchy-depth-comp}. Then the solutions of \eqref{eq:microscopic-recursion-comp}--\eqref{eq:side-profile-inflow-decay-comp} have the following properties.

\begin{enumerate}[label=\textup{(\roman*)},ref=\textup{(\roman*)},
leftmargin=*,itemsep=0.8em]

\item
\label{item:hierarchy-canonical-comp}
\emph{Existence and uniqueness.} For $0\le k\le\Nexp$, equations \eqref{eq:microscopic-recursion-comp}--\eqref{eq:side-profile-inflow-decay-comp} uniquely determine $\widehat{\Uint}_k$, $e_{j,k}$, $\rho_k$, $\Uint_k$ and $\Uside_{j,k}$. The function $\rho_0$ is the weak harmonic solution, while $\rho_k$ for $k\ge1$ belongs to the class with prescribed singular part in Definition~\ref{def:finite-part-solution-comp}. The side layer $\Uside_{j,k}$ is the unique decaying mild Milne solution and belongs to $X_{\kappa_k}$ together with the tangential derivatives in \eqref{eq:parameterized-complete-profile-norm}.

The solutions, their vertex expansion coefficients and their remainders depend continuously and linearly on $\{g_{j,\ell}:0\le\ell\le1\}$ in the closed linear subspace defined by \eqref{eq:leading-compatibility-comp}.

\item
\label{item:hierarchy-ledger-comp}
\emph{Vertex expansions.} The decompositions \eqref{eq:weighted-interior-decomposition-comp}--\eqref{eq:weighted-endstate-decomposition-comp} hold at every vertex and incident ray. Their indices belong to $\mathscr I_m(\mathsf d)$ in Definition~\ref{def:singular-exponent-ledger}. The set depends only on the geometry and the rules corresponding to \eqref{eq:microscopic-recursion-comp}--\eqref{eq:side-profile-inflow-decay-comp}, and is independent of the values of the data. Some coefficients may be zero. By Lemma~\ref{lem:ledger-finite-closed}, $\mathscr I_m(\mathsf d)$ is finite and closed under these rules. Every combined exponent satisfies
\begin{align}
 \alpha\ge0.
 \label{eq:generated-exponents-nonnegative}
\end{align}

\item \emph{Weighted estimates.} For the truncation at $\alpha<\mathsf d$, we have the following bounds.
\begin{enumerate}[label=\textup{(\alph*)},ref=\textup{(\alph*)},
leftmargin=*,itemsep=0.5em]

\item
\label{item:hierarchy-complete-norm-comp}
With the norms defined in \eqref{eq:parameterized-complete-profile-norm} and \eqref{eq:hierarchy-primitive-data-norm}, for $0\le k\le\Nexp$,
\begin{align}
 \mathfrak N_k^\sharp(\mathsf d) \leq C_{\mathsf d}\mathcal G,
 \label{eq:generic-hierarchy-bound-comp}
\end{align}
where $C_{\mathsf d}$ depends only on the polygon, the fixed coordinate charts and radii, the singular-part cutoffs and harmonic-lift normalization, the parameters in \eqref{eq:N-choice-comp}--\eqref{eq:Kreg-choice-comp}, $\mathsf d$, the H\"older exponent and the decay rates in \eqref{eq:milne-rate-ladder-comp}. It is independent of $\e$ and the in-flow data. Since $\kappa_\ast\le\kappa_k$, the coefficients and remainders in the side-layer expansions also satisfy this bound in $X_{\kappa_\ast}$.

\item
\label{item:hierarchy-local-remainders-comp}
The remainders satisfy the following estimates. In the first estimate, $a+c\le q_k$, and in the other three, $0\le q\le q_k$. We take $0<\bulkdistat{m}<2\ellvtx$ in the first two, $0<\sidedistat{m}{i}<4\ellvtx$ in the third and $0<\sidedistat{m}{i}<2\ellvtx$ in the fourth:
\begin{align}
 \nm{\partial_{\bulkdistat{m}}^a\partial_{\theta_m}^c \mathcal R^{U,\mathsf d}_{m,k}(\bulkdistat{m},\theta_m,\cdot)}_{L^\infty_w} &\leq C\mathcal G\, \bulkdistat{m}^{\mathsf d-k-a} \bigl(1+\abs{\ln \bulkdistat{m}}\bigr)^{\Llog},
 \label{eq:outer-ledger-remainder-comp}\\
 \nm{D_x^q\mathcal R^{U,\mathsf d}_{m,k}(x,\cdot)}_{L^\infty_w} &\leq C\mathcal G\, \bulkdistat{m}^{\mathsf d-k-q} \bigl(1+\abs{\ln \bulkdistat{m}}\bigr)^{\Llog},
 \label{eq:outer-cartesian-remainder-comp}\\
 \abs{\partial_{\sidedistat{m}{i}}^q \mathcal R^{e,\mathsf d}_{m,i,k}(\sidedistat{m}{i})} &\leq C\mathcal G\, \sidedistat{m}{i}^{\mathsf d-k-q} \bigl(1+\abs{\ln \sidedistat{m}{i}}\bigr)^{\Llog},
 \label{eq:endstate-ledger-remainder-comp}\\
 \nm{ \ue^{\kappa_\ast\eta}\partial_{\sidedistat{m}{i}}^q \mathcal R^{B,\mathsf d}_{m,i,k}(\sidedistat{m}{i},\eta,\cdot) }_{L^\infty_{\eta,w}} &\leq C\mathcal G\, \sidedistat{m}{i}^{\mathsf d-k-q} \bigl(1+\abs{\ln \sidedistat{m}{i}}\bigr)^{\Llog}.
 \label{eq:side-ledger-remainder-comp}
\end{align}
The estimates hold at every spatial point in the indicated charts. The two interior estimates also hold for $\mathcal R^{\rho,\mathsf d}_{m,k}$, with absolute values instead of velocity norms.

\item
\label{item:hierarchy-terminal-bounds-comp}
At order $\Nexp$, the interior solution $\Uint_\Nexp$ and the side layer $\Uside_{m,i,\Nexp}$ satisfy, for $0<\bulkdistat{m},\sidedistat{m}{i}<2\ellvtx$,
\begin{align}
 \bnm{w\cdot\nabla\Uint_\Nexp(x,w)}_{L^\infty_w} &\leq C\mathcal G\, \bulkdistat{m}^{-\Nexp-1} \bigl(1+\abs{\ln \bulkdistat{m}}\bigr)^{\Llog},
 \label{eq:terminal-complete-interior-comp}\\
 \nm{ \ue^{\kappa_\ast\eta}\partial_{\sidedistat{m}{i}}\Uside_{m,i,\Nexp} (\sidedistat{m}{i},\eta,\cdot) }_{L^\infty_{\eta,w}} &\leq C\mathcal G\, \sidedistat{m}{i}^{-\Nexp-1} \bigl(1+\abs{\ln \sidedistat{m}{i}}\bigr)^{\Llog}.
 \label{eq:terminal-complete-side-comp}
\end{align}
These estimates precede the one-sided ray extensions and the regularizations at the vertex scale in Section~\ref{subsec:fixed-local-continuations-v6}.
\end{enumerate}

\item
\label{item:hierarchy-splitting-comp}
\emph{Equations for the finite sums and remainders.} For $0\le k\le\Nexp$, the functions in \eqref{eq:exact-construction-projection-comp}--\eqref{eq:exact-side-construction-projection} satisfy
\begin{align}
 &\qk C^U_{m,k}+w\cdot\nabla_xC^U_{m,k-1}=0,\qquad \qk\mathcal R^U_{m,k}+w\cdot\nabla_x\mathcal R^U_{m,k-1}=0,
 \label{eq:exact-interior-displayed-recurrence-comp}\\
 &\big(\mu_{m,i}\partial_\eta+\qk\big)C^\Uside_{m,i,k} +\tau_{m,i}\partial_{\sidedistat{m}{i}}C^\Uside_{m,i,k-1}=0,
 \label{eq:exact-side-displayed-recurrence-comp}\\
 &\big(\mu_{m,i}\partial_\eta+\qk\big)\mathcal R^\Uside_{m,i,k} +\tau_{m,i}\partial_{\sidedistat{m}{i}}\mathcal R^\Uside_{m,i,k-1}=0.
 \label{eq:exact-side-projected-recurrence-comp}
\end{align}
Thus the remainders give no transport source term below order $\Nexp$, before multiplication by the localization cutoffs and regularization at the vertex scale.
\end{enumerate}
\end{proposition}

\begin{proof}
We prove the result by induction on $k$. At each order, we first estimate the Milne problem and then the harmonic Dirichlet problem.

\paragraph{\underline{Step 1: solvability and construction}} Away from the vertices, Lemma~\ref{lem:angular-moments-comp} gives, for $1\le k\le\Nexp-1$,
\begin{align}
 \pk(w\cdot\nabla\Uint_k)&=-\sum_{q=1}^{\lfloor\frac{k+1}{2}\rfloor}c_q\Delta^q\rho_{k+1-2q}=-\frac{1}{2}\Delta\rho_{k-1}=0.
 \label{eq:hierarchy-microscopic-solvability}
\end{align}
Here \eqref{eq:explicit-interior-profile-comp} and \eqref{eq:even-angular-moment-comp} give the first equality, and harmonicity makes the terms with $q\ge2$ vanish. At $k=0$, the velocity average vanishes by oddness. Thus $-w\cdot\nabla\Uint_{k-1}$ is microscopic for every $1\le k\le\Nexp$, and \eqref{eq:microscopic-recursion-comp}--\eqref{eq:total-interior-recursion-comp} imply $\qk\Uint_k+w\cdot\nabla\Uint_{k-1}=0$. The solvability condition first uses $\Delta\rho_j=0$ at order $j+2$, for $0\le j\le\Nexp-2$. Lemma~\ref{lem:finite-part-remainder-trace} supplies the trace regularity for the Dirichlet problems.

At $k=0$, use the Milne inputs $(g_{j,0},0)$. For $k\ge1$, suppose $\rho_\ell$, $\Uint_\ell$ and $\Uside_{j,\ell}$ have been constructed for $0\le\ell\le k-1$. We then know $\widehat{\Uint}_k=-w\cdot\nabla\Uint_{k-1}$ and $-\tau_j\partial_s\Uside_{j,k-1}$. The maps $\mathcal E_{\mu_j}$ and $\mathfrak M_{\mu_j}$ in \eqref{eq:milne-maps-comp} give $e_{j,k}$ and $\Uside_{j,k}$ from the input pair in \eqref{eq:endstate-recursion-comp}. We solve the weak Dirichlet problem at $k=0$ and the Dirichlet problem with prescribed singular part at $k\ge1$ to obtain $\rho_k$. Equation \eqref{eq:total-interior-recursion-comp} gives the interior solution $\Uint_k$. Since $\rho_k\big|_{E_j}=e_{j,k}$, the side layer already obtained from $\mathfrak M_{\mu_j}$ has precisely the in-flow boundary trace in \eqref{eq:side-profile-inflow-decay-comp}.

The Milne and Dirichlet problems are unique in their stated classes, and their solution maps are bounded and linear. Steps~2 and~3 verify that the input data at each order satisfy their hypotheses and prove the quantitative bounds. Since the index sets are fixed independently of the data, the coefficients in the vertex expansions and their remainders then depend continuously and linearly on $g_{j,0}$ and $g_{j,1}$. This proves part~\ref{item:hierarchy-canonical-comp} once the induction estimates are established.

\paragraph{\underline{Step 2: construction at order zero}}
\label{step:hierarchy-base-comp}
For $0\le k\le\Nexp$, put $\sigma_k:=\mathsf d-k$, $J_k^{\rm con}:=\lfloor\sigma_k\rfloor$ and $J^{\rm con}_{m,i,k}(\sidedistat{m}{i},w):=\sum_{n=0}^{J_k^{\rm con}}\frac{\sidedistat{m}{i}^n}{n!}\,\partial_{\sidedistat{m}{i}}^n\gaux_{m,i,k}(0,w)$. By \eqref{eq:auxiliary-endpoint-pullbacks}, $\partial_{\sidedistat{m}{i}}^n\gaux_{m,i,k}(0,w)=\varsigma_{m,i}^{\,n}\partial_s^n\gaux_{j(m,i),k}(s_{m,i}^{\rm ve}(0),w)$. This gives the change of sign for odd derivatives at $s=L_j$.

Since $\mathsf d\notin\mathscr E_\Nexp$, we have $\sigma_k\notin\mathbb N_0$. For $0\le a\le J_k^{\rm con}$, Taylor's theorem gives
\begin{align}
 \nm{\partial_{\sidedistat{m}{i}}^a\big(\gaux_{m,i,k}-J^{\rm con}_{m,i,k}\big)}_{L^\infty(\{\mu_{m,i}>0\})}=O\big(\sidedistat{m}{i}^{J_k^{\rm con}+1-a}\big)=O(\sidedistat{m}{i}^{\sigma_k-a}).
 \label{eq:hierarchy-taylor-remainder}
\end{align}
For $J_k^{\rm con}+1\le a\le q_k$, the same velocity norm is bounded by $C\mathcal G$ using \eqref{eq:hierarchy-primitive-data-norm}. Since $\sigma_k-a<0$, this also gives the required weight on $0<\sidedistat{m}{i}<1$. To obtain the H\"older seminorm, write $\sidedistat{m}{i}=\varrho t$ on $\frac{1}{2}<t<2$. For $k=0,1$, the regularity assumptions give $\Kreg-k\ge\max\{J_k^{\rm con}+1,q_k+1\}$. We may therefore apply the same derivative estimates through $a=q_k+1$. After rescaling, every derivative through this order of $(\gaux_{m,i,k}-J^{\rm con}_{m,i,k})(\varrho t,\cdot)$ is bounded in $L^\infty(\{\mu_{m,i}>0\})$ by $C\mathcal G\varrho^{\sigma_k}$. The mean value estimate for the derivative of order $q_k$ gives its $\gamma$-H\"older seminorm on the fixed interval. This proves the weighted $\mathcal C^{q_k,\gamma}_{\sigma_k,0}$ estimate in \eqref{eq:weighted-boundary-norm-comp}. For $k\ge2$, the estimates hold immediately because $\gaux_{j,k}=0$.

At $k=0$, Proposition~\ref{prop:flat-milne-comp}, with weights $\kappa_{-1}$ and $\kappa_0$, gives the Milne end state $e_{j,0}$. By the leading compatibility condition, the sidewise values agree continuously at the vertices. Their Taylor polynomials contain $\bulkdistat{m}^n$ for $0\le n\le J_0^{\rm con}$. Apply Proposition~\ref{prop:mellin-mapping-comp} with these Taylor terms and weight $\mathsf d$. We obtain the prescribed harmonic lifts, the modes $Z^+_{m,n}$ with $\lambda_{m,n}<\mathsf d$, and the remainder.

We verify that the resulting solution is $\rho_0\in H^1(\Om)$. Every harmonic lift of positive degree belongs locally to $H^1$, including the terms with logarithms. By \eqref{eq:leading-compatibility-comp}, the lift of degree zero is the same constant on the two incident rays and also has finite energy. Lemma~\ref{lem:finite-part-uniqueness-comp} therefore allows us to include all these cutoff lifts in the regular correction without changing the solution. Hence we recover the weak harmonic solution and may take $H_0^{\rm sing}=0$, $v_0^{\rm reg}=\rho_0$ in \eqref{eq:parameterized-complete-profile-norm}.

We have $\Uint_0=\rho_0$. By Proposition~\ref{prop:flat-milne-comp}, the side layer satisfies $\Uside_{j,0}\in X_{\kappa_0}$ and has the expansion \eqref{eq:weighted-side-decomposition-comp}.

The expansions of $e_{m,i,0}$ and $\Uside_{m,i,0}$ have logarithmic degree at most $L_0^{\rm ray}=0$, and the expansion of $\rho_0=\Uint_0$ has degree at most $L_0^{\rm int}=1$. Combining the variational estimate with the Schauder estimates on the nested sets below gives $\mathfrak N_0^\sharp(\mathsf d)\le C_0\mathcal G$.

\paragraph{\underline{Step 3: estimates at order $k$}}
\label{step:hierarchy-induction-comp}
Let $1\le k\le\Nexp$ and assume the estimates hold through order $k-1$ at the fixed exponent $\mathsf d$. We first bound the two inputs of \eqref{eq:endstate-recursion-comp} by $\mathcal G+\mathfrak N_{k-1}^\sharp(\mathsf d)$. The Milne maps then estimate both $e_{j,k}$ and $\Uside_{j,k}$. The coefficients and remainders of $e_{j,k}$ supply the Dirichlet data, to which we apply the variational and elliptic expansion estimates for $\rho_k$.

\emph{Milne estimates.} Since $q_{k-1}=q_k+1$, the derivatives $-w\cdot\nabla\Uint_{k-1}$ and $-\tau_j\partial_s\Uside_{j,k-1}$ retain $q_k$ spatial or tangential derivatives. For the remainders, differentiation changes the weight from $\mathsf d-(k-1)$ to $\mathsf d-k$. The angular functions $A^\rho_{m,k-1,\alpha,b}$ and $A^U_{m,k-1,\alpha,b}$ lie in fixed finite-dimensional spaces generated by the angular factors of the harmonic lifts and $Z^+_{m,n}$, their spatial derivatives and multiplication by velocity components. Norm equivalence on these spaces bounds their differentiated H\"older seminorms by the coefficient norms in \eqref{eq:weighted-interior-germ-norm-comp}.

By \eqref{eq:order-dependent-log-budgets}, $L_{k-1}^{\rm int}=k=L_k^{\rm ray}$. Differentiation preserves or lowers the logarithmic degree. Thus $\partial$ maps $\mathcal C^{q_{k-1},\gamma}_{\mathsf d-(k-1),L_{k-1}^{\rm int}}$ into $\mathcal C^{q_k,\gamma}_{\mathsf d-k,L_k^{\rm ray}}$.

On each dyadic interval, write $s=\varrho t$ with $\frac{1}{2}<t<2$ and apply Corollary~\ref{cor:milne-polyhom-comp} with $(\kappa_{\rm out},\kappa_{\rm in})=(\kappa_k,\kappa_{k-1})$. The Milne maps preserve the logarithmic degree and give the weighted H\"older remainder estimate uniformly in $\varrho$. Their source and solution spaces are $X_{\kappa_{k-1}}$ and $X_{\kappa_k}$, respectively.

It follows that the expansions of $e_{j,k}$ and $\Uside_{j,k}$ have logarithmic degree at most $L_k^{\rm ray}=k$. The harmonic lift adds at most one logarithm, so $L_k^{\rm int}=k+1$. Both degree bounds are at most $\Llog=\Nexp+1$.

\emph{Dirichlet data.} Consider the singular part $H_k^{\rm sing}$ and the source term $f_k^{\rm cut}$ for the Dirichlet problem. The source term is supported on fixed annuli away from the vertices. To apply the elliptic theorem, we need the ray remainder on $(0,4\ellvtx)$. The local norms of $\rho_{k-1},\Uint_{k-1}$ are given on $K_{\omega_m}\cap B_{2\ellvtx}$, and those of $\Uside_{m,i,k-1}$ on $(0,2\ellvtx)$. We therefore also estimate $e_{j,k}$ between radii $2\ellvtx$ and $4\ellvtx$.

For $2\ellvtx\le\sidedistat{m}{i}\le4\ellvtx$, the weights $\sidedistat{m}{i}^{-\sigma_k}\bigl(1+\abs{\ln\sidedistat{m}{i}}\bigr)^{-L}$ are bounded above and below. The finite expansion is smooth there and is bounded by its coefficients. In the side coordinate, this interval lies in $E_{j,\rm away}^{(k-1)}$ because $\ell_{k-1}^{\rm away}<\frac{\ellvtx}{4}$ and $L_j\ge16\ellvtx$. Applying the parameter estimate of Proposition~\ref{prop:flat-milne-comp} on this set gives the $C^{q_k+1,\gamma}$ bound for $e_{j,k}$. The right-hand side is controlled by $\mathfrak F_{k-1,\kappa_{k-1}}$ and $\gaux_{j,k}$.

We next estimate the global Dirichlet data. The ray coefficients determine $H_k^{\rm sing}$ and $f_k^{\rm cut}$. Since $\sigma_k>0$, Lemma~\ref{lem:finite-part-remainder-trace} gives an $H^{\frac{1}{2}}$ trace for the ray remainders. The norms at the preceding order control the boundary data away from the endpoints. Combining these bounds with the Milne coefficient estimates, we obtain
\begin{align}
 \sum_m\sum_{i=1}^2\left\{\mathfrak N^{e,\mathsf d}_{m,i,k}(e_{m,i,k};q_k,k) +\mathfrak N^{B,\mathsf d}_{m,i,k,\kappa_k}(\Uside_{m,i,k};q_k,k)\right\} &\le C_k\big\{\mathcal G+\mathfrak N_{k-1}^\sharp(\mathsf d)\big\},
 \label{eq:hierarchy-Milne-data-estimate}\\
 \sum_j\nm{e_{j,k}}_{C^{q_k+1,\gamma}(E_{j,\rm away}^{(k-1)})} +\bnm{e_k-H_k^{\rm sing}}_{H^{\frac{1}{2}}(\partial\Om)} +\nm{f_k^{\rm cut}}_{C^{q_k-1,\gamma}(\overline\Om)} &\le C_k\big\{\mathcal G+\mathfrak N_{k-1}^\sharp(\mathsf d)\big\}.
 \label{eq:hierarchy-global-data-estimate}
\end{align}
For \eqref{eq:hierarchy-Milne-data-estimate}, use \eqref{eq:milne-mapping-estimate-comp} on the coefficients and Corollary~\ref{cor:milne-polyhom-comp} on the remainders, with $(\kappa_{\rm out},\kappa_{\rm in})=(\kappa_k,\kappa_{k-1})$, $\sigma=\sigma_k$, $q=q_k$ and $L=k$. Together with \eqref{eq:hierarchy-global-data-estimate}, this bounds $\EllData_{m,\sigma_k}(e_k)$ in \eqref{eq:elliptic-total-data-size}.

The derivative orders in \eqref{eq:hierarchy-global-data-estimate} also suffice for the Schauder estimate below. Indeed, $\mathfrak F_{k-1,\kappa_{k-1}}$ contains $q_{k-1}+1=q_k+2$ derivatives. Differentiating once in \eqref{eq:microscopic-recursion-comp} and \eqref{eq:endstate-recursion-comp} therefore gives $e_{j,k}$ in $C^{q_k+1,\gamma}(E_{j,\rm away}^{(k-1)})$. The source term $f_k^{\rm cut}$ is a finite sum formed from smooth harmonic lifts and cutoff derivatives on annuli away from the vertices. Its $C^{q_k-1,\gamma}$ norm is consequently bounded by the ray coefficients.

\emph{Elliptic estimates near the vertices.} By Definition~\ref{def:admissible-hierarchy-depth-comp}, the weight $\sigma_k=\mathsf d-k$ is positive and noncritical. We apply Proposition~\ref{prop:mellin-mapping-comp} at this weight to obtain $\rho_k$ and its remainder.

We check the hypotheses of Proposition~\ref{prop:mellin-mapping-comp} with $q=q_k$. By \eqref{eq:hierarchy-Milne-data-estimate}, $\mathcal R^{e,\mathsf d}_{m,i,k}\in\mathcal C^{q_k,\gamma}_{\sigma_k,L_k^{\rm ray}}(0,4\ellvtx)$, as required in \textup{(A2)}. The condition $q\ge2$ in \textup{(A1)} follows from $q_k\ge q_\Nexp=\max\{\pacc,3\}\ge3$. The norm $\mathcal N_{\rm glob}$ in \eqref{eq:explicit-Mellin-far-norm} requires $e_k\in C^{q_k,\gamma}$ away from the vertex balls and $f_k^{\rm cut}\in C^{q_k-2,\gamma}(\overline\Om)$. Both follow from \eqref{eq:hierarchy-global-data-estimate}, which gives one more derivative. The source term is supported where a derivative of $\chi_m^{\rm sing}$ is nonzero. By \textup{(A1)} at radius $4\ellvtx$, $\chi_m^{\rm sing}=1$ on $B_{4\ellvtx}(V_m)$, so this support is contained in $\Om_{\rm away}$. All hypotheses are therefore satisfied. By \eqref{eq:mellin-remainder-bound-comp}, $\mathcal R^{\rho,\mathsf d}_{m,k}\in\mathcal C^{q_k,\gamma}_{\sigma_k,L_k^{\rm ray}+1}$, with the derivative order and logarithmic bound used in \eqref{eq:weighted-interior-germ-norm-comp}.

The harmonic lifts depend linearly on the ray coefficients. Also, \eqref{eq:mellin-remainder-bound-comp} bounds every coefficient $c_{m,n}$ multiplying $Z^+_{m,n}$ in the global solution. Summing these finitely many bounds gives the coefficient estimate for $\rho_k$. Since $\Uint_k=\rho_k+\widehat{\Uint}_k$ and $\widehat{\Uint}_k=-w\cdot\nabla\Uint_{k-1}$ is already estimated, we obtain both the coefficient and remainder bounds for the vertex expansion of $\Uint_k$. The Milne estimate gives the corresponding bounds for $\Uside_{m,i,k}$. Thus all local terms in $\mathfrak N_k^\sharp(\mathsf d)$ are controlled.

The additional indices belong to the sets specified by \textup{(L2)} and \textup{(L5)}. A harmonic lift preserves the exponent in the boundary data and adds at most one logarithm. The homogeneous harmonic function $Z^+_{m,n}=\bulkdistat{m}^{\frac{n\pi}{\omega_m}}\Phi_{m,n}$ has $\alpha=k+\frac{n\pi}{\omega_m}>0$, and its coefficient is determined by the global Dirichlet problem.

\emph{Estimates away from the vertices.} It remains to bound $\mathfrak F_{k,\kappa_k}$. Estimate \eqref{eq:hierarchy-global-data-estimate} gives the boundary and source norms on the larger sets in \eqref{eq:hierarchy-buffered-domain-set}. We apply the variational estimate and local Schauder estimates to obtain $q_k+1$ derivatives on $\Om_{\rm away}^{(k)}$.

Cover $\Om_{\rm away}^{(k)}$ by finitely many interior balls and half-balls at flat sides, with doubled neighborhoods contained in $\Om_{\rm away}^{(k-1)}$ and disjoint from the vertices. The lower norm in the local Schauder estimate is bounded by the global $H^1$ norm. Hence
\begin{align}
 \nm{v_k^{\rm reg}}_{C^{q_k+1,\gamma}(\Om_{\rm away}^{(k)})}&\le C_k\left(\nm{v_k^{\rm reg}}_{H^1(\Om)}+\nm{f_k^{\rm cut}}_{C^{q_k-1,\gamma}(\overline\Om)}\right)
 +C_k\sum_j\nm{e_{j,k}-H_k^{\rm sing}}_{C^{q_k+1,\gamma}(E_{j,\rm away}^{(k-1)})}.
 \label{eq:hierarchy-buffered-Schauder}
\end{align}
By the variational estimate, the boundary remainder and source bounds above control $\nm{v_k^{\rm reg}}_{H^1(\Om)}$. The explicit harmonic lifts bound $H_k^{\rm sing}$ and its derivatives on $\Om_{\rm away}^{(k)}$ and $E_{j,\rm away}^{(k)}$. Thus $\rho_k=H_k^{\rm sing}+v_k^{\rm reg}$ has $q_k+1$ derivatives away from the vertices. Using $\Uint_k=\rho_k-w\cdot\nabla\Uint_{k-1}$ and $q_{k-1}+1=q_k+2$ gives the same estimate for $\Uint_k$.

For the side layer, apply $\mathfrak M_{\mu_j}$ on $E_{j,\rm away}^{(k-1)}$ to the inputs in \eqref{eq:endstate-recursion-comp}: the in-flow boundary function $\gaux_{j,k}-\widehat{\Uint}_k\big|_{E_j}$ and the source $-\tau_j\partial_s\Uside_{j,k-1}$. Their $C^{q_k+1,\gamma}$ norms with values in $L^\infty(\{\mu_j>0\})$ and $X_{\kappa_{k-1}}$, respectively, follow from the preceding order. The parameter estimate in Proposition~\ref{prop:flat-milne-comp} therefore bounds the decaying output $\Uside_{j,k}$ in $C^{q_k+1,\gamma}(E_{j,\rm away}^{(k-1)};X_{\kappa_k})$. Restricting to $E_{j,\rm away}^{(k)}$ gives the final term in $\mathfrak F_{k,\kappa_k}$. We conclude
\begin{align}
 \mathfrak F_{k,\kappa_k} &\leq C_k\left\{\mathcal G+\mathfrak N_{k-1}^\sharp(\mathsf d)\right\}.
 \label{eq:hierarchy-away-induction}
\end{align}
At $k=0$, the same argument uses the smooth in-flow boundary data $g_{j,0}$ on $E_{j,\rm away}^{(-1)}$. The local vertex estimates use only $q_k$ derivatives.

Combining the coefficient estimates, the weighted remainder estimates and \eqref{eq:hierarchy-away-induction}, we obtain $\mathfrak N_k^\sharp(\mathsf d)\le C_k\{\mathcal G+\mathfrak N_{k-1}^\sharp(\mathsf d)\}$, with $\mathfrak N_{-1}^\sharp(\mathsf d):=0$. Induction gives $\mathfrak N_k^\sharp(\mathsf d)\le C_{\mathsf d}\mathcal G$ for $0\le k\le\Nexp$. Since $X_{\kappa_k}\hookrightarrow X_{\kappa_\ast}$, the same bounds hold at weight $\kappa_\ast$. This proves part~\ref{item:hierarchy-complete-norm-comp} and the quantitative assertions in part~\ref{item:hierarchy-canonical-comp}.

\paragraph{\underline{Step 4: the exponents}} The exponents introduced by \textup{(L1)} and \textup{(L2)} satisfy $\alpha=k+n\ge0$ and $\alpha=k+\frac{n\pi}{\omega_m}>0$, respectively. Rule \textup{(L3)} preserves $\alpha$ by \eqref{eq:recurrence-degree-invariance-comp}, and \textup{(L4)}--\textup{(L7)} preserve both $k$ and $\lambda$. Hence \eqref{eq:generated-exponents-nonnegative} holds.

By Lemma~\ref{lem:ledger-finite-closed}, the index set is finite and closed under the stated rules. Together with $\alpha\ge0$, this proves part~\ref{item:hierarchy-ledger-comp}.

We also examine $\alpha=0$. By leading compatibility, the corresponding term of $\rho_0$ is a constant common to the two rays. The corresponding term of $\Uside_{m,i,0}$ is independent of the tangential variable. Their derivatives therefore vanish at the next order. All remaining prescribed Taylor terms have $\alpha\ge1$, and convexity gives $\alpha>1$ for every $Z^+_{m,n}$ introduced at any order. Induction yields
\begin{align}
 A^Q_{m,k,0,b}=e_{m,i,k,0,b}=B^{\rm sl}_{m,i,k,0,b}&=0, \qquad k+b>0,\quad Q\in\{\rho,U\}.
 \label{eq:hierarchy-zero-exponent-coefficients}
\end{align}
We retain these indices with zero coefficients when comparing expansions. Apart from the leading terms with $k=b=0$ and $\alpha=0$, every nonzero coefficient has $\alpha\ge1$.

\paragraph{\underline{Step 5: estimates for the remainders and the order-$\Nexp$ solutions}} Apply \eqref{eq:weighted-boundary-pointwise-comp}--\eqref{eq:weighted-sector-pointwise-comp} at each positive radius. Since $L=\Llog$ bounds both logarithmic degrees, these estimates give \eqref{eq:outer-cartesian-remainder-comp}, \eqref{eq:endstate-ledger-remainder-comp} and \eqref{eq:side-ledger-remainder-comp}. For the side remainder, we also use $X_{\kappa_k}\hookrightarrow X_{\kappa_\ast}$.

To estimate the polar derivatives, use $\partial_{\theta_m}=\bigl(x-V_m\bigr)^\perp\cdot\nabla_x$. Iterating this identity, with $\mathcal R=\mathcal R^{U,\mathsf d}_{m,k}$, gives
\begin{align}
 \bulkdistat{m}^{a}\nm{\partial_{\bulkdistat{m}}^a\partial_{\theta_m}^c\mathcal R}_{L^\infty_w} &\le C_{q_k}\sum_{j=0}^{a+c}\bulkdistat{m}^{j}\nm{D_x^j\mathcal R}_{L^\infty_w}, \qquad a+c\le q_k.
 \label{eq:polar-cartesian-weighted-derivative-comparison}
\end{align}
By \eqref{eq:outer-cartesian-remainder-comp}, each summand is bounded by $C\mathcal G\bulkdistat{m}^{\mathsf d-k}\bigl(1+\abs{\ln\bulkdistat{m}}\bigr)^{\Llog}$. Dividing by $\bulkdistat{m}^a$ proves \eqref{eq:outer-ledger-remainder-comp}. In particular, angular differentiation causes no additional radial loss. Expressing Cartesian derivatives in the fixed polar frame gives the reverse inequality, so the two weighted norms are equivalent on each rescaled annulus.

For a term in the vertex expansion at order $\Nexp$, \eqref{eq:generated-exponents-nonnegative} gives $\lambda=\alpha-\Nexp\ge-\Nexp$. One spatial or tangential derivative therefore has degree at least $-\Nexp-1$. The coefficient norms in \eqref{eq:weighted-interior-germ-norm-comp}--\eqref{eq:weighted-side-germ-norm-comp} give the required $L^\infty_w$ and $X_{\kappa_\ast}$ bounds for each term, with $\kappa_\ast=\kappa_\Nexp$. By \eqref{eq:outer-cartesian-remainder-comp} and \eqref{eq:side-ledger-remainder-comp}, the derivatives $w\cdot\nabla_x\mathcal R^{U,\mathsf d}_{m,\Nexp}$ and $\partial_{\sidedistat{m}{i}}\mathcal R^{B,\mathsf d}_{m,i,\Nexp}$ have weight $\mathsf d-\Nexp-1>-\Nexp-1$, since $\mathsf d>0$. They therefore satisfy the same bounds as the differentiated finite sums. Adding the corresponding remainder derivative proves \eqref{eq:terminal-complete-interior-comp} and \eqref{eq:terminal-complete-side-comp}. This completes parts~\ref{item:hierarchy-local-remainders-comp}--\ref{item:hierarchy-terminal-bounds-comp}.

\paragraph{\underline{Step 6: equations for the finite sums and remainders}}
\label{step:hierarchy-splitting-comp}
By the construction, the interior solution and the side layer satisfy
\begin{align}
 \qk\Uint_k+w\cdot\nabla_x\Uint_{k-1}&=0,
 \label{eq:full-local-recurrences-construction-split}\\
 \big(\mu_{m,i}\partial_\eta+\qk\big)\Uside_{m,i,k} +\tau_{m,i}\partial_{\sidedistat{m}{i}}\Uside_{m,i,k-1}&=0, \qquad 0\leq k\leq \Nexp,
 \label{eq:full-side-recurrence-construction-split}
\end{align}
where the quantities with index $-1$ are zero.

Insert \eqref{eq:exact-construction-projection-comp} into \eqref{eq:full-local-recurrences-construction-split}. The terms containing $C^U_{m,k}$ and $C^U_{m,k-1}$ form a finite sum of powers and logarithms of $\bulkdistat{m}$, with exponents $\alpha-k<\mathsf d-k$. By part~\ref{item:hierarchy-local-remainders-comp}, the terms containing $\mathcal R^U_{m,k}$ and $\mathcal R^U_{m,k-1}$ are $O(\bulkdistat{m}^{\mathsf d-k}\bigl(1+\abs{\ln\bulkdistat{m}}\bigr)^{\Llog})$ in $C([0,\omega_m];L^\infty(\Sone))$. Thus every exponent in the finite sum is below the remainder weight.

Collect the terms with each radial power. Their coefficients are polynomials in $\ln\bulkdistat{m}$ with values in $C([0,\omega_m];L^\infty(\Sone))$. Since \eqref{eq:full-local-recurrences-construction-split} is zero, Lemma~\ref{lem:finite-polyhom-uniqueness} with $d_{\rm rem}=\mathsf d-k$ implies that every polynomial vanishes. This gives the first identity in \eqref{eq:exact-interior-displayed-recurrence-comp}. Subtracting it from \eqref{eq:full-local-recurrences-construction-split} gives the second.

For the side layer, apply Corollary~\ref{cor:milne-polyhom-comp} with the remainder estimates from Step~3. Each coefficient of $C^\Uside_{m,i,k}$ satisfies the corresponding Milne equation. By \eqref{eq:recurrence-degree-invariance-comp}, differentiating a term of $C^\Uside_{m,i,k-1}$ preserves $\alpha<\mathsf d$. Its index therefore belongs to $\mathcal J^\Uside_{m,i,k}(\mathsf d)$ by \textup{(L3)}. Summing the coefficient equations gives \eqref{eq:exact-side-displayed-recurrence-comp}. Subtracting this identity from \eqref{eq:full-side-recurrence-construction-split} gives \eqref{eq:exact-side-projected-recurrence-comp}.

Multiplying the remainder identity in \eqref{eq:exact-interior-displayed-recurrence-comp} by $\e^{k-1}$ and summing, with $\mathcal R^U_{m,-1}=0$, gives
\begin{align}
 \Le\Bigl[\sum_{k=0}^\Nexp\e^k\mathcal R^U_{m,k}\Bigr]&=\e^\Nexp w\cdot\nabla_x\mathcal R^U_{m,\Nexp}.
 \label{eq:residual-telescope-comp}
\end{align}
Applying the same argument to \eqref{eq:exact-side-projected-recurrence-comp} gives the corresponding identity for the side remainders. In both cases, the only source term is the transport derivative at order $\Nexp$.

These identities hold on the punctured chart before localization and regularization at the vertices. The cutoff commutators are estimated in Section~\ref{sec:matching-realization}. This proves part~\ref{item:hierarchy-splitting-comp}.
\end{proof}

\paragraph*{\underline{Consistency across construction depths}}
Fix the polygon, data, $\Nexp$, cutoff choices and harmonic-lift normalization, and let $\Nexp<\mathsf d_1<\mathsf d_2$ be two admissible construction depths. Restriction to $\alpha<\mathsf d_1$ commutes with the index rules. At order zero, both constructions use the same Milne problem and weak harmonic Dirichlet problem. Inductively, equality of the preceding profiles gives identical order-$k$ Milne inputs, hence identical $e_{j,k}$ and $\Uside_{j,k}$. Lemma~\ref{lem:finite-polyhom-uniqueness} identifies their common ray coefficients. The additional prescribed harmonic lifts have $\mathsf d_1\le\alpha<\mathsf d_2$ and spatial degree $\alpha-k\ge\mathsf d_1-k>0$. Their cutoff realizations belong to $H^1(\Om)$, including the logarithmic terms. Lemma~\ref{lem:finite-part-uniqueness-comp} therefore transfers them between the prescribed part and regular correction without changing $\rho_k$, so $\Uint_k$ also agrees. Coefficient uniqueness and the fixed lifting convention identify the common interior coefficients. Thus the profiles and coefficients with $\alpha<\mathsf d_1$ are unchanged; each remainder at depth $\mathsf d_1$ equals its counterpart at depth $\mathsf d_2$ plus the intervening finite sum. The matching coefficients are then obtained by restricting the construction at $\Tcon$ to $\alpha<\Tmatch$.

\section{Wedge Layer}
\label{sec:wedge-layer}

Near a vertex, the side layers on the two incident sides interact with the singular parts of the interior solutions. We study this interaction in the rescaled sector $K_\omega$. For a prescribed $M$, we construct a wedge layer $D\in L^\infty(K_\omega\times\Sone)$ such that $\Uwedge=M+D$ solves the homogeneous transport equation with the prescribed in-flow boundary condition. The estimates require algebraically weighted bounds for $\Kw M$ and the difference between the prescribed in-flow boundary data and the trace of $M$. Section~\ref{sec:matching-realization} defines $M$ and verifies these bounds.

\subsection{Wedge Geometry}
\label{subsec:exact-wedge-blowup}

Fix a vertex $V_m$ and use the rotation and coordinates in \eqref{eq:wedge-sector-def}--\eqref{eq:corner-scaling-full}. Write $\omega=\omega_m$, $Y=Y_m$, $v=v_m$ and $\scaledbulk=\abs Y$. The coordinate $\eta_i$ is the distance to the line containing ray $i$, and $\sigma_i$ is the coordinate of the orthogonal projection onto that line. This projection lies on the ray if and only if $\sigma_i\ge0$.

If $F(x,w)=\widehat F(Y,v)$, then the chain rule and rotation invariance of the velocity average give
\begin{align}
 \e\Le F &=v\cdot\nabla_Y\widehat F+\widehat F-\pk\widehat F =:\Kw\widehat F.
 \label{eq:exact-blowup}
\end{align}
Since both components of $x-V_m$ are divided by $\e$, both spatial derivatives occur in $\Kw$.

Take ray $1$ as the positive horizontal axis. For $Y=\sigma_1\twedge{m}{1}+\eta_1\nwedge{m}{1}$, the coordinates associated with ray $2$ are $\sigma_2=\sigma_1\cos\omega+\eta_1\sin\omega$ and $\eta_2=\sigma_1\sin\omega-\eta_1\cos\omega$. If $\omega>\frac\pi2$, then $\sigma_2$ can be negative near ray $1$. We therefore extend the side layers from $\sigma_i>0$ as in Section~\ref{sec:matching-realization}.

For $0<\omega<\pi$, we have $\eta_1+\eta_2\ge\scaledbulk\sin\omega$. Hence the intersection of the two strips $0<\eta_i<c$ is contained in $\scaledbulk<\frac{2c}{\sin\omega}$. The wedge layer accounts for both side layers in this neighborhood of the vertex.

\subsection{Formulation of the Wedge Layer}
\label{subsec:M-W-D-precise}

Substitute $\bulkdist=\e\scaledbulk$ into the vertex expansions of the interior solutions and side layers. We collect the coefficients of each $\e^\alpha(\ln\e)^b$ and match them with a solution of the wedge equation, as in Section~\ref{subsec:physical-exponent-comp}. Under the assumptions of Lemma~\ref{lem:leading-field-convex-corner}, $\nabla\rho_0$ belongs to $L^\infty$ near $V_m$. The convex opening also permits $1<\betaw<\frac{\pi}{\omega_m}$, which is needed in the weighted estimate below.

Let $0<\omega<\pi$ and let $K_\omega$ be the sector in \eqref{eq:wedge-sector-def}. Denote the directed unit tangent and inward unit normal on ray $i$ by $t_i^\omega$ and $n_i^\omega$. The in-flow boundary on this ray is $\gamma_{i,-}:=\{(\sigma t_i^\omega,v):\sigma>0,\ v\cdot n_i^\omega>0\}$. We write its boundary data as functions of $(\sigma,v)$. All boundary traces and essential suprema use the measure $\abs{v\cdot n_i^\omega}\ud\sigma\ud\nu(v)$. At $V_m$, these definitions agree with $\omega=\omega_m$, $t_i^\omega=\twedge{m}{i}$ and $n_i^\omega=\nwedge{m}{i}$.

Let $M\in L^\infty_{\rm loc}(K_\omega\times\Sone)$. Assume that its restriction to almost every characteristic is locally absolutely continuous and has a finite limit at each finite endpoint on the in-flow boundary. Denote these limits by $M\big|_{\gamma_{i,-}}$ and define $\Kw M=v\cdot\nabla_YM+\qk M$ using the characteristic derivative. In our application, $M=M^{\rm con}_{m,\alpha,b}$ is given by \eqref{eq:Mcon-coefficient-expanded}. It can grow as a power of $\scaledbulk$ multiplied by a polynomial in $\ln\scaledbulk$. Let $G_i$ be the prescribed in-flow boundary data for $\Uwedge$.

Set $D=\Uwedge-M$. To obtain the homogeneous equation and the prescribed in-flow boundary condition for $\Uwedge$, we define
\begin{align}
 F&:=-\Kw M,\qquad h_i:=G_i-M\big|_{\gamma_{i,-}}.
 \label{eq:wedge-reference-source}
\end{align}
Here $h=h_i$ on $\gamma_{i,-}$. The vertex has zero kinetic boundary measure.

We seek the wedge layer $D$ satisfying
\begin{align}
 \Kw D&=F\quad\text{in }K_\omega\times\Sone,\qquad D\big|_{\gamma_{i,-}}=h_i,\quad i=1,2,
 \label{eq:defect-wedge-problem}
\end{align}
with $D\in L^\infty(K_\omega\times\Sone)$. We will prove the algebraic decay estimate \eqref{eq:wedge-weighted-solution-estimate}.

By the definitions of $F$ and $h_i$, adding the equations for $M$ and $D$ gives
\begin{align}
 \Kw\Uwedge&=0\quad\text{in }K_\omega\times\Sone,\qquad \Uwedge\big|_{\gamma_{i,-}}=G_i,\quad i=1,2.
 \label{eq:total-wedge-problem}
\end{align}
The decay of $D$ means that $\Uwedge$ matches $M$ as $\scaledbulk\to\infty$. In particular, $\Uwedge$ may inherit the polynomial and logarithmic growth of $M$. For the fixed $M$, we interpret \eqref{eq:total-wedge-problem} in the class $M+L^\infty(K_\omega\times\Sone)$. More precisely, $\widetilde\Uwedge$ is a solution in this class if $\widetilde\Uwedge-M$ satisfies the mild formula \eqref{eq:wedge-mild} with data \eqref{eq:wedge-reference-source}. Uniqueness below is understood in this sense.

The function $\Ubase^{\e,[\Nexp]}$ already contains the rescaled interior solutions and side layers that define $M$. The additional contribution at each collected order is therefore the wedge layer $D=\Uwedge-M$. We add these wedge layers with the cutoffs specified in Section~\ref{subsec:full-realization} to define $\Uapp^{\e,[\Nexp]}(\delta)$.

\subsection{Well-Posedness of the Wedge Layer Equation}

We prove existence, uniqueness, and algebraic decay for \eqref{eq:defect-wedge-problem} under the weighted data assumption \eqref{eq:abstract-wedge-data-assumption-v6}. The difficulty is the estimate of $\rho_D:=\pk D$ on the unbounded sector. Since $\qk$ vanishes on functions independent of $v$, the collision term does not directly control $\rho_D$. When we integrate along characteristics, $\rho_D$ enters as a source at every point of the characteristic. Averaging this formula over velocity gives a closed integral equation for $\rho_D$.

We estimate this equation by the barriers of Lemma~\ref{lem:sector-barrier}. The decaying barrier bounds the sum of the absolute values of the Neumann iterates. This gives a solution and its weighted $L^\infty$ estimate. To prove uniqueness among all $L^\infty$ mild solutions, we use a growing barrier, which forces every homogeneous $L^\infty$ solution to vanish without prescribing a limit at infinity. The decaying barrier requires $\betaw<\frac{\pi}{\omega}$. Proposition~\ref{prop:wedge-wellposedness} gives the resulting existence, uniqueness, and decay statement. Remark~\ref{rem:exponential-data-algebraic-wedge-decay} explains the relation to the expected diffusive behavior as $\scaledbulk\to\infty$.

\subsubsection{An integral equation for the velocity average}
\label{subsec:collision-skeleton}

For $(Y,v)\in K_\omega\times\Sone$, define $t_-(Y,v):=\inf\{t>0:Y-tv\notin K_\omega\}\in(0,\infty]$ and set $\ue^{-\infty}=0$. Let $F\in L^\infty(K_\omega\times\Sone)$ and $h_i\in L^\infty(\gamma_{i,-})$. We call $D\in L^\infty(K_\omega\times\Sone)$ a mild solution of \eqref{eq:defect-wedge-problem} if it satisfies, for almost every $(Y,v)$,
\begin{align}
 D(Y,v) &=\mathds 1_{\{t_-<\infty\}}\ue^{-t_-} h(Y-t_-v,v)+\int_0^{t_-}\ue^{-\ell} \big\{\rho_D(Y-\ell v)+F(Y-\ell v,v)\big\}\ud\ell.
 \label{eq:wedge-mild}
\end{align}
If $t_-<\infty$, we use $h$ on the ray reached by the backward characteristic. If $t_-=\infty$, the boundary term is zero. For each $Y$, the velocities tangent to a ray or directed through the vertex form a $\nu$-null set. Formula \eqref{eq:wedge-mild} also gives the prescribed characteristic trace almost everywhere on the in-flow boundary. Indeed, for almost every $(y,v)$ there, $t_-(y+\tau v,v)=\tau$ for small $\tau>0$. The integral has absolute value at most $\tau(\nm{\rho_D}_{L^\infty(K_\omega)}+\nm F_{L^\infty(K_\omega\times\Sone)})$ for almost every such $\tau$, and hence tends to zero as $\tau\downarrow0$.

The only unknown on the right-hand side of \eqref{eq:wedge-mild} is $\rho_D$. Averaging over $v$ therefore gives an equation for this velocity average alone:
\begin{align}
 \rho_D&=\Pw\rho_D+H.
 \label{eq:rho-integral-equation}
\end{align}
The function $H$ collects the prescribed boundary and volume data:
\begin{align}
 H(Y)&:=\int_{\Sone}\mathds 1_{\{t_-<\infty\}}\ue^{-t_-} h(Y-t_-v,v)\ud\nu(v) +\int_{\Sone}\int_0^{t_-}\ue^{-\ell} F(Y-\ell v,v)\ud\ell\ud\nu(v).
 \label{eq:wedge-averaged-data}
\end{align}
The contribution of $\rho_D$ is the integral operator
\begin{align}
 (\Pw f)(Y)&:=\int_{K_\omega}q(Y-Z)f(Z)\ud Z,\qquad q(X):=\frac{\ue^{-\abs X}}{2\pi\abs X}.
 \label{eq:restricted-convolution-operator}
\end{align}

To derive \eqref{eq:restricted-convolution-operator}, set $X=\ell v$ and use $\ud X=\ell\ud\ell\ud\varphi$. Then $\ue^{-\ell}\ud\ell\,\frac{\ud\varphi}{2\pi}=q(X)\ud X$. For fixed $Y\in K_\omega$, put $Z=Y-X$. The condition $0<\ell<t_-(Y,v)$ ensures that $Z$ lies in the sector. Conversely, every $Z\in K_\omega$ with $Z\ne Y$ can be reached in this way, since convexity keeps the segment from $Y$ to $Z$ inside $K_\omega$. This gives the restriction of the convolution to $K_\omega$ in \eqref{eq:restricted-convolution-operator}.

Let $\mathds 1$ denote the constant function one on $K_\omega$. For every finite $Y$, we have $\Pw\mathds 1(Y)<1$ by \eqref{eq:strict-mass-deficit}. The quantity $1-\Pw\mathds 1(Y)$ is the probability that a characteristic reaches the boundary before the next collision. This probability can be arbitrarily small far from both rays. Indeed, let $e_{\rm bis}$ be the unit bisector. The ball with center $\scaledbulk e_{\rm bis}$ and radius $\scaledbulk\sin\bigl(\frac{\omega}{2}\bigr)$ is contained in $K_\omega$. Integrating the exponential tail of $q$ gives
\begin{align}
 0<1-\Pw\mathds 1(\scaledbulk e_{\rm bis})&\le\ue^{-\scaledbulk\sin\left(\frac{\omega}{2}\right)}\longrightarrow0\quad\text{as }\scaledbulk\to\infty.
 \label{eq:wedge-escape-probability}
\end{align}
Hence $\nm{\Pw}_{L^\infty\to L^\infty}=1$. The contraction argument on the bounded domain $\Om$ therefore does not apply. We instead estimate the iterates of $\Pw$ by a positive barrier.

The barrier construction uses the positivity of $q$ and the moment identities
\begin{align}
 \int_{\mathbb R^2}q(X)\,\ud X&=1,\qquad \int_{\mathbb R^2}X_iq(X)\,\ud X=0,\qquad \int_{\mathbb R^2}X_iX_jq(X)\,\ud X=\delta_{ij},
 \label{eq:q-moments}
\end{align}
for $i,j\in\{1,2\}$. Polar integration also gives $\int_{\mathbb R^2}\abs X^kq(X)\,\ud X<\infty$ for every $k\ge0$.

For $s>0$, define the weighted norms
\begin{align}
 \nm F_{\infty,s}&:=\esssup_{Y\in K_\omega,\,v\in\Sone} \jbr{\scaledbulk}^{s}\abs{F(Y,v)},
 \label{eq:wedge-data-weighted-norms} \\
 \abs h_{\infty,s}&:=\max_{i=1,2}\esssup_{ (\sigma t_i^\omega,v)\in\gamma_{i,-}} \jbr{\sigma}^{s}\abs{h_i(\sigma,v)}.
 \label{eq:wedge-incoming-weighted-norm}
\end{align}

\begin{lemma}
\label{lem:H-decay}
Let $s>0$. If
\begin{align}
 \mathcal A_{s} :=\nm F_{\infty,s}+\abs h_{\infty,s}<\infty,
 \label{eq:abstract-wedge-data-assumption-v6}
\end{align}
then
\begin{align}
 \abs{H(Y)}\le C_s\mathcal A_{s}\jbr{\scaledbulk}^{-s}\quad\text{for almost every }Y\in K_\omega.
 \label{eq:wedge-averaged-data-decay}
\end{align}
\end{lemma}

\begin{proof}
For the volume source, split the characteristic integral at $\ell=\frac{\scaledbulk}{2}$. If $\ell<\frac{\scaledbulk}{2}$, then $\abs{Y-\ell v}\ge\scaledbulk-\ell>\frac{\scaledbulk}{2}$. The weighted bound for $F$ gives a contribution at most $C\nm F_{\infty,s}\jbr{\scaledbulk}^{-s}$. The remaining integral is bounded by $\ue^{-\frac{\scaledbulk}{2}}\nm F_{\infty,s}\le C_s\jbr{\scaledbulk}^{-s}\nm F_{\infty,s}$.

For the boundary term, let $t_-<\infty$ and set $Z=Y-t_-v$. If $t_-\ge\frac{\scaledbulk}{2}$, then $\ue^{-t_-}\le\ue^{-\frac{\scaledbulk}{2}}$. Otherwise, $\abs Z\ge\abs Y-t_->\frac{\scaledbulk}{2}$, and we apply the weighted bound for $h_i(\abs Z,v)$. Integrating in $v$ proves \eqref{eq:wedge-averaged-data-decay}.
\end{proof}

\subsubsection{A shifted superharmonic barrier}
\label{subsec:shifted-barrier}

\begin{figure}[!htbp]
\centering
\begin{tikzpicture}[scale=1.05,>=latex]
  \def\Rmax{3.7}
  \def\lsh{1.5}
  \fill[blue!8] (0,0) -- (40:\Rmax) arc (40:140:\Rmax) -- cycle;
  \draw[gray!65,dashed] (0,-\lsh) -- ++(40:5.5);
  \draw[gray!65,dashed] (0,-\lsh) -- ++(140:5.5);
  \draw[gray!40,dotted] (0,-\lsh) -- ++(24:5.8);
  \draw[gray!40,dotted] (0,-\lsh) -- ++(156:5.8);
  \draw[very thick] (40:\Rmax) -- (0,0) -- (140:\Rmax);
  \draw[->,thin,gray!80] (0,-\lsh) -- (0,\Rmax+0.3) node[right,black]{$e_{\rm bis}$};
  \draw[|<->|,thin] (-0.3,-\lsh) -- (-0.3,0) node[midway,left]{$\ell_0$};
  \coordinate (Y) at (70:2.55);
  \draw[very thick,red!65!black] (0,0) -- (Y);
  \node[red!55!black] at (0.30,1.55) {$\scaledbulk$};
  \draw[very thick,blue!55!black] (0,-\lsh) -- (Y);
  \node[blue!45!black] at (1.02,1.25) {$\varrho$};
  \draw[thin] (40:0.8) arc (40:140:0.8);
  \node at (0,1.02) {$\omega$};
  \draw[thin] (0,-\lsh) ++(77:1.75) arc (77:90:1.75);
  \node at (0.34,0.53) {$\vartheta$};
  \fill (Y) circle (1.7pt) node[above right]{$Y$};
  \fill (0,0) circle (1.7pt) node[below right,xshift=1pt]{$V_m$};
  \fill (0,-\lsh) circle (1.7pt);
  \node[below] at (0,-\lsh-0.08) {$-\ell_0e_{\rm bis}$};
\end{tikzpicture}
\caption{The shifted barrier on $K_\omega$. Its pole $-\ell_0e_{\rm bis}$ lies outside $\overline{K_\omega}$. On the sector, $\varrho\ge\ell_0$ and $\abs\vartheta\le\frac\omega2$, so $\varrho^s\cos(a\vartheta)$ is smooth and positive for $0<a<\frac\pi\omega$. The dashed cone has aperture $\omega$. The dotted cone has aperture $\omega'$ and contains the angular support of the extension. The cutoffs change the barrier only outside the closed sector.}
\label{fig:shifted-barrier-geometry}
\end{figure}

By \eqref{eq:q-moments}, Taylor expansion against $q$ gives the function itself and the second-order term $\frac12\Delta$. The first-order term vanishes. We therefore construct a positive superharmonic function and estimate the Taylor remainder. The decaying barrier will satisfy $B_{\rm dec}\simeq\jbr{\scaledbulk}^{-\betaw}$ and $(I-\Pw)B_{\rm dec}\ge c\jbr{\scaledbulk}^{-\betaw-2}$, as stated in \eqref{eq:global-barrier}--\eqref{eq:decaying-barrier-defect}. This gives the two additional powers required of the source term.

We move the pole outside the sector along the negative bisector so that the negative power remains smooth at $Y=0$. The shifted radius is then bounded below by the shift length. Increasing this length makes the third-order Taylor remainder uniformly smaller than the superharmonic second-order term. We also extend the barrier to a nonnegative function on the whole plane. Its full-plane convolution bounds the restricted convolution from above, allowing us to use the moment identities for $q$.

\FloatBarrier

\begin{lemma}
\label{lem:sector-barrier}
Let $0<\omega<\pi$, and let $s\in\mathbb R$ and $a>0$ satisfy
\begin{align}
 \abs s<a<\frac\pi\omega .
 \label{eq:barrier-exponent-range}
\end{align}
For $\ell_0\ge1$, let $(\varrho,\vartheta)$ be polar coordinates centered at the pole $-\ell_0e_{\rm bis}$, with angle measured from the unit bisector $e_{\rm bis}$. There is $\ell_0=\ell_0(\omega,s,a)\ge1$ such that
\begin{align}
 B(Y):=\varrho^{\,s}\cos(a\vartheta)
 \label{eq:sector-barrier-def}
\end{align}
is positive and smooth on $\overline{K_\omega}$ and satisfies, for constants $c,C>0$ depending only on $\omega,s,a$,
\begin{align}
 c\jbr{\scaledbulk}^{\,s}&\le B(Y)\le C\jbr{\scaledbulk}^{\,s},
 \label{eq:sector-barrier-estimates} \\
 (I-\Pw)B(Y)&\ge c\jbr{\scaledbulk}^{\,s-2},\qquad Y\in\overline{K_\omega}.
 \label{eq:sector-barrier-defect}
\end{align}
\end{lemma}

\begin{proof}\leavevmode

\paragraph{\underline{Step 1: geometry}} For $Y\in K_\omega$, adding $\ell_0e_{\rm bis}$ decreases the angle with the bisector. Since the angle of $Y$ is at most $\frac{\omega}{2}$, we obtain
\begin{align}
 \abs\vartheta&\le\frac\omega2,\qquad \varrho=\abs{Y+\ell_0e_{\rm bis}}\ge\ell_0,\qquad \varrho\simeq\ell_0+\scaledbulk\quad\text{on }\overline{K_\omega}.
 \label{eq:shift-geometry}
\end{align}
The last comparison follows from $\varrho^2=\scaledbulk^2+\ell_0^2+2\ell_0\scaledbulk\cos\angle(Y,e_{\rm bis})$ and $\cos\angle(Y,e_{\rm bis})\ge\cos\bigl(\frac{\omega}{2}\bigr)>0$. In particular, the pole lies outside $\overline{K_\omega}$, and $\varrho^s$ is smooth and finite at the vertex. See Figure~\ref{fig:shifted-barrier-geometry}.

\paragraph{\underline{Step 2: superharmonicity}} By \eqref{eq:barrier-exponent-range}, $\frac{a\omega}{2}<\frac\pi2$. Together with \eqref{eq:shift-geometry}, this gives $\cos(a\vartheta)\ge\cos\bigl(\frac{a\omega}{2}\bigr)>0$ and hence $B>0$ on $\overline{K_\omega}$. Direct differentiation gives
\begin{align}
 -\Delta B&=(a^2-s^2)\varrho^{\,s-2}\cos(a\vartheta)\ge c_\ast\varrho^{\,s-2},
 \label{eq:sector-barrier-laplacian}
\end{align}
where $c_\ast=(a^2-s^2)\cos\bigl(\frac{a\omega}{2}\bigr)>0$. Both strict inequalities in \eqref{eq:barrier-exponent-range} are used here. If $a=\abs s$, then $\Delta B=0$. If $a=\frac\pi\omega$, then $\cos(a\vartheta)$ tends to zero along either boundary ray as $\scaledbulk\to\infty$.

\paragraph{\underline{Step 3: extension to $\mathbb R^2$}} To apply the full-plane moment identities, we extend $B$ to a nonnegative smooth function on $\mathbb R^2$. Choose $\omega<\omega'<\min\bigl\{\pi,\frac\pi a\bigr\}$. Let $\chi$ be smooth, with $0\le\chi\le1$, equal to one on $\abs\vartheta\le\frac\omega2$ and zero on $\abs\vartheta\ge\frac{\omega'}2$. Let $\zeta$ be smooth, with $0\le\zeta\le1$, zero on $[0,\frac14]$ and one on $[\frac12,\infty)$. Set
\begin{align}
 \widetilde B(Z):=\zeta\left(\frac{\varrho(Z)}{\ell_0}\right)\varrho\bigl(Z\bigr)^{\,s}\chi\bigl(\vartheta(Z)\bigr)\cos\bigl(a\vartheta(Z)\bigr),
 \label{eq:sector-barrier-extension}
\end{align}
and set $\widetilde B=0$ near the pole. The angular support avoids the polar branch cut and lies where the cosine is positive. Hence $\widetilde B$ is nonnegative and smooth on $\mathbb R^2$, and $\widetilde B=B$ on $\overline{K_\omega}$.

\paragraph{\underline{Step 4: derivative estimates}} Differentiating \eqref{eq:sector-barrier-extension}, we obtain
\begin{align}
 \abs{\nabla^j\widetilde B(Z)}\le C\bigl(\ell_0+\varrho(Z)\bigr)^{\,s-j}, \qquad 0\le j\le3,
 \label{eq:barrier-extension-derivatives}
\end{align}
with $C$ independent of $\ell_0$. Indeed, the angular cutoff does not depend on $\ell_0$. Each derivative of the radial cutoff contributes $\ell_0^{-1}$, and $\ell_0^{-1}\simeq\bigl(\ell_0+\varrho\bigr)^{-1}$ on its support.

\paragraph{\underline{Step 5: Taylor expansion}} For $\ell_0\ge1$, Taylor's formula and \eqref{eq:q-moments} give, uniformly for $Y\in\overline{K_\omega}$,
\begin{align}
 \int_{\mathbb R^2}q(X)\widetilde B(Y-X)\,\ud X&=B(Y)+\frac12\Delta B(Y)+E(Y),
 \label{eq:barrier-taylor} \\
 \abs{E(Y)}&\le C_\ast\,\ell_0^{-1}\varrho^{\,s-2},
 \label{eq:barrier-taylor-remainder}
\end{align}
where $C_\ast$ is independent of $\ell_0\ge1$. Write $\varrho=\varrho(Y)$. If $\abs X\le\frac\varrho2$, the shifted radius on the segment from $Y$ to $Y-X$ lies between $\frac\varrho2$ and $\frac{3\varrho}2$. By \eqref{eq:barrier-extension-derivatives}, the third-order remainder is at most $C\abs X^3\varrho^{s-3}$. Its integral against $q$ is therefore bounded by $C\varrho^{s-3}$.

If $\abs X>\frac\varrho2$, then $\varrho\ge\ell_0\ge1$. For $s\ge0$, we have $\ell_0+\varrho(Y-X)\le\ell_0+\varrho+\abs X\le5\abs X$, so \eqref{eq:barrier-extension-derivatives} gives $\babs{\widetilde B(Y-X)}\le C\abs X^s$. For $s<0$, it gives $\babs{\widetilde B(Y-X)}\le C\ell_0^s\le C$. The integral over this region is bounded by $C\ue^{-\frac\varrho4}$. The tails of the first three Taylor terms satisfy the same bound after their polynomial factors are absorbed by the exponential. Since $\ue^{-\frac\varrho4}\le C\varrho^{s-3}$ for $\varrho\ge1$, we obtain $\abs E\le C\varrho^{s-3}\le C\ell_0^{-1}\varrho^{s-2}$. This proves \eqref{eq:barrier-taylor-remainder}.

\paragraph{\underline{Step 6: choice of $\ell_0$}} By \eqref{eq:barrier-taylor-remainder} and \eqref{eq:sector-barrier-laplacian}, choose
\begin{align}
 \ell_0=\ell_0(\omega,s,a):=\max\Big\{1,\frac{4C_\ast}{c_\ast}\Big\},
 \label{eq:shift-length-choice}
\end{align}
Then $\abs E\le\frac{c_\ast}4\varrho^{s-2}$. Since $\widetilde B\ge0$ and $\widetilde B=B$ on the sector, \eqref{eq:barrier-taylor} gives
\begin{align}
 (I-\Pw)B(Y)&\ge B(Y)-\int_{\mathbb R^2}q(X)\widetilde B(Y-X)\,\ud X \ge-\frac12\Delta B(Y)-\abs{E(Y)}\\*
 &\ge\frac{c_\ast}2\varrho^{s-2}-\frac{c_\ast}4\varrho^{s-2} =\frac{c_\ast}4\varrho^{s-2}.\notag
\end{align}
For this fixed $\ell_0$, \eqref{eq:shift-geometry} gives $\varrho\simeq\jbr{\scaledbulk}$. Combining the last inequality with $\cos(a\vartheta)\ge\cos\bigl(\frac{a\omega}2\bigr)>0$ proves \eqref{eq:sector-barrier-estimates}--\eqref{eq:sector-barrier-defect}.
\end{proof}

For existence, apply Lemma~\ref{lem:sector-barrier} with
\begin{align}
 (s,a)&=\left(-\betaw,\frac12\Bigl(\betaw+\frac\pi\omega\Bigr)\right),
 \label{eq:decaying-barrier-parameters}
\end{align}
The condition \eqref{eq:barrier-exponent-range} holds because $\betaw<\frac\pi\omega$. Thus
\begin{align}
 c_1\jbr{\scaledbulk}^{-\betaw}&\le B_{\rm dec}(Y)\le C_1\jbr{\scaledbulk}^{-\betaw},
 \label{eq:global-barrier} \\
 (I-\Pw)B_{\rm dec}(Y)&\ge c_1\jbr{\scaledbulk}^{-\betaw-2},
 \label{eq:decaying-barrier-defect}
\end{align}
where $c_1,C_1$ depend only on $\omega,\betaw$. If the data exponent $s$ in \eqref{eq:abstract-wedge-data-assumption-v6} satisfies $s\ge\betaw+2$, Lemma~\ref{lem:H-decay} bounds $\abs H$ by a constant multiple of $(I-\Pw)B_{\rm dec}$. Applying the positive operator $\Pw$ repeatedly then bounds the Neumann partial sums by a telescoping sum. Thus the data decay with at least two more powers than the solution estimate produced by this barrier.

For uniqueness, let $\gamma_{\rm g}=\frac12\min\bigl\{1,\frac\pi{2\omega}\bigr\}$ and apply Lemma~\ref{lem:sector-barrier} with a separate shift and
\begin{align}
 (s,a)&=\left(\gamma_{\rm g},\frac{\pi}{2\omega}\right),
 \label{eq:growing-barrier-parameters}
\end{align}
Since $0<\gamma_{\rm g}<\frac\pi{2\omega}<\frac\pi\omega$, these parameters satisfy \eqref{eq:barrier-exponent-range}. We obtain
\begin{align}
 c\jbr{\scaledbulk}^{\gamma_{\rm g}}&\le B_{\rm grow}(Y)\le C\jbr{\scaledbulk}^{\gamma_{\rm g}},
 \label{eq:growing-wedge-barrier} \\
 (I-\Pw)B_{\rm grow}(Y)&\ge c\jbr{\scaledbulk}^{\gamma_{\rm g}-2}\ge0,
 \label{eq:growing-barrier-defect}
\end{align}
where $c,C$ depend only on $\omega$. A homogeneous $L^\infty$ solution of the averaged equation has a continuous representative. Subtracting a positive multiple of $B_{\rm grow}$ from this representative gives a function tending to $-\infty$ as $\scaledbulk\to\infty$. The resulting difference attains its maximum at a finite point, where the strict inequality $\Pw\mathds 1<1$ gives the comparison argument below.

\subsubsection{Existence, uniqueness, and decay of the wedge layer}

\begin{lemma}
\label{lem:growing-wedge-barrier}
Let $0<\omega<\pi$. Every $u\in L^\infty(K_\omega)$ satisfying
\begin{align}
 u&=\Pw u \quad\hbox{almost everywhere in }K_\omega
 \label{eq:wedge-homogeneous-scalar}
\end{align}
vanishes almost everywhere.
\end{lemma}

\begin{proof}
Extend $u$ by zero outside $K_\omega$ and set $u^\sharp=q*(u\mathds 1_{K_\omega})$. Since $q\in L^1(\mathbb R^2)$ and $u\mathds 1_{K_\omega}\in L^\infty(\mathbb R^2)$, the function $u^\sharp$ is continuous on $\mathbb R^2$. By \eqref{eq:wedge-homogeneous-scalar}, it agrees with $u$ almost everywhere in $K_\omega$. Hence $\Pw u^\sharp=u^\sharp$ almost everywhere in $K_\omega$, and continuity extends this identity to every point of $\overline{K_\omega}$. Also $\abs{u^\sharp(Y)}\le\nm u_{L^\infty(K_\omega)}$ for every $Y$, since $q$ has mass one. It suffices to treat real-valued $u$, because the same argument applies to its real and imaginary parts.

Fix $\delta>0$ and $\iota\in\{-1,1\}$, and set $\phi_\delta=\iota u^\sharp-\delta B_{\rm grow}$. Since $u^\sharp$ is bounded and $B_{\rm grow}(Y)\ge c\jbr{\scaledbulk}^{\gamma_{\rm g}}$, we have $\phi_\delta(Y)\to-\infty$ as $\scaledbulk\to\infty$. Thus $\phi_\delta$ attains its maximum $m_\delta=\max_{\overline{K_\omega}}\phi_\delta$ at some finite $Y_\delta\in\overline{K_\omega}$.

The convolution $\Pw\phi_\delta(Y)$ is absolutely convergent at every finite $Y$. Indeed, $B_{\rm grow}(Z)\le C\jbr{\abs Z}^{\gamma_{\rm g}}$ and
\begin{align}
 \int_{K_\omega}q(Y-Z)\jbr{\abs Z}^{\gamma_{\rm g}}\,\ud Z&\le C\jbr{\abs Y}^{\gamma_{\rm g}}\int_{\mathbb R^2}q(X)\jbr{\abs X}^{\gamma_{\rm g}}\,\ud X<\infty.
 \label{eq:growing-comparison-integrability}
\end{align}
Thus linearity and positivity of $\Pw$ apply to $\phi_\delta$.

Suppose $m_\delta>0$. By \eqref{eq:wedge-homogeneous-scalar} and \eqref{eq:growing-barrier-defect}, $\phi_\delta=\Pw\phi_\delta-\delta(I-\Pw)B_{\rm grow}\le\Pw\phi_\delta$. Since $\phi_\delta\le m_\delta$ and $\Pw$ is positive, we obtain $\Pw\phi_\delta(Y_\delta)\le m_\delta\Pw\mathds 1(Y_\delta)$.

For every finite $Y\in\overline{K_\omega}$, we have
\begin{align}
 \Pw\mathds 1(Y)=\int_{K_\omega}q(Y-Z)\,\ud Z<1.
 \label{eq:strict-mass-deficit}
\end{align}
Indeed, $\mathbb R^2\setminus K_\omega$ has positive Lebesgue measure and $q>0$ almost everywhere. Therefore $m_\delta\le m_\delta\Pw\mathds 1(Y_\delta)<m_\delta$, a contradiction. Only the strict inequality at the finite point $Y_\delta$ is needed.

Hence $\phi_\delta\le0$. For each fixed $Y\in\overline{K_\omega}$, let $\delta\downarrow0$ to obtain $\iota u^\sharp(Y)\le0$. Taking both signs gives $u^\sharp=0$, so $u=0$ almost everywhere.
\end{proof}

\begin{proposition}[Well-posedness of wedge layer]
\label{prop:wedge-wellposedness}
Let $0<\omega<\pi$, $1<\betaw<\frac{\pi}{\omega}$ and $s\ge\betaw+2$. Suppose the measurable source and in-flow boundary data in \eqref{eq:defect-wedge-problem} satisfy $\mathcal A_{s}<\infty$ in \eqref{eq:abstract-wedge-data-assumption-v6}. Then this problem has a unique mild solution $D\in L^\infty(K_\omega\times\Sone)$, and
\begin{align}
 \esssup_{Y,v}\jbr{\scaledbulk}^\betaw \big(\abs{D(Y,v)}+\abs{\pk D(Y)}\big) \le C(\omega,\betaw,s)\mathcal A_{s}.
 \label{eq:wedge-weighted-solution-estimate}
\end{align}
The solution map is linear and continuous from the data space with norm $\mathcal A_s$ into the space with norm $\nm{\jbr{\scaledbulk}^{\betaw}D}_{L^\infty(K_\omega\times\Sone)}$.
\end{proposition}

\begin{proof}\leavevmode

\paragraph{\underline{Step 1: the averaged equation}} By Lemma~\ref{lem:H-decay} and $s\ge\betaw+2$,
\begin{align}
 \abs{H(Y)}\le C\mathcal A_{s}\jbr{\scaledbulk}^{-s} \le C\mathcal A_{s}\jbr{\scaledbulk}^{-\betaw-2}\quad\text{for almost every }Y\in K_\omega.
 \label{eq:H-two-orders-down}
\end{align}
The right-hand side has the same weight as $(I-\Pw)B_{\rm dec}$ in \eqref{eq:decaying-barrier-defect}. The estimate includes $s=\betaw+2$.

Combining \eqref{eq:H-two-orders-down} with \eqref{eq:decaying-barrier-defect} gives $\abs H\le C\mathcal A_{s}(I-\Pw)B_{\rm dec}$. Positivity of $\Pw$ then gives, for every $J\ge0$,
\begin{align}
 \sum_{n=0}^{J}\Pw^n\abs H&\le C\mathcal A_{s}\sum_{n=0}^{J}\Pw^n(I-\Pw)B_{\rm dec}
 \label{eq:neumann-telescope} =C\mathcal A_{s}\bigl(B_{\rm dec}-\Pw^{J+1}B_{\rm dec}\bigr)\le C\mathcal A_{s}B_{\rm dec},
\end{align}
Since $\Pw^{J+1}B_{\rm dec}\ge0$, $\rho^\sharp:=\sum_{n\ge0}\Pw^nH$ converges absolutely for almost every $Y\in K_\omega$. Its sum belongs to $L^\infty(K_\omega)$ and satisfies, almost everywhere,
\begin{align}
 \abs{\rho^\sharp(Y)}&\le\sum_{n\ge0}\abs{\Pw^nH(Y)}\le C\mathcal A_s B_{\rm dec}(Y)\le C\mathcal A_s\jbr{\scaledbulk}^{-\betaw}.
 \label{eq:rho-sharp-bound}
\end{align}
By \eqref{eq:rho-sharp-bound}, the partial sums are bounded almost everywhere by $C\mathcal A_sB_{\rm dec}$. Also, $\Pw B_{\rm dec}(Y)\le B_{\rm dec}(Y)<\infty$ for each $Y\in K_\omega$. Dominated convergence in \eqref{eq:restricted-convolution-operator} therefore permits us to interchange $\Pw$ and the series. Thus, almost everywhere,
\begin{align}
 \rho^\sharp&=H+\sum_{n\ge1}\Pw^nH=H+\Pw\rho^\sharp.
 \label{eq:rho-sharp-fixed-point}
\end{align}

\paragraph{\underline{Step 2: construction and decay of $D$}} Define $D$ by the right-hand side of \eqref{eq:wedge-mild}, with $\rho_D$ replaced by $\rho^\sharp$. Split the characteristic integrals as in Lemma~\ref{lem:H-decay}, now with exponent $\betaw$. By \eqref{eq:rho-sharp-bound} and $s\ge\betaw+2$, the same estimates apply to $\rho^\sharp$, $F$ and $h$ before velocity averaging. For almost every $Y\in K_\omega$, we obtain
\begin{align}
 \nm{D(Y,\cdot)}_{L^\infty(\Sone)}\le C\mathcal A_{s}\jbr{\scaledbulk}^{-\betaw}.
 \label{eq:D-pointwise-bound}
\end{align}
Averaging the defining characteristic formula for $D$ over velocity gives $\pk D=\Pw\rho^\sharp+H=\rho^\sharp$ almost everywhere, by \eqref{eq:rho-sharp-fixed-point}. Thus the function inserted into the characteristic formula is the velocity average of the resulting $D$, and $D\in L^\infty(K_\omega\times\Sone)$ satisfies \eqref{eq:wedge-mild} with $\rho_D=\pk D$. Its in-flow boundary trace is $h$ by the argument following \eqref{eq:wedge-mild}. Combining \eqref{eq:rho-sharp-bound} and \eqref{eq:D-pointwise-bound} proves \eqref{eq:wedge-weighted-solution-estimate}.

The constants in \eqref{eq:global-barrier} depend only on $\omega,\betaw$, and the characteristic integral estimates also depend on $s$. Hence the constant in \eqref{eq:wedge-weighted-solution-estimate} has the stated dependence.

\paragraph{\underline{Step 3: uniqueness and dependence on the data}} Let $D\in L^\infty(K_\omega\times\Sone)$ solve the homogeneous problem with $F=h=0$. Then $\rho_D=\pk D\in L^\infty(K_\omega)$ satisfies \eqref{eq:wedge-homogeneous-scalar}. Lemma~\ref{lem:growing-wedge-barrier} gives $\rho_D=0$ almost everywhere, and \eqref{eq:wedge-mild} gives $D=0$ almost everywhere. Linearity and \eqref{eq:wedge-weighted-solution-estimate} give the stated continuity of the solution map.
\end{proof}

\begin{corollary}
\label{cor:affine-total-wedge-field-v6}
Let $M$ and $G_i$ be prescribed as in Section~\ref{subsec:M-W-D-precise}, and suppose the data \eqref{eq:wedge-reference-source} satisfy the hypotheses of Proposition~\ref{prop:wedge-wellposedness}. Then $\Uwedge=M+D$ is the unique total wedge solution of \eqref{eq:total-wedge-problem} in $M+L^\infty(K_\omega\times\Sone)$, with the mild interpretation specified in that subsection. The difference $\Uwedge-M$ satisfies the weighted $L^\infty$ estimate \eqref{eq:wedge-weighted-solution-estimate}.
\end{corollary}

\begin{proof}
Subtracting $M$ reduces the stated mild problem to \eqref{eq:defect-wedge-problem}. Proposition~\ref{prop:wedge-wellposedness} gives existence, uniqueness, and the weighted estimate for $D=\Uwedge-M$.
\end{proof}

\begin{remark}
\label{rem:exponential-data-algebraic-wedge-decay}
Suppose $\ue^{c\scaledbulk}F\in L^\infty(K_\omega\times\Sone)$ and $\ue^{c\sigma}h_i\in L^\infty(\gamma_{i,-})$ for some $c>0$. Splitting the characteristic integrals as in Lemma~\ref{lem:H-decay} gives $\ue^{c'\scaledbulk}H\in L^\infty(K_\omega)$ for some $c'>0$. Proposition~\ref{prop:wedge-wellposedness} then applies for every $1<\betaw<\frac\pi\omega$. To compare this range with diffusion, note that the full-plane kernel has Fourier transform
\begin{align}
 \widehat q(\xi)&=\bigl(1+\abs{\xi}^2\bigr)^{-\frac12},\qquad 1-\widehat q(\xi)=\frac12\abs{\xi}^2+O(\abs{\xi}^4)\quad\text{as }\xi\to0.
 \label{eq:wedge-diffusive-low-frequency-mode}
\end{align}
Hence $I-q*$ has the same quadratic behavior at $\xi=0$ as $-\frac12\Delta$. The first decaying harmonic mode with zero trace on both rays is $\scaledbulk^{-\frac\pi\omega}\sin\bigl(\frac{\pi\theta}{\omega}\bigr)$. Its radial exponent $\frac\pi\omega$ is the square root of the first angular Dirichlet eigenvalue. This explains the exponent in the barrier condition \eqref{eq:barrier-exponent-range}.

To prove a sharp decay rate for given data, one would need an asymptotic expansion of the restricted resolvent and its first nonzero coefficient. Special data can give faster decay. For example, every nonzero $D\in C_c^\infty(K_\omega\times\Sone)$ is the unique solution with source $F=\Kw D$ and zero in-flow boundary data.

\end{remark}
\section{Wedge Matching and Construction of the Approximate Solution}
\label{sec:matching-realization}

We construct the wedge layer from the interior solution $\Uint_k$ and the side layer $\Uside_{j,k}$ in Section~\ref{sec:interior-side-hierarchy-comp}. For each matching pair $(\alpha,b)$ with $\alpha<\Tmatch$, we define $M_{m,\alpha,b}$ and apply Proposition~\ref{prop:wedge-wellposedness} to obtain $D_{m,\alpha,b}$. The sum $M_{m,\alpha,b}+D_{m,\alpha,b}$ then solves the homogeneous wedge equation with in-flow boundary data $G_{m,i,\alpha,b}$. We use the same $M_{m,\alpha,b}$ in \eqref{eq:Ubase-finite-expanded} to obtain exact cancellation in the approximate solution.

We first rescale the finite vertex expansions and collect equal powers of $\e$ and $\ln\e$. The grouping retains all contributions at the same combined exponent, so the consecutive transport and collision terms can cancel. We extend the finite expansion of each side layer from its ray to the wedge and multiply the resulting sum by the vertex cutoff $\chi_{\rm vtx}$. This gives $M_{m,\alpha,b}$ with finite $L^\infty$ norm near the vertex. The extension preserves the trace on its ray, and its cutoff derivatives are supported where Milne decay controls the resulting source. The wedge theorem then applies to finitely many problems independent of $\e$. After estimating their source and in-flow boundary data, we construct $D_{m,\alpha,b}$ and add the finite sum of wedge layers to $\Ubase^{\e,[\Nexp]}$.

\subsection{Rescaling and Truncation Exponents}
\label{subsec:physical-exponent-comp}

We use the coordinates \eqref{eq:signed-endpoint-coordinate}--\eqref{eq:bulk-signed-side-coordinate-relation} and \eqref{eq:wedge-side-coordinates}--\eqref{eq:corner-scaling-full}. In an abstract sector the radial variables are $\bulkdist$, $\scaledbulk$, with indices $\bulkdistat{m}$, $\scaledbulkat{m}$ near $V_m$.

At order $k$, substitute $\bulkdistat{m}=\e\scaledbulkat{m}$ in the vertex expansion of $\Uint_k$ and $\sidedistat{m}{i}=\e\sigma_i$, $d_{m,i}=\e\eta_i$ in that of $\Uside_{j,k}$. Writing $b_0$ for the logarithmic degree before rescaling, we obtain
\begin{align}
 \e^k\bulkdistat{m}^\lambda\bigl(\ln \bulkdistat{m}\bigr)^{b_0} &= \e^{k+\lambda}\scaledbulkat{m}^\lambda\sum_{c=0}^{b_0}\binom{b_0}{c}(\ln\e)^{b_0-c}\bigl(\ln\scaledbulkat{m}\bigr)^c,
 \label{eq:block-binomial-rule}\\
 \e^k\sidedistat{m}{i}^\lambda\bigl(\ln\sidedistat{m}{i}\bigr)^{b_0} &= \e^{k+\lambda}\sigma_i^\lambda\sum_{c=0}^{b_0}\binom{b_0}{c}(\ln\e)^{b_0-c}\bigl(\ln\sigma_i\bigr)^c.
 \label{eq:side-block-binomial-rule}
\end{align}
Write $r=\bulkdistat{m}$ and $s=\sidedistat{m}{i}$. For a summand $r^\lambda(\ln r)^{b_0}A(\theta_m,w)$ or $s^\lambda(\ln s)^{b_0}F(\eta_i,w)$, the angular factor becomes $A(\theta_m,O_mv)$ and the Milne factor becomes $F(\eta_i,O_mv)$ under the velocity pullback $w=O_mv$. Multiplying these factors by \eqref{eq:block-binomial-rule} or \eqref{eq:side-block-binomial-rule} gives terms with combined exponent $\alpha=k+\lambda$ and spatial degree $\alpha-k$.

The coefficients on the wedge use the velocity pullback $w=O_mv$, explicit in \eqref{eq:primitive-matching-coefficient-comp}. Under this pullback, the components $\tau_{m,i},\mu_{m,i}$ in \eqref{eq:local-ray-tangent-sign} become $\tau_i^Y,\mu_i^Y$ in \eqref{eq:wedge-frame-velocity-components}.

The recurrences preserve $\alpha\ge0$ by \eqref{eq:recurrence-degree-invariance-comp} and \eqref{eq:generated-exponents-nonnegative}, while repeated differentiation can make the spatial degree negative. When $\bulkdistat{m}\simeq\e$, every term generated with the same $\alpha$ contributes at the same power of $\e$. Matching therefore includes all such terms through order $\Nexp$. Their consecutive transport and collision terms cancel, leaving the decaying source needed for the wedge problem.

For example, the zero-trace mode $Z^+_{m,1}=\bulkdistat{m}^{\frac{3}{2}}\Phi_{m,1}(\theta_m)$ at an opening $\omega_m=\frac{2\pi}{3}$ has $Z^+_{m,1},\nabla Z^+_{m,1}\in L^\infty$ on a fixed vertex chart. The terms obtained from its successive derivatives in the Hilbert expansion have size $\e^\ell\bulkdistat{m}^{\frac{3}{2}-\ell}$, including a negative spatial power from $\ell=2$ onward. Successive derivative terms differ by $\frac{\e}{\bulkdistat{m}}$, so their ordering is asymptotic for $\bulkdistat{m}\gg\e$ and becomes nonuniform when $\bulkdistat{m}\simeq\e$. All have combined exponent $\alpha=\frac{3}{2}$ and contribute to the coefficient of $\e^{\frac{3}{2}}$.

We use $\Tcon$ to truncate the vertex expansions of $\Uint_k$ and $\Uside_{j,k}$, and $\Tmatch$ to determine which powers receive a wedge layer. The first truncation retains $\alpha<\Tcon$ and leaves a remainder of spatial degree $\Tcon-k$. This degree must remain positive after differentiating the order-$\Nexp$ remainder. We add wedge layers only for $\alpha<\Tmatch$. This second truncation determines the error in the in-flow boundary condition near each vertex. Keeping the two truncations separate allows us to construct every harmonic function $\rho_k$ with the positive Mellin weight $\Tcon-k$, and then restrict its finite expansion to the terms required for matching. In particular, matching to exponent $\Tmatch$ still uses contributions from all expansion orders $0\le k\le\Nexp$.

The choices in \eqref{eq:depth-choices-comp} give the inequalities
\begin{align}
 \Tcon-\Nexp&>\pacc+1,\qquad \Tcon>\Tmatch>\pacc+3,\qquad \Tmatch,\Tcon\notin\mathscr E_\Nexp,
 \label{eq:two-depth-separation-comp}
\end{align}
so $\Tcon$ is an admissible truncation exponent (Definition~\ref{def:admissible-hierarchy-depth-comp}). We construct $\rho_k$, $\Uint_k$ and $\Uside_{j,k}$ once, with $\mathsf d=\Tcon$ in Proposition~\ref{prop:data-generated-hierarchy-comp}.

We use three index sets. The construction set $\mathscr I_m^{\rm con}$ consists of the indices with $\alpha<\Tcon$. Its matching restriction $\mathscr I_m^{\rm mat}$ retains the entries whose combined $\e$-exponent is below $\Tmatch$. The intermediate set $\mathscr I_m^{\rm mid}$ collects the entries between the two truncation exponents:
\begin{align}
 \mathscr I_m^{\rm con} &:=\mathscr I_m(\Tcon),
 \label{eq:construction-ledger-comp}\\
 \mathscr I_m^{\rm mat} &:=\big\{\mathfrak l\in\mathscr I_m^{\rm con}: \alpha(\mathfrak l)<\Tmatch\big\},
 \label{eq:two-ledgers-comp}\\
 \mathscr I_m^{\rm mid} &:=\big\{\mathfrak l\in\mathscr I_m^{\rm con}: \Tmatch\le\alpha(\mathfrak l)<\Tcon\big\}.
 \label{eq:intermediate-ledger-comp}
\end{align}
The matching index set is closed under the operations in \eqref{eq:recurrence-degree-invariance-comp} because they preserve $\alpha$. Its coefficients inherit \eqref{eq:generic-hierarchy-bound-comp}. Binomial regrouping sends these index tuples to the pairs $(\alpha,b)$ defined below. Intermediate terms, $\Tmatch\le\alpha<\Tcon$, remain in $\Ubase^{\e,[\Nexp]}$ without a wedge layer. On each fixed region $\scaledbulk\le R_0$, their regularized contributions $\e^\alpha(\ln\e)^bM^{\rm con}_{m,\alpha,b}$ to $\Ubase^{\e,[\Nexp]}$ are $O(\e^\alpha\Lameps)$ in $L^\infty_{Y,v}$. Section~\ref{sec:main-proof} estimates their growth across the full chart in $x$ after summing the identities \eqref{eq:completed-string-coefficient-recurrence}--\eqref{eq:completed-side-coefficient-recurrence}.

All harmonic functions $\rho_k$ are constructed with positive Mellin weights $\Tcon-k$, including the orders $k>\Tmatch$. Matching restricts their finite expansions algebraically to $\alpha<\Tmatch$ and retains all terms with that combined exponent through order $\Nexp$.

Each term retained in $\mathscr I_m^{\rm mid}$ has $\alpha>\Tmatch$. At a fixed expansion order $k$, its spatial power $\bulkdistat{m}^{\alpha-k}$ or $\sidedistat{m}{i}^{\alpha-k}$ gains the factor $\bulkdistat{m}^{\alpha-\Tmatch}$ or $\sidedistat{m}{i}^{\alpha-\Tmatch}$ relative to a remainder of degree $\Tmatch-k$. Finiteness and $\Tmatch\notin\mathscr E_\Nexp$ give a fixed positive gap that absorbs its logarithms. This gap will also control the intermediate contribution at order $\delta^{\Tmatch}$ in the endpoint trace estimate.

The following bounds on the order-$\Nexp$ remainders supply the tangential and spatial derivatives used in Section~\ref{sec:main-proof}.

\begin{corollary}
\label{cor:terminal-derivative-bounds}
Let $\mathsf d=\Tcon$ in Proposition~\ref{prop:data-generated-hierarchy-comp}, so that $\mathsf d-\Nexp>\pacc+1$ as required by Definition~\ref{def:admissible-hierarchy-depth-comp}. Then, on every fixed vertex chart,
\begin{align}
 \nm{\ue^{\kappa_\ast\eta}\partial_{\sidedistat{m}{i}}^q \mathcal R^{B,\Tcon}_{m,i,\Nexp}}_{L^\infty_{\eta,w}} &\le C\mathcal G\,\sidedistat{m}{i}^{\pacc-q}, \qquad 0\le q\le \pacc,
 \label{eq:terminal-regular-remainder-comp}\\
 \bnm{D_x^q\mathcal R^{U,\Tcon}_{m,\Nexp}}_{L^\infty_w} &\le C\mathcal G\,\bulkdistat{m}^{\pacc-q}, \qquad 0\le q\le \pacc.
 \label{eq:terminal-interior-remainder-comp}
\end{align}
\end{corollary}

\begin{proof}
Apply \eqref{eq:outer-cartesian-remainder-comp} and \eqref{eq:side-ledger-remainder-comp} at $\mathsf d=\Tcon$ and $k=\Nexp$. For $0\le q\le \pacc$ they give $\bulkdistat{m}^{\Tcon-\Nexp-q}\bigl(1+\abs{\ln \bulkdistat{m}}\bigr)^{\Llog}$ and $\sidedistat{m}{i}^{\Tcon-\Nexp-q}\bigl(1+\abs{\ln \sidedistat{m}{i}}\bigr)^{\Llog}$. Since $\Tcon-\Nexp-q>(\pacc-q)+1$, the additional positive power absorbs the logarithm on the fixed chart.
\end{proof}

From here on, $C^U_{m,k}$, $\mathcal R^U_{m,k}$, $C^\Uside_{m,i,k}$ and $\mathcal R^\Uside_{m,i,k}$ denote the finite sums and remainders with truncation exponent $\Tcon$, as defined in \eqref{eq:exact-construction-projection-comp}--\eqref{eq:exact-side-construction-projection}. In particular, $\mathcal R^U_{m,k}=\mathcal R^{U,\Tcon}_{m,k}$ and $\mathcal R^\Uside_{m,i,k}=\mathcal R^{B,\Tcon}_{m,i,k}$. Thus the side remainder in \eqref{eq:terminal-regular-remainder-comp} is the same function that appears below as $\mathcal R^\Uside_{m,i,\Nexp}$.

\subsection{Vertex-Scale Regularization}
\label{subsec:fixed-local-continuations-v6}

For fixed $\alpha$, differentiation in \eqref{eq:microscopic-recursion-comp}--\eqref{eq:side-profile-inflow-decay-comp} makes $\lambda=\alpha-k$ negative at sufficiently large $k$, as in the example in Subsection~\ref{subsec:physical-exponent-comp}. The function $\Cvertex$ in \eqref{eq:C-coefficient-expanded} is therefore a finite sum on the punctured wedge, with possible singularities at $\scaledbulk=0$. To use it in Proposition~\ref{prop:wedge-wellposedness}, we choose $M$ so that the volume source $F=-\Kw M$ and trace mismatch $h_i=G_i-M\big|_{\gamma_{i,-}}$ satisfy the global weighted bounds.

We modify $\Cvertex$ on $\scaledbulk<2$, corresponding to $\bulkdist<2\e$, so the wedge problem is independent of $\e$. The resulting local errors belong to $L^\infty_{Y,v}$. On $\scaledbulk\ge2$, the identities \eqref{eq:completed-string-coefficient-recurrence}--\eqref{eq:completed-side-coefficient-recurrence} and Milne decay control the source and the ray-extension commutator.

We multiply the finite vertex expansion of $\Uint_k$ directly by $\chi_{\rm vtx}$. For $\Uside_{j,k}$, we first extend its expansion from the corresponding ray to the wedge. As described in Subsection~\ref{subsec:exact-wedge-blowup}, its tangential coordinate $\sigma_i$ is signed away from the ray. At an obtuse opening, the orthogonal foot of a point near one ray can lie behind the vertex on the other ray's line. The powers and logarithms in the corresponding side layer are defined for a positive tangential argument, so the extension is set to zero for $\sigma_i\le0$.

On $\sigma_i>0$, the ray cutoff below confines the side layer to a cone adjacent to its own ray. It preserves the trace at $\eta_i=0$ and makes the two incident extensions compatible with the wedge geometry, playing the localization role of the cutoff in the disk construction of \cite{Wu.Guo2015}.

Choose $\cfs>0$ so small that $2\cfs<\tan\omega_m$ at every acute opening, and choose a smooth $\Theta\colon[0,\infty)\to[0,1]$ such that
\begin{align}
 \Theta(z)&=\begin{cases}
 1,&0\le z\le\cfs,\\
 0,&z\ge2\cfs,
 \end{cases}
 \label{eq:Theta-choice}
\end{align}
For a function $\mathsf F\colon(0,\infty)_\sigma\times[0,\infty)_\eta\times\Sone\to\mathbb R$, define its one-sided extension by
\begin{align}
 \mathfrak E_i\mathsf F(\sigma_i,\eta_i,v)&:= \begin{cases}
 \Theta\big(\frac{\eta_i}{\sigma_i}\big)\mathsf F(\sigma_i,\eta_i,v),&\sigma_i>0,\\
 0,&\sigma_i\le0.
 \end{cases}
 \label{eq:forward-sector-extension}
\end{align}

As $\Theta(0)=1$, the trace on the in-flow boundary on ray $i$ is preserved. The support cone has angular width at most $\arctan(2\cfs)<\omega_m$, so this extension vanishes in a neighborhood of the opposite ray away from the vertex, including at a right opening. Derivatives of $\Theta\big(\frac{\eta_i}{\sigma_i}\big)$ are supported in $\cfs\sigma_i\le\eta_i\le2\cfs\sigma_i$, where $\eta_i\simeq\sigma_i\simeq\scaledbulk$ and $\abs{\nabla_{\sigma_i,\eta_i}\Theta}\lesssim\scaledbulk^{-1}$.

Hence exponential Milne decay in $\eta_i$ becomes exponential decay in $\scaledbulk$: on $\supp\Theta'\cap\{\scaledbulk\ge1\}$, for every power-logarithmic $P$, every exponentially decaying Milne factor $F_{\rm sl}$, every multi-index $a$ and every $L>0$,
\begin{align}
 \nm{\partial_Y^a\Theta\big(\tfrac{\eta_i}{\sigma_i}\big)P(\sigma_i)F_{\rm sl}(\eta_i,\cdot)}_{L^\infty_v} \le C_{a,L}\jbr{\scaledbulk}^{-L}.
 \label{eq:forward-commutator-rapid}
\end{align}
Here and below, coefficient estimates at fixed $Y$ or $\sigma$ hold for almost every spatial point or ray coordinate, and velocity norms are essential suprema. Every use of \eqref{eq:forward-commutator-rapid} below includes the vertex cutoff, which vanishes on $\scaledbulk\le1$.

We regularize the vertex with a cutoff $\chi_{\rm vtx}\in C^\infty([0,\infty);[0,1])$, independent of $\e$, satisfying
\begin{align}
 \chi_{\rm vtx}(\scaledbulk)&=\begin{cases}
 0,&0\le\scaledbulk\le1,\\
 1,&\scaledbulk\ge2.
 \end{cases}
 \label{eq:vertex-cutoff-choice}
\end{align}
The vertex-scale regularization of a function $C(Y,v)$ on the punctured wedge is the product $\chi_{\rm vtx}(\scaledbulk)\,C(Y,v)$, interpreted as zero on $\scaledbulk\le1$ before a logarithm or negative power is evaluated. Its $L^\infty_{Y,v}$ norm on $\scaledbulk\le2$ is at most the $L^\infty_{Y,v}$ norm of $C$ on $1\le\scaledbulk\le2$, and it equals $C$ for $\scaledbulk\ge2$.

The extension \eqref{eq:forward-sector-extension} vanishes identically near $\{\sigma_i=0,\ \eta_i>0\}$, since $\frac{\eta_i}{\sigma_i}\to\infty$ there. It is therefore smooth across that interface. Multiplication by $\chi_{\rm vtx}$ removes the vertex singularity.

Lemma~\ref{lem:localized-wedge-data} combines the cutoff bounds with \eqref{eq:completed-string-coefficient-recurrence}--\eqref{eq:completed-side-coefficient-recurrence} to verify the weighted hypotheses for both the volume source and the in-flow boundary discrepancy.

\subsection{Construction of Wedge Data}
\label{subsec:coefficient-data-localization-v6}

We now define the source and in-flow boundary data for the wedge equation. For each $(\alpha,b)$, we collect every term multiplying $\e^\alpha(\ln\e)^b$ after rescaling. All logarithmic powers produced by the binomial formula are included. The resulting functions satisfy the identities in Step~1 of Lemma~\ref{lem:localized-wedge-data} and are independent of $\e$.

The coefficient pairs produced by the binomial formulas are obtained from the index set \eqref{eq:construction-ledger-comp} by
\begin{align}
 \mathfrak J_m^{\rm con}&:=\bigl\{(k+\lambda,b):(m,k,\lambda,b_0;\mathsf t)\in\mathscr I_m^{\rm con},\ 0\le b\le b_0\bigr\}.
 \label{eq:construction-coefficient-pairs}
\end{align}
This finite set is partitioned into
\begin{align}
 \mathfrak J_m^{\rm mat} &:=\big\{(\alpha,b)\in\mathfrak J_m^{\rm con}:\alpha<\Tmatch\big\},
 \label{eq:finite-index-split-expanded}\\
 \mathfrak J_m^{\rm mid} &:=\big\{(\alpha,b)\in\mathfrak J_m^{\rm con}: \Tmatch\le\alpha<\Tcon\big\}.
 \label{eq:intermediate-coefficient-pairs}
\end{align}

For each $(\alpha,b)\in\mathfrak J_m^{\rm con}$, collect the terms multiplying $\e^\alpha(\ln\e)^b$ in the vertex expansions of $\sum_{k=0}^\Nexp\e^k\Uint_k$ and $\sum_{k=0}^\Nexp\e^k\Uside_{j,k}$. Let $\mathcal I^I_{m,\alpha,b}$ and $\mathcal I^B_{m,i,\alpha,b}$ index the resulting sums from $\Uint_k$ and the side layer on ray $i$. For an index $q$ coming from expansion order $k_q$ and original logarithmic degree $b_{0,q}$, put $\lambda_q=\alpha-k_q$ and $c_q=b_{0,q}-b$. The factors $A_q$ and $\mathsf F_q$ include the binomial coefficient $\binom{b_{0,q}}b$ and the velocity pullback. Thus $b+c_q=b_{0,q}\le\Llog$. With $\scaledbulk=\abs Y$ and $\theta$ the polar angle in $K_{\omega_m}$, define
\begin{align}
 C^I_{m,\alpha,b}(Y,v) &:={} \sum_{q\in\mathcal I^I_{m,\alpha,b}} \scaledbulk^{\lambda_q} \bigl(\ln \scaledbulk\bigr)^{c_q}A_q(\theta,v),
 \label{eq:raw-coefficient-blocks-v6}\\
 C^B_{m,i,\alpha,b}(\sigma_i,\eta_i,v) &:={} \sum_{q\in\mathcal I^B_{m,i,\alpha,b}} \sigma_i^{\lambda_q}\bigl(\ln\sigma_i\bigr)^{c_q}\mathsf F_q(\eta_i,v),
 \label{eq:raw-side-coefficient-sum}
\end{align}
and define $\Cvertex_{m,\alpha,b}$ by the following sum, with each side term extended by \eqref{eq:forward-sector-extension} and an empty sum interpreted as zero:
\begin{align}
 \Cvertex_{m,\alpha,b}(Y,v) &:=C^I_{m,\alpha,b}(Y,v) +\sum_{i=1}^2\mathfrak E_i \bigl[C^B_{m,i,\alpha,b}\bigr](Y,v).
 \label{eq:C-coefficient-expanded}
\end{align}

Regularize each construction coefficient by
\begin{align}
 M^{\rm con}_{m,\alpha,b}(Y,v) &:=\chi_{\rm vtx}(\scaledbulk)\,\Cvertex_{m,\alpha,b}(Y,v).
 \label{eq:Mcon-coefficient-expanded}
\end{align}

The restriction of $M^{\rm con}_{m,\alpha,b}$ to $\scaledbulk\le2$ belongs to $L^\infty_{Y,v}$. Indeed, on $1\le\scaledbulk\le2$ the side cutoff forces $\eta_i\le2\cfs\sigma_i$ and hence $\sigma_i\ge\frac{\scaledbulk}{\sqrt{1+4\cfs^2}}$. Each retained power and logarithm has absolute value at most a fixed constant on this annulus, and the Milne factors belong to $L^\infty_{\eta,v}$. The zero convention removes the remaining vertex singularity.

For each construction pair, $M^{\rm con}_{m,\alpha,b}$ has the tangential regularity in Proposition~\ref{prop:data-generated-hierarchy-comp} and is locally absolutely continuous along almost every nongrazing characteristic. The combination $\mu_i^Y\partial_{\eta_i}+\qk$ on each Milne factor is evaluated using its equation. The function $M^{\rm con}_{m,\alpha,b}$ is independent of $\e$, equals $\Cvertex_{m,\alpha,b}$ for $\scaledbulk\ge2$, and uses the same ray cutoff as $\Ubase^{\e,[\Nexp]}$.

For a matching index, set $M_{m,\alpha,b}:=M^{\rm con}_{m,\alpha,b}$, $(\alpha,b)\in\mathfrak J_m^{\rm mat}$. The definition of $\Ubase^{\e,[\Nexp]}$ uses this same $M_{m,\alpha,b}$. Solving the wedge problem with source $-\mathscr K_{\omega_m}M_{m,\alpha,b}$ will therefore cancel the matched source, including its regularization errors, wherever the overlap cutoff equals one.

For $0\le k\le1$ set the Taylor polynomial degree $q_k^{\Tmatch}:=\max\{n\in\mathbb N_0:k+n<\Tmatch\}$ and the endpoint Taylor polynomial
\begin{align}
 J_{m,i,k}^{\Tmatch}[\gaux](\sidedistat{m}{i},w)&:= \sum_{n=0}^{q_k^{\Tmatch}} \frac{\sidedistat{m}{i}^{n}}{n!}\, \partial_{\sidedistat{m}{i}}^n\gaux_{m,i,k}(0,w).
 \label{eq:taylor-jet-def}
\end{align}
The maximum exists and $k+q_k^{\Tmatch}<\Tmatch<k+q_k^{\Tmatch}+1$, since $\Tmatch\notin\mathbb N_0$. This polynomial retains precisely the Taylor terms with $k+n<\Tmatch$.

At vertex $m$, write $\Kw=\mathscr K_{\omega_m}$ for the wedge operator of Section~\ref{sec:wedge-layer}. For $(\alpha,b)\in\mathfrak J_m^{\rm mat}$, the Taylor coefficient of the in-flow data on ray $i$ is
\begin{align}
 G_{m,i,\alpha,b}(\sigma,v) := \begin{cases}
 \sum_{0\le k\le1,\ n\in\mathbb N_0,\ k+n=\alpha} \frac{\sigma^n}{n!}\, \partial_{\sidedistat{m}{i}}^n\gaux_{m,i,k}(0,O_mv), &b=0,\\
 0,&b\ge1,
 \end{cases}
 \label{eq:primitive-matching-coefficient-comp}
\end{align}
The data for the wedge layer are
\begin{align}
 F_{m,\alpha,b}&:=-\Kw M_{m,\alpha,b},\qquad h_{m,i,\alpha,b}:=G_{m,i,\alpha,b}-M_{m,\alpha,b}\big|_{\gamma_{i,-}}.
 \label{eq:generated-wedge-data}
\end{align}

Only $0\le k\le1$ occurs in \eqref{eq:primitive-matching-coefficient-comp}, since \eqref{eq:auxiliary-inflow-coefficients} sets the auxiliary in-flow data to zero for $2\le k\le \Nexp$.

It remains to estimate $F=-\Kw M$ and $h_i=G_i-M\big|_{\gamma_{i,-}}$. Although $G_i$ and $M\big|_{\gamma_{i,-}}$ may grow along the ray, their difference is supported near the vertex. We prove this by using the in-flow boundary condition for $\Uint_k+\Uside_{j,k}$ at each order. We also use the cancellation of consecutive transport and collision terms to estimate $\Kw M$.

\begin{lemma}
\label{lem:localized-wedge-data}
Fix $(\alpha,b)\in\mathfrak J_m^{\rm mat}$ and retain all terms through order $\Nexp$ in \eqref{eq:raw-coefficient-blocks-v6}--\eqref{eq:raw-side-coefficient-sum}. Set $\sdata=\betaw+2$ as in \eqref{eq:N-choice-comp}. For the regularization \eqref{eq:Mcon-coefficient-expanded}, the data \eqref{eq:generated-wedge-data} have the following properties.
\begin{enumerate}[label=\textup{(\roman*)}]
\item The source is well defined in the characteristic and distributional senses and has a measurable representative satisfying
\begin{align}
 \nm{F_{m,\alpha,b}(Y,\cdot)}_{L^\infty_v} &\le C\jbr{\scaledbulk}^{-\sdata}.
 \label{eq:localized-data-decay}
\end{align}
\item The in-flow boundary data have measurable representatives supported in $0\le\sigma\le2$ and satisfy
\begin{align}
 \nm{h_{m,i,\alpha,b}(\sigma,\cdot)}_{L^\infty(\{\mu_i^Y>0\})} &\le C\jbr{\sigma}^{-\sdata}, \qquad i=1,2.
 \label{eq:localized-incoming-decay}
\end{align}
\end{enumerate}
Here $Y=\sigma\twedge{m}{i}$ and $\abs Y=\sigma$ on ray $i$. The constants depend continuously on finitely many norms in Assumption~\ref{ass:primitive-data-comp}. Thus the data satisfy \eqref{eq:abstract-wedge-data-assumption-v6} with $s=\sdata$.
\end{lemma}

\begin{proof}
The source estimate follows by canceling consecutive transport and collision terms, leaving the order-$\Nexp$ derivatives and two cutoff commutators. For the boundary data, we use coefficient uniqueness in the in-flow boundary condition at each fixed order, then regroup by $\alpha=k+\lambda$.

Fix $(\alpha,b)\in\mathfrak J_m^{\rm mat}$ and suppress $m$. Let $C_k^I$ and $C_{i,k}^B$ be the order-$k$ terms in \eqref{eq:raw-coefficient-blocks-v6}--\eqref{eq:raw-side-coefficient-sum} after binomial regrouping, with an absent term read as zero. Thus $C^I_{\alpha,b}=\sum_{k=0}^\Nexp C_k^I$ and $C^B_{i,\alpha,b}=\sum_{k=0}^\Nexp C_{i,k}^B$.

\paragraph{\underline{Step 1: coefficient recurrences before regularization}} For $0\le k\le\Nexp$, we first show that
\begin{align}
 \qk C_k^I+v\cdot\nabla_YC_{k-1}^I&=0,
 \label{eq:completed-string-coefficient-recurrence}\\
 \big(\mu_i^Y\partial_{\eta_i}+\qk\big)C_{i,k}^B +\tau_i^Y\partial_{\sigma_i}C_{i,k-1}^B&=0,
 \label{eq:completed-side-coefficient-recurrence}
\end{align}
with index $-1$ equal to zero. These are the coefficient forms of the identities at each expansion order $k$ \eqref{eq:exact-interior-displayed-recurrence-comp} and \eqref{eq:exact-side-displayed-recurrence-comp} for the displayed construction parts. By \eqref{eq:block-binomial-rule}--\eqref{eq:side-block-binomial-rule}, every rescaled term carries the factor $\e^\alpha(\ln\e)^\ell$ with $\alpha=k+\lambda$. The collision operator preserves $(k,\alpha-k)$. The transport derivative pairs an order-$(k-1)$ term of degree $\alpha-k+1$ with an order-$k$ term of degree $\alpha-k$. Both preserve the combined $\e$-exponent $\alpha$.

Multiply each order-$k$ identity by $\e^k$ and use $w\cdot\nabla_x=\e^{-1}v\cdot\nabla_Y$ and $\partial_{\sidedistat{m}{i}}=\e^{-1}\partial_{\sigma_i}$. By the binomial formulas, the coefficient of $\e^\alpha(\ln\e)^b$ is the corresponding left-hand side of \eqref{eq:completed-string-coefficient-recurrence} or \eqref{eq:completed-side-coefficient-recurrence}. The finite sum vanishes for every small $\e$. Linear independence of the functions $\e^\alpha(\ln\e)^b$ gives the two identities coefficientwise.

The extraction uses only $\alpha<\Tcon$, so the identities \eqref{eq:completed-string-coefficient-recurrence}--\eqref{eq:completed-side-coefficient-recurrence} hold for every regrouped pair $(\alpha,b)\in\mathfrak J_m^{\rm con}$, matched or intermediate.

\paragraph{\underline{Step 2: the regularized source}} By \eqref{eq:Mcon-coefficient-expanded}, $M_{\alpha,b}=\chi_{\rm vtx}\Cvertex_{\alpha,b}$. Apply $\Kw$ to \eqref{eq:C-coefficient-expanded} on the punctured wedge using the characteristic or distributional product rule. Summing \eqref{eq:completed-string-coefficient-recurrence}--\eqref{eq:completed-side-coefficient-recurrence} over $k$ cancels consecutive transport and collision terms. The remaining terms are
\begin{align}
 \Kw\Cvertex_{\alpha,b} ={}&v\cdot\nabla_YC_\Nexp^I+\sum_{i=1}^2\mathfrak E_i\bigl[\tau_i^Y\partial_{\sigma_i}C_{i,\Nexp}^B\bigr]
 +\sum_{i=1}^2\sum_{k=0}^\Nexp\bigl(v\cdot\nabla_Y\Theta\big(\tfrac{\eta_i}{\sigma_i}\big)\bigr)C_{i,k}^B.
 \label{eq:raw-coefficient-wedge-operator}
\end{align}
The last sum is extended by zero for $\sigma_i\le0$, as in \eqref{eq:forward-sector-extension}.

Multiplication by $\chi_{\rm vtx}$ creates one additional radial commutator. We separate the resulting four terms because the truncation order $\Nexp$ determines the spatial degree of the first two derivatives, while the cutoff supports control the last two terms:
\begin{align}
 \Kw M_{\alpha,b} ={}&\mathcal E^{I,\mathrm{loc}}_{\alpha,b} +\mathcal E^{B,\mathrm{loc}}_{\alpha,b} +\mathcal E^{F,\mathrm{loc}}_{\alpha,b} +\mathcal E^{\rm vtx,\mathrm{loc}}_{\alpha,b},
 \label{eq:completed-string-error-split}
\end{align}
The derivatives of $C_\Nexp^I$ and $C_{i,\Nexp}^B$ and the two cutoff commutators give
\begin{align}
 \mathcal E^{I,\mathrm{loc}}_{\alpha,b} &:=\chi_{\rm vtx}\,v\cdot\nabla_YC_\Nexp^I,
 \label{eq:local-interior-source}\\
 \mathcal E^{B,\mathrm{loc}}_{\alpha,b} &:=\chi_{\rm vtx}\sum_{i=1}^2 \mathfrak E_i\!\left[ \tau_i^Y\partial_{\sigma_i}C_{i,\Nexp}^B \right],
 \label{eq:local-side-source}\\
 \mathcal E^{F,\mathrm{loc}}_{\alpha,b} &:=\chi_{\rm vtx}\sum_{i=1}^2\sum_{k=0}^\Nexp \bigl(v\cdot\nabla_Y\Theta\big(\tfrac{\eta_i}{\sigma_i}\big)\bigr) C_{i,k}^B,
 \label{eq:local-ray-commutator}\\
 \mathcal E^{\rm vtx,\mathrm{loc}}_{\alpha,b} &:=\bigl(v\cdot\nabla_Y\chi_{\rm vtx}(\scaledbulk)\bigr)\Cvertex_{\alpha,b}.
 \label{eq:local-vertex-commutator}
\end{align}
The radial commutator is zero on $\scaledbulk\le1$ with its value assigned before the singular coefficient is evaluated. Since $\chi_{\rm vtx}'$ is supported in $[1,2]$,
\begin{align}
 \supp \mathcal E^{\rm vtx,\mathrm{loc}}_{\alpha,b} \subset\{1\le\scaledbulk\le2\}.
 \label{eq:localized-core-error-support}
\end{align}

\paragraph{\underline{Step 3: cutoff commutators}} The vertex commutator belongs to $L^\infty_{Y,v}$ and has support \eqref{eq:localized-core-error-support}. The ray commutator is controlled by \eqref{eq:forward-commutator-rapid}, since its support has $\eta_i\simeq\sigma_i\simeq\scaledbulk$. Thus, for every $L>0$,
\begin{align}
 \bnm{\mathcal E^{F,\mathrm{loc}}_{\alpha,b}(Y,\cdot)}_{L^\infty_v} \le C_L\jbr{\scaledbulk}^{-L}.
 \label{eq:completed-forward-error-bound}
\end{align}
The constants depend on the fixed interior solutions $\Uint_k$ and side layers $\Uside_{j,k}$ for $0\le k\le\Nexp$.

\paragraph{\underline{Step 4: derivatives at order $\Nexp$}} For a coefficient with combined exponent $\alpha$, each order-$\Nexp$ term has spatial degree $\alpha-\Nexp$. For $C_\Nexp^I$, Cartesian differentiation of the radial power, logarithm, and angular factor lowers the degree by one. The angular derivatives are controlled by Proposition~\ref{prop:data-generated-hierarchy-comp}. For $C_{i,\Nexp}^B$, the tangential derivative is taken at fixed $(\eta_i,v)$ and acts on the power and logarithm, while the ray cutoff gives $\sigma_i\simeq\scaledbulk$ on its support. Together with the Milne estimate, this gives
\begin{align}
 \bnm{\mathcal E^{I,\mathrm{loc}}_{\alpha,b}(Y,\cdot)}_{L^\infty_v} +\bnm{\mathcal E^{B,\mathrm{loc}}_{\alpha,b}(Y,\cdot)}_{L^\infty_v} \le C\scaledbulk^{\alpha-\Nexp-1}\bigl(1+\ln\scaledbulk\bigr)^{\Llog} \quad(\scaledbulk\ge3).
 \label{eq:completed-terminal-error-bound}
\end{align}
The parameter inequality $\Nexp>\Tmatch+\sdata$ in \eqref{eq:N-dominates-matching-decay}, together with $\alpha<\Tmatch$, gives $\alpha-\Nexp-1<-\sdata-1$. The positive gap absorbs the logarithm, since $\scaledbulk^{-1}\bigl(1+\ln\scaledbulk\bigr)^{\Llog}\le C$ for $\scaledbulk\ge3$. Hence the $L^\infty_v$ norms in \eqref{eq:completed-terminal-error-bound} are $O(\scaledbulk^{-\sdata})$.

\paragraph{\underline{Step 5: the volume source}} On $\scaledbulk\le3$, the regularized coefficients and the four terms in \eqref{eq:completed-string-error-split} belong to $L^\infty_{Y,v}$. For a side summand, the Milne equation evaluates $\mu_i^Y\partial_{\eta_i}+\qk$ as a whole and gives an $L^\infty_v$ estimate without a grazing singularity. Summing the compactly supported vertex term, \eqref{eq:completed-forward-error-bound}, and \eqref{eq:completed-terminal-error-bound} proves \eqref{eq:localized-data-decay} for $F_{\alpha,b}=-\Kw M_{\alpha,b}$.

\paragraph{\underline{Step 6: in-flow boundary condition}} On ray $i$, the cutoff in \eqref{eq:forward-sector-extension} equals one for the coefficient from that side, and the extension from the other side vanishes by the choice of $\cfs$. Hence the trace of $\Cvertex_{\alpha,b}$ is the sum of the traces of $C^I_{\alpha,b}$ and $C^B_{i,\alpha,b}$. All boundary identities below use the almost-everywhere characteristic traces. We determine their sum at each expansion order.

For $\sigma\ge2$, the vertex cutoff and the ray cutoff for side $i$ both equal one. Set $j=j(m,i)$. We determine the trace coefficients there from the exact identity at each fixed $0\le k\le \Nexp$,
\begin{align}
 \Uint_k\big|_{\gamma_{j,-}}\bigl(x_j(s_{m,i}^{\rm ve}(s)),O_mv\bigr) +\Uside_{m,i,k}(s,0,O_mv) &=\gaux_{m,i,k}(s,O_mv), \qquad \mu_i^Y(v)>0.
 \label{eq:own-side-order-trace}
\end{align}
Here $s=\sidedistat{m}{i}\in(0,\ellvtx)$ is distance from $V_m$, and $x_j(s_{m,i}^{\rm ve}(s))$ is the corresponding point of $E_j$. The side layer is evaluated at its characteristic trace on the in-flow boundary $\eta_i=0$. Insert the decompositions truncated at $\alpha<\Tcon$ in \eqref{eq:weighted-interior-decomposition-comp}--\eqref{eq:weighted-side-decomposition-comp} and the Taylor expansion of the prescribed in-flow coefficients from Step~2 of Proposition~\ref{prop:data-generated-hierarchy-comp}. At order $k$ this gives
\begin{align}
 \sum_{(\lambda,b)} s^{\lambda}(\ln s)^b\, \mathfrak h_{i,k,\lambda,b}(v) +\mathcal R^{\rm tel}_{i,k}(s,v) &=0,
 \label{eq:telescope-with-residual}
\end{align}
where $\mathfrak h_{i,k,\lambda,b}$ is the prescribed Taylor coefficient minus the coefficients of $s^\lambda(\ln s)^b$ in $\Uint_k$ and $\Uside_{m,i,k}$, with zero padding as in clause \textup{(L10)}. The function $\mathcal R^{\rm tel}_{i,k}$ is the Taylor remainder of $\gaux_{m,i,k}$ minus the two trace remainders of $\Uint_k$ and $\Uside_{m,i,k}$. For $k\ge2$, its Taylor-remainder contribution is zero. Every displayed exponent satisfies $\lambda<\Tcon-k$. The bounds for these three remainders give
\begin{align}
 \nm{\mathcal R^{\rm tel}_{i,k}(s,\cdot)}_{L^\infty(\{\mu_i^Y>0\})} &\le C s^{\Tcon-k}\bigl(1+\abs{\ln s}\bigr)^{\Llog},
 \label{eq:telescope-residual-bound}
\end{align}
by \eqref{eq:outer-ledger-remainder-comp}, \eqref{eq:side-ledger-remainder-comp}, and the Taylor remainder bound. The retained order-dependent weight $\Tcon-k$ lies strictly above every displayed spatial exponent and permits coefficient uniqueness at each order.

\paragraph{\underline{Step 7: coefficientwise ray cancellation}} Collect equal spatial powers in \eqref{eq:telescope-with-residual} into $L^\infty(\{\mu_i^Y>0\})$-valued polynomials in $\ln s$. Every such power is strictly below $\Tcon-k$, so Lemma~\ref{lem:finite-polyhom-uniqueness}, applied with $d_{\rm rem}=\Tcon-k$ and \eqref{eq:telescope-residual-bound}, gives $\mathfrak h_{i,k,\lambda,b}=0$ for every displayed index.

For $(\alpha,b)\in\mathfrak J_m^{\rm mat}$, every order-$k$ exponent satisfies $\lambda=\alpha-k<\Tmatch-k<\Tcon-k$ and is included in this cancellation. The coefficient of each $\sigma^\lambda\bigl(\ln\sigma\bigr)^{c-b}$ in the rescaled trace difference is $\binom cb\mathfrak h_{i,k,\lambda,c}=0$, with $k+\lambda=\alpha$ and $c\ge b$. Therefore
\begin{align}
 C^I_{\alpha,b}(\sigma\twedge{m}{i},v)+C^B_{i,\alpha,b}(\sigma,0,v)&=G_{i,\alpha,b}(\sigma,v),\qquad\sigma>0,
 \label{eq:coefficientwise-incoming-cancellation}
\end{align}
for in-flow velocities. Coefficient uniqueness was applied as $s\to0$ within the vertex chart. Once every coefficient vanishes, the resulting finite power-logarithmic identity holds for all $\sigma>0$, as stated in \eqref{eq:coefficientwise-incoming-cancellation}. The Taylor coefficients in \eqref{eq:primitive-matching-coefficient-comp} are precisely those retained in \eqref{eq:taylor-jet-def}.

Since $\alpha<\Tmatch<\Tcon$, every retained coefficient is included in this cancellation. Thus $\Tcon>\Tmatch$ suffices to match the in-flow boundary data. The inequality $\Tcon-\Nexp>\pacc+1$ is used separately to estimate derivatives of the construction remainders in Section~\ref{sec:main-proof}.

\paragraph{\underline{Step 8: in-flow boundary discrepancy}} Combining \eqref{eq:coefficientwise-incoming-cancellation} with the product definition \eqref{eq:Mcon-coefficient-expanded} gives
\begin{align}
 h_{i,\alpha,b}(\sigma,v)&=\bigl(1-\chi_{\rm vtx}(\sigma)\bigr)G_{i,\alpha,b}(\sigma,v),\qquad \supp h_{i,\alpha,b}\subset\{0\le\sigma\le2\}.
 \label{eq:completed-ray-mismatch-split}
\end{align}
The polynomial $G_{i,\alpha,b}$ belongs to $L^\infty((0,2)\times\{\mu_i^Y>0\})$, which proves \eqref{eq:localized-incoming-decay}. The finite coefficient bounds \eqref{eq:generic-hierarchy-bound-comp} at $\mathsf d=\Tcon$ control both the source and the trace data continuously in the prescribed in-flow coefficients $g_{j,0},g_{j,1}$. This completes the proof.
\end{proof}

The recurrences \eqref{eq:completed-string-coefficient-recurrence}--\eqref{eq:completed-side-coefficient-recurrence}, the exact split \eqref{eq:completed-string-error-split}, the commutator estimates in Steps~2--3, and the $L^\infty_{Y,v}$ estimate on $\scaledbulk\le3$ in Step~5 apply to every $(\alpha,b)\in\mathfrak J_m^{\rm con}$. So does the order-$\Nexp$ derivative estimate \eqref{eq:completed-terminal-error-bound}. The conclusion $\bnm{\mathscr K_{\omega_m}M^{\rm con}_{m,\alpha,b}(Y,\cdot)}_{L^\infty_v}\lesssim\jbr{\scaledbulk}^{-\sdata}$ additionally uses $\alpha<\Tmatch$. Section~\ref{sec:main-proof} uses the exact split and the degree $\alpha-\Nexp-1$ for intermediate pairs.

\subsection{Application of Wedge Theorem}
\label{subsec:apply-abstract-wedge-v6}

For every vertex, $1<\betaw<\lambdaMin\le\frac{\pi}{\omega_m}$ and $\sdata=\betaw+2$. Lemma~\ref{lem:localized-wedge-data} therefore verifies the source and in-flow boundary data bounds of Proposition~\ref{prop:wedge-wellposedness} with $s=\sdata$.

For $(\alpha,b)\in\mathfrak J_m^{\rm mat}$, let $D_{m,\alpha,b}$ be the unique $L^\infty$ mild solution of \eqref{eq:defect-wedge-problem} with data \eqref{eq:generated-wedge-data}, and set
\begin{align}
 \Uwedge_{m,\alpha,b} &:= M_{m,\alpha,b}+D_{m,\alpha,b}.
 \label{eq:coefficientwise-total-wedge-definition}
\end{align}
By \eqref{eq:generated-wedge-data}, $\Kw\Uwedge_{m,\alpha,b}=0$ in $K_{\omega_m}\times\Sone$ and $\Uwedge_{m,\alpha,b}\big|_{\gamma_{i,-}}=G_{m,i,\alpha,b}$, $i=1,2$. Corollary~\ref{cor:affine-total-wedge-field-v6} gives uniqueness among solutions $\Uwedge_{m,\alpha,b}$ with $\Uwedge_{m,\alpha,b}-M_{m,\alpha,b}\in L^\infty(K_{\omega_m}\times\Sone)$, with the mild interpretation specified in Subsection~\ref{subsec:M-W-D-precise}.

Recall the finite sums $\mathcal M_m^\e$ and $\mathcal D_m^\e$ from \eqref{eq:intro-summed-reference-field}--\eqref{eq:intro-summed-wedge-fields}:
\begin{align}
 \mathcal M_m^\e &= \sum_{(\alpha,b)\in\mathfrak J_m^{\rm mat}} \e^\alpha(\ln\e)^bM_{m,\alpha,b},
 \label{eq:Mless-finite-expanded}\\
 \mathcal D_m^\e &= \sum_{(\alpha,b)\in\mathfrak J_m^{\rm mat}} \e^\alpha(\ln\e)^bD_{m,\alpha,b},
 \label{eq:summed-wedge-corrector-definition}
\end{align}
Thus $\Uwedge_m^\e=\mathcal M_m^\e+\mathcal D_m^\e$ is the summed wedge solution,
\begin{align}
 \Uwedge_m^\e &=\sum_{(\alpha,b)\in\mathfrak J_m^{\rm mat}} \e^\alpha(\ln\e)^b\Uwedge_{m,\alpha,b}.
 \label{eq:summed-total-wedge-definition}
\end{align}

Summing the in-flow boundary conditions $\Uwedge_{m,\alpha,b}\big|_{\gamma_{i,-}}=G_{m,i,\alpha,b}$ gives
\begin{align}
 \Uwedge_m^\e\big|_{\gamma_{i,-}}(\sigma,v) &= \sum_{(\alpha,b)\in\mathfrak J_m^{\rm mat}} \e^\alpha(\ln\e)^b G_{m,i,\alpha,b}(\sigma,v)
 \label{eq:summed-wedge-incoming-trace}= \sum_{k=0}^1\e^k J_{m,i,k}^{\Tmatch}[\gaux](\e\sigma,O_mv),
\end{align}
the second equality in \eqref{eq:summed-wedge-incoming-trace} being the finite regrouping of Taylor coefficients of the prescribed in-flow data by $\alpha=k+n$.

The sum $\mathcal D_m^\e$ obeys the decay bound
\begin{align}
 \nm{\mathcal D_m^\e(Y,\cdot)}_{L^\infty_v} &\le C\Lameps\jbr{\scaledbulk}^{-\betaw},
 \label{eq:summed-wedge-decay}
\end{align}
Indeed, Proposition~\ref{prop:wedge-wellposedness} gives $\nm{D_{m,\alpha,b}(Y,\cdot)}_{L^\infty_v}\le C\jbr{\scaledbulk}^{-\betaw}$ for almost every $Y$. Every generated exponent has $\alpha\ge0$ by \eqref{eq:generated-exponents-nonnegative}, and every logarithmic degree has $b\le\Llog=\Nexp+1$ by \eqref{eq:delta-lambda-comp}. Summing the finite family proves \eqref{eq:summed-wedge-decay}.

The wedge layer satisfies $\jbr{\scaledbulk}^{\betaw}D_{m,\alpha,b}\in L^\infty_{Y,v}$. The total wedge solution $\Uwedge_{m,\alpha,b}$ retains any growth of $M_{m,\alpha,b}$. Intermediate coefficients stay regularized in $\Ubase^{\e,[\Nexp]}$ and have no wedge layer.

\subsection{Extension of Side Layers}
\label{subsec:global-side-realization}

We now extend the side layers to the polygon. Near each vertex, the extension agrees with $\mathfrak E_i$. Along the middle of each side, it equals the side layer with a cutoff in the normal direction. The extension preserves the in-flow boundary condition on that side and gives the same $M_{m,\alpha,b}$ as in the wedge problem.

\paragraph{\underline{Fixed length scales}} Let $\ellsep$ be the separation of a vertex from every nonincident supporting line, defined through the affine hull of the side by
\begin{align}
 \ellsep &:= \min\left\{ \dist(V_m,\operatorname{aff}(E_j)): V_m\notin\overline E_j \right\} >0.
 \label{eq:nonincident-side-separation}
\end{align}
Convexity and maximality of the sides give $\overline\Om\cap\operatorname{aff}(E_j)=\overline E_j$. Each nonincident vertex therefore has positive distance from this line, and the finite minimum in \eqref{eq:nonincident-side-separation} is positive. Choose $0<\ellnor<\frac{\ellsep}{2}$ and $\ellend>0$ so small that the intervals of radius $2\ellend$ about the two endpoints of every side are disjoint.

The fixed vertex-chart radius was chosen so that the balls $\{\bulkdistat{m}<4\ellvtx\}$ lie in the larger analytic vertex neighborhoods reserved before \eqref{eq:weighted-interior-decomposition-comp}, are disjoint, agree with their tangent sectors, and satisfy
\begin{align}
 4\ellvtx&<\min\{1,\ellnor,\ellend\},\qquad 2\ellnor+4\ellvtx<\ellsep.
 \label{eq:fixed-realization-separation}
\end{align}

All four lengths are fixed independently of $\e$. The endpoint scale separates the two endpoint transitions of each side, while the normal scale and \eqref{eq:fixed-realization-separation} keep nonincident sides out of the vertex balls.

\paragraph{\underline{Side layers and endpoint cutoffs}} Let $B_j=B_j(s,\eta,w)$, $0<s<L_j$, $\eta\ge0$, denote a side layer, a finite sum of side layers, a tangential derivative, or a construction remainder. The notation allows dependence on $\e$.

In the extensions below we set $(s,\eta)=\big(s_j(x),\frac{d_j(x)}{\e}\big)$. By \eqref{eq:global-affine-side-coordinate} only $s_j(x)$ is extended off the side, and $B_j$ itself is never evaluated outside $0\le s\le L_j$.

For an incident pair $(m,i)$ with $j=j(m,i)$, the signed endpoint coordinate is $\signedsideat{m}{i}=s_j$ or $L_j-s_j$ according to \eqref{eq:signed-endpoint-coordinate}. We use $d_{m,i}=d_j$, $w\cdot\nabla_x\signedsideat{m}{i}=\tau_{m,i}(w)$, and the identity for tangential derivatives \eqref{eq:tangential-operator-frame-invariance}.

For each endpoint $V_m\in\partial E_j$, let $i=i(m,j)$ and choose $\chi_{m,i}^{\rm end}\in C_c^\infty\big((-2\ellend,2\ellend);[0,1]\big)$ with $\chi_{m,i}^{\rm end}=1$ on $[-\ellend,\ellend]$. These cutoffs are fixed independently of $\e$ and $\delta$. Write $\mathcal I(j):=\{(m,i):j(m,i)=j\}$. Its two endpoint supports are disjoint, and $0<s_j(x)<L_j$ gives $\signedsideat{m}{i}(x)>0$ for both incident pairs.

In the endpoint chart of $V_m$, with $\sigma_i=\frac{\signedsideat{m}{i}}{\e}$ and $\eta_i=\frac{d_j}{\e}$,
\begin{align}
 \Theta\left(\frac{d_j}{\signedsideat{m}{i}}\right) &=\Theta\left(\frac{\eta_i}{\sigma_i}\right),
 \label{eq:physical-wedge-ratio-identity}
\end{align}
so the constructions in $(x,w)$ and $(Y,v)$ use the same scale-invariant ray cutoff.

On the inward half-strip $0<s_j(x)<L_j$, $d_j(x)\ge0$, set
\begin{align}
 a_j^{\rm end}(x)&:=1+\sum_{(m,i)\in\mathcal I(j)}\chi_{m,i}^{\rm end}(\signedsideat{m}{i})\left[\Theta\left(\frac{d_j}{\signedsideat{m}{i}}\right)-1\right].
 \label{eq:endpoint-assembly-factor}
\end{align}
At most one endpoint term is active. The factor equals $\Theta$ near that endpoint and equals one in the middle of the side. On $d_j=0$ it equals one everywhere, so it preserves the trace on the in-flow boundary.

\paragraph{\underline{Side extension and normal localization}} Define the extension along side $E_j$ by
\begin{align}
 (\mathfrak E_j^{\rm glob}B_j)(x,w):= \begin{cases}
 a_j^{\rm end}(x) B_j\left(s_j(x),\frac{d_j(x)}{\e},w\right), &0<s_j(x)<L_j,\ d_j(x)\ge0,\\
 0,&\text{otherwise},
 \end{cases}
 \label{eq:explicit-global-forward-glue}
\end{align}

In the commutator below, $\Le$ acting on $B_j(s,\eta,w)$ means $\tau_j\partial_s+\e^{-1}(\mu_j\partial_\eta+\qk)$. The derivative $\partial_s$ is taken at fixed $(\eta,w)$, before the substitution $(s,\eta)=(s_j(x),\frac{d_j(x)}{\e})$. Since $a_j^{\rm end}$ is velocity independent, the collision operator commutes with it. On the open half-strip we obtain
\begin{align}
 \label{eq:explicit-global-forward-commutator} [\Le,\mathfrak E_j^{\rm glob}]B_j ={}& \sum_{(m,i)\in\mathcal I(j)}\chi_{m,i}^{\rm end} \left(w\cdot\nabla_x \Theta\left(\frac{d_j}{\signedsideat{m}{i}}\right)\right) B_j\left(s_j,\frac{d_j}{\e},w\right) \\
 &\quad+ \sum_{(m,i)\in\mathcal I(j)} \tau_{m,i}(w)(\chi_{m,i}^{\rm end})'(\signedsideat{m}{i}) \left[ \Theta\left(\frac{d_j}{\signedsideat{m}{i}}\right)-1 \right] B_j\left(s_j,\frac{d_j}{\e},w\right).\notag
\end{align}

Finally, the normal cutoff $\zeta_j\in C_c^\infty([0,\infty);[0,1])$, applied last to keep the side chart inside the fixed collar $0\le d_j<2\ellnor$, satisfies
\begin{align}
 \zeta_j(d)&=\begin{cases}
 1,&0\le d\le\ellnor,\\
 0,&d\ge2\ellnor.
 \end{cases}
 \label{eq:normal-cutoff-definition}
\end{align}

For fixed $d_j>0$ and $\signedsideat{m}{i}>0$ small, $\chi_{m,i}^{\rm end}=1$ and $\Theta\big(\frac{d_j}{\signedsideat{m}{i}}\big)=0$. Thus the half-strip branch vanishes in a neighborhood of each endpoint interface, and its extension by zero preserves the distributional product rule. The single point $(\signedsideat{m}{i},d_j)=(0,0)$ may be assigned arbitrarily.

The commutator vanishes outside $0<s_j<L_j$, and \eqref{eq:explicit-global-forward-commutator} holds in distributions on $\Om$, including across the endpoint interfaces $s_j=0,L_j$ that meet the open polygon. After rescaling and multiplication by $\chi_{\rm vtx}$, the ray-cutoff derivatives are controlled by \eqref{eq:forward-commutator-rapid}. For the terms containing $(\chi_{m,i}^{\rm end})'$, the support has $\signedsideat{m}{i}\simeq\ellend$. A nonzero bracket forces $d_j\gtrsim\signedsideat{m}{i}$ and hence $\frac{d_j}{\e}\gtrsim\e^{-1}$. Milne decay makes these terms $O(\ue^{-\frac{c}{\e}})$ with all required tangential derivatives.

Finally, define
\begin{align}
 (\Ext_j^\e B_j)(x,w) :=\zeta_j(d_j(x))(\mathfrak E_j^{\rm glob}B_j)(x,w).
 \label{eq:physical-side-assembly}
\end{align}

\begin{lemma}
\label{lem:local-global-side-consistency}
Fix a vertex $V_m$. If $\bulkdistat{m}<2\ellvtx$, then every $\Ext_j^\e B_j$ from a side not incident to $V_m$ vanishes identically. For an incident side, set $j=j(m,i)$. Away from the single point $\sigma_i=\eta_i=0$,
\begin{align}
 (\Ext_j^\e B_j)\bigl(V_m+\e O_mY,O_mv\bigr)=\begin{cases}
 \Theta\big(\frac{\eta_i}{\sigma_i}\big)B_j\bigl(s_{m,i}^{\rm ve}(\e\sigma_i),\eta_i,O_mv\bigr),&\sigma_i>0,\\
 0,&\sigma_i\le0.
 \end{cases}
 \label{eq:local-global-side-consistency}
\end{align}
Thus coefficient extraction from the finite sum of extended side layers uses the local operator $\mathfrak E_i$, with the factor $\Theta$ occurring once.
\end{lemma}

\begin{proof}
In $\bulkdistat{m}<2\ellvtx$, condition \eqref{eq:fixed-realization-separation} gives $\abs{\signedsideat{m}{i}}<\ellend$ and $d_j<\ellnor$, hence $\chi_{m,i}^{\rm end}=\zeta_j=1$, while the cutoff at the other endpoint of the same global side vanishes there because $\ellend$ was chosen to make the endpoint supports disjoint. Then \eqref{eq:physical-wedge-ratio-identity} reduces $\Ext_j^\e$ to the right-hand side of \eqref{eq:local-global-side-consistency}. If $V_m$ is not incident to $E_j$, then \eqref{eq:nonincident-side-separation} and \eqref{eq:fixed-realization-separation} give $d_j\ge\ellsep-\bulkdistat{m}>\ellsep-2\ellvtx>2\ellnor$, so $\zeta_j(d_j)=0$ and $\Ext_j^\e B_j$ vanishes identically on $\{\bulkdistat{m}<2\ellvtx\}$.
\end{proof}

\subsection{The Approximate Solution on the Polygon}
\label{subsec:full-realization}

We first define $\Ubase^{\e,[\Nexp]}$ by replacing $\Cvertex$ with $M^{\rm con}=\chi_{\rm vtx}\Cvertex$ in the finite sum near each vertex. We then add the wedge layers, multiplied by the cutoff at radius $\delta$, to obtain $\Uapp^{\e,[\Nexp]}(\delta)$.

\subsubsection{Regularization of the finite expansion}

Recall the sum $\mathscr H_\Nexp^\e$ of the interior solution and extended side layers defined in \eqref{eq:raw-hierarchy-expanded}.

Near $V_m$, on $0<\bulkdistat{m}<2\ellvtx$, the expansions truncated at $\alpha<\Tcon$ in \eqref{eq:weighted-interior-decomposition-comp}--\eqref{eq:weighted-side-decomposition-comp} of Proposition~\ref{prop:data-generated-hierarchy-comp}, the regrouping \eqref{eq:block-binomial-rule}--\eqref{eq:side-block-binomial-rule} and Lemma~\ref{lem:local-global-side-consistency} give the exact finite decomposition
\begin{align}
 \mathscr H_\Nexp^\e &=\sum_{(\alpha,b)\in\mathfrak J_m^{\rm con}} \e^\alpha(\ln\e)^b\Cvertex_{m,\alpha,b}(Y_m,v_m) +\mathcal R_{m,\rm con}^\e,
 \label{eq:raw-finite-decomposition-expanded}
\end{align}
Here $\mathcal R_{m,\rm con}^\e$ is the sum over $0\le k\le\Nexp$, with the factor $\e^k$ at order $k$, of $\mathcal R^U_{m,k}(x,w)$ and the two extended remainders $\mathfrak E_i[\mathcal R^\Uside_{m,i,k}(\e\sigma_i,\eta_i,O_mv)](Y_m,v_m)$. The truncation exponent in all three remainders is $\Tcon$, as fixed after Corollary~\ref{cor:terminal-derivative-bounds}. The ray extension acts on the side remainders as well as on their finite expansions.

Before adding wedge layers, define $\Ubase^{\e,[\Nexp]}$ on an overlapping cover. On each vertex chart we regularize the finite sum and retain its construction remainder:
\begin{align}
 \Ubase^{\e,[\Nexp]} &:={} \begin{cases}
 \sum_{(\alpha,b)\in\mathfrak J_m^{\rm con}} \e^\alpha(\ln\e)^b M^{\rm con}_{m,\alpha,b}(Y_m,v_m) +\mathcal R_{m,\rm con}^\e, &\bulkdistat{m}<2\ellvtx,\\
 \mathscr H_\Nexp^\e, &\bulkdistat{m}>2\e\quad\text{for every }m.
 \end{cases}
 \label{eq:Ubase-finite-expanded}
\end{align}

The local vertex neighborhoods in \eqref{eq:Ubase-finite-expanded} are disjoint. On $2\e<\bulkdistat{m}<2\ellvtx$ the two branches agree exactly, because $M^{\rm con}_{m,\alpha,b}=\Cvertex_{m,\alpha,b}$ for $\scaledbulkat{m}\ge2$ and Lemma~\ref{lem:local-global-side-consistency} identifies their side contributions. They therefore define one function in $L^\infty(\Om\times\Sone)$ for each fixed $\e$. The agreement on an open overlap also permits either formula to be used when applying the transport operator there.

The replacement acts on the finite expansion at $\bulkdistat{m}=O(\e)$ and leaves the construction remainder unchanged. The incident-side cutoffs $\chi_{m,i}^{\rm end}$ and $\zeta_j$ equal one in the corresponding vertex chart. Their derivatives occur at a fixed positive distance from the vertex, where Milne decay controls the resulting commutators beyond all algebraic orders.

The finite sum in \eqref{eq:Ubase-finite-expanded} has coefficient $M^{\rm con}_{m,\alpha,b}$ at $\e^\alpha(\ln\e)^b$. Its matching subfamily $(\alpha,b)\in\mathfrak J_m^{\rm mat}$ is exactly the family $M_{m,\alpha,b}$ used in the wedge problems.

\subsubsection{Shrinking wedge overlap and coefficientwise cancellation}
\label{subsec:matched-composite-realization}

Let $\chi_{\rm ov}\in C_c^\infty([0,\infty))$ equal one on $[0,1]$ and zero on $[2,\infty)$, and let the shrinking vertex cutoff at radius $\bulkdistat{m}=\abs{x-V_m}$ be
\begin{align}
 \chi_m^\delta(x)&:=\chi_{\rm ov}\left(\frac{\bulkdistat{m}}{\delta}\right).
 \label{eq:corner-cutoff-physical}
\end{align}
For the rest of the construction we assume
\begin{align}
 0&<3\e<\delta,\qquad 2\delta<\ellvtx.
 \label{eq:realization-scale-chain}
\end{align}

Together with \eqref{eq:fixed-realization-separation}, these inequalities place the vertex regularization inside the overlap region and the overlap region inside the fixed vertex chart. In particular, the supports of the $\chi_m^\delta$ are disjoint.

Recall from \eqref{eq:forward-composite-definition} that the approximate solution is
\begin{align}
 \Uapp^{\e,[\Nexp]}(\delta)=\Ubase^{\e,[\Nexp]}+\sum_m\chi_m^\delta\mathcal D_m^\e(Y_m,v_m).
\end{align}
The added term $\chi_m^\delta\mathcal D_m^\e(Y_m,v_m)$ equals $\mathcal D_m^\e(Y_m,v_m)$ in $\bulkdistat{m}<\delta$ and vanishes for $\bulkdistat{m}\ge2\delta$. The derivative of $\chi_m^\delta$ is supported in the transition annulus $\delta\le\bulkdistat{m}\le2\delta$.

In $\bulkdistat{m}<\delta$, $\chi_m^\delta=1$ and $\Ubase+\mathcal D_m^\e=(\Ubase-\mathcal M_m^\e)+\Uwedge_m^\e$. Thus adding $\mathcal D_m^\e$ replaces $\mathcal M_m^\e$ in $\Ubase$ by the total wedge solution $\Uwedge_m^\e$. The construction remainders remain unchanged in $\Ubase-\mathcal M_m^\e$.

Intermediate coefficients of $\mathfrak J_m^{\rm mid}$ also remain in $\Ubase$. Their larger combined exponents supply the smallness needed for their uncanceled source. Section~\ref{sec:main-proof} proves an $O(\Lameps\e^{\Tmatch-1})$ bound in $L^\infty(\Om\times\Sone)$ after $\Le$ acts, using the estimates on $\scaledbulk\le3$ and the degree $\alpha-\Nexp-1$ of the derivatives of the order-$\Nexp$ coefficients.

The scaling \eqref{eq:corner-scaling-full} gives $\e\Le=v\cdot\nabla_Y+\qk=\mathscr K_{\omega_m}$ on a pulled-back coefficient at $V_m$. Hence \eqref{eq:generated-wedge-data} and the wedge layer equation imply
\begin{align}
 \Le\!\left[\e^\alpha(\ln\e)^b(M_{m,\alpha,b}+D_{m,\alpha,b})(Y_m,v_m)\right]&=0.
 \label{eq:physical-matched-coefficient-equation}
\end{align}
Its trace on the in-flow boundary is $\e^\alpha(\ln\e)^bG_{m,i,\alpha,b}$ by \eqref{eq:generated-wedge-data} and \eqref{eq:defect-wedge-problem}. These identities give the required equation and trace cancellation for each matching coefficient.

\begin{remark}
\label{rem:cutoff-ledger}
The vertex regularization acts where $\bulkdistat{m}=O(\e)$, and the overlap cutoff acts where $\bulkdistat{m}=O(\delta)$. On the latter scale the exact product rule is
\begin{align}
 \Le\bigl(\chi_m^\delta\mathcal D_m^\e(Y_m,v_m)\bigr) &= \e^{-1}\chi_m^\delta\bigl(\mathscr K_{\omega_m}\mathcal D_m^\e\bigr)(Y_m,v_m)
 +(w\cdot\nabla\chi_m^\delta)\mathcal D_m^\e(Y_m,v_m).
 \label{eq:shrinking-cutoff-exact-product-rule}
\end{align}
Every cutoff is independent of velocity, so it commutes with $\qk$.
\end{remark}

The construction gives the following well-posedness and dependence statement.

\begin{proposition}[Construction of the approximate solution]
\label{prop:profile-wellposedness-comp}
Assume \eqref{eq:strict-convexity} and Assumption~\ref{ass:primitive-data-comp}, with the parameters of Section~\ref{subsec:parameters-data}. Fix the lifting normalization of Section~\ref{subsec:harmonic-lifts-comp} and the cutoffs, extensions and length scales of Section~\ref{sec:matching-realization}. Then the following are uniquely defined:
\begin{enumerate}[label=\textup{(\roman*)}]
\item the harmonic functions $\rho_k$, interior solutions $\Uint_k$ and decaying side layers $\Uside_{j,k}$ for $0\le k\le\Nexp$ generated by \eqref{eq:microscopic-recursion-comp}--\eqref{eq:side-profile-inflow-decay-comp}, with $\rho_k$, $k\ge1$, as in Definition~\ref{def:finite-part-solution-comp}.
\item the construction expansion truncated at $\Tcon$ and the matching subexpansion truncated at $\Tmatch$, as in \eqref{eq:construction-ledger-comp}--\eqref{eq:two-ledgers-comp}.
\item every wedge layer $D_{m,\alpha,b}$, $(\alpha,b)\in\mathfrak J_m^{\rm mat}$, as the unique $L^\infty$ mild solution of \eqref{eq:defect-wedge-problem}.
\item with the fixed cutoffs and extensions above, $\Ubase^{\e,[\Nexp]}$ and the approximate solution $\Uapp^{\e,[\Nexp]}(\delta)$.
\end{enumerate}
Moreover, $\mathcal D_m^\e$ obeys \eqref{eq:summed-wedge-decay}. All the functions above, together with the coefficients in the finite vertex expansions of $\Uint_k$ and $\Uside_{j,k}$, depend continuously and linearly on $g_{j,0},g_{j,1}$.
\end{proposition}

\begin{proof}
For \textup{(i)}--\textup{(ii)}, Proposition~\ref{prop:data-generated-hierarchy-comp} constructs $\rho_k$, $\Uint_k$ and $\Uside_{j,k}$ with truncation exponent $\Tcon$ and gives their bounds and continuous linear dependence. The matching subset is obtained by restricting the finite expansion as in \eqref{eq:two-ledgers-comp}.

For \textup{(iii)}, Lemma~\ref{lem:localized-wedge-data} verifies the weighted data hypotheses of Proposition~\ref{prop:wedge-wellposedness}. That proposition gives the unique wedge layers and their continuous linear dependence. Their finite sum obeys \eqref{eq:summed-wedge-decay}.

For \textup{(iv)}, use \eqref{eq:Ubase-finite-expanded} and \eqref{eq:forward-composite-definition}. The two definitions of $\Ubase^{\e,[\Nexp]}$ agree on the open overlap, as shown after \eqref{eq:Ubase-finite-expanded}. All cutoffs and extensions are fixed, so these finite sums preserve continuous linear dependence on the data.
\end{proof}

Let $\Ext_j^\e$ denote the side extension with the fixed cutoffs defined in \eqref{eq:physical-side-assembly}. For the terms with $\alpha=0$, use the unregularized sum $\Cvertex_{m,0,0}$ and its regularization $M^{\rm con}_{m,0,0}$ from \eqref{eq:C-coefficient-expanded}--\eqref{eq:Mcon-coefficient-expanded}. Set
\begin{align}
 \mathscr H_0^\e &:=\rho_0+\sum_j\Ext_j^\e\bigl[\Uside_{j,0}\bigr],
 \label{eq:leading-raw-hierarchy}\\
 \Ubase^{\e,[0]} &:=\mathscr H_0^\e +\sum_m\bigl(M^{\rm con}_{m,0,0}-\Cvertex_{m,0,0}\bigr)(Y_m,v_m),
 \label{eq:leading-base-definition}\\
 \Uapp^{\e,[0]}(\delta) &:=\Ubase^{\e,[0]} +\sum_m\chi_m^\delta(x)D_{m,0,0}(Y_m,v_m).
 \label{eq:leading-composite-def}
\end{align}
Each regularization difference belongs to $L^\infty(K_{\omega_m}\times\Sone)$, is evaluated in the corresponding vertex coordinates, and vanishes for $\scaledbulkat{m}\ge2$. The terms with $\alpha=0$ are the endpoint value of $e_0$ and the tangentially constant terms in $\Uside_{j,0}$. 

We use the higher-order interior solutions and side layers to make $\Le\Uapp^{\e,[\Nexp]}$ small enough for Proposition~\ref{prop:physical-stability-comp}. Their contribution to the approximate solution is estimated as follows. We compare the sum of wedge layers first. We then compare $\Ubase^{\e,[\Nexp]}$ with $\Ubase^{\e,[0]}$ both outside and inside $\bulkdistat{m}=2\e$. This last comparison includes the different regularizations of the positive-degree terms at expansion order zero.

\begin{lemma}[Leading approximation]
\label{lem:leading-composite-gap}
Under the assumptions of Proposition~\ref{prop:profile-wellposedness-comp}, the leading approximation defined in \eqref{eq:leading-composite-def} satisfies
\begin{align}
 \nm{\Uapp^{\e,[\Nexp]}(\delta)-\Uapp^{\e,[0]}(\delta)}_\infty &\le C\e\Lameps
 \label{eq:leading-composite-gap-bound}
\end{align}
uniformly for the radii in \eqref{eq:realization-scale-chain}.
\end{lemma}

\begin{proof}\leavevmode

\paragraph{\underline{Step 1: wedge layers}} By \eqref{eq:hierarchy-zero-exponent-coefficients}, the only nonzero coefficients in the interior and side expansions with $\alpha=0$ have $k=b=0$. Thus for every formal pair $(0,b)$ with $b\ge1$, $M_{m,0,b}=0$ and $G_{m,i,0,b}=0$, and uniqueness in Proposition~\ref{prop:wedge-wellposedness} gives $D_{m,0,b}=0$. The same coefficient classification in Proposition~\ref{prop:data-generated-hierarchy-comp} gives $\alpha\ge1$ for every other nonzero coefficient.

For each remaining pair, $\e^\alpha\abs{\ln\e}^b\le\e\Lameps$, since $\alpha\ge1$ and $b\le\Llog$ by \eqref{eq:delta-lambda-comp}. The finite sum \eqref{eq:summed-wedge-corrector-definition} and Proposition~\ref{prop:wedge-wellposedness} therefore give
\begin{align}
 \nm{\sum_m\chi_m^\delta\bigl(\mathcal D_m^\e-D_{m,0,0}\bigr)(Y_m,v_m)}_\infty&\le C\e\Lameps.
 \label{eq:leading-gap-wedge}
\end{align}
The constants include $\nm{\chi_{\rm ov}}_{L^\infty(0,\infty)}$ and are independent of $\delta$.

\paragraph{\underline{Step 2: $2\e\le r<2\ellvtx$}} Write $r=\bulkdistat{m}$ on a fixed vertex chart. In this region, $\Ubase^{\e,[\Nexp]}-\Ubase^{\e,[0]}$ is the sum of $\e^k\Uint_k$ and the extended side layers for $1\le k\le\Nexp$. A nonzero power-logarithmic term in the vertex expansion of $\Uint_k$ has $\alpha\ge1$ and $L^\infty_w$ norm at most $C\e^k r^{\alpha-k}\bigl(1+\abs{\ln r}\bigr)^{\Llog}$. If $\alpha-k<0$, then $r\ge2\e$ gives
\begin{align}
 \e^k r^{\alpha-k} &=\e^\alpha\left(\frac{r}{\e}\right)^{\alpha-k}\le\e^\alpha\le\e.
 \label{eq:leading-gap-negative-degree}
\end{align}
If $\alpha-k\ge0$, the fixed upper bound $r<1$ gives $\e^k r^{\alpha-k}\le\e$. Also $1+\abs{\ln r}\le C(1+\abs{\ln\e})$ on this region. Both signs of the spatial degree therefore give the same $O(\e\Lameps)$ bound in $L^\infty_{x,w}$ on this annulus.

For a side term, Lemma~\ref{lem:local-global-side-consistency} identifies the ray-cutoff support, where the positive tangential argument $s_i=\e\sigma_i$ is comparable to $r$ because $\eta_i\le2\cfs\sigma_i$. Thus $\frac{s_i}{\e}\ge c>0$, and the same two power estimates hold with a constant depending on the cutoff and the finite exponents. The Milne factor belongs to $L^\infty_{\eta,w}$. The construction remainder at order $k$ has $L^\infty_w$ norm at most $Cr^{\Tcon-k}\bigl(1+\abs{\ln r}\bigr)^{\Llog}$, or the corresponding bound with $s_i$ in the side term, by \eqref{eq:outer-cartesian-remainder-comp} and \eqref{eq:side-ledger-remainder-comp} at $\mathsf d=\Tcon$. Since $\Tcon-k>0$, its $L^\infty_{x,w}$ norm on the chart is finite and its contribution after multiplication by $\e^k$ is $O(\e)$ for $k\ge1$.

\paragraph{\underline{Step 3: $r<2\e$}} The higher-order displayed terms in $\Ubase^{\e,[\Nexp]}$ are regularized. Their coefficients have finite $L^\infty_{Y,v}$ norms on $\scaledbulkat{m}<2$, and $\alpha\ge1$ again gives an $L^\infty_{x,w}$ bound $C\e\Lameps$ after summation. The higher-order construction remainders obey the same bound as above.

At order zero, $\Ubase^{\e,[0]}$ in \eqref{eq:leading-base-definition} regularizes the coefficient with combined exponent zero. The function $\Ubase^{\e,[\Nexp]}$ also regularizes the remaining terms of the finite expansion. This difference is supported in $r<2\e$, where every nonzero remaining term has $\alpha\ge1$ and hence
\begin{align}
 \sup_{0<r<2\e}r^\alpha\bigl(1+\abs{\ln r}\bigr)^b &\le C\e^\alpha\bigl(1+\abs{\ln\e}\bigr)^b \le C\e\Lameps, \qquad b\le \Llog.
 \label{eq:leading-gap-order-zero}
\end{align}
This follows by setting $r=\e R$ and using $\sup_{0<R<2}R^\alpha\bigl(1+\abs{\ln R}\bigr)^b<\infty$. For a side term, the same power-logarithmic estimate applies to its positive tangential argument, and its Milne factor belongs to $L^\infty_{\eta,w}$. The order-zero construction remainder is identical in $\Ubase^{\e,[\Nexp]}$ and $\Ubase^{\e,[0]}$ and cancels.

\paragraph{\underline{Step 4: away from the vertices}} On the complement of the vertex charts, the estimates for $\Uint_k$ and $\Uside_{j,k}$ on $\Om_{\rm away}$ and the fixed cutoffs give $L^\infty_{x,w}$ bounds independent of $\e$ for $\Uint_k$ and $\Ext_j^\e\Uside_{j,k}$, $1\le k\le\Nexp$. The functions $\Ubase^{\e,[\Nexp]}$ and $\Ubase^{\e,[0]}$ therefore differ there by $O(\e)$ in $L^\infty_{x,w}$. Combining the finitely many chart estimates with \eqref{eq:leading-gap-wedge} proves \eqref{eq:leading-composite-gap-bound}.
\end{proof}
\section{\texorpdfstring{$L^\infty$}{} Diffusive Limit}
\label{sec:main-proof}

We now estimate the remainder between $\Uphys^\e$ and the approximate solution $\Uapp^{\e,[\Nexp]}(\delta)$ constructed in Section~\ref{sec:matching-realization}. The equation is satisfied exactly by each matched sum $M_{m,\alpha,b}+D_{m,\alpha,b}$ near its vertex. The source that remains comes from the order-$\Nexp$ derivatives outside the overlap, the intermediate exponents, the construction remainders, and the cutoff derivatives. We estimate these sources and the error in the in-flow boundary condition, then apply Proposition~\ref{prop:physical-stability-comp}. The radius $\delta$ will be chosen after these estimates.

Define the remainder
\begin{align}
 \Err^\e(\delta):=\Uphys^\e-\Uapp^{\e,[\Nexp]}(\delta)
\end{align}
Since $\Le\Uphys^\e=0$ and $\Uphys^\e\big|_{\gamma_-}=g^\e$, setting $S_{\rm app}^\e:=\Le\Uapp^{\e,[\Nexp]}(\delta)$ gives
\begin{align}
 \Le\Err^\e(\delta) &= -S_{\rm app}^\e,
 \label{eq:error-equation}\\
 \Err^\e(\delta)\big|_{\gamma_-} &= g^\e- \Uapp^{\e,[\Nexp]}(\delta)\big|_{\gamma_-}.
 \label{eq:error-incoming-data}
\end{align}
We verify the hypotheses of Lemma~\ref{lem:admissible-mild-comp} for each fixed $\e,\delta$.
\begin{enumerate}
\item In a vertex chart, the finite coefficients of \eqref{eq:Ubase-finite-expanded} are $M^{\rm con}_{m,\alpha,b}=\chi_{\rm vtx}\Cvertex_{m,\alpha,b}$ from \eqref{eq:Mcon-coefficient-expanded}. They vanish for $\scaledbulkat{m}\le1$, so their unregularized powers and derivatives are evaluated away from the vertex. For every remainder order $0\le k\le\Nexp$, \eqref{eq:depth-choices-comp} gives $\Tcon-k-1\ge\Tcon-\Nexp-1>\pacc>0$. This positive exponent absorbs the logarithms after one spatial derivative, including the factor $\bulkdistat{m}^{-1}$ from a differentiated ratio cutoff. Outside the charts, the estimates away from the vertices and the positive endpoint distance on each active side branch place $\Uint_k$, the extended side layers, and their tangential derivatives in $L^\infty_{x,w}$. The two formulas for $\Ubase^{\e,[\Nexp]}$ agree on their overlap, and $\Uint_k$ has a continuous restriction on each open side.

\item For the side factors, the mild Milne formulas give absolute continuity in the normal variable for almost every nongrazing velocity. Combined with the tangential regularity of the side layers, this yields the weak streaming identity
\begin{align}
 w\cdot\nabla_x\left[\mathsf F\Big(s_j,\frac{d_j}{\e},w\Big)\right] &=\tau_j(w)\partial_{s_j}\mathsf F+\e^{-1}\bigl(\mathsf S-\qk\mathsf F\bigr),
\end{align}
Here $\mathsf F(s,\eta,w)$ denotes the side factor under consideration, and $\mathsf S$ is its Milne source, so $(\mu_j\partial_\eta+\qk)\mathsf F=\mathsf S$. The derivative $\partial_{s_j}\mathsf F$ is taken at fixed $(\eta,w)$. After this differentiation, every function on the right-hand side is evaluated at $(s_j(x),\frac{d_j(x)}{\e},w)$. On the branches above, this grouped identity gives an $L^\infty_{x,w}$ estimate for the streaming derivative at fixed $\e$. The ray-cutoff support satisfies $\sigma_i\simeq\scaledbulkat{m}$, so the terms containing derivatives of $\Theta$ belong to $L^\infty_{x,w}$ after regularization or multiplication by the positive remainder weight. For $\eta_i>0$, the ray factor vanishes near $\sigma_i=0$, and the global formulas agree across their interfaces, so extension by zero creates no distributional interface source. On its own side, a Milne factor has the prescribed trace on the in-flow boundary by \eqref{eq:milne-char-plus-comp}. A contribution from another side is evaluated at its interior normal coordinate through the characteristic representative. On each compact tangential interval, the tangential derivative bound and $\mu_j(w)\partial_\eta\mathsf F=\mathsf S-\qk\mathsf F$ give a locally Lipschitz representative in $(s,\eta)$ for almost every fixed nongrazing velocity, extending to $\eta=0$. Its weighted essential value bound therefore holds after restriction to another side's coordinate curve, by continuity. Division by $\mu_j(w)$ is used only to identify this representative at a fixed velocity and does not enter the uniform value bound.

\item Each wedge layer satisfies $v\cdot\nabla_YD_{m,\alpha,b}=F_{m,\alpha,b}-\qk D_{m,\alpha,b}\in L^\infty_{Y,v}$ along almost every characteristic. The substitution $Y=\e^{-1}O_m^T(x-V_m)$ preserves membership in $L^\infty$ for fixed $\e$, and multiplication by $\chi_m^\delta$ gives a function satisfying \eqref{eq:admissible-transport-identity} inside the wedge chart. The uncut wedge layer has trace on the in-flow boundary $h_{m,i,\alpha,b}$ from \eqref{eq:generated-wedge-data}, and the product has trace on the in-flow boundary $\chi_m^\delta h_{m,i,\alpha,b}$ after rescaling.
\end{enumerate}
Thus each local expression for $\Uapp^{\e,[\Nexp]}(\delta)$ satisfies \eqref{eq:admissible-transport-identity}, and its trace is the one computed in Step~5 below. Lemma~\ref{lem:admissible-mild-comp} and linearity give the mild remainder problem \eqref{eq:error-equation}--\eqref{eq:error-incoming-data}. Proposition~\ref{prop:physical-stability-comp} consequently gives
\begin{align}
 \nm{\Err^\e(\delta)}_\infty &\le \babs{g^\e- \Uapp^{\e,[\Nexp]}(\delta)}_{\infty,-} + C\e^{-1}\nm{S_{\rm app}^\e}_\infty.
 \label{eq:physical-stability-recall}
\end{align}
Here $S_{\rm app}^\e$ includes both streaming and collision. Thus the streaming derivative is $S_{\rm app}^\e-\e^{-1}\qk\Uapp^{\e,[\Nexp]}(\delta)$. All boundary values in \eqref{eq:error-incoming-data} are characteristic limits almost everywhere. By \eqref{eq:physical-stability-recall}, it remains to estimate $S_{\rm app}^\e$ and the error in the in-flow boundary condition. The source estimate loses one power of $\e$, whereas the boundary estimate enters directly. Thus the exact cancellation of the matched source must be used before estimating its separate contributions. In boundary norms such as $\abs{\,\cdot\,}_{\infty,-}$, the arguments are understood as their characteristic traces.

Proposition~\ref{prop:full-residual} gives both bounds for every radius satisfying \eqref{eq:realization-scale-chain}. The endpoint Taylor-remainder bound $\delta^{\Tmatch}$ increases with $\delta$, while the bound $\e^{\betaw-1}\delta^{-\betaw-1}$ from the wedge cutoff decreases. We choose $\delta$ by balancing these two powers, with the logarithmic factor $\Lameps$ included in the final estimate.

\subsection{Source, Trace, and Remainder Estimates}
\label{subsec:remainder-equation}

Applying the product rule in \eqref{eq:shrinking-cutoff-exact-product-rule} gives
\begin{align}
 S_{\rm app}^\e &= \Le\Ubase^{\e,[\Nexp]} -\e^{-1}\sum_m\chi_m^\delta \bigl(\mathscr K_{\omega_m}\mathcal M_m^\e\bigr)(Y_m,v_m) +\sum_m (w\cdot\nabla_x\chi_m^\delta) \mathcal D_m^\e(Y_m,v_m),
 \label{eq:composite-source-identity}
\end{align}
where the middle sum has been rewritten from $\e^{-1}\sum_m\chi_m^\delta(\mathscr K_{\omega_m}\mathcal D_m^\e)(Y_m,v_m)$ by $\mathscr K_{\omega_m}\mathcal D_m^\e=-\mathscr K_{\omega_m}\mathcal M_m^\e$, and the last sum is the only place a derivative of the shrinking cutoff occurs.

\begin{proposition}[Source, trace, and remainder estimates]
\label{prop:full-residual}
Assume \eqref{eq:realization-scale-chain}, and use $\Ubase^{\e,[\Nexp]}$ and $\mathcal M_m^\e$ from \eqref{eq:Ubase-finite-expanded} and \eqref{eq:Mless-finite-expanded}. Then the following estimates hold.
\begin{enumerate}
\item The source $S_{\rm app}^\e$ has the decomposition \eqref{eq:source-pieces-expanded} and satisfies
\begin{align}
 \nm{S_{\rm app}^\e}_\infty &\le C\Lameps \left\{\e^\Nexp\delta^{-\Nexp-1}+\e^{\Tmatch-1} +\delta^{-1}\left(\frac{\e}{\delta}\right)^\betaw\right\},
 \label{eq:S-app-est}
\end{align}
\item The error in the in-flow boundary condition satisfies
\begin{align}
 \abs{\Uapp^{\e,[\Nexp]}(\delta)-g^\e}_{\infty,-} &\le C\left\{\delta^{\Tmatch}+\e^2+ \Lameps\Big[\left(\frac{\e}{\delta}\right)^\betaw+\ue^{-c\frac{\delta}{\e}}\Big]\right\}.
 \label{eq:trace-est}
\end{align}
\item The remainder satisfies
\begin{align}
 \nm{\Err^\e(\delta)}_\infty \le{}& C\bigl(\delta^{\Tmatch}+\e^2\bigr)+C\Lameps\left[\left(\frac{\e}{\delta}\right)^\betaw+\ue^{-c\frac{\delta}{\e}}+\e^{\Nexp-1}\delta^{-\Nexp-1}+\e^{\Tmatch-2}+\e^{\betaw-1}\delta^{-\betaw-1}\right].
 \label{eq:preoptimized-error}
\end{align}
\end{enumerate}
The constants depend only on the fixed polygon, parameter choices, cutoffs, extension and regularization norms, and the finite data norms in Assumption~\ref{ass:primitive-data-comp}. They are independent of $\e,\delta$.
\end{proposition}

\begin{proof}
Step~1 identifies the source exactly, and Steps~2 to~4 estimate its terms in $L^\infty(\Om\times\Sone)$. Step~5 estimates the trace on the in-flow boundary in $L^\infty(\gamma_-)$. We use $\mathscr K_{\omega_m}(\mathcal M_m^\e+\mathcal D_m^\e)=0$ before taking these norms. All $O(\cdot)$ assertions for source terms below refer to $L^\infty(\Om\times\Sone)$. Summing \eqref{eq:completed-string-coefficient-recurrence}--\eqref{eq:completed-side-coefficient-recurrence} leaves the order-$\Nexp$ derivatives in \eqref{eq:local-interior-source}--\eqref{eq:local-side-source}.

\paragraph{\underline{Step 1: Exact decomposition of $S_{\rm app}^\e$}}

We first fix a partition of unity on $\Om$ adapted to the vertices. Choose disjointly supported $\psi_m\in C_c^\infty(\{\bulkdistat{m}<2\ellvtx\})$, equal to one on $\{\bulkdistat{m}<\ellvtx\}$, and set $\psi_0:=1-\sum_m\psi_m$, so that $\psi_0+\sum_m\psi_m=1$. We take $2\delta<\ellvtx$, so each localized contribution $\chi_m^\delta(x)\mathcal D_m^\e(Y_m,v_m)$ is supported where $\psi_m=1$. We use this partition to decompose the computed source $\Le\Ubase$, so the factors $\psi_m$ occur without derivatives.

In each vertex chart, we use the decomposition \eqref{eq:Ubase-finite-expanded} into the regularized finite sum and the construction remainder. We then combine the source of each $M_{m,\alpha,b}$ with $(\alpha,b)\in\mathfrak J_m^{\rm mat}$ with that of $D_{m,\alpha,b}$.

\subparagraph{Interior and side derivatives} Summing \eqref{eq:microscopic-recursion-comp}--\eqref{eq:side-profile-inflow-decay-comp} leaves derivatives of $\Uint_\Nexp$ and $\Uside_{j,\Nexp}$, together with the commutator from extending the side layers. Define $\mathscr S_{j,\Nexp}^{\e,{\rm unc}}:=\mathfrak E_j^{\rm glob}\big[\sum_{k=0}^\Nexp\e^k\Uside_{j,k}\big]$ using \eqref{eq:explicit-global-forward-glue}, with the normal cutoff omitted, and set
\begin{align}
 \mathcal T_I^\e&:=\Le\Big[\sum_{k=0}^\Nexp\e^k\Uint_k\Big]=\e^\Nexp\big(w\cdot\nabla_x\Uint_\Nexp\big),
 \label{eq:global-interior-terminal-expanded}\\
 \mathcal T_B^\e&:=\sum_j\mathfrak E_j^{\rm glob}\!\left[ \e^\Nexp\tau_j(w)\partial_{s_j}\Uside_{j,\Nexp}\right],
 \label{eq:global-side-terminal-expanded}\\
 \mathcal C_{\rm fwd}^\e&:=\sum_j\left\{ \Le\mathscr S_{j,\Nexp}^{\e,{\rm unc}} -\mathfrak E_j^{\rm glob}\!\left[ \e^\Nexp\tau_j(w)\partial_{s_j}\Uside_{j,\Nexp}\right] \right\}.
 \label{eq:global-forward-commutator-expanded}
\end{align}

On $\supp\psi_0$, $\Ubase^{\e,[\Nexp]}$ equals the unregularized expansion $\mathscr H_\Nexp^\e$ of \eqref{eq:raw-hierarchy-expanded}, so these formulas apply there. In the coordinates \eqref{eq:wedge-side-coordinates}--\eqref{eq:corner-scaling-full}, the order-$k$ endpoint contribution to \eqref{eq:global-forward-commutator-expanded} is
\begin{align}
 \label{eq:forward-commutator-chart-expanded} {}&\e^{k-1}\chi_{m,i}^{\rm end}(\e\sigma_{m,i}) \left[\Big(\tau_{m,i}(w)\partial_{\sigma_{m,i}} +\mu_{m,i}(w)\partial_{\eta_{m,i}}\Big) \Theta\Big(\frac{\eta_{m,i}}{\sigma_{m,i}}\Big)\right] \Uside_{m,i,k}(\e\sigma_{m,i},\eta_{m,i},w) \\
 &\quad+\e^k\tau_{m,i}(w) \left(\chi_{m,i}^{\rm end}\right)'(\e\sigma_{m,i}) \Big[\Theta\Big(\frac{\eta_{m,i}}{\sigma_{m,i}}\Big)-1\Big] \Uside_{m,i,k}(\e\sigma_{m,i},\eta_{m,i},w).\notag
\end{align}
On $\supp\psi_m$, the identities $\chi_{m,i}^{\rm end}=1$ and $(\chi_{m,i}^{\rm end})'=0$ leave only the ray-cutoff term containing the derivative of $\Theta$. This term decays rapidly on the wedge scale. The term containing $\bigl(\chi_{m,i}^{\rm end}\bigr)^{\prime}$ is supported outside the vertex charts, where Milne decay at a fixed positive normal distance makes it $O(\ue^{-\frac{c}{\e}})$.

Before using $M^{\rm con}_{m,\alpha,b}$, remove the normal cutoff through the exact identity $\Le(\zeta_j\mathscr S_{j,\Nexp}^{\e,{\rm unc}})=\Le\mathscr S_{j,\Nexp}^{\e,{\rm unc}}+(\zeta_j-1)\Le\mathscr S_{j,\Nexp}^{\e,{\rm unc}}+(w\cdot\nabla\zeta_j)\mathscr S_{j,\Nexp}^{\e,{\rm unc}}$. Summed over $j$, the first term is $\mathcal T_B^\e+\mathcal C_{\rm fwd}^\e$ by the definition of $\mathcal C_{\rm fwd}^\e$. The last two terms form the normal-cutoff source $S_{\rm nor}$ in \eqref{eq:Snor-expanded}.

\subparagraph{Source in a vertex chart} Inside the chart of $V_m$, each pair $(\alpha,b)\in\mathfrak J_m^{\rm con}$ satisfies the rescaling identity
\begin{align}
 \Le\!\left[\e^\alpha(\ln\e)^bM^{\rm con}_{m,\alpha,b}(Y_m,v_m)\right] &=\e^{\alpha-1}(\ln\e)^b\bigl(\mathscr K_{\omega_m}M^{\rm con}_{m,\alpha,b}\bigr)(Y_m,v_m).
 \label{eq:local-source-rescaling}
\end{align}
The decomposition \eqref{eq:completed-string-error-split} gives, for every $(\alpha,b)\in\mathfrak J_m^{\rm con}$,
\begin{align}
 \mathscr K_{\omega_m} M^{\rm con}_{m,\alpha,b} &=E^I_{m,\alpha,b}+E^B_{m,\alpha,b}+E^{\rm fwd}_{m,\alpha,b}+E^{\rm vtx}_{m,\alpha,b}.
 \label{eq:EC-def-identity}
\end{align}
Here $E^I$, $E^B$, $E^{\rm fwd}$ and $E^{\rm vtx}$ denote $\mathcal E^{I,\mathrm{loc}}$, $\mathcal E^{B,\mathrm{loc}}$, $\mathcal E^{F,\mathrm{loc}}$ and $\mathcal E^{\rm vtx,\mathrm{loc}}$ from \eqref{eq:local-interior-source}--\eqref{eq:local-vertex-commutator}, respectively, with the vertex index $m$ restored. In the source sums below these functions are evaluated at $(Y_m(x),v_m(w))$, and $\psi_m,\chi_m^\delta$ at $x$.

The terms $E^I_{m,\alpha,b}$ and $E^B_{m,\alpha,b}$ are $\chi_{\rm vtx}(\scaledbulkat{m})$ times $v\cdot\nabla_YC^I_\Nexp$ and $\sum_i\mathfrak E_i[\tau_i^Y\partial_{\sigma_i}C^B_{i,\Nexp}]$, respectively. The ray cutoff in the order-$\Nexp$ side term remains undifferentiated. Its derivative belongs to $E^{\rm fwd}_{m,\alpha,b}$, which is the ray part of \eqref{eq:forward-commutator-chart-expanded} regrouped by $\e^{\alpha-1}(\ln\e)^b$. The grouped Milne equation gives the $L^\infty_v$ estimate for the side operator.

All four coefficients belong to $L^\infty_{Y,v}$ on $\scaledbulkat{m}\le3$ and vanish on $\scaledbulkat{m}\le1$, where their value is set to zero before the unregularized formulas are evaluated. The factor $v\cdot\nabla_Y\chi_{\rm vtx}$ in the regularization source gives
\begin{align}
 \supp E^{\rm vtx}_{m,\alpha,b} &\subset\{1\le\scaledbulkat{m}\le2\}.
 \label{eq:EC-support-expanded}
\end{align}
On $\scaledbulkat{m}\ge2$, $\chi_{\rm vtx}=1$ and $M^{\rm con}_{m,\alpha,b}=\Cvertex_{m,\alpha,b}$, so the first three terms are the unregularized order-$\Nexp$ derivatives and ray-cutoff terms.

With $\chi_m^{\delta}$ from \eqref{eq:corner-cutoff-physical}, define $\Xi_{m,\alpha}^\delta:=1-\chi_m^\delta$ for $\alpha<\Tmatch$ and $\Xi_{m,\alpha}^\delta:=1$ for $\alpha\ge\Tmatch$. For a matched pair, evaluate the wedge equation at $(Y_m,v_m)$ and multiply it by $\chi_m^\delta(x)$. This gives
\begin{align}
 \mathscr K_{\omega_m}M_{m,\alpha,b}+\chi_m^\delta\mathscr K_{\omega_m}D_{m,\alpha,b} &=(1-\chi_m^\delta)\mathscr K_{\omega_m}M_{m,\alpha,b}.
 \label{eq:matched-coefficient-source-cancellation}
\end{align}
Consequently all four terms of \eqref{eq:EC-def-identity} acquire the same factor $\Xi_{m,\alpha}^\delta$. For intermediate pairs, that factor is one. The derivative of $\chi_m^\delta$ remains the last term of \eqref{eq:composite-source-identity}.

By \eqref{eq:EC-support-expanded} and \eqref{eq:realization-scale-chain}, $\supp E^{\rm vtx}_{m,\alpha,b}\subset\{\bulkdistat{m}\le2\e\}\subset\{\chi_m^\delta=1\}$, so $\Xi_{m,\alpha}^\delta=0$ there and
\begin{align}
 \Xi_{m,\alpha}^\delta E^{\rm vtx}_{m,\alpha,b}=0, \qquad\alpha<\Tmatch.
 \label{eq:wedge-cutoff-cancellation}
\end{align}
The matched regularization source cancels exactly. Moreover, the matched order-$\Nexp$ derivative sources and ray-cutoff sources are supported after cancellation in $\bulkdistat{m}\ge\delta$. Their $L^\infty_{x,w}$ estimates therefore use only this region, giving the factor $\delta^{-\Nexp-1}$ in the order-$\Nexp$ source bound below.

The construction remainders with truncation exponent $\Tcon$, fixed after Proposition~\ref{prop:data-generated-hierarchy-comp}, are unchanged by \eqref{eq:Ubase-finite-expanded} and contain no factor $\chi_{\rm vtx}$. Their exact projected recurrences \eqref{eq:exact-interior-displayed-recurrence-comp} and \eqref{eq:exact-side-projected-recurrence-comp} cancel every forcing below order $\Nexp$. Since $\Le$ acts on the side remainder as $\tau_{m,i}(w)\partial_{\sidedistat{m}{i}}+\e^{-1}\big(\mu_{m,i}(w)\partial_{\eta_{m,i}}+\qk\big)$, this gives, before the derivative hits the ray cutoff,
\begin{align}
 \Le\Big[\sum_{k=0}^\Nexp\e^k\mathcal R^U_{m,k}\Big] &=\e^\Nexp\left(w\cdot\nabla_x\mathcal R^U_{m,\Nexp}\right),\\
 \Le\Big[\sum_{k=0}^\Nexp\e^k \mathcal R^\Uside_{m,i,k} (\sidedistat{m}{i},\eta_{m,i},w)\Big] &=\e^\Nexp\tau_{m,i}(w) \partial_{\sidedistat{m}{i}}\mathcal R^\Uside_{m,i,\Nexp} (\sidedistat{m}{i},\eta_{m,i},w),
\end{align}
\subparagraph{Summation of the sources} Set $S_{I,\rm far}^{[\Nexp]}:=\psi_0\mathcal T_I^\e$ and $S_{B,\rm far}^{[\Nexp]}:=\psi_0\mathcal T_B^\e$ for the contributions supported in $\supp\psi_0$. The sources from the order-$\Nexp$ interior and side derivatives are
\begin{align}
 \label{eq:SIterm-expanded} S_I^{[\Nexp]} &:=S_{I,\rm far}^{[\Nexp]} +\sum_m\psi_m \sum_{(\alpha,b)\in\mathfrak J_m^{\rm con}} \e^{\alpha-1}(\ln\e)^b \Xi_{m,\alpha}^\delta E^I_{m,\alpha,b},\\
 \label{eq:SBterm-expanded} S_B^{[\Nexp]} &:=S_{B,\rm far}^{[\Nexp]} +\sum_m\psi_m \sum_{(\alpha,b)\in\mathfrak J_m^{\rm con}} \e^{\alpha-1}(\ln\e)^b \Xi_{m,\alpha}^\delta E^B_{m,\alpha,b}.
\end{align}
The normal cutoff removed above contributes the source supported in $\supp\psi_0$
\begin{align}
 S_{\rm nor} &:=\psi_0\sum_j\left\{ (\zeta_j-1)\Le\mathscr S_{j,\Nexp}^{\e,{\rm unc}} +(w\cdot\nabla\zeta_j)\mathscr S_{j,\Nexp}^{\e,{\rm unc}} \right\}.
 \label{eq:Snor-expanded}
\end{align}
In each vertex chart, $\Ubase^{\e,[\Nexp]}$ is computed from its local coefficients. The incident normal cutoffs equal one, and nonincident extended layers vanish by Lemma~\ref{lem:local-global-side-consistency}. Thus the normal-cutoff correction in \eqref{eq:Snor-expanded} occurs only in the $\psi_0$ branch of the definition.

By \eqref{eq:wedge-cutoff-cancellation}, all terms multiplying $E^{\rm vtx}_{m,\alpha,b}$ with $\alpha<\Tmatch$ vanish. The remaining contribution is
\begin{align}
 S_{\rm vtx}^{\rm mid} &:=\sum_m\psi_m \sum_{(\alpha,b)\in\mathfrak J_m^{\rm mid}} \e^{\alpha-1}(\ln\e)^bE^{\rm vtx}_{m,\alpha,b}.
\end{align}
The two identities for the construction remainders give the contribution
\begin{align}
 S_{\rm rem}^{[\Nexp]} &:=\sum_m\psi_m\e^\Nexp\bigg\{ w\cdot\nabla_x\mathcal R^U_{m,\Nexp} +\sum_{i=1}^2\mathds 1_{\{\sigma_{m,i}>0\}} \Theta\Big(\frac{\eta_{m,i}}{\sigma_{m,i}}\Big) \tau_{m,i}(w)\partial_{\sidedistat{m}{i}} \mathcal R^\Uside_{m,i,\Nexp} (\e\sigma_{m,i},\eta_{m,i},w) \bigg\}.
 \label{eq:Scon-exact-expanded}
\end{align}

The remaining source collects the commutators containing derivatives of $\Theta$ and $\chi_{m,i}^{\rm end}$. Put $S_{{\rm fwd},\rm far}:=\psi_0\mathcal C_{\rm fwd}^\e$, and let $S_{\rm fwd}^{\rm rem}$ be obtained from the ray-cutoff contribution in \eqref{eq:forward-commutator-chart-expanded}, summed over $m$, over $i=1,2$ and over $0\le k\le \Nexp$, by inserting the factor $\psi_m\mathds 1_{\{\sigma_{m,i}>0\}}$ and replacing $\Uside_{m,i,k}$ by $\mathcal R^\Uside_{m,i,k}$.

The total commutator source is
\begin{align}
 S_{\rm fwd} &:=S_{{\rm fwd},\rm far}+S_{\rm fwd}^{\rm rem} +\sum_m\psi_m \sum_{(\alpha,b)\in\mathfrak J_m^{\rm con}} \e^{\alpha-1}(\ln\e)^b \Xi_{m,\alpha}^\delta E^{\rm fwd}_{m,\alpha,b}.
\end{align}
Direct substitution and \eqref{eq:composite-source-identity} now give the exact finite decomposition announced in the proposition,
\begin{align}
 S_{\rm app}^\e &=S_I^{[\Nexp]}+S_B^{[\Nexp]}+S_{\rm nor}+S_{\rm vtx}^{\rm mid} +S_{\rm rem}^{[\Nexp]}+S_{\rm fwd} +\sum_m(w\cdot\nabla_x\chi_m^\delta)\mathcal D_m^\e(Y_m,v_m),
 \label{eq:source-pieces-expanded}
\end{align}
all six named sources depending on the suppressed arguments $(\e,\delta)$.

\paragraph{\underline{Step 2: Estimates of $S_I^{[\Nexp]}+S_B^{[\Nexp]}$}}

Fix $V_m$ and abbreviate $r:=\bulkdistat{m}$, $R:=\scaledbulkat{m}=\frac{r}{\e}$, $(Y,v):=(Y_m,v_m)$, $(\sigma_i,\eta_i):=(\sigma_{m,i},\eta_{m,i})$, and $s_i:=\e\sigma_i$. Tangential derivatives $\partial_{s_i}$ of a side layer or its remainder are taken at fixed $(\eta_i,w)$ before these coordinates are substituted. On the ray-cutoff support, $s_i$ is the positive tangential argument at which the side layer is evaluated. Put $\tau_i:=\tau_{m,i}(w)=\tau_i^Y(v)$ after the velocity pullback, and suppress $m$ from $E^\bullet_{\alpha,b}$ and $\Xi_\alpha^\delta$. Constants are uniform over the finite vertex set.

By \eqref{eq:global-interior-terminal-expanded}, \eqref{eq:global-side-terminal-expanded} and \eqref{eq:local-global-side-consistency}, the side layers are extended by $\mathfrak E_i$ on $\supp\psi_m$. The finite coefficient sums contribute $\e^\Nexp w\cdot\nabla_xC^U_{m,\Nexp}$ and $\mathfrak E_i[\e^\Nexp\tau_i\partial_{s_i}C^\Uside_{m,i,\Nexp}]$. Their construction remainders contribute exactly $S_{\rm rem}^{[\Nexp]}$ in \eqref{eq:Scon-exact-expanded}.

Rescaling to the kinetic variables gives the chart coefficients. Since $r=\e R$ and $w\cdot\nabla_x=\e^{-1}v\cdot\nabla_Y$, the term $r^{\alpha-\Nexp}(\ln r)^{b_0}A(\theta_m,v)$ in the vertex expansion of $\Uint_\Nexp$, with $0\le b_0\le\Llog$, contributes to $\mathcal T_I^\e$
\begin{align}
 \e^\Nexp w\cdot\nabla_x \left[r^{\alpha-\Nexp}(\ln r)^{b_0}A\right] &=\e^{\alpha-1}\sum_{b=0}^{b_0}\binom{b_0}{b}(\ln\e)^b v\cdot\nabla_Y\left[R^{\alpha-\Nexp}(\ln R)^{b_0-b}A\right],
\end{align}
where $A=A(\theta_m,v)$. A side term $s_i^{\alpha-\Nexp}(\ln s_i)^{b_0}F_i(\eta_i,v)$ with $\nm{F_i(\eta_i,\cdot)}_{L^\infty_v}\le C\ue^{-\kappa_\ast\eta_i}$ has the same rescaling, with $v\cdot\nabla_Y[R^{\alpha-\Nexp}(\ln R)^{b_0-b}A]$ replaced by $\Theta\big(\frac{\eta_i}{\sigma_i}\big)\tau_i\partial_{\sigma_i}[\sigma_i^{\alpha-\Nexp}(\ln\sigma_i)^{b_0-b}F_i]$.

The spatial derivative of an order-$\Nexp$ term lowers its radial degree by one, giving $\lambda_\Nexp=\lambda_\Nexp(\alpha):=\alpha-\Nexp-1$. Derivatives of $\Theta$ are included in $S_{\rm fwd}$. Collecting equal powers therefore gives $v\cdot\nabla_YC^I_\Nexp$ and $\sum_i\mathfrak E_i[\tau_i^Y\partial_{\sigma_i}C^B_{i,\Nexp}]$. These equal $E^I_{\alpha,b}$ and $E^B_{\alpha,b}$ on $R\ge2$ by \eqref{eq:EC-def-identity}.

For $R\ge3$, $\chi_{\rm vtx}=1$, so differentiation of $C_\Nexp^I$ gives the first bound below. For $C_{i,\Nexp}^B$, the support of the ray cutoff satisfies $\sigma_i>0$, $0\le\eta_i<2\cfs\sigma_i$ and $R\simeq\sigma_i$. Hence the Milne estimate gives the second bound:
\begin{align}
 \nm{E^I_{\alpha,b}(Y,\cdot)}_{L^\infty_v} &\le C R^{\lambda_\Nexp} \bigl(1+\ln R\bigr)^{\Llog-b},
 \label{eq:EI-degree-expanded}\\
 \nm{E^B_{\alpha,b}(Y,\cdot)}_{L^\infty_v} &\le C R^{\lambda_\Nexp} \bigl(1+\ln R\bigr)^{\Llog-b} \sum_{i=1}^2 \mathds 1_{\{\sigma_i>0,\ 0\le\eta_i<2\cfs\sigma_i\}} \ue^{-\kappa_\ast\eta_i}.
 \label{eq:EB-degree-expanded}
\end{align}
These estimates hold for almost every $Y$, with essential suprema in velocity. Both families vanish on $R\le1$ and belong to $L^\infty_{Y,v}$ on $1<R<3$.

Every product $(\ln\e)^b(\ln R)^c$ comes from a logarithm in $r$ or $s_i$ of degree at most $\Llog$. Thus, for $3\le R\le\frac{2\ellvtx}{\e}$,
\begin{align}
 b+c&\le\Llog,\qquad \abs{\ln\e}^b\bigl(1+\ln R\bigr)^c\le C\Lameps.
 \label{eq:paired-log-budget-expanded}
\end{align}
For $R\le3$, the regularized coefficients have finite $L^\infty_{Y,v}$ norms, so their factors $\abs{\ln\e}^b$ are also bounded by $C\Lameps$. We estimate the remaining powers separately for matched and intermediate coefficients.

\subparagraph{\underline{Matched coefficients, $\alpha<\Tmatch$}} The scaling identity is $\e^{\alpha-1}R^{\lambda_\Nexp}=\e^\Nexp r^{\lambda_\Nexp}$. Since $\Nexp>\Tmatch+\sdata$, we have $\lambda_\Nexp<\Tmatch-\Nexp-1<0$. The factor $\Xi_\alpha^\delta=1-\chi_m^\delta$ restricts this coefficient to $r\ge\delta>3\e$, where \eqref{eq:EI-degree-expanded}--\eqref{eq:EB-degree-expanded} apply. Also $\alpha\ge0$ gives $\lambda_\Nexp\ge-\Nexp-1$, and therefore
\begin{align}
 \e^\Nexp r^{\lambda_\Nexp} &\le \e^\Nexp\delta^{\lambda_\Nexp} \le \e^\Nexp\delta^{-\Nexp-1}.
 \label{eq:dnegative-outer-expanded}
\end{align}

With \eqref{eq:paired-log-budget-expanded}, every matched local contribution is $O(\Lameps\e^\Nexp\delta^{-\Nexp-1})$.

\subparagraph{\underline{Intermediate coefficients, $\Tmatch\le\alpha<\Tcon$}} Here $\Xi_\alpha^\delta=1$, since no wedge layer is added. If $\lambda_\Nexp<0$, polynomial decay dominates the logarithm for $R\ge3$ while the cutoff coefficient belongs to $L^\infty_{Y,v}$ on $R<3$. Including its prefactor gives a bound $C\Lameps\e^{\alpha-1}\le C\Lameps\e^{\Tmatch-1}$. If $\lambda_\Nexp\ge0$, then $\alpha\ge \Nexp+1$ and $r<2\ellvtx<1$, so
\begin{align}
 \e^{\alpha-1}R^{\lambda_\Nexp} &= \e^\Nexp r^{\lambda_\Nexp} \le C\e^\Nexp \le C\e^{\Tmatch-1}
 \label{eq:dpositive-expanded}
\end{align}
by $\Nexp>\Tmatch-1$. Every intermediate term is $O(\Lameps\e^{\Tmatch-1})$.

Summing the finite families in \eqref{eq:SIterm-expanded}--\eqref{eq:SBterm-expanded} therefore gives $O\big(\Lameps\{\e^\Nexp\delta^{-\Nexp-1}+\e^{\Tmatch-1}\}\big)$ for their local contributions.

The terms $S_{I,\rm far}^{[\Nexp]}$ and $S_{B,\rm far}^{[\Nexp]}$ contain $\psi_0$, with $\supp\psi_0\subset\Om_{\rm away}$. Where $\chi_{m,i}^{\rm end}=1$, a nonzero ray factor requires $0\le d_j<2\cfs\signedsideat{m}{i}$. The identity \eqref{eq:bulk-signed-side-coordinate-relation} therefore gives
\begin{align}
 \signedsideat{m}{i}&\ge\frac{\bulkdistat{m}}{\sqrt{1+4\cfs^2}}\ge\frac{\ellvtx}{\sqrt{1+4\cfs^2}}>0\qquad\text{on }\supp\psi_0.
 \label{eq:far-side-endpoint-separation}
\end{align}
The factors on the endpoint transition and on the middle of the side are supported a fixed positive tangential distance from that endpoint. Applying this argument at both endpoints keeps every active tangential coefficient away from its singularities.

The estimates for $\Uint_k$ and $\Uside_{j,k}$ in Proposition~\ref{prop:data-generated-hierarchy-comp}, using the local expansions on fixed annuli as well as the norms $\mathfrak F_{k,\kappa}$, consequently give $\bnm{S_{I,\rm far}^{[\Nexp]}}_\infty+\bnm{S_{B,\rm far}^{[\Nexp]}}_\infty\le C\e^\Nexp$, including the construction remainders. Their local contributions form $S_{\rm rem}^{[\Nexp]}$ and are treated in Step~3. Combining these estimates gives
\begin{align}
 \bnm{S_I^{[\Nexp]}}_\infty + \bnm{S_B^{[\Nexp]}}_\infty \le C\Lameps \left\{ \e^\Nexp\delta^{-\Nexp-1}+\e^{\Tmatch-1} \right\}.
 \label{eq:terminal-total-bound-expanded}
\end{align}

\paragraph{\underline{Step 3: Estimates of $S_{\rm nor}+S_{\rm vtx}^{\rm mid} +S_{\rm rem}^{[\Nexp]}+S_{\rm fwd}$}}

By \eqref{eq:wedge-cutoff-cancellation}, only $(\alpha,b)\in\mathfrak J_m^{\rm mid}$ contribute to $S_{\rm vtx}^{\rm mid}$. Their coefficients are independent of $\e$ and have finite $L^\infty_{Y,v}$ norms on $1\le\scaledbulkat{m}\le2$, so
\begin{align}
 \nm{S_{\rm vtx}^{\rm mid}}_\infty &\le C\sum_m \sum_{(\alpha,b)\in\mathfrak J_m^{\rm mid}} \e^{\alpha-1}\bigl(1+\abs{\ln\e}\bigr)^b \le C\Lameps\e^{\Tmatch-1}.
 \label{eq:unmatched-core-bound-expanded}
\end{align}
\subparagraph{\underline{Construction remainders}} The positive construction weight survives one derivative and controls the source throughout the chart.

On $\supp\psi_m$, where $0<s_i\le\bulkdistat{m}<2\ellvtx$ on the local ray-cutoff support, \eqref{eq:terminal-interior-remainder-comp} and \eqref{eq:terminal-regular-remainder-comp} with $q=1$ give $\e^\Nexp\nm{\nabla_x\mathcal R^U_{m,\Nexp}(x,\cdot)}_{L^\infty_w}\le C\e^\Nexp\bulkdistat{m}^{\pacc-1}$ and $\e^\Nexp\nm{\ue^{\kappa_\ast\eta}\partial_{s_i}\mathcal R^\Uside_{m,i,\Nexp}(s_i,\eta,\cdot)}_{L^\infty_{\eta,w}}\le C\e^\Nexp s_i^{\pacc-1}$. Both right-hand sides are at most $C\e^\Nexp$ because $\pacc\ge1$. Since $\psi_m$ is applied after $\Le$ and $\Nexp>\Tmatch-1$,
\begin{align}
 \bnm{S_{\rm rem}^{[\Nexp]}}_\infty \le C\e^\Nexp \le C\e^{\Tmatch-1}.
 \label{eq:Scon-bound-expanded}
\end{align}
The construction remainders on $\supp\psi_0$ contribute to $\psi_0\mathcal T_I^\e$, $\psi_0\mathcal T_B^\e$ and were estimated in \eqref{eq:terminal-total-bound-expanded}.

\subparagraph{\underline{Ray-cutoff commutators}} On the ray-cutoff transition, $\eta_{m,i}\simeq\sigma_{m,i}\simeq\scaledbulkat{m}$, so normal Milne decay absorbs the power-logarithmic factors. Thus \eqref{eq:completed-forward-error-bound} gives $\bnm{E^{\rm fwd}_{m,\alpha,b}(Y,\cdot)}_{L^\infty_v}\le C_L\jbr{\scaledbulkat{m}}^{-L}$ for every $L>0$.

For $\alpha<\Tmatch$ and $\scaledbulkat{m}<\frac{\delta}{\e}$, the factor $\Xi_{m,\alpha}^\delta$ vanishes. The terms with $\alpha\ge\Tmatch$ are $O(\Lameps\e^{\Tmatch-1})$. In the remaining region $\scaledbulkat{m}\ge\frac{\delta}{\e}>1$, taking $L>\Nexp+1$ and using $\alpha\ge0$ from \eqref{eq:generated-exponents-nonnegative} gives
\begin{align}
 \e^{\alpha-1}\abs{\ln\e}^b\bnm{\Xi_{m,\alpha}^\delta E^{\rm fwd}_{m,\alpha,b}(Y,\cdot)}_{L^\infty_v}&\le C\Lameps\e^{-1}\scaledbulkat{m}^{-\Nexp-1}\le C\Lameps\e^\Nexp\delta^{-\Nexp-1}.
\end{align}

We estimate the contribution of the construction remainder to this commutator directly. For $S_{\rm fwd}^{\rm rem}$, use the same support relation and $\abs{\nabla_{(\sigma_{m,i},\eta_{m,i})}\Theta}\le C\scaledbulkat{m}^{-1}$. Apply \eqref{eq:side-ledger-remainder-comp} with truncation exponent $\Tcon$, $q=0$ and $s_i=\e\sigma_{m,i}$. Abbreviating $\sigma=\sigma_{m,i}$ and $\eta=\eta_{m,i}$ gives, for $0\le k\le\Nexp$,
\begin{align}
 \e^{k-1} \abs{ \nabla_{(\sigma,\eta)} \Theta\!\left(\frac{\eta}{\sigma}\right) } \nm{\mathcal R^\Uside_{m,i,k}(\e\sigma,\eta,\cdot)}_{L^\infty_w} &\le C\e^{\Tcon-1} \sigma^{\Tcon-k-1} \bigl(1+\abs{\ln(\e\sigma)}\bigr)^{\Llog} \ue^{-\kappa_\ast\eta} \\
 &\le C\e^{\Tcon-1} \scaledbulkat{m}^{\Tcon-k-1} \bigl(1+\abs{\ln(\e\scaledbulkat{m})}\bigr)^{\Llog} \ue^{-c_0\scaledbulkat{m}} \le C\Lameps\e^{\Tcon-1},\notag
\end{align}
with $c_0>0$ depending only on $\kappa_\ast$ and the fixed ray cutoff. Here $\Tcon-k-1\ge\Tcon-\Nexp-1>\pacc>0$, so the positive power controls the logarithm at the origin and the exponential controls infinity. Since $\Tcon>\Tmatch$, this gives $\nm{S_{\rm fwd}^{\rm rem}}_\infty\le C\Lameps\e^{\Tcon-1}\le C\Lameps\e^{\Tmatch-1}$.

For the term in $S_{{\rm fwd},\rm far}$ containing a derivative of $\Theta$, $\psi_0\ne0$ keeps the point a fixed distance from the vertex. Together with $\eta_{m,i}\simeq\sigma_{m,i}$, this gives $\eta_{m,i}\ge\frac{c}{\e}$. The term containing $\bigl(\chi_{m,i}^{\rm end}\bigr)^{\prime}$ has $\signedsideat{m}{i}\simeq\ellend$ and $d_j\ge c\ellend$. Milne decay therefore gives $\nm{S_{{\rm fwd},\rm far}}_\infty\le C_L\e^L$ for every $L$. Altogether
\begin{align}
 \nm{S_{\rm fwd}}_\infty \le C\Lameps \left\{ \e^\Nexp\delta^{-\Nexp-1}+\e^{\Tmatch-1} \right\}.
 \label{eq:Sfwd-bound-expanded}
\end{align}

\subparagraph{\underline{Macroscopic normal-cutoff commutators}} On the support of $\zeta_j-1$ or of $\nabla\zeta_j$ we have $d_j\ge\ellnor$, hence $\eta_j=\frac{d_j}{\e}\ge\frac{\ellnor}{\e}$. We also need $L^\infty$ estimates for the tangential coefficients on this support. If the cutoff belonging to an endpoint $(m,i)\in\mathcal I(j)$ equals one, a nonzero extended side layer requires $\Theta\big(\frac{d_j}{\signedsideat{m}{i}}\big)\ne0$ and hence $\signedsideat{m}{i}\ge \frac{d_j}{2\cfs}\ge\frac{\ellnor}{2\cfs}$. The supports of $\bigl(\chi_{m,i}^{\rm end}\bigr)^{\prime}$ and of the portion where all endpoint cutoffs vanish are likewise a fixed tangential distance from the endpoints. The power-logarithmic coefficients and the derivatives used here therefore belong to $L^\infty$ on this support.

Equation \eqref{eq:side-profile-recursion-comp} controls $\mu_j(w)\partial_{\eta_j}\Uside_{j,k}+\qk\Uside_{j,k}$ as a whole at grazing velocities. Using it in $\Le\big[\sum_{k=0}^\Nexp\e^k\Uside_{j,k}\big]$ cancels consecutive collision and transport terms. The remaining terms have $L^\infty_{x,w}$ norms of order $\ue^{-\frac{\kappa_\ast\ellnor}{\e}}$ before any factor $\e^{-1}$ from $\nabla_x=\e^{-1}O_m\nabla_Y$. Allowing that factor, we obtain the bound
\begin{align}
 \nm{S_{\rm nor}}_\infty \le C\e^{-1}\ue^{-\kappa_\ast\frac{\ellnor}{\e}} \le C_L\e^L
 \label{eq:Snor-bound-expanded}
\end{align}
for every fixed $L>0$. Taking $L>\Nexp+2$ absorbs this term.

The preceding estimates give
\begin{align}
 \nm{S_I^{[\Nexp]}+S_B^{[\Nexp]}+S_{\rm nor}+S_{\rm vtx}^{\rm mid}+S_{\rm rem}^{[\Nexp]}+S_{\rm fwd}}_\infty &\le C\Lameps\big\{\e^\Nexp\delta^{-\Nexp-1}+\e^{\Tmatch-1}\big\}.
 \label{eq:six-source-terms-bound}
\end{align}

\paragraph{\underline{Step 4: commutator with $\chi_m^\delta$}}

By \eqref{eq:shrinking-cutoff-exact-product-rule}, the last term of \eqref{eq:source-pieces-expanded} contains the derivative of $\chi_m^\delta$. Since $\chi_m^\delta=\chi_{\rm ov}\left(\frac{\bulkdistat{m}}{\delta}\right)$, $\supp\nabla\chi_m^\delta\subset\{\delta\le\bulkdistat{m}\le2\delta\}$ and $\abs{\nabla\chi_m^\delta}\le C\delta^{-1}$. Thus the source is a derivative of size $\delta^{-1}$ multiplied by the wedge layer at wedge radius at least $\frac{\delta}{\e}$.

These annuli are disjoint because $2\delta<\ellvtx$ and satisfy $\scaledbulkat{m}\ge\frac{\delta}{\e}>3$ by \eqref{eq:realization-scale-chain}. The decay \eqref{eq:summed-wedge-decay} and $\abs w=1$ therefore bound the source by $C\delta^{-1}\Lameps\jbr{\frac{\bulkdistat{m}}{\e}}^{-\betaw}\le C\Lameps\delta^{-1}\left(\frac{\e}{\delta}\right)^\betaw$. With Step~3 this proves \eqref{eq:S-app-est}. The stability estimate \eqref{eq:physical-stability-recall} will multiply this bound by a further $\e^{-1}$, producing the overlap error balanced against the endpoint trace below.

\paragraph{\underline{Step 5: In-flow boundary trace}}

For measurable $A\subset\partial\Om$ put $\abs{f}_{\infty,-;A}:=\esssup\{\abs{f(x,w)}:(x,w)\in\gamma_-,\ x\in A\}$. On an open side $E_j$, $f\big|_{\gamma_{j,-}}(x_j(s),w)$ denotes the characteristic trace on the in-flow boundary at global side coordinate $s\in(0,L_j)$ and velocity $\mu_j(w)>0$. The auxiliary in-flow boundary data are $g_{\rm aux}^\e(x_j(s),w):=\sum_{k=0}^1\e^k\gaux_{j,k}(s,w)$, since \eqref{eq:auxiliary-inflow-coefficients} sets $\gaux_{j,k}=0$ for $k\ge2$.

Near an endpoint, the total wedge solution supplies the retained Taylor trace. Subtracting the matched sum leaves the intermediate coefficients and the construction remainders, whose boundary norms must be estimated separately. In the overlap transition, the prescribed wedge layer trace and the extension from the other ray vanish. In the middle of a side, the sum of the interior solution and the layer on side $j$ reproduces the auxiliary in-flow data, while any contribution from another side decays exponentially. Write $\dist_{{\rm end},j}(x_j(s)):=\min\{s,L_j-s\}$ for the distance along $E_j$ to its nearer endpoint. We partition $E_j$ into
\begin{align}
 A_{j,1}^\delta&:=\big\{x\in E_j:0<\dist_{{\rm end},j}(x)\le\delta\big\},\qquad A_{j,2}^\delta:=\big\{x\in E_j:\delta<\dist_{{\rm end},j}(x)<2\delta\big\},\\
 A_{j,3}^\delta&:=\big\{x\in E_j:\dist_{{\rm end},j}(x)\ge2\delta\big\}.
\end{align} Since $2\delta<\ellvtx$, every point of $A_{j,1}^\delta\cup A_{j,2}^\delta$ has a unique nearby endpoint pair $(m,i)\in\mathcal I(j)$.

The coefficient identity \eqref{eq:coefficientwise-incoming-cancellation}, obtained from \eqref{eq:own-side-order-trace} in Steps~6--7 of Lemma~\ref{lem:localized-wedge-data}, leaves only the Taylor coefficients of $g_{j,0}$ and $g_{j,1}$ at orders $k=0,1$ below $\Tmatch$. The remaining coefficients in the vertex expansions of $\Uint_k$ and the incident side layer cancel, giving \eqref{eq:summed-wedge-incoming-trace}.

\subparagraph{\underline{Region (a): the endpoint regions $A_{j,1}^\delta$}} Let $s=\sidedistat{m}{i}=\dist_{{\rm end},j}(x)$ for the unique nearby pair, so $x=x_j(s_{m,i}^{\rm ve}(s))$, and put $(\sigma,v)=\left(\frac{s}{\e},O_m^Tw\right)$. Here $\chi_m^\delta=1$ while every other wedge cutoff vanishes. Evaluating \eqref{eq:forward-composite-definition} in the characteristic trace on the in-flow boundary gives $\Uapp^{\e,[\Nexp]}(\delta)\big|_{\gamma_{j,-}}(x,w)=\Ubase^{\e,[\Nexp]}\big|_{\gamma_{j,-}}(x,w)-\mathcal M_m^\e\big|_{\gamma_{i,-}}(\sigma,v)+\Uwedge_m^\e\big|_{\gamma_{i,-}}(\sigma,v)$.

By \eqref{eq:Ubase-finite-expanded} the sum of the matched coefficients in $\Ubase^{\e,[\Nexp]}$ is exactly $\mathcal M_m^\e$, so their subtraction is exact. The remaining total wedge trace is prescribed by \eqref{eq:summed-wedge-incoming-trace}. Thus the trace on the in-flow boundary of $\Uapp^{\e,[\Nexp]}(\delta)$ consists of the retained endpoint Taylor polynomial and a complement,
\begin{align}
 \Uapp^{\e,[\Nexp]}(\delta)\big|_{\gamma_{j,-}}(x,w) &= \sum_{k=0}^1 \e^kJ_{m,i,k}^{\Tmatch}[\gaux](s,w) + \mathcal R_{m,i}^{\mathrm{tr},\e,\Tmatch}(s,w),
\end{align}
the complement being
\begin{align}
 \mathcal R_{m,i}^{\mathrm{tr},\e,\Tmatch}(s,w) &:= \sum_{(\alpha,b)\in\mathfrak J_m^{\rm mid}} \e^\alpha(\ln\e)^b \left. M^{\rm con}_{m,\alpha,b} \right|_{\gamma_{i,-}} \left(\frac{s}{\e},O_m^Tw\right)+ \mathcal R_{m,\rm con}^\e\big|_{\gamma_{j,-}}(x,w).
\end{align}
Each coefficient $M^{\rm con}_{m,\alpha,b}$ contains the interior and both side contributions. On ray $i$, the extension from the other side vanishes by \eqref{eq:Theta-choice}--\eqref{eq:forward-sector-extension}, and subtraction of $\mathcal M_m^\e$ cancels the full matched trace.

We estimate the retained Taylor polynomial $J_{m,i,k}^{\Tmatch}[\gaux]$ from \eqref{eq:taylor-jet-def}. For $0\le k\le1$, the relation $\Tmatch\in(\pacc+3,\pacc+4)$ gives $k+q_k^{\Tmatch}+1=\pacc+4>\Tmatch$. Taylor's formula in $L^\infty(\{w:\mu_{m,i}(w)>0\})$ and Assumption~\ref{ass:primitive-data-comp}, together with $s\le\delta$ and $\e\le\delta$, give
\begin{align}
 \e^k\nm{\gaux_{m,i,k}(s,\cdot)-J_{m,i,k}^{\Tmatch}[\gaux](s,\cdot)}_{L^\infty(\{\mu_{m,i}>0\})} &\le C\e^ks^{\pacc+4-k}\le C\delta^{\pacc+4}\le C\delta^{\Tmatch}.
 \label{eq:primitive-full-jet-remainder-bound-expanded}
\end{align}

For $(\alpha,b)\in\mathfrak J_m^{\rm mid}$ and $\frac{s}{\e}\ge2$, we have $M^{\rm con}_{m,\alpha,b}=\Cvertex_{m,\alpha,b}$. Before taking norms, we sum the binomial contributions from each term of the original vertex expansion. If the original logarithmic degree is $b_0$ and $\alpha=k+\lambda$, then $\e^\alpha\sum_{c=0}^{b_0}\binom{b_0}{c}(\ln\e)^{b_0-c}\left(\frac{s}{\e}\right)^\lambda\bigl(\ln\left(\frac{s}{\e}\right)\bigr)^c=\e^ks^\lambda(\ln s)^{b_0}$. This restores the logarithm of the boundary distance $s$. The positive gap above $\Tmatch$ then absorbs the logarithm for $s\le\delta$.

For $s<2\e$, the regularized trace on the in-flow boundary belongs to $L^\infty((0,2)\times\{v:\mu_i^Y(v)>0\})$ in wedge variables. These two regimes give
\begin{align}
 2\e\le s\le\delta:\qquad& \e^ks^{\alpha-k} \bigl(1+\abs{\ln s}\bigr)^{\Llog} \le s^\alpha\bigl(1+\abs{\ln s}\bigr)^{\Llog},
 \label{eq:trace-poly-sgeeps-expanded}\\
 0<s<2\e:\qquad& \e^\alpha\abs{\ln\e}^b\nm{\left.M^{\rm con}_{m,\alpha,b}\right|_{\gamma_{i,-}}\left(\frac{s}{\e},\cdot\right)}_{L^\infty(\{\mu_i^Y>0\})}\le C\e^\alpha\bigl(1+\abs{\ln\e}\bigr)^{\Llog}.
 \label{eq:trace-poly-sleeps-expanded}
\end{align}

Every intermediate exponent satisfies $\alpha>\Tmatch$ by \eqref{eq:generated-exponents-in-exceptional-set} and $\Tmatch\notin\mathscr E_\Nexp$. Since the index sets are finite, choose $\Delta_{\rm tr}>0$ such that $\alpha\ge\Tmatch+\Delta_{\rm tr}$ for every $(\alpha,b)\in\mathfrak J_m^{\rm mid}$. The inequality $\sup_{0<z<1}z^{\Delta_{\rm tr}}\bigl(1+\abs{\ln z}\bigr)^{\Llog}<\infty$ then bounds both \eqref{eq:trace-poly-sgeeps-expanded} and \eqref{eq:trace-poly-sleeps-expanded} by $C\delta^{\Tmatch}$.

At tangential distance $s$ on the ray, the sum of the $L^\infty$ norms in in-flow velocity of the construction remainders for $\Uint_k$ and $\Uside_{m,i,k}$, multiplied by $\e^k$, is at most $C\e^ks^{\Tcon-k}\bigl(1+\abs{\ln s}\bigr)^{\Llog}$ by \eqref{eq:outer-cartesian-remainder-comp} and \eqref{eq:side-ledger-remainder-comp}. We estimate this bound separately for $s\ge\e$ and $s\le\e$. For $s\ge\e$, the bound is $Cs^{\Tcon}\bigl(1+\abs{\ln s}\bigr)^{\Llog}$. For $s\le\e$, set $s=\e z$ and use $\Tcon-k>0$ to absorb $\abs{\ln z}$ into $z^{\Tcon-k}$. The resulting bound is $C\e^{\Tcon}\bigl(1+\abs{\ln\e}\bigr)^{\Llog}$.

Both bounds are $O(\delta^{\Tmatch})$ since $\Tcon>\Tmatch$ and $\e<\delta$. The remainder from the other side vanishes on this ray under the same extension. Thus the intermediate coefficients and construction remainders satisfy
\begin{align}
 \bnm{\mathcal R_{m,i}^{\mathrm{tr},\e,\Tmatch}}_{ L^\infty((0,\delta)\times \{w:\mu_{m,i}(w)>0\})} &\le C\delta^{\Tmatch}.
 \label{eq:trace-complement-bound-expanded}
\end{align}
The terms in \eqref{eq:trace-complement-bound-expanded} come from the approximate solution, while \eqref{eq:primitive-full-jet-remainder-bound-expanded} estimates the difference between its Taylor polynomial and $g_{\rm aux}^\e$. Adding these two bounds gives the discrepancy in the in-flow boundary condition on the endpoint region,
\begin{align}
 \abs{ \Uapp^{\e,[\Nexp]}(\delta) - g_{\rm aux}^\e }_{\infty,-;A_{j,1}^\delta} &\le C\delta^{\Tmatch}.
 \label{eq:trace-core-expanded}
\end{align}

\subparagraph{\underline{Region (b): the endpoint transition annuli $A_{j,2}^\delta$}} Here $s=\sidedistat{m}{i}>\delta>3\e$, so the regularized coefficients equal their unregularized counterparts and \eqref{eq:profile-boundary-telescope-comp} gives $g_{\rm aux}^\e$ from the sum of the interior solution and the layer on side $j$. Contributions from nonincident sides vanish by Lemma~\ref{lem:local-global-side-consistency}. The other incident side also contributes zero: its coordinates are $\sigma_{\rm opp}=\left(\frac{s}{\e}\right)\cos\omega_m$ and $\eta_{\rm opp}=\left(\frac{s}{\e}\right)\sin\omega_m$, so either $\sigma_{\rm opp}\le0$ or $\frac{\eta_{\rm opp}}{\sigma_{\rm opp}}=\tan\omega_m>2\cfs$. In both cases \eqref{eq:forward-sector-extension} vanishes.

The wedge layer satisfies the in-flow boundary condition \eqref{eq:defect-wedge-problem}. Since \eqref{eq:completed-ray-mismatch-split} gives $h_{m,i,\alpha,b}=0$ for $\sigma>2$ and $\frac{s}{\e}>3$, every matched wedge layer has zero trace on the in-flow boundary here. Hence
\begin{align}
 \abs{\Uapp^{\e,[\Nexp]}(\delta)-g_{\rm aux}^\e}_{\infty,-;A_{j,2}^\delta}&=0.
 \label{eq:trace-transition-expanded}
\end{align}

\subparagraph{\underline{Region (c): away from the endpoints, $A_{j,3}^\delta$}} Both wedge cutoffs vanish here. Since $2\delta>2\e$, the regularized coefficients equal their unregularized counterparts. Also $a_j^{\rm end}=\zeta_j=1$ on $d_j=0$, so the sum of the interior solution and the layer on side $j$ again gives $g_{\rm aux}^\e$. A different side sharing the endpoint $V_m$ has normal coordinate $d_{\rm opp}(x)=\sidedistat{m}{i}(x)\sin\omega_m\ge2\delta\sin\omega_m$. Its contribution vanishes for $\bulkdistat{m}<2\ellvtx$ by the local side identity and the ray-cutoff condition. On the remaining support, the endpoint separation argument of Step~2 bounds its tangential coefficients and remainders at both endpoints independently of $\e$. For each such side $E_\ell$, Milne decay bounds $\Ext_\ell^\e\big[\sum_{k=0}^\Nexp\e^k\Uside_{\ell,k}\big]$ by $C\exp\big(-\frac{c\dist_{{\rm end},j}(x)}{\e}\big)$.

For a side sharing neither endpoint of $E_j$, the signed distance to its supporting line is affine along $E_j$. By \eqref{eq:nonincident-side-separation} it is at least $\ellsep$ at both endpoints, hence throughout $E_j$. Its normal cutoff therefore vanishes because $\ellsep>2\ellnor$. Summing the at most two incident-side contributions gives
\begin{align}
 \abs{\Uapp^{\e,[\Nexp]}(\delta)-g_{\rm aux}^\e}_{\infty,-;A_{j,3}^\delta}&\le C\Lameps\ue^{-\frac{2c\delta}{\e}}.
 \label{eq:trace-middle-expanded}
\end{align}

\subparagraph{\underline{Return to the prescribed in-flow data}} The three regions partition each open side, and the vertices have zero boundary measure. Estimates \eqref{eq:trace-core-expanded}, \eqref{eq:trace-transition-expanded}, and \eqref{eq:trace-middle-expanded} give $C\delta^{\Tmatch}+C\Lameps\ue^{-\frac{c\delta}{\e}}$. Adding the nonnegative term $C\Lameps\bigl(\frac{\e}{\delta}\bigr)^\betaw$ gives the form stated in the proposition:
\begin{align}
 \abs{\Uapp^{\e,[\Nexp]}(\delta)-g_{\rm aux}^\e}_{\infty,-} \le C\left\{\delta^{\Tmatch} +\Lameps\Big[\left(\frac{\e}{\delta}\right)^\betaw+\ue^{-\frac{c\delta}{\e}}\Big]\right\}.
\end{align}
Assumption~\ref{ass:primitive-data-comp} and \eqref{eq:data-expansion-comp} also give
\begin{align}
 \max_j\Bnm{g_j^\e-\sum_{k=0}^1\e^k\gaux_{j,k}} _{L^\infty((0,L_j)\times\{w:\mu_j(w)>0\})} \le C\e^2.
\end{align}
This last bound is the only change from the auxiliary in-flow data to the prescribed data. The triangle inequality now proves \eqref{eq:trace-est}.

Finally, inserting the source bound \eqref{eq:S-app-est} and the boundary bound \eqref{eq:trace-est} into \eqref{eq:physical-stability-recall} gives \eqref{eq:preoptimized-error}.
\end{proof}

All derivative orders used in the proof are at most $\Kreg$ in \eqref{eq:Kreg-choice-comp}.

\begin{remark}[Role of wedge layer]
\label{rem:source-mechanisms}
For each $\alpha<\Tmatch$, matching the sums in \eqref{eq:raw-coefficient-blocks-v6}--\eqref{eq:raw-side-coefficient-sum} cancels the corresponding source on $\bulkdistat{m}<\delta$. It also supplies the retained Taylor coefficients of the in-flow boundary data. In $L^\infty(\Om\times\Sone)$, the unmatched source is $O(\Lameps\e^{\Tmatch-1})$ and the remaining matched order-$\Nexp$ source is $O(\Lameps\e^\Nexp\delta^{-\Nexp-1})$. The endpoint trace error is $O(\delta^{\Tmatch})$ in the boundary norm defined in Step~5. The stability estimate \eqref{eq:physical-stability-recall} multiplies the source bounds by $\e^{-1}$.
\end{remark}

\subsection{Proof of \texorpdfstring{$L^\infty$}{L-infinity} Theorem}

\begin{proof}[Proof of Theorem~\ref{thm:main-comp}]\leavevmode
\paragraph{\underline{Step 1: Choice of the overlap radius}} Put $\delta=\e^a$, $0<a<1$. The endpoint trace discrepancy $\delta^{\Tmatch}$ has exponent $a\Tmatch$, while the wedge-cutoff term after applying \eqref{eq:physical-stability-recall} $\e^{\betaw-1}\delta^{-\betaw-1}=\e^{\betaw-1-a(\betaw+1)}$ has exponent $\betaw-1-a(\betaw+1)$. Increasing $a$ improves the endpoint accuracy by shrinking the endpoint region. It also moves the cutoff toward the wedge vertex, where the wedge layer is larger, and increases the cutoff derivative.

The first exponent increases and the second decreases in $a$, so equating them, $a\Tmatch=\betaw-1-a(\betaw+1)$, maximizes the smaller of the two. This recovers the parameters fixed in \eqref{eq:overlap-rate-parameters-comp}, namely $a=\aoverlap$ and $\Tmatch \aoverlap=\qrate$. We next verify that the other errors in \eqref{eq:preoptimized-error} are smaller at this choice.

Since $\betaw>1$, we have $0<\aoverlap<1$. Set $\delta=\deleps=\e^{\aoverlap}$ as in \eqref{eq:delta-lambda-comp}. Then $\deleps\to0$ and $\frac{\deleps}{\e}=\e^{\aoverlap-1}\to\infty$, so $3\e<\deleps$ and $2\deleps<\ellvtx$ for sufficiently small $\e$. Thus \eqref{eq:realization-scale-chain} hold.

\paragraph{\underline{Step 2: Comparison of the error terms}} The balancing identity gives
\begin{align}
 \deleps^{\Tmatch}&=\e^{\betaw-1}\deleps^{-\betaw-1}=\e^{\qrate}.
\end{align}
The order-$\Nexp$ source and the wedge trace have strictly larger exponents:
\begin{align}
 (1-\aoverlap)\Nexp-1-\aoverlap-\qrate&=(1-\aoverlap)(\Nexp-\betaw)>0,\qquad \betaw(1-\aoverlap)-\qrate=1+\aoverlap>0.
 \label{eq:terminal-rate-gap-expanded}
\end{align}
The first identity is \eqref{eq:N-rate-redundancy-comp}, and $\Nexp>\betaw$ follows from \eqref{eq:N-choice-comp}. The remaining algebraic errors have exponents $\Tmatch-2$ and $2$, both larger than $\qrate$, since $\qrate<\betaw-1<2$ and $\Tmatch-2>\pacc+1\ge2$. Finally, $\aoverlap<1$ gives, for every $L>0$,
\begin{align}
 \ue^{-\frac{c\deleps}{\e}}&=\ue^{-c\e^{\aoverlap-1}}=o(\e^L).
 \label{eq:superalg-term}
\end{align}

Substituting these comparisons into \eqref{eq:preoptimized-error} and absorbing the wedge-trace and superalgebraic terms into $\e^{\qrate}$,
\begin{align}
 \nm{\Err^\e(\deleps)}_\infty \le C\Lameps \left[ \e^{\qrate} + \e^{(1-\aoverlap)\Nexp-1-\aoverlap} + \e^{\Tmatch-2} \right] + C\e^2,
\end{align}
which, as $\Err^\e(\deleps)=\Uphys^\e-\Uapp^{\e,[\Nexp]}(\deleps)$, is exactly \eqref{eq:main-estimate-comp}.

For every $0<q<\qrate$, the positive gap absorbs the logarithm, since $\sup_{0<\e\le1}\e^{\qrate-q}\bigl(1+\abs{\ln\e}\bigr)^{\Llog}<\infty$. The other gaps are larger and $\e^2\le\e^q$ because $q<\qrate<2$. This proves \eqref{eq:strict-rate-comp}.

\end{proof}

\begin{proof}[Proof of Corollary~\ref{cor:interior-convergence}]

Let $K_0\Subset\Om$ and $d_{K_0}:=\dist(K_0,\partial\Om)>0$. Since $\dist(K_0,V_m)\ge d_{K_0}$, all wedge cutoffs vanish on $K_0$ once $2\deleps<d_{K_0}$.

For every nonzero extended side layer on $K_0$, $\eta_j=\frac{d_j(x)}{\e}\ge \frac{d_{K_0}}{\e}$. The endpoint and ray cutoffs exclude a simultaneous endpoint singularity. Indeed, near an endpoint an active ray cutoff gives $\signedsideat{m}{i}\ge\frac{d_j(x)}{2\cfs}\ge\frac{d_{K_0}}{2\cfs}$, while the endpoint transition and the middle of a side already have fixed positive tangential distance from that endpoint. The finite sum of remaining power-logarithmic coefficients has $L^\infty$ norm at most $C\e^{-A}$ for some fixed $A>0$, so Milne decay gives
\begin{align}
 \Bnm{\sum_j\Ext_j^\e\Big[\sum_{k=0}^\Nexp\e^k\Uside_{j,k}\Big]}_{L^\infty(K_0\times\Sone)}&\le C\e^{-A}\ue^{-\frac{cd_{K_0}}{\e}}=o(1).
\end{align}

For small $\e$, \eqref{eq:Ubase-finite-expanded} gives $\Ubase^{\e,[\Nexp]}=\mathscr H_\Nexp^\e$ on $K_0$. Since $\Uint_k\in L^\infty(K_0\times\Sone)$ and $\Uint_0=\rho_0$, it follows that $\Uapp^{\e,[\Nexp]}(\deleps)=\rho_0+\sum_{k=1}^\Nexp\e^k\Uint_k+o(1)=\rho_0+o(1)$ in $L^\infty(K_0\times\Sone)$. The estimate \eqref{eq:strict-rate-comp} from Theorem~\ref{thm:main-comp} and the triangle inequality prove \eqref{eq:interior-limit-comp}.
\end{proof}

\begin{remark}
\label{rem:order-two-source-failure}
In smooth domains it is often enough to truncate the interior expansion at order $\Nexp=2$. A polygonal vertex mode shows why the source of an order-two interior truncation need not belong to $L^\infty$ near a vertex, already due to $\Uint_0=\rho_0$.

\begin{enumerate}
\item Suppose $\rho_0$ contains $cZ^+_{m,n}$ with $c\ne0$, where $Z^+_{m,n}$ is a zero-trace harmonic function with positive radial exponent from \eqref{eq:intro-zero-ray-pairs-comp} at a vertex of opening angle $\omega$. Write it as $Z_\lambda(\bulkdist,\theta)=\bulkdist^\lambda\sin(\lambda\theta)$, with $\lambda=\frac{n\pi}{\omega}$ and $n\in\mathbb N$. Assume $\lambda<3$ and $\lambda\ne2$. Convexity forces $\lambda>1$, so the range is $1<\lambda<3$ with $\lambda\ne2$. Such exponents occur for convex polygons, for example $\omega=\frac{2\pi}{3}$ gives $\lambda=\frac{3}{2}$. The global Dirichlet problem can determine a nonzero coefficient of this mode.

\item Apply \eqref{eq:microscopic-recursion-comp}--\eqref{eq:total-interior-recursion-comp}, writing $\bullet^{[H]}$ for the interior solutions in the Hilbert expansion. Homogeneity of noninteger degree below three gives $\abs{\nabla^3Z_\lambda}\simeq\bulkdist^{\lambda-3}$ on a nonempty angular subsector. The first three interior solutions are $\Uint_0^{[H]}=\rho_0$, $\Uint_1^{[H]}=-(w\cdot\nabla_x)\rho_0+\rho_1^{[H]}$, and $\Uint_2^{[H]}=\bigl(w\cdot\nabla_x\bigr)^2\rho_0-(w\cdot\nabla_x)\rho_1^{[H]}+\rho_2^{[H]}$. Their order-two source contains $c\e^2\bigl(w\cdot\nabla_x\bigr)^3Z_\lambda$.

Since $Z_\lambda$ is harmonic, $D^3Z_\lambda$ is trace-free, and its contraction with $w^{\otimes3}$ is a third angular harmonic in velocity. Let $O_m^Tw=(\cos\varphi,\sin\varphi)$ in the same local Cartesian frame that defines the spatial angle $\theta$. Direct differentiation gives
\begin{align}
 \bigl(w\cdot\nabla_x\bigr)^3Z_\lambda &=\lambda(\lambda-1)(\lambda-2)\bulkdist^{\lambda-3}\sin\bigl(3\varphi+(\lambda-3)\theta\bigr).
\end{align}
The $\rho_1^{[H]}$ and $\rho_2^{[H]}$ contributions have velocity degrees at most two and one, so they cannot cancel this nonzero third harmonic. The construction remainder satisfies $D_x^3\mathcal R^U_{m,0}\in L^\infty$ on the vertex chart by \eqref{eq:outer-cartesian-remainder-comp} and $\Tcon>3$. After three derivatives, the finite expansion contains distinct radial powers and their logarithmic companions. Lemma~\ref{lem:finite-polyhom-uniqueness} shows that membership of the projected source in $L^\infty_{x,w}$ would force every negative-degree coefficient to vanish, contradicting the coefficient $c$. Thus the order-two interior source does not belong to $L^\infty_{x,w}$ on any neighborhood of $V_m$. For $\bulkdist=\e\scaledbulk$, the component with coefficient $c$ has $L^\infty_w$ norm proportional to $\e^{\lambda-1}\scaledbulk^{\lambda-3}$. The excluded exponent $\lambda=2$ gives a harmonic polynomial whose third derivatives vanish.

\item For each matched exponent, \eqref{eq:raw-coefficient-blocks-v6}--\eqref{eq:raw-side-coefficient-sum} retain all terms through the fixed order $\Nexp$. Adding $D_{m,\alpha,b}$ cancels the source of the corresponding $M_{m,\alpha,b}$ on $\bulkdist<\deleps$. The remaining matched order-$\Nexp$ source is supported in $\bulkdist\ge\deleps$, where Proposition~\ref{prop:full-residual} bounds it by $C\Lameps\e^\Nexp\deleps^{-\Nexp-1}$. Its exponent after applying \eqref{eq:physical-stability-recall} exceeds $\qrate$ by \eqref{eq:terminal-rate-gap-expanded}. Steps~2--3 of that proposition give the separate $L^\infty$ source estimates for the intermediate coefficients and construction remainders.
\end{enumerate}
\end{remark}
\appendix
\section{\texorpdfstring{$L^2$}{L2} Diffusive Limit in General Polygonal Domains}
\label{sec:velocity-dependent-l2}

We prove Theorem~\ref{thm:vd-l2-wedge-free} under Assumption~\ref{ass:vd-l2-data} on a polygon satisfying \eqref{eq:vd-l2-angle-assumption}. The approximate solution is the sum of the harmonic function $\rho_0$ and the leading side layer, cut off within distance $O(\e)$ of each endpoint. The cutoff introduces an error in the in-flow boundary condition with $O(1)$ $L^\infty(\gamma_-)$ norm on intervals of total length $O(\e)$. Its kinetic $L^2$ norm is therefore $O(\e^{\frac12})$. We will also prove that the source term of the cutoff side layer is $O(\e)$ in $L^1(\Om\times\Sone)$. In the energy identity, this source pairs with the remainder, whose $L^\infty$ norm is bounded uniformly in $\e$. The resulting bound, together with the factor $\e^{-1}$ in the collision term, gives an $O(\e)$ estimate for the microscopic part in $L^2$.

For the kernel estimate, let $\phi^\e$ solve the Dirichlet Poisson problem \eqref{eq:vd-l2-poisson-remainder}, and set $f^\e=\Uphys^\e-\rho_0$. Testing its equation with $w\cdot\nabla_x\phi^\e$ introduces the term $\e^{-1}\int_{\Om\times\Sone}\qk f^\e(w\cdot\nabla_x\phi^\e)\ud x\ud\nu$. The identity obtained with the test function $\phi^\e$, multiplied by $\e^{-1}$, cancels this term exactly. This combination of test identities, as in \cite{Guo.Wu2025}, removes the factor $\e^{-1}$ from the microscopic contribution to the kernel estimate. The source involving $\rho_0$ also vanishes by harmonicity. We then estimate the Poisson Hessian in $\Om$ and its gradient on $\partial\Om$, taking account of the singularities at reentrant vertices. The second Milne moment \eqref{eq:vd-l2-localized-second-moment} removes the side layer from the boundary integral, leaving the remainder traces controlled by the energy estimate and the in-flow boundary data.

We use the equation \eqref{eq:physical-model}, the projections $\pk,\qk$, and the coordinates \eqref{eq:local-global-side-coordinate}--\eqref{eq:endpoint-pullbacks}. All geometric constants are chosen for the present polygon. The leading in-flow boundary data $g_{j,0}$ and the norm $\mathcal G_{\rm vd}$ are given in Assumption~\ref{ass:vd-l2-data}. For each $s$, the Milne end state is $e_{j,0}(s)=\mathcal E_{\mu_j}(g_{j,0}(s,\cdot),0)$, where $\mathcal E_{\mu_j}$ is defined in \eqref{eq:milne-maps-comp}. It determines the boundary value $\rho_0\big|_{E_j}=e_{j,0}$. The harmonic function $\rho_0\in H^1(\Om)\cap L^\infty(\Om)$ satisfies \eqref{eq:vd-l2-harmonic-limit}--\eqref{eq:vd-l2-harmonic-estimate}. We use $d_j$ only in the collar \eqref{eq:vd-l2-one-sided-collar}, where $d_j>0$.

We first define the cutoff side layer and estimate its source term. We then prove the elliptic estimates needed for the remainder equation. Until the last step of the proof, assume $0<\e\le\e_0(\Om)$ so that the two endpoint cutoff intervals on each side are disjoint. The case $\e_0(\Om)\le\e\le1$ follows from the energy and comparison estimates.

Section~\ref{sec:physical-stability} applies to the present polygon because each characteristic segment before the first exit has length at most $\operatorname{diam}(\Om)$. Use the boundaries $\gamma_\pm$ and the kinetic measure and norms in \eqref{eq:l2-linfty-boundary-norms-comp}, with $p=2$. Set $\gamma=\gamma_-\cup\gamma_+$. The vertices and grazing set have zero kinetic boundary measure. Define
\begin{align}
 \mathscr W_{\rm kin}^2 :=\big\{F\in L^2(\Om\times\Sone):{}&w\cdot\nabla_xF\in L^2(\Om\times\Sone), \; F\big|_{\gamma_\pm}\!\!\in L^2(\gamma_\pm;\ud\gamma_{\rm kin})\big\}.
 \label{eq:vd-l2-graph-space}
\end{align}
In \eqref{eq:vd-l2-graph-space}, the boundary values are the one-sided limits along almost every characteristic interval. For the kinetic solution, they are the mild traces of Section~\ref{sec:physical-stability}. All boundary integrals below use these traces. Since a line can intersect a nonconvex polygon in several intervals, Green's identity includes both endpoints of each interval.

For $H^1(\Om)$ functions, these characteristic limits agree almost everywhere with the Sobolev traces. The same holds for the Poisson gradient locally on each open side, by the $H^2$ regularity away from vertices in \eqref{eq:vd-l2-poisson-singular-decomposition}. Its piecewise $L^2$ trace is bounded by \eqref{eq:vd-l2-poisson-decomposition-bound}. The traces of the cutoff side layers at $d=0$ follow from their definition in Subsection~\ref{subsec:vd-l2-side-localization}.

\begin{lemma}
\label{lem:vd-l2-green-identity}
For every $F,K\in\mathscr W_{\rm kin}^2$,
\begin{align}
 &\int_{\Om\times\Sone} \big(w\cdot\nabla_xF\big)K\,\ud x\ud\nu + \int_{\Om\times\Sone} F\big(w\cdot\nabla_xK\big)\,\ud x\ud\nu = \int_{\gamma_+}FK\,\ud\gamma_{\rm kin} - \int_{\gamma_-}FK\,\ud\gamma_{\rm kin}.
 \label{eq:vd-l2-green-identity}
\end{align}
\end{lemma}

\begin{proof}
For almost every $w$, write $x=y+tw$ with $y\in w^\perp$. For almost every $y$, the set $\{t:y+tw\in\Om\}$ is a finite union of disjoint open intervals. By Fubini's theorem, the distributional streaming derivatives restrict to the weak derivatives on these intervals. Hence $F$ and $K$ belong to $H^1$ in $t$ and have absolutely continuous representatives, with derivatives $w\cdot\nabla_xF$ and $w\cdot\nabla_xK$. Their endpoint values are the traces in \eqref{eq:vd-l2-graph-space}.

Integrate by parts on each interval, sum the endpoint terms, and integrate in $y$ and $w$. The projection of a side onto $w^\perp$ has Jacobian $\abs{w\cdot n}$, so the resulting boundary integral is $\int_{\partial\Om\times\Sone}FK(w\cdot n)\ud\sigma\ud\nu$. The right endpoints lie on the out-flow boundary and the left endpoints on the in-flow boundary. Splitting by the sign of $w\cdot n$ proves \eqref{eq:vd-l2-green-identity}.
\end{proof}

\subsection{Construction of the Side Layer}
\label{subsec:vd-l2-side-localization}

Fix $0<\kappaLtwo<1$. By Proposition~\ref{prop:flat-milne-comp} with $\kappa_{\rm out}=\kappaLtwo$, $\kappa_{\rm in}=\frac{1+\kappaLtwo}{2}$, and zero source, the leading side layer is
\begin{align}
 \Uside_{j,0}(s,\eta,w) &:={} \mathfrak M_{\mu_j} \bigl(g_{j,0}(s,\cdot),0\bigr)(\eta,w),
 \label{eq:vd-l2-leading-side-profile}
\end{align}
and satisfies
\begin{align}
 \big(\mu_j\partial_\eta+\qk\big)\Uside_{j,0} &=0,
 \label{eq:vd-l2-leading-side-equation}\\
 \Uside_{j,0}\big|_{\eta=0}(s,w)&=g_{j,0}(s,w)-e_{j,0}(s),\qquad\mu_j(w)>0,
 \label{eq:vd-l2-leading-side-incoming}\\
 \nm{\Uside_{j,0}(s,\eta,\cdot)}_{L^\infty(\Sone)}&\rightarrow0\quad\text{as }\eta\to\infty.
 \label{eq:vd-l2-leading-side-limit}
\end{align}
Here $\Uside_{j,0}$ is the Milne solution minus its end state $e_{j,0}(s)$. The limit at infinity is taken in $L^\infty(\Sone)$ for the characteristic representative in \eqref{eq:milne-char-plus-comp}--\eqref{eq:milne-char-minus-comp}. By the Milne estimate \eqref{eq:milne-mapping-estimate-comp}, the commutation identity \eqref{eq:milne-parameter-commutation-comp}, and the regularity \eqref{eq:vd-l2-leading-regularity}, we have, for every $s\in[0,L_j]$ and almost every $\eta>0$,
\begin{align}
 \sum_{a=0}^1\bnm{\partial_s^a\Uside_{j,0}(s,\eta,\cdot)}_{L^\infty(\Sone)} &\leq C\mathcal G_{\rm vd}\ue^{-\kappaLtwo\eta}.
 \label{eq:vd-l2-leading-side-estimate}
\end{align}
At $\eta=0$, the in-flow boundary condition \eqref{eq:vd-l2-leading-side-incoming} holds in $L^\infty(\{w:\mu_j(w)>0\},\ud\nu)$. The trace for $\mu_j(w)<0$ is given by \eqref{eq:milne-char-minus-comp}. Together with \eqref{eq:vd-l2-leading-side-estimate}, this gives
\begin{align}
 \bnm{\Uside_{j,0}\big|_{\eta=0}(s,\cdot)}_{L^\infty(\Sone)}&\leq C\mathcal G_{\rm vd}.
 \label{eq:vd-l2-leading-side-trace}
\end{align}

All decay rates, cutoffs, and auxiliary exponents in this appendix are fixed independently of $\e$. Their dependence is included in $C_\Om$.

Integrating the flat Milne equation gives the two moment identities below. The first is conservation of normal flux. The second removes the side layer from the boundary term in the kernel estimate.

\begin{lemma}
\label{lem:vd-l2-milne-moments}
For every side $E_j$, $s\in[0,L_j]$, and $\eta\geq0$,
\begin{align}
 \int_{\Sone}\mu_j(w)^a\,\Uside_{j,0}(s,\eta,w)\,\ud\nu(w)&=0,\qquad a\in\{1,2\}.
 \label{eq:vd-l2-milne-moments}
\end{align}
Both identities remain valid after multiplication by a spatial cutoff independent of $w$.
\end{lemma}

\begin{proof}
Fix $j,s$ and write $\mu=\mu_j(w)$. We first derive the identities in the sense of distributions in $\eta$. The Milne equation and \eqref{eq:vd-l2-leading-side-estimate} give absolutely continuous representatives of the moments. At $\eta=0$, we use the values determined by the characteristic traces.

Integrating the Milne equation in $w$ and using $\int_{\Sone}\qk\,\Uside_{j,0}\,\ud\nu=0$, we obtain
\begin{align}
 \frac{\ud}{\ud\eta} \int_{\Sone} \mu\,\Uside_{j,0}(\eta,w)\,\ud\nu(w) &=0.
 \label{eq:vd-l2-flux-constant}
\end{align}
Thus the normal flux is constant in $\eta$. By the decay estimate \eqref{eq:vd-l2-leading-side-estimate}, this constant is zero.

Multiplying the Milne equation by $\mu$ and integrating in $w$, we use $\int_{\Sone}\mu\,\ud\nu=0$ and the first identity to obtain
\begin{align}
 \frac{\ud}{\ud\eta}\int_{\Sone}\mu^2\,\Uside_{j,0}(\eta,w)\,\ud\nu(w) =-\int_{\Sone}\mu\,\qk\,\Uside_{j,0}\,\ud\nu =-\int_{\Sone}\mu\,\Uside_{j,0}\,\ud\nu=0.
 \label{eq:vd-l2-second-moment-constant}
\end{align}
The second moment is therefore constant and vanishes by decay at infinity. Any cutoff independent of $w$ factors out of both integrals and preserves the identities.
\end{proof}

We define the side layer near $E_j$ by multiplying $\Uside_{j,0}\bigl(s,\frac d\e,w\bigr)$ by a cutoff $a_{j,\e}(s,d)$ and extending by zero. By \eqref{eq:vd-l2-leading-side-estimate}, this function decays exponentially in $\frac d\e$. We choose the support inside a collar that narrows at each endpoint. Then the layer lies inside $\Om$ and has zero trace on every other side.

The collar cutoffs produce derivatives of size $O(s^{-1})$ and $O\big((L_j-s)^{-1}\big)$. Their derivatives are supported where the normal distance is comparable to the distance from an endpoint, so the Milne decay controls these terms. We also cut off the layer within tangential distance $O(\e)$ of each endpoint. The derivatives of this endpoint cutoff are $O(\e^{-1})$ on regions of area $O(\e^2)$. These two observations give an $O(\e)$ source estimate in $L^1$ and an $O(1)$ estimate in $L^2$. The first will be used in the energy inequality, and the second justifies Green's identity for the cutoff layer.

Choose $\chi_{\rm end}\in C^\infty([0,\infty);[0,1])$ with $\chi_{\rm end}=0$ on $[0,1]$ and $\chi_{\rm end}=1$ on $[2,\infty)$. We choose the constants $\cfs$ and $\ellnor$ using the following geometric lemma.

\begin{lemma}
\label{lem:vd-l2-separated-side-collars}
There exist $\cfs,\ellnor>0$, depending only on $\Om$, such that for every side $E_j$ the collar
\begin{align}
 \mathsf C_j^{\rm col} &:={} \big\{ x_j(s)-dn_j: 0<s<L_j,\; 0<d<2\ellnor,\; d<2\cfs s,\; d<2\cfs(L_j-s) \big\}
 \label{eq:vd-l2-one-sided-collar}
\end{align}
is contained in $\Om$ and satisfies
\begin{align}
 \overline{\mathsf C_j^{\rm col}}\cap\partial\Om &\subset \overline E_j.
 \label{eq:vd-l2-collar-boundary-separation}
\end{align}
\end{lemma}

\begin{proof}
Let $V_m=x_j(0)$. Near $V_m$, the polygon agrees with the sector of opening $\omega_m$, and $-n_j$ points into the sector along $E_j$. Since the other ray makes a positive angle with $E_j$, there exist $c_{m,j},s_{m,j}>0$ such that $x_j(s)-dn_j\in\Om$ for $0<s<s_{m,j}$ and $0<d<2c_{m,j}s$. The closure of this cone meets the adjacent side only at $V_m$. This argument holds for every $0<\omega_m<2\pi$, including a reentrant vertex. At $x_j(L_j)$, the same argument applies with $L_j-s$.

The remaining middle portion of $E_j$ is compact and has positive distance from the other sides. A sufficiently thin normal strip over this portion therefore lies in $\Om$ and meets its boundary only on $E_j$. Taking the minimum of the finitely many allowable constants gives $\cfs,\ellnor$ satisfying \eqref{eq:vd-l2-one-sided-collar}--\eqref{eq:vd-l2-collar-boundary-separation}.
\end{proof}

Fix $\cfs,\ellnor$ as in Lemma~\ref{lem:vd-l2-separated-side-collars}, and choose $\Theta$ satisfying \eqref{eq:Theta-choice} with this value of $\cfs$. Decrease $\ellvtx$ so that $B_{2\ellvtx}(V_m)$ meets only the collars of the two incident sides. This is possible by \eqref{eq:vd-l2-collar-boundary-separation}, since $V_m$ has positive distance from the closed collars of every other side.

Choose $\zeta_j\in C_c^\infty([0,\infty);[0,1])$ with $\zeta_j=1$ on $[0,\ellnor]$ and $\zeta_j=0$ on $[2\ellnor,\infty)$. Fix $0<\e_0(\Om)\le\min\{1,\ellvtx\}$ so that $4\e_0<\min_jL_j$. For $0<s<L_j$ and $d>0$, set
\begin{align}
 &a_{j,\e}(s,d) :={} \chi_{\rm end}\left(\frac{s}{\e}\right) \chi_{\rm end}\Big(\frac{L_j-s}{\e}\Big) \zeta_j(d)\, \Theta\Big(\frac{d}{s}\Big) \Theta\Big(\frac{d}{L_j-s}\Big).
 \label{eq:vd-l2-side-cutoff}
\end{align}
By the supports of $\zeta_j$ and $\Theta$, the cutoff lies inside $\mathsf C_j^{\rm col}$. Define
\begin{align}
 &\UsideZeroAt{j}^\e(x,w) :={} \begin{cases}
 a_{j,\e}(s,d)\, \Uside_{j,0}\left(s,\frac{d}{\e},w\right), &x=x_j(s)-dn_j\in\mathsf C_j^{\rm col},\\
 0,&x\notin\mathsf C_j^{\rm col},
 \end{cases}.
 \label{eq:vd-l2-localized-side-layer}
\end{align}
Their sum is $\UsideZero^\e$ in \eqref{eq:vd-l2-global-leading-layer}. By Lemma~\ref{lem:vd-l2-separated-side-collars}, the trace on every other open side vanishes:
\begin{align}
 &\UsideZeroAt{j}^\e\big|_{\gamma}=0 \quad\hbox{on }E_i\times\Sone \qquad(i\neq j).
 \label{eq:vd-l2-cross-side-trace-zero}
\end{align}

The cutoff is zero near the tangential endpoints and vanishes smoothly at the other artificial interfaces. Thus extending by zero preserves the transport equation in the sense of distributions. Since all cutoffs are independent of $w$, the two Milne moment identities remain valid.

\begin{lemma}
\label{lem:vd-l2-localized-leading-layer}
Let $S_{\rm sl}^\e=\Le\UsideZero^\e$. The cutoff side layer has the following properties.
\begin{enumerate}
\item For $0<\e\le1$ and $1\le q<\infty$,
\begin{align}
 \nm{\UsideZero^\e}_{L^\infty(\Om\times\Sone)} &\leq C_\Om\mathcal G_{\rm vd},
 \label{eq:vd-l2-layer-Linfty}\\
 \nm{\UsideZero^\e}_{L^q(\Om\times\Sone)} &\leq C_{\Om,q}\mathcal G_{\rm vd}\e^{\frac{1}{q}}.
 \label{eq:vd-l2-layer-Lr}
\end{align}
\item For $0<\e\le\e_0(\Om)$, the source term belongs to $L^1(\Om\times\Sone)\cap L^2(\Om\times\Sone)$ and satisfies
\begin{align}
 \nm{S_{\rm sl}^\e}_{L^1(\Om\times\Sone)} &\leq C_\Om\mathcal G_{\rm vd}\e,
 \label{eq:vd-l2-layer-residual-norms}\\
 \nm{S_{\rm sl}^\e}_{L^2(\Om\times\Sone)} &\leq C_\Om\mathcal G_{\rm vd}.
 \label{eq:vd-l2-layer-residual-L2}
\end{align}
\item On each open side $E_j$,
\begin{align}
 \int_{\Sone}\bigl(w\cdot n_j\bigr)^2 \UsideZero^\e(x,w)\,\ud\nu(w) &=0 \quad(x\in E_j).
 \label{eq:vd-l2-localized-second-moment}
\end{align}
\item For $0<\e\le\e_0(\Om)$, the in-flow boundary data $h^\e=g^\e-(\rho_0+\UsideZero^\e)\big|_{\gamma_-}$ satisfy
\begin{align}
 \abs{h^\e}_{2,-} &\leq C_\Om\mathcal G_{\rm vd}\e^{\frac12}.
 \label{eq:vd-l2-incoming-error-bound}
\end{align}
\end{enumerate}
\end{lemma}

\begin{proof}\leavevmode
\paragraph{\underline{Norms of the side layer}} Evaluating \eqref{eq:vd-l2-leading-side-estimate} at $\eta=\frac d\e$, we have, for every $s\in[0,L_j]$ and almost every $d>0$,
\begin{align}
 \sum_{a=0}^1\nm{\partial_s^a\Uside_{j,0}\Big(s,\frac{d}{\e},\cdot\Big)}_{L^\infty(\Sone)} &\leq C\mathcal G_{\rm vd}\ue^{-\frac{\kappaLtwo d}{\e}},
 \label{eq:vd-l2-leading-side-physical-decay}
\end{align}
For $1\le q<\infty$, integration in $s,d$ gives
\begin{align}
 \int_0^{L_j}\int_0^\infty \ue^{-\frac{q\kappaLtwo d}{\e}}\,\ud d\ud s &\leq C_{q,\kappaLtwo}L_j\e.
 \label{eq:vd-l2-normal-scaling}
\end{align}
Since $0\le a_{j,\e}\le1$ and $\nu(\Sone)=1$, summing over the sides gives \eqref{eq:vd-l2-layer-Lr}. The estimate \eqref{eq:vd-l2-layer-Linfty} follows directly from \eqref{eq:vd-l2-leading-side-physical-decay}.

\paragraph{\underline{Source estimates}} In the side coordinates, $w\cdot\nabla_x=\tau_j\partial_s+\mu_j\partial_d$. The cutoff $a_{j,\e}$ commutes with $\qk$, and the Milne equation \eqref{eq:vd-l2-leading-side-equation} cancels the terms containing $\e^{-1}$. The distributional product rule therefore gives
\begin{align}
 \Le\UsideZeroAt{j}^\e &= \bigl(\tau_j\partial_sa_{j,\e}+\mu_j\partial_da_{j,\e}\bigr)\Uside_{j,0}+a_{j,\e}\tau_j\partial_s\Uside_{j,0},
 \label{eq:vd-l2-exact-side-residual}
\end{align}
where $\Uside_{j,0}$ and its derivatives are evaluated at $\bigl(s,\frac d\e,w\bigr)$. It remains to estimate the cutoff derivatives and the last tangential derivative. We use \eqref{eq:vd-l2-leading-side-physical-decay} for $\Uside_{j,0}$ and $\partial_s\Uside_{j,0}$.

The derivative of $\chi_{\rm end}\bigl(\frac s\e\bigr)$ is bounded by $C\e^{-1}$ and supported on $\e\le s\le2\e$. There $\Theta\bigl(\frac ds\bigr)\ne0$ implies $d\le4\cfs\e$. Thus this source term has $L^\infty$ norm at most $C\mathcal G_{\rm vd}\e^{-1}$ and support of area $O(\e^2)$, giving bounds $C\mathcal G_{\rm vd}\e$ in $L^1$ and $C\mathcal G_{\rm vd}$ in $L^2$. The cutoff at $s=L_j$ gives the same bounds.

By \eqref{eq:Theta-choice}, a derivative of $\Theta\bigl(\frac ds\bigr)$ is supported where $\cfs s\le d\le2\cfs s$. The factor $\chi_{\rm end}\bigl(\frac s\e\bigr)$ also gives $s\ge\e$. On this set, $\abs{\nabla_{s,d}\Theta\bigl(\frac ds\bigr)}\le Cs^{-1}$. Integrating first in $d$ and then setting $s=\e u$, we obtain, for $q\in\{1,2\}$,
\begin{align}
 \int_\e^{L_j}\!\int_{\cfs s}^{2\cfs s}s^{-q}\ue^{-\frac{q\kappaLtwo d}{\e}}\ud d\ud s &\leq \cfs\e^{2-q}\int_1^\infty u^{1-q}\ue^{-q\kappaLtwo\cfs u}\ud u\leq C\e^{2-q}.
\end{align}
Replacing $s$ by $L_j-s$ proves the same estimates at the other endpoint. Thus these source terms are bounded by $C\mathcal G_{\rm vd}\e$ in $L^1$ and $C\mathcal G_{\rm vd}$ in $L^2$.

A derivative of $\zeta_j$ is supported on $d\ge\ellnor$. By exponential decay, its contribution is $O\big(\mathcal G_{\rm vd}\ue^{-\frac{\kappaLtwo\ellnor}\e}\big)$ in both norms. The last term in \eqref{eq:vd-l2-exact-side-residual} has $L^1$ norm $O(\mathcal G_{\rm vd}\e)$ and $L^2$ norm $O(\mathcal G_{\rm vd}\e^{\frac12})$, by \eqref{eq:vd-l2-leading-side-physical-decay}--\eqref{eq:vd-l2-normal-scaling}. Summing over the sides proves \eqref{eq:vd-l2-layer-residual-norms}--\eqref{eq:vd-l2-layer-residual-L2}.

\paragraph{\underline{Boundary moment}} On $E_j$, only $\UsideZeroAt{j}^\e$ has a nonzero trace, by \eqref{eq:vd-l2-cross-side-trace-zero}. At $d=0$, the normal and ray cutoffs equal one. The trace is therefore $\chi_{\rm end}\bigl(\frac s\e\bigr)\chi_{\rm end} \bigl(\frac{L_j-s}\e\bigr)\Uside_{j,0}(s,0,w)$. Both cutoff factors are independent of $w$. Lemma~\ref{lem:vd-l2-milne-moments} and $(w\cdot n_j)^2=\mu_j(w)^2$ give \eqref{eq:vd-l2-localized-second-moment}.

\paragraph{\underline{In-flow boundary data}} On the in-flow boundary over $E_j$, use $\rho_0\big|_{E_j}=e_{j,0}$ and \eqref{eq:vd-l2-leading-side-incoming} to obtain
\begin{align}
 h^\e &=(g^\e-g_0)+\left[1-\chi_{\rm end}\Big(\frac{s}{\e}\Big)\chi_{\rm end}\Big(\frac{L_j-s}{\e}\Big)\right](g_{j,0}-e_{j,0}).
 \label{eq:vd-l2-incoming-discrepancy}
\end{align}
By \eqref{eq:vd-l2-data-norm}, the first term is bounded by $\mathcal G_{\rm vd}\e^{\frac12}$ in $L^2(\gamma_-)$. The second term has $L^\infty(\gamma_-)$ norm at most $C\mathcal G_{\rm vd}$ and is supported on $(0,2\e)\cup(L_j-2\e,L_j)$. Its kinetic $L^2$ norm is therefore at most $C\mathcal G_{\rm vd}\e^{\frac12}$. Summing over the sides proves \eqref{eq:vd-l2-incoming-error-bound}.
\end{proof}

\subsection{Elliptic Estimates}
\label{subsec:vd-l2-poisson-pairing}

The kernel estimate uses a Poisson equation whose source is the velocity average of the remainder. Testing transport by the directional derivative of the Poisson solution gives its Hessian in the volume integral and its gradient in the boundary integral. We first bound these derivatives in terms of the $L^2$ norm of the Poisson source.

We decompose the Poisson solution into an $H^2$ function and the first singular mode at each reentrant vertex. Since every such exponent exceeds $\frac12$, the boundary gradient belongs to $L^2$ and the Hessian belongs to $L^{p_{\rm P}}$ for some $p_{\rm P}>\frac43$. This range is sufficient for the interpolation in Subsection~\ref{subsec:vd-l2-energy} and the $O(\e^{\frac12})$ remainder estimate.

\begin{lemma}
\label{lem:vd-l2-polygonal-poisson}
Let $\varrho\in L^2(\Om)$, and let $\phi\in H^1_0(\Om)$ be the unique weak solution of the Poisson problem
\begin{align}
 &-\Delta\phi=\varrho \quad\hbox{in }\Om.
 \label{eq:vd-l2-poisson-problem}
\end{align}
For each reentrant vertex $V_m$, recall from \eqref{eq:intro-zero-ray-pairs-comp} that $\lambda_{m,1}=\frac{\pi}{\omega_m}\in(\frac12,1)$. Fix cutoffs $\chi_m^{\rm D}\in C_c^\infty(B_{2\ellvtx}(V_m))$ equal to one in $B_{\ellvtx}(V_m)$. Their supports lie in pairwise disjoint vertex neighborhoods and meet only the two side collars incident to $V_m$. Then the following statements hold.
\begin{enumerate}
\item There are coefficients $c_m\in\mathbb R$ such that
\begin{align}
 &\phi = \phi_{\rm reg} +\sum_{m:\,\omega_m>\pi} c_m\chi_m^{\rm D} \bulkdistat{m}^{\lambda_{m,1}} \sin(\lambda_{m,1}\theta_m), \qquad \phi_{\rm reg}\in H^2(\Om).
 \label{eq:vd-l2-poisson-singular-decomposition}
\end{align}
With $\nm{\nabla_x\phi}_{L^2_{\rm pw}(\partial\Om)}^2:=\sum_j\nm{\nabla_x\phi}_{L^2(E_j)}^2$, the regular part, the singular coefficients, and the boundary gradient satisfy
\begin{align}
 &\nm{\phi_{\rm reg}}_{H^2(\Om)} +\sum_{m:\,\omega_m>\pi}\abs{c_m} +\nm{\nabla_x\phi}_{L^2_{\rm pw}(\partial\Om)} \leq C_\Om\nm{\varrho}_{L^2(\Om)}.
 \label{eq:vd-l2-poisson-decomposition-bound}
\end{align}
\item The range of exponents
\begin{align}
 &\frac43<p_{\rm P}\leq2, \qquad p_{\rm P}<\frac{2}{2-\lambda_{m,1}} \quad\hbox{for every }m\hbox{ with }\omega_m>\pi, \qquad p_{\rm P}':=\frac{p_{\rm P}}{p_{\rm P}-1},
 \label{eq:vd-l2-poisson-p-range}
\end{align}
is nonempty, and for every $p_{\rm P}$ in it the Hessian obeys
\begin{align}
 &\nm{D_x^2\phi}_{L^{p_{\rm P}}(\Om)} \leq C_{\Om,p_{\rm P}}\nm{\varrho}_{L^2(\Om)},
 \label{eq:vd-l2-poisson-W2p}
\end{align}
and the conjugate exponent satisfies
\begin{align}
 &\frac{2}{p_{\rm P}'}=2-\frac{2}{p_{\rm P}} >\frac12.
 \label{eq:vd-l2-poisson-interpolation-gain}
\end{align}
\end{enumerate}
\end{lemma}

\begin{proof}\leavevmode
\paragraph{\underline{Singular decomposition and boundary gradient}} The decomposition \eqref{eq:vd-l2-poisson-singular-decomposition} and the bounds for $\phi_{\rm reg}$ and $c_m$ follow from the Dirichlet theory in a plane polygon \cite[Theorem~4.4.3.7]{Grisvard1985}, \cite[Chapters~2--3]{Dauge1988}. Its hypotheses hold here: $\Om$ is a bounded plane polygon, $\varrho\in L^2(\Om)$, and the Dirichlet data vanish on $\partial\Om$.

The singular functions have exponents $\frac{n\pi}{\omega_m}\in(0,1)$. Since $0<\omega_m<2\pi$, this range contains exactly $\lambda_{m,1}=\frac\pi{\omega_m}$ at each reentrant vertex and contains no exponent at a convex vertex. The threshold $\frac{n\pi}{\omega_m}=1$ would require $n=1$ and $\omega_m=\pi$, which is excluded by \eqref{eq:vd-l2-angle-assumption} after collinear sides are merged. Thus subtracting the finite sum in \eqref{eq:vd-l2-poisson-singular-decomposition} gives the stated $H^2$ regular part.

The decomposition is unique: a nonzero singular mode with exponent $\lambda_{m,1}<1$ does not belong to $H^2$, and the disjoint vertex neighborhoods prevent cancellation between these modes. Hence the coefficient maps and $R\colon\varrho\mapsto\phi_{\rm reg}$ are linear. To obtain the quantitative bound for $R$, suppose $\varrho_n\to\varrho$ in $L^2(\Om)$ and $R\varrho_n\to v$ in $H^2(\Om)$. The energy estimate gives convergence of the weak solutions $\phi_n\to\phi$ in $H^1_0(\Om)$, and continuity of the singular coefficients gives $c_m(\varrho_n)\to c_m(\varrho)$. Since the fixed cutoff singular functions belong to $H^1(\Om)$, \eqref{eq:vd-l2-poisson-singular-decomposition} then gives $R\varrho_n\to R\varrho$ in $H^1(\Om)$, so $v=R\varrho$. The closed graph theorem yields $\nm{\phi_{\rm reg}}_{H^2(\Om)}\le C_\Om\nm{\varrho}_{L^2(\Om)}$. Together with the coefficient bound, this proves the first two terms of \eqref{eq:vd-l2-poisson-decomposition-bound}; the constant includes the fixed cutoff norms.

To estimate the boundary gradient and the Hessian, use
\begin{align}
 \abs{ \nabla_x\left( \bulkdistat{m}^{\lambda_{m,1}} \sin(\lambda_{m,1}\theta_m) \right) } &\leq C\bulkdistat{m}^{\lambda_{m,1}-1},
 \label{eq:vd-l2-singular-gradient-size}\\
 \abs{ D_x^2\left( \bulkdistat{m}^{\lambda_{m,1}} \sin(\lambda_{m,1}\theta_m) \right) } &\leq C\bulkdistat{m}^{\lambda_{m,1}-2}.
 \label{eq:vd-l2-singular-hessian-size}
\end{align}
The boundary integral $\int_0^{\ellvtx}\bulkdistat{m}^{2\lambda_{m,1}-2}\ud\bulkdistat{m}$ is finite exactly when $\lambda_{m,1}>\frac12$. Combining this with the $H^2$ trace estimate for $\phi_{\rm reg}$ on each open side proves the last term in \eqref{eq:vd-l2-poisson-decomposition-bound}.

\paragraph{\underline{Hessian estimate}} The area element $\bulkdistat{m}\ud\bulkdistat{m}\ud\theta_m$ gives local $L^{p_{\rm P}}$ integrability exactly when $p_{\rm P}(\lambda_{m,1}-2)+1>-1$, or $p_{\rm P}<\frac2{2-\lambda_{m,1}}$. Terms containing derivatives of $\chi_m^{\rm D}$ are supported away from $V_m$ and have $L^{p_{\rm P}}$ norm at most $C_{\Om,p_{\rm P}}\abs{c_m}$. Their sum is bounded by \eqref{eq:vd-l2-poisson-decomposition-bound}. Since $\Om$ has finite measure and $p_{\rm P}\le2$, the $H^2$ bound for $\phi_{\rm reg}$ gives its Hessian bound in $L^{p_{\rm P}}$. This proves \eqref{eq:vd-l2-poisson-W2p}.

For every reentrant vertex, $\lambda_{m,1}>\frac12$ implies $\frac2{2-\lambda_{m,1}}>\frac43$. Since there are finitely many vertices, we can choose one $p_{\rm P}>\frac43$ satisfying all the strict upper bounds. The conjugate exponent then satisfies \eqref{eq:vd-l2-poisson-interpolation-gain}.

As $\omega_m\to2\pi$, the exponent $\lambda_{m,1}$ tends to $\frac12$, and the upper bound for $p_{\rm P}$ tends to $\frac43$. We therefore choose $p_{\rm P}$ after fixing the polygon.
\end{proof}

For the regular Hessian, Cauchy's inequality uses the $O(\e^{\frac12})$ norm of the side layer in $L^2$. For a singular Hessian, the distance to the vertex is comparable to the tangential coordinate in each incident collar. Integration in the normal coordinate gives a factor $\e$ by Milne decay. The endpoint cutoff restricts the tangential integral to distances at least $\e$, giving $\e^{\lambda_{m,1}}$. Since $\lambda_{m,1}>\frac12$, this is sufficient.

\begin{lemma}
\label{lem:vd-l2-layer-hessian}
Let $\phi$ solve \eqref{eq:vd-l2-poisson-problem} and let $0<\e\leq\e_0(\Om)$. Then
\begin{align}
 \abs{ \int_{\Om\times\Sone} \UsideZero^\e (w\otimes w):D_x^2\phi \,\ud x\ud\nu } &\leq C_\Om\mathcal G_{\rm vd} \Big( \e^{\frac12} +\sum_{m:\,\omega_m>\pi} \e^{\lambda_{m,1}} \Big) \nm{\varrho}_{L^2(\Om)} \leq C_\Om\mathcal G_{\rm vd}\e^{\frac12} \nm{\varrho}_{L^2(\Om)}.
 \label{eq:vd-l2-layer-hessian-bound}
\end{align}
\end{lemma}

\begin{proof}
Cauchy's inequality bounds the contribution of $\phi_{\rm reg}$ by $C_\Om\mathcal G_{\rm vd}\e^{\frac12}\nm{\varrho}_{L^2(\Om)}$. Here we use \eqref{eq:vd-l2-layer-Lr} with $q=2$ and \eqref{eq:vd-l2-poisson-decomposition-bound}.

Fix a reentrant vertex $V_m$ and write $\lambda=\lambda_{m,1}$. Only its two incident collars meet $\supp\chi_m^{\rm D}$ near the vertex. In the $i$-th collar, the cutoff gives $\sidedistat{m}{i}\ge\e$ and $0\le d_{m,i}\le2\cfs\sidedistat{m}{i}$, so $\bulkdistat{m}\simeq\sidedistat{m}{i}$. By \eqref{eq:vd-l2-leading-side-physical-decay}, \eqref{eq:vd-l2-singular-hessian-size}, and normal integration as in \eqref{eq:vd-l2-normal-scaling},
\begin{align}
 \int_{(B_{2\ellvtx}(V_m)\cap\Om)\times\Sone}\babs{\UsideZeroAt{j(m,i)}^\e}\abs{D_x^2\bigl(\bulkdistat{m}^\lambda\sin(\lambda\theta_m)\bigr)}\,\ud x\ud\nu &\le C\mathcal G_{\rm vd}\e\int_\e^{2\ellvtx}(\sidedistat{m}{i})^{\lambda-2}\,\ud\sidedistat{m}{i}\le C\mathcal G_{\rm vd}\e^\lambda,
\end{align}
where the last inequality uses $\lambda-2<-1$ and an integral of size $\frac{\e^{\lambda-1}}{1-\lambda}$.

The other incident collar gives the same estimate. Terms with a derivative of $\chi_m^{\rm D}$ are supported a fixed distance from $V_m$ and belong to $L^2$. By Cauchy's inequality and \eqref{eq:vd-l2-layer-Lr}, their contribution is $O(\mathcal G_{\rm vd}\e^{\frac12})$. Including $c_m$, the total contribution from $V_m$ is at most $C\abs{c_m}\mathcal G_{\rm vd}(\e^{\frac12}+\e^{\lambda_{m,1}})$. Sum over the reentrant vertices and apply \eqref{eq:vd-l2-poisson-decomposition-bound} to obtain the first inequality in \eqref{eq:vd-l2-layer-hessian-bound}. The second follows from $\lambda_{m,1}>\frac12$ and $0<\e\le1$.
\end{proof}

\subsection{Energy and Kernel Estimates}
\label{subsec:vd-l2-energy}

We first estimate the microscopic part and the out-flow boundary trace by the energy identity. These estimates will then be used in the Poisson argument for the velocity average, which is the component in the kernel of $\qk$. Define
\begin{align}
 \Err^\e &:={} \Uphys^\e-\rho_0-\UsideZero^\e,\qquad \rho_{\Err}^\e:=\pk\Err^\e,\qquad \Err^{\e,\perp}:=\qk\Err^\e.
 \label{eq:vd-l2-residual-decomposition}
\end{align}
Since $\qk\rho_0=0$ and $\Le\UsideZero^\e=S_{\rm sl}^\e$, the remainder satisfies
\begin{align}
 \Le\Err^\e &=-w\cdot\nabla_x\rho_0-S_{\rm sl}^\e,\qquad \Err^\e\big|_{\gamma_-}=h^\e.
 \label{eq:vd-l2-residual-equation}
\end{align}
By Lemma~\ref{lem:physical-comparison-comp}, $\nm{\Uphys^\e}_{L^\infty(\Om\times\Sone)}\le\abs{g^\e}_{\infty,-}$. For fixed $\e>0$, the equation gives $w\cdot\nabla_x\Uphys^\e=-\e^{-1}\qk\Uphys^\e\in L^2(\Om\times\Sone)$. The characteristic formula also gives $\abs{\Uphys^\e}_{\infty,\pm}\le\abs{g^\e}_{\infty,-}$. Since the kinetic boundary measure is finite, $\Uphys^\e\in\mathscr W_{\rm kin}^2$.

For the side layer, \eqref{eq:vd-l2-layer-Lr} and \eqref{eq:vd-l2-layer-residual-L2} give $w\cdot\nabla_x\UsideZero^\e=S_{\rm sl}^\e-\e^{-1}\qk\,\UsideZero^\e\in L^2(\Om\times\Sone)$. By \eqref{eq:vd-l2-leading-side-trace} and the cutoff definition, $\babs{\UsideZero^\e}_{\infty,\pm}\le C_\Om\mathcal G_{\rm vd}$. Also, $\rho_0\in H^1(\Om)$ has trace $\rho_0\big|_{\partial\Om}=e_0\in L^\infty(\partial\Om)$ by \eqref{eq:vd-l2-harmonic-limit}--\eqref{eq:vd-l2-harmonic-estimate}. Thus $\rho_0,\UsideZero^\e,\Err^\e\in\mathscr W_{\rm kin}^2$.

Combining the comparison estimate for $\Uphys^\e$, \eqref{eq:vd-l2-harmonic-estimate}, and \eqref{eq:vd-l2-layer-Linfty}, we obtain, uniformly for $0<\e\le1$,
\begin{align}
 &\nm{\Err^\e}_{L^\infty(\Om\times\Sone)} +\nm{\Err^{\e,\perp}}_{L^\infty(\Om\times\Sone)} \leq C_\Om\mathcal G_{\rm vd}.
 \label{eq:vd-l2-uniform-residual-bound}
\end{align}
\paragraph{\underline{Basic energy estimate}} Test \eqref{eq:vd-l2-residual-equation} with $\Err^\e$. Since $\int_{\Om\times\Sone}(\qk\Err^\e)\Err^\e\ud x\ud\nu=\nm{\Err^{\e,\perp}}_{L^2(\Om\times\Sone)}^2$, Green's identity and the in-flow boundary condition \eqref{eq:vd-l2-residual-equation} give
\begin{align}
 \frac12\abs{\Err^\e}_{2,+}^2+\e^{-1}\nm{\Err^{\e,\perp}}_{L^2(\Om\times\Sone)}^2&=-\int_{\Om\times\Sone}(w\cdot\nabla_x\rho_0)\Err^{\e,\perp}\,\ud x\ud\nu
 \label{eq:vd-l2-energy-identity}\\
 &\quad-\int_{\Om\times\Sone}S_{\rm sl}^\e\Err^\e\,\ud x\ud\nu+\frac12\abs{h^\e}_{2,-}^2.\notag
\end{align}
Here the term containing $\rho_{\Err}^\e$ vanishes because $\int_{\Sone}w\ud\nu=0$.

By Young's inequality, the first integral is bounded by $\frac1{2\e}\nm{\Err^{\e,\perp}}_{L^2}^2+C\e\bnm{\nabla_x\rho_0}_{L^2(\Om)}^2$. Use \eqref{eq:vd-l2-harmonic-estimate} to bound $\bnm{\nabla_x\rho_0}_{L^2(\Om)}$ by $C_\Om\mathcal G_{\rm vd}$.

For the side-layer source term, \eqref{eq:vd-l2-layer-residual-norms} and \eqref{eq:vd-l2-uniform-residual-bound} give $\babs{\int_{\Om\times\Sone}S_{\rm sl}^\e\Err^\e\ud x\ud\nu}\le\nm{S_{\rm sl}^\e}_{L^1}\nm{\Err^\e}_{L^\infty}\le C_\Om\mathcal G_{\rm vd}^2\e$. The boundary term is bounded by $\frac12\abs{h^\e}_{2,-}^2\le C_\Om\mathcal G_{\rm vd}^2\e$, by \eqref{eq:vd-l2-incoming-error-bound}.

Insert these estimates into \eqref{eq:vd-l2-energy-identity} and absorb half of the collision term. We obtain
\begin{align}
 \nm{\Err^{\e,\perp}}_{L^2(\Om\times\Sone)} &\leq C_\Om\mathcal G_{\rm vd}\e,
 \label{eq:vd-l2-energy-consequences}\\
 \abs{\Err^\e}_{2,+} &\leq C_\Om\mathcal G_{\rm vd}\e^{\frac12}.
 \label{eq:vd-l2-energy-outgoing}
\end{align}

\paragraph{\underline{Kernel estimate}} Apply Lemma~\ref{lem:vd-l2-polygonal-poisson} with $\varrho=\rho_{\Err}^\e$. Denote the Poisson solution by $\phi^\e$, so that
\begin{align}
 -\Delta\phi^\e&=\rho_{\Err}^\e \quad\hbox{in }\Om,\qquad \phi^\e\in H^1_0(\Om).
 \label{eq:vd-l2-poisson-remainder}
\end{align}
Let $\phi^\e_{\rm reg}$ be the regular part in \eqref{eq:vd-l2-poisson-singular-decomposition}. Set $f^\e=\Uphys^\e-\rho_0=\Err^\e+\UsideZero^\e$. Adding $\UsideZero^\e$ to the remainder cancels $S_{\rm sl}^\e$ and gives
\begin{align}
 \Le f^\e=-w\cdot\nabla_x\rho_0.
 \label{eq:vd-l2-total-error-equation}
\end{align}

By Lemma~\ref{lem:physical-comparison-comp} and \eqref{eq:vd-l2-harmonic-estimate},
\begin{align}
 \nm{f^\e}_{L^\infty(\Om\times\Sone)}+\nm{f^\e}_{L^\infty(\gamma;\ud\gamma_{\rm kin})}&\leq C_\Om\mathcal G_{\rm vd}.
 \label{eq:vd-l2-total-error-Linfty}
\end{align}
We combine the weak identities for \eqref{eq:vd-l2-total-error-equation} with test functions $\phi^\e$ and $w\cdot\nabla_x\phi^\e$. The latter contains $\e^{-1}\int_{\Om\times\Sone}\qk f^\e(w\cdot\nabla_x\phi^\e)\ud x\ud\nu$, whose direct estimate would contribute $\e^{-1}\nm{\qk f^\e}_{L^2}\nm{\nabla_x\phi^\e}_{L^2}$. Adding $\e^{-1}$ times the former identity cancels this contribution exactly, as in \cite{Guo.Wu2025}. We then substitute $f^\e=\rho_{\Err}^\e+\Err^{\e,\perp}+\UsideZero^\e$ and estimate the remaining terms.

\paragraph{\underline{Testing with $\phi^\e$}} The function $\phi^\e$ belongs to $H^1_0(\Om)$ and is independent of $w$. It therefore belongs to $\mathscr W_{\rm kin}^2$ with zero kinetic trace and can be used in Green's identity.

Take $\phi^\e$ as the test function in \eqref{eq:vd-l2-total-error-equation}. Its zero trace eliminates the boundary term. The identities $\int_{\Sone}\qk f^\e\ud\nu=0$ and $\int_{\Sone}w\ud\nu=0$ eliminate the collision and source terms, respectively. Integration by parts therefore gives
\begin{align}
 -\int_{\Om\times\Sone}f^\e\big(w\cdot\nabla_x\phi^\e\big)\ud x\ud\nu =-\int_{\Om\times\Sone}\qk f^\e\big(w\cdot\nabla_x\phi^\e\big)\ud x\ud\nu=0.
 \label{eq:vd-l2-first-test-cancellation}
\end{align}
The first equality again uses $\int_{\Sone}w\ud\nu=0$.

\paragraph{\underline{Testing with $w\cdot\nabla_x\phi^\e$}} The streaming derivative of this test function is $(w\otimes w):D_x^2\phi^\e$, whose velocity average is $\frac12\Delta\phi^\e=-\frac12\rho_{\Err}^\e$. Thus its pairing with $\rho_{\Err}^\e$ will give the square of the desired $L^2$ norm. The Hessian may fail to belong to $L^2$ at a reentrant vertex, so we first justify this test away from the vertices and then remove the artificial boundaries. Fix $\e>0$ and set $\Om_{\delta_{\rm exc}}=\Om\setminus\bigcup_m\overline{B_{\delta_{\rm exc}}(V_m)}$.

For small $\delta_{\rm exc}>0$, the domain $\Om_{\delta_{\rm exc}}$ is bounded and Lipschitz. By \eqref{eq:vd-l2-poisson-singular-decomposition}, $\phi^\e\in H^2(\Om_{\delta_{\rm exc}})$, so $w\cdot\nabla_x\phi^\e$ and its streaming derivative belong to $L^2(\Om_{\delta_{\rm exc}}\times\Sone)$. Equation \eqref{eq:vd-l2-total-error-equation} gives $w\cdot\nabla_xf^\e\in L^2(\Om\times\Sone)$ for fixed $\e$. Thus, for almost every $w$ and almost every line parallel to $w$, $f^\e$ has an $H^1$ representative on each connected characteristic interval. Its bulk bound \eqref{eq:vd-l2-total-error-Linfty} extends by continuity to this representative and to the endpoints of its subintervals. The characteristic proof of Lemma~\ref{lem:vd-l2-green-identity} applies also to the circular arcs: projection onto $w^\perp$ has Jacobian $\abs{w\cdot n_{\delta_{\rm exc}}}$, so exceptional lines and tangencies have zero kinetic boundary measure. The induced arc traces of $f^\e$ are therefore bounded by $C_\Om\mathcal G_{\rm vd}$ independently of $\delta_{\rm exc}$. The trace of $w\cdot\nabla_x\phi^\e$ belongs to $L^2$ by $H^2$ regularity. Both traces are therefore in kinetic $L^2$. Testing on $\Om_{\delta_{\rm exc}}$, we obtain
\begin{align}
 &\int_{(\partial\Om\cap\partial\Om_{\delta_{\rm exc}})\times\Sone} f^\e(w\cdot\nabla_x\phi^\e)(w\cdot n) \,\ud\sigma\ud\nu+ \int_{(\partial\Om_{\delta_{\rm exc}}\setminus\partial\Om)\times\Sone} f^\e(w\cdot\nabla_x\phi^\e)(w\cdot n_{\delta_{\rm exc}}) \,\ud\sigma\ud\nu
 \label{eq:vd-l2-punctured-second-test}\\
 &\quad- \int_{\Om_{\delta_{\rm exc}}\times\Sone} f^\e (w\otimes w):D_x^2\phi^\e\,\ud x\ud\nu +\e^{-1}\int_{\Om_{\delta_{\rm exc}}\times\Sone} \qk f^\e\,w\cdot\nabla_x\phi^\e\,\ud x\ud\nu \notag\\
 &={} -\int_{\Om_{\delta_{\rm exc}}\times\Sone} (w\cdot\nabla_x\rho_0)(w\cdot\nabla_x\phi^\e) \,\ud x\ud\nu,\notag
\end{align}
where $n_{\delta_{\rm exc}}$ is the outward normal of $\Om_{\delta_{\rm exc}}$ on the artificial circular arcs.

\paragraph{\underline{Passage to the polygon}} We choose radii along which the contribution of $\nabla_x\phi^\e_{\rm reg}$ on the circular arcs tends to zero. Set $A(t)=\sum_m\int_{\partial B_t(V_m)\cap\Om}\abs{\nabla_x\phi^\e_{\rm reg}}^2\ud\sigma$ and $\mathfrak m(\delta)=\sum_m\int_{B_\delta(V_m)\cap\Om}\abs{\nabla_x\phi^\e_{\rm reg}}^2\ud x$. The coarea formula gives $\mathfrak m(\delta)=\int_0^\delta A(t)\ud t\to0$. Choose the radii from the full-measure set on which these slices agree with the Sobolev traces, and start with $n$ large enough that the vertex balls are disjoint and contained in their fixed vertex neighborhoods. Averaging over $(2^{-n-1},2^{-n})$ then gives a radius $\delta_n$ in this interval with $A(\delta_n)\le2^{n+1}\mathfrak m(2^{-n})$. Thus $\delta_n\downarrow0$ and $\delta_nA(\delta_n)\le2\mathfrak m(2^{-n})\to0$.

Each arc has length at most $C\delta_n$. By Cauchy's inequality,
\begin{align}
 &\bigg|\sum_m\int_{(\partial B_{\delta_n}(V_m)\cap\Om)\times\Sone} f^\e(w\cdot\nabla_x\phi^\e_{\rm reg}) (w\cdot n_{\delta_n})\,\ud\sigma\ud\nu\bigg|\le C\nm{f^\e}_{L^\infty} \bigl(\delta_n A(\delta_n)\bigr)^{\frac12}\to0.
\end{align}
For a singular mode of exponent $\lambda_{m,1}$, \eqref{eq:vd-l2-singular-gradient-size} gives a gradient of size $O\big(\delta_n^{\lambda_{m,1}-1}\big)$ on an arc of length $O(\delta_n)$. The $L^\infty$ bound for the characteristic trace of $f^\e$ therefore gives a contribution $O(\delta_n^{\lambda_{m,1}})\to0$.

On $\partial\Om\cap\partial\Om_{\delta_n}$, dominated convergence gives the integral over all of $\partial\Om$. We use $f^\e\big|_\gamma\in L^\infty(\gamma;\ud\gamma_{\rm kin})$ from \eqref{eq:vd-l2-total-error-Linfty} and $\nabla_x\phi^\e\in L^2_{\rm pw}(\partial\Om)$ from \eqref{eq:vd-l2-poisson-decomposition-bound}.

For fixed $\e$, the collision integrand belongs to $L^1$ because $\qk f^\e\in L^2(\Om\times\Sone)$ and $\nabla_x\phi^\e\in L^2(\Om)$. Its limit is the full integral $\e^{-1}\int_{\Om\times\Sone}\qk f^\e(w\cdot\nabla_x\phi^\e)\ud x\ud\nu$. Adding $\e^{-1}$ times \eqref{eq:vd-l2-first-test-cancellation} cancels this term exactly.

The source integrand also belongs to $L^1$, since $\nabla_x\rho_0,\nabla_x\phi^\e\in L^2(\Om)$. Its integral vanishes:
\begin{align}
 \int_{\Om\times\Sone} (w\cdot\nabla_x\rho_0)(w\cdot\nabla_x\phi^\e) \,\ud x\ud\nu &= \frac12\int_\Om \nabla_x\rho_0\cdot\nabla_x\phi^\e\,\ud x =0.
 \label{eq:vd-l2-harmonic-two-test-cancellation}
\end{align}
Here we use $\int_{\Sone}w\otimes w\ud\nu=\frac12I$ and test the weak harmonic equation \eqref{eq:vd-l2-harmonic-limit} by $\phi^\e\in H^1_0(\Om)$.

In the Hessian term, write $f^\e=\rho_{\Err}^\e+\Err^{\e,\perp}+\UsideZero^\e$ on $\Om_{\delta_{\rm exc}}$. We first integrate the term containing $\rho_{\Err}^\e$ in velocity. This is justified by $\rho_{\Err}^\e\in L^\infty(\Om)$ and $D_x^2\phi^\e\in L^2(\Om_{\delta_{\rm exc}})$. Fubini's theorem and $\int_{\Sone}w\otimes w\ud\nu=\frac12I$ give
\begin{align}
 \int_{\Om_{\delta_{\rm exc}}\times\Sone} \rho_{\Err}^\e (w\otimes w):D_x^2\phi^\e\,\ud x\ud\nu &= \frac12\int_{\Om_{\delta_{\rm exc}}} \rho_{\Err}^\e\Delta\phi^\e\,\ud x =-\frac12\int_{\Om_{\delta_{\rm exc}}} \abs{\rho_{\Err}^\e}^2\,\ud x,
 \label{eq:vd-l2-scalar-contracted-pairing}
\end{align}
by \eqref{eq:vd-l2-poisson-remainder}. The integrable bound $\frac12\abs{\rho_{\Err}^\e}^2$ now permits $\delta_{\rm exc}\downarrow0$ by dominated convergence. The limit is $-\frac12\nm{\rho_{\Err}^\e}_{L^2(\Om)}^2$.

For the other two Hessian terms, \eqref{eq:vd-l2-uniform-residual-bound} and \eqref{eq:vd-l2-layer-Linfty} give
\begin{align}
 \bigl(\abs{\Err^{\e,\perp}}+\abs{\UsideZero^\e}\bigr)\abs{(w\otimes w):D_x^2\phi^\e} &\leq C_\Om\mathcal G_{\rm vd}\abs{D_x^2\phi^\e}\in L^1(\Om\times\Sone),
 \label{eq:vd-l2-hessian-integrable-majorant}
\end{align}
since \eqref{eq:vd-l2-poisson-W2p} gives $D_x^2\phi^\e\in L^{p_{\rm P}}(\Om)\subset L^1(\Om)$. For fixed $\e$, this bound is independent of $\delta_{\rm exc}$, so dominated convergence applies. We will estimate the limits by the $L^{p_{\rm P}'}$--$L^{p_{\rm P}}$ pairing in \eqref{eq:vd-l2-microscopic-interpolation} and \eqref{eq:vd-l2-poisson-W2p}, and by Lemma~\ref{lem:vd-l2-layer-hessian}.

Set $\delta_{\rm exc}=\delta_n$ in \eqref{eq:vd-l2-punctured-second-test} and let $n\to\infty$. Using \eqref{eq:vd-l2-first-test-cancellation} and \eqref{eq:vd-l2-harmonic-two-test-cancellation}, we obtain
\begin{align}
 &\int_{\partial\Om\times\Sone} f^\e (w\cdot\nabla_x\phi^\e)(w\cdot n) \,\ud\sigma\ud\nu +\frac12\nm{\rho_{\Err}^\e}_{L^2(\Om)}^2 = \int_{\Om\times\Sone} \bigl(\Err^{\e,\perp}+\UsideZero^\e\bigr) (w\otimes w):D_x^2\phi^\e \,\ud x\ud\nu.
 \label{eq:vd-l2-global-second-test-identity}
\end{align}

\paragraph{\underline{Boundary cancellation and estimates}} The second Milne moment removes the side layer from the boundary integral in \eqref{eq:vd-l2-global-second-test-identity}. On each open side $E_j$, the zero trace of $\phi^\e$ and its local $H^2$ regularity imply $t_j\cdot\nabla_x\phi^\e\big|_{E_j}=0$. Together with \eqref{eq:vd-l2-poisson-decomposition-bound}, this gives
\begin{align}
 \nabla_x\phi^\e\big|_{E_j}&=(\partial_{n_j}\phi^\e)n_j, \qquad \partial_{n_j}\phi^\e\in L^2(E_j).
 \label{eq:vd-l2-normal-gradient-trace}
\end{align}
Thus $w\cdot\nabla_x\phi^\e=(w\cdot n_j)\partial_{n_j}\phi^\e$ on $E_j$. By \eqref{eq:vd-l2-cross-side-trace-zero} and \eqref{eq:vd-l2-localized-second-moment}, the contribution of $\UsideZero^\e$ to the boundary integral is zero. Substituting $f^\e=\Err^\e+\UsideZero^\e$ into that integral gives
\begin{align}
 &\int_{\partial\Om\times\Sone}f^\e(w\cdot\nabla_x\phi^\e)(w\cdot n)\,\ud\sigma\ud\nu
 \label{eq:vd-l2-remainder-boundary-functional}=\sum_j\int_{E_j}\partial_{n_j}\phi^\e\int_{\Sone}(w\cdot n_j)^2\Err^\e\,\ud\nu\,\ud\sigma=:\mathcal T_{\Err}^\e.
\end{align}
Combining \eqref{eq:vd-l2-remainder-boundary-functional} with \eqref{eq:vd-l2-global-second-test-identity}, we obtain
\begin{align}
 &\frac12\nm{\rho_{\Err}^\e}_{L^2(\Om)}^2 = \int_{\Om\times\Sone} \bigl(\Err^{\e,\perp}+\UsideZero^\e\bigr) (w\otimes w):D_x^2\phi^\e\,\ud x\ud\nu -\mathcal T_{\Err}^\e.
 \label{eq:vd-l2-exact-macro-identity}
\end{align}
It remains to bound the right-hand side of \eqref{eq:vd-l2-exact-macro-identity} by $C_\Om\mathcal G_{\rm vd}\e^{\frac12}\nm{\rho_{\Err}^\e}_{L^2(\Om)}$. For the side-layer term, this follows from Lemma~\ref{lem:vd-l2-layer-hessian} with $\varrho=\rho_{\Err}^\e$.

Choose $p_{\rm P}$ in \eqref{eq:vd-l2-poisson-p-range}. Since $p_{\rm P}'\ge2$, interpolation between \eqref{eq:vd-l2-energy-consequences} and \eqref{eq:vd-l2-uniform-residual-bound} gives
\begin{align}
 &\nm{\Err^{\e,\perp}}_{L^{p_{\rm P}'}(\Om\times\Sone)} \leq \nm{\Err^{\e,\perp}}_{L^2(\Om\times\Sone)}^{\frac{2}{p_{\rm P}'}} \nm{\Err^{\e,\perp}}_{L^\infty(\Om\times\Sone)}^{1-\frac{2}{p_{\rm P}'}} \leq C_\Om\mathcal G_{\rm vd}\e^{\frac{2}{p_{\rm P}'}},
 \label{eq:vd-l2-microscopic-interpolation}
\end{align}
By \eqref{eq:vd-l2-poisson-W2p}, the microscopic term in \eqref{eq:vd-l2-exact-macro-identity} is therefore bounded by $C_{\Om,p_{\rm P}}\mathcal G_{\rm vd}\e^{\frac2{p_{\rm P}'}}\nm{\rho_{\Err}^\e}_{L^2(\Om)}$.

For the boundary term, apply Cauchy's inequality in the kinetic boundary measure and use $\Err^\e\big|_{\gamma_-}=h^\e$. Then
\begin{align}
 \abs{\mathcal T_{\Err}^\e} &\leq C\big( \abs{\Err^\e}_{2,+} +\abs{h^\e}_{2,-} \big) \bnm{\nabla_x\phi^\e}_{L^2_{\rm pw}(\partial\Om)}\leq C_\Om\mathcal G_{\rm vd}\e^{\frac12} \nm{\rho_{\Err}^\e}_{L^2(\Om)},
\end{align}
by \eqref{eq:vd-l2-energy-outgoing}, \eqref{eq:vd-l2-incoming-error-bound}, and \eqref{eq:vd-l2-poisson-decomposition-bound}. The remaining velocity integral is the constant $\bigl(\int_{\Sone}\abs{w\cdot n_j}^3\ud\nu\bigr)^{\frac12}$, independent of $j$.

By \eqref{eq:vd-l2-poisson-interpolation-gain}, $\frac2{p_{\rm P}'}>\frac12$, so $\e^{\frac2{p_{\rm P}'}}\le\e^{\frac12}$ for $0<\e\le1$. Insert these estimates into \eqref{eq:vd-l2-exact-macro-identity}. The right-hand side is bounded by $C_\Om\mathcal G_{\rm vd}\e^{\frac12}\nm{\rho_{\Err}^\e}_{L^2(\Om)}$. If this norm is nonzero, divide by it. If it is zero, the estimate is immediate. Thus
\begin{align}
 \nm{\rho_{\Err}^\e}_{L^2(\Om)} \leq C_\Om\mathcal G_{\rm vd}\e^{\frac12}.
 \label{eq:vd-l2-macroscopic-residual-bound}
\end{align}

\begin{proof}[Proof of Theorem~\ref{thm:vd-l2-wedge-free}]
For $0<\e\le\e_0(\Om)$, combine \eqref{eq:vd-l2-energy-consequences} and \eqref{eq:vd-l2-macroscopic-residual-bound} to obtain
\begin{align}
 \nm{\Err^\e}_{L^2(\Om\times\Sone)} &\leq \nm{\rho_{\Err}^\e}_{L^2(\Om)} +\nm{\Err^{\e,\perp}}_{L^2(\Om\times\Sone)} \leq C_\Om\mathcal G_{\rm vd}\e^{\frac12}.
 \label{eq:vd-l2-full-residual-bound}
\end{align}

Together with \eqref{eq:vd-l2-layer-Lr} at $q=2$, this proves \eqref{eq:vd-l2-full-residual-rate}. The microscopic and out-flow boundary estimates \eqref{eq:vd-l2-residual-rates}--\eqref{eq:vd-l2-outgoing-residual-rate} follow from \eqref{eq:vd-l2-energy-consequences}--\eqref{eq:vd-l2-energy-outgoing}. Finally, write $\Uphys^\e-\rho_0=\Err^\e+\UsideZero^\e$ and $\pk\Uphys^\e-\rho_0=\rho_{\Err}^\e+\pk\UsideZero^\e$. The $L^2$ contraction of $\pk$ and the preceding bounds give \eqref{eq:vd-l2-main-rate}.

For $\e_0(\Om)\le\e\le1$, take $F=K=\Uphys^\e$ in Lemma~\ref{lem:vd-l2-green-identity} and use \eqref{eq:physical-model}. We obtain
\begin{align}
 \frac12\abs{\Uphys^\e}_{2,+}^2+\e^{-1}\nm{\qk\Uphys^\e}_{L^2(\Om\times\Sone)}^2&=\frac12\abs{g^\e}_{2,-}^2.
 \label{eq:vd-l2-physical-energy}
\end{align}

The layer estimates \eqref{eq:vd-l2-layer-Linfty}--\eqref{eq:vd-l2-layer-Lr} hold for $0<\e\le1$. Its traces satisfy $\babs{\UsideZero^\e}_{\infty,\pm}\le C_\Om\mathcal G_{\rm vd}$ by \eqref{eq:vd-l2-leading-side-trace} and $\abs{a_{j,\e}}\le1$. Combining these bounds with \eqref{eq:vd-l2-physical-energy}, Lemma~\ref{lem:physical-comparison-comp}, and \eqref{eq:vd-l2-harmonic-estimate}, every norm in the theorem is bounded by $C_\Om\mathcal G_{\rm vd}$. The kinetic boundary measure is finite. Since $\e\ge\e_0(\Om)$ in this case, increasing $C_\Om$ proves all the estimates on this interval as well.
\end{proof}
\section{Numerical Illustration of Vertex Singularity}
\label{sec:numerical-vertex-singularity}

This appendix illustrates the possible singularity of the leading interior solution $\rho_0$ discussed in Subsection~\ref{subsec:intro-difficulties}. Even with smooth Dirichlet data on a convex polygon, its Hessian can fail to belong to $L^\infty$ near a vertex. These derivatives enter the higher-order interior solution and the data for the side layers, which motivates regularizing their finite vertex expansions near each vertex. The regularization introduces errors in the transport equation and the in-flow boundary condition. The wedge layer compensates these errors together with those arising from the interaction of the two incident side layers. The comparison below illustrates the elliptic singularity underlying this construction.

\paragraph{\underline{Dirichlet problem}} Let $\Omega_{\mathrm{hex}}$ be the regular hexagon of side length one with vertices $\left(\cos\frac{k\pi}{3},\sin\frac{k\pi}{3}\right)$, $k=0,\ldots,5$. Let $\Omega_{\mathrm{disk}}=B_1(0)$. On both domains $\Omega_\star\in\{\Omega_{\mathrm{hex}},\Omega_{\mathrm{disk}}\}$, prescribe the restriction of $g_{\rm num}(x_1,x_2)=\sin(3x_1)$ and solve
\begin{align}
 -\Delta \rho_0^\star &=0 \quad\text{in }\Omega_\star, \qquad \rho_0^\star\big|_{\partial\Omega_\star}=g_{\rm num}\big|_{\partial\Omega_\star}.
 \label{eq:numerical-laplace-problem}
\end{align}
On each side $E_j$, the tangential derivatives of the data satisfy $\partial_s^a(g_{\rm num}|_{E_j})=\big[\bigl(t_j\cdot\nabla\bigr)^ag_{\rm num}\big]\big|_{E_j}$ for $a\in\mathbb N_0$. Endpoint derivatives are taken with the tangent of the corresponding side.

For the hexagon, use conforming quadratic ($P_2$) elements on a uniform triangular mesh with $N_{\rm mesh}$ subdivisions per side. At $N_{\rm mesh}=80$, there are $6N_{\rm mesh}^2=38{,}400$ elements. We impose the Dirichlet data by nodal interpolation of $g_{\rm num}$. For the disk, use the Jacobi--Anger series
\begin{align}
 \rho_0^{\mathrm{disk}}(r,\theta)&=2\sum_{k\geq0}\bigl(-1\bigr)^kJ_{2k+1}(3)\,r^{2k+1}\cos\bigl((2k+1)\theta\bigr).
 \label{eq:numerical-disk-series}
\end{align}
Here $J_n$ is the Bessel function of the first kind. We retain $24$ odd modes in the numerical evaluation.

We compare $D_1^\star=\max_i\abs{\partial_{x_i}\rho_0^\star}$ and $D_2^\star=\max_{i,j}\abs{\partial_{x_ix_j}\rho_0^\star}$. On the disk, the derivatives are computed from the series. On the hexagon, we average the constant element Hessians at each incident finite-element node, using the element areas as weights. We then take the largest absolute component.

\begin{figure}[ht]
    \centering
    \includegraphics[width=0.75\textwidth]
{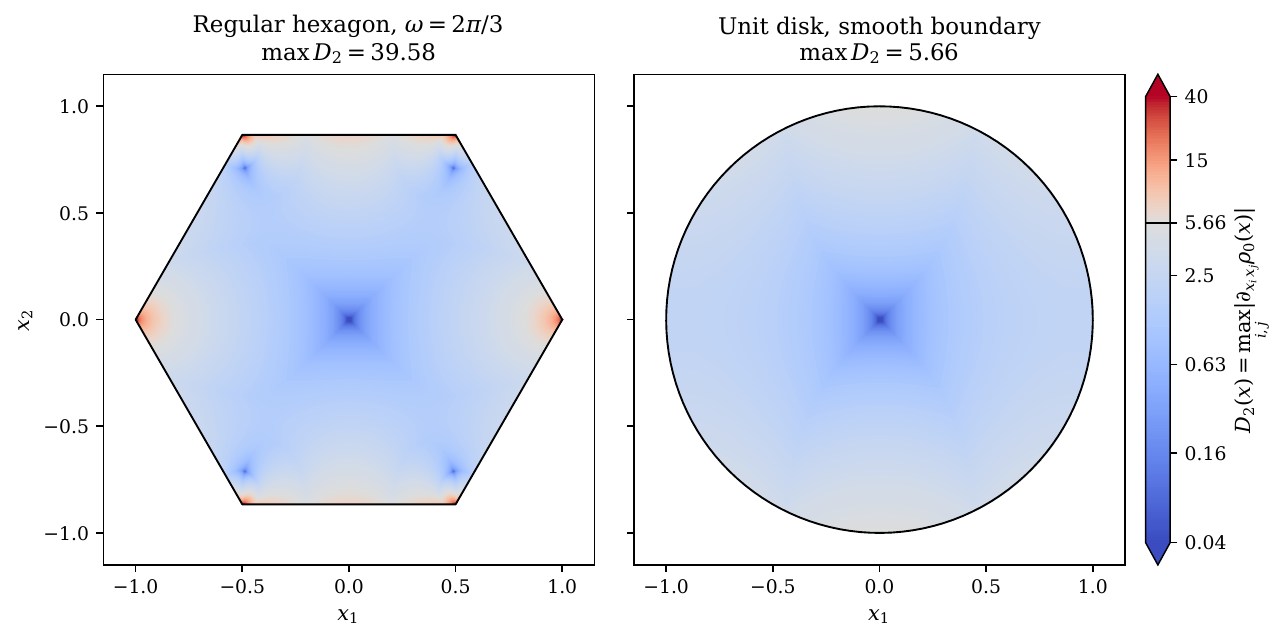}
    \caption{Comparison of $D_2^\star=\max_{i,j}\abs{\partial_{x_i x_j}\rho_0^\star}$ for the harmonic solutions of \eqref{eq:numerical-laplace-problem}}
    \label{fig:hexagon-disk-D2}
\end{figure}

\begin{figure}[ht]
    \centering
    \includegraphics[width=0.70\textwidth]
{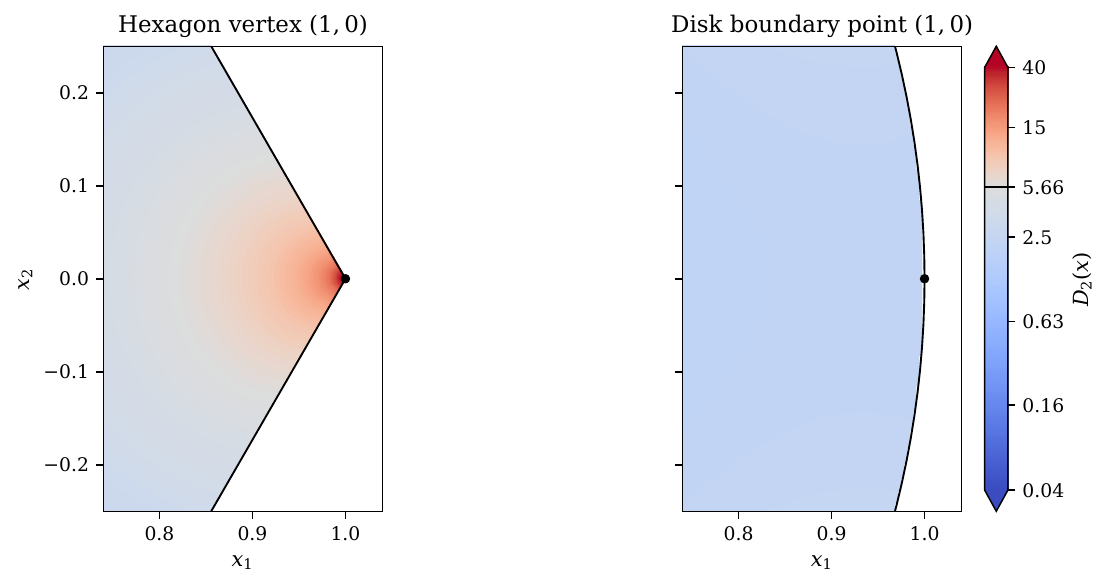}
    \caption{Matched close-ups of $D_2^\star$ near $(1,0)$}
    \label{fig:hexagon-disk-D2-zoom}
\end{figure}

\paragraph{\underline{Vertex singularity}} The smooth function $g_{\rm num}$ is defined on the whole plane and provides a lift of the Dirichlet data. The difference $\rho_0^{\mathrm{hex}}-g_{\rm num}$ has zero boundary values and satisfies $-\Delta(\rho_0^{\mathrm{hex}}-g_{\rm num})=\Delta g_{\rm num}$ in the hexagon. Its source is smooth, so the Dirichlet vertex expansion applies. At each vertex, $\omega=\frac{2\pi}{3}$ and the first singular exponent is $\frac\pi\omega=\frac32$. In polar coordinates $(r,\phi)$ centered at a vertex $P$, with $0<\phi<\omega$ between the incident sides, we have
\begin{align}
 \rho_0^{\mathrm{hex}} &= \rho_{\mathrm{reg},P} + c_P\, r^{\frac{3}{2}} \sin\Big(\frac{3\phi}{2}\Big) + \mathcal R_P,
 \label{eq:numerical-vertex-expansion}
\end{align}
where $\rho_{\mathrm{reg},P}$ is smooth near $P$ and $D^q\mathcal R_P=o(r^{\frac32-q})$ for $q=0,1,2$. If $c_P\ne0$, then the gradient belongs to $L^\infty$ near $P$ and $\abs{D^2\rho_0^{\mathrm{hex}}}\simeq r^{-\frac12}$ as $r\to0$. In this case, the Hessian belongs locally to $L^p$ for $1\le p<4$, but does not belong to $L^\infty$ near $P$. In particular, the solution is locally in $H^2$. The disk solution is smooth up to the boundary.

The global Dirichlet problem determines $c_P$, as in Remark~\ref{rem:zero-ray-mode}. The concentration near the vertices in Figure~\ref{fig:hexagon-disk-D2} is consistent with a nonzero value of this coefficient.

The computed maxima of $D_1^{\mathrm{disk}}$ and $D_1^{\mathrm{hex}}$ are about $3$. The maximum of $D_2^{\mathrm{disk}}$ is about $5.66$. On the hexagon, the recovered finite-element Hessian has maximum about $39.58$ at $N_{\rm mesh}=80$. If $c_P\ne0$, then \eqref{eq:numerical-vertex-expansion} gives $D_2^{\mathrm{hex}}\simeq N_{\rm mesh}^{\frac12}$ at distance $r\simeq N_{\rm mesh}^{-1}$ from $P$. Figures~\ref{fig:hexagon-disk-D2} and~\ref{fig:hexagon-disk-D2-zoom} show the two domains and the neighborhoods of $(1,0)$ with the same logarithmic color scale $[0.04,40]$.

\bibliographystyle{siam}
\bibliography{Reference}

\end{document}